\documentclass[10pt,reqno]{amsart}  
\usepackage[dvipdfmx]{graphicx}
\usepackage{amssymb,color}                 
\allowdisplaybreaks[1]
\newtheorem{thm}{Theorem}[section]
\newtheorem{cor}[thm]{Corollary}
\newtheorem{prop}[thm]{Proposition}
\newtheorem{lem}[thm]{Lemma}

\newcommand{\sect}[1]{\section{#1}%
\setcounter{equation}{0}}
\renewcommand{\theequation}{\thesection.\arabic{equation}}
\newenvironment{prf}[1]
   {{\noindent \bf Proof of {#1}.}}{\hfill \qed}

\newcommand{\nn}{\nonumber}  
\newcommand{\wh}{\widehat}
\renewcommand{\r}{\rho}                
\newcommand{\pt}{\partial}             
\renewcommand{\th}{\theta}

\newcommand{\D}{\mathcal{D}}
           
\newcommand{\F}{\mathcal {F}}  
            
\renewcommand{\S}{\mathcal {S}}

\newcommand{\Z}{\mathbb {Z}}

\newcommand{\re}{\mathbb R}

\newcommand{\Nt}{\mathbb N}

\newcommand{\Om}{\Omega}
\newcommand{\al}{\alpha}

\newcommand{\gm}{\gamma}
\newcommand{\Gm}{\Gamma}
\renewcommand{\L}{\mathcal{L}}

\newcommand{\ep}{\varepsilon}
\newcommand{\lam}{\lambda}
\newcommand{\om}{\omega}
\newcommand{\del}{\delta}
\newcommand{\Del}{\Delta}

\newcommand{\sg}{\sigma}
\newcommand{\s}{\sigma}
\newcommand{\x}{\xi}

\renewcommand{\r}{\rho}
\renewcommand{\t}{\tau}
\newcommand{\bd}{\widetilde{d}}

\newcommand{\teta}{\tilde{\eta}}

\newcommand{\tu}{\tilde{u}}

\newcommand{\tp}{\tilde{p}}

\newcommand{\td}{\tilde{d}}

\newcommand{\bu}{\bar{u}}
\newcommand{\bp}{\bar{p}}

\newcommand{\bet}{\bar{\eta}}

\newcommand{\N}{\nabla }

\newcommand{\dsp}{\displaystyle}
\renewcommand{\div}{{\rm {div\, }}}
\newcommand{\rot}{{\rm {rot\, }}}

\newcommand{\cof}{{\rm {cof\, }}}
\newcommand{\dB}{\dot{B}}
\newcommand{\dF}{\dot{F}}
\newcommand{\dH}{\dot{H}}

\def\<{\langle }

\renewcommand{\>}{\rangle }

\newcommand{\gr}{\Green}

\newcommand{\rd}{\color{red}}
\newcommand{\crd}{\Crimson}
\newcommand{\bl}{\color{blue}}

\newcommand{\mt}{\color{magenta}}
\newcommand{\ppl}{\Purple}
\newcommand{\bk}{\color{black}}
\renewcommand{\gr}{\color{black}} 
\renewcommand{\bl}{\color{black}} %
\renewcommand{\mt}{\color{black}} %

\renewcommand{\crd}{\color{black}} %
\renewcommand{\rd}{\color{black}}

\renewcommand{\ppl}{\color{black}}

\newcommand{\supp}{\text{ supp }}

\renewcommand{\qed}{\qquad\kern1pt   
   \vbox{\hrule height 0.6pt      
         \hbox{\vrule width 0.6pt 
               \vbox{\vskip 6pt}  
               \hskip 3pt
              \vrule width 1.3pt} 
         \hrule depth 1.3pt}     
   \kern1pt}

\newcommand{\beq}{\begin{equation}}
\newcommand{\eeq}{\end{equation}}
\newcommand{\bea}{\begin{eqnarray}}
\newcommand{\eea}{\end{eqnarray}}
\newcommand{\eq}[1]{{\begin{equation}#1\end{equation}}}
\newcommand{\spl}[1]{{\begin{aligned}#1\end{aligned}}}
\newcommand{\eqn}[1]{\begin{equation*}#1\end{equation*}}

\newcommand{\algn}[1]{\begin{align*}{#1}\end{align*}} 
\newcommand{\alg}[1]{\begin{align}{#1}\end{align}} 
\newcommand{\eqntag}{\addtocounter{equation}{1}\tag{\theequation}} 

\begin{document}
\baselineskip 4.1mm

\title{Maximal $L^1$-regularity and  
free boundary problems for the incompressible 
Navier--Stokes equations \\
in critical spaces 
}

\maketitle
\vskip-1mm
 \begin{center}
{\footnotesize
{\large\sf Takayoshi Ogawa${}^{*}$}  and
{\large\sf Senjo Shimizu${}^{\ddag}$}
\vskip3mm
\vbox{
        Mathematical Institute${}^{*}$\\
        Tohoku University  \\
        Sendai  980-8578, Japan  \\
        takayshi.ogawa.c8@tohoku.ac.jp\\
      }
\hskip5mm\noindent
\vbox{
    {\gr
        Department of Mathematics${}^{\ddag}$,\\
        Faculty of Science, 
    }
        Kyoto University  \\
        Kyoto 606-850{\gr 2}, Japan  \\
        shimizu.senjo.5s@kyoto-u.ac.jp
      }
}
\end{center}        
\begin{abstract}
Time-dependent free surface problem for the incompressible 
Navier--Stokes equations which 
describes the motion of viscous incompressible fluid nearly half-space 
are considered. 
We obtain global well-posedness of the problem for a small initial data 
in scale invariant critical Besov spaces. 
Our proof is based on maximal $L^1$-regularity of the corresponding 
Stokes problem in the half-space and special structures of the 
quasi-linear term appearing from the Lagrangian transform of the 
coordinate. 
\end{abstract}      
        
\section{Introduction and Main Results}\label{Sec1}
We consider a time-dependent free surface problem for 
the Navier--Stokes equations which describes the motion of 
viscous incompressible fluid. 
The domain $\Omega_t\subset \re^n$ ($n\ge2$) is occupied 
by the fluid and the velocity of fluid $\bu(t,y)$ and   
the pressure $\bp(t,y)$ for $y\in\Omega_t$ satisfy the   
incompressible Navier--Stokes equations:
\begin{equation} \label{eqn;FNS}
   \left\{
   \begin{split}
     \pt_t \bu + &(\bu\cdot \N)  \bu-\div T(\bu,\bp) =0,
        \quad&\,  &t>0,\ \ y\in\Omega_t,  \\
     &\qquad \div \bu=0,
        \quad&\,  &t>0,\ \  y\in \Omega_t, \\
     &\quad T(\bu,\bp) \nu_t =0,
        \quad&\,  &t>0,\ \  y\in \pt\Om_t,\\
     &\qquad \bu(0,y) =u_0(y),
        \quad&\,  &\phantom{t>0,}\ \  y\in \Om_0.
     \end{split}
     \right.  \hfill
\end{equation}
Here, $\pt\Om_t$ denotes the boundary of $\Omega_t$, 
$\nu_t=\nu_t(y)$ is the unit outward normal at a point $y\in \pt\Om_t$, 
$T(\bu,\bp)$ is the stress tensor defined by 
$T(\bu,\bp)=(\nabla\bu+(\nabla\bu)^{\sf T})-\bp I$, where $I$ is 
the $n\times n$ identity matrix, 
$(\nabla_y\bu)_{ij}=\frac{\pt \bu_j}{\pt y_i}$, and 
$(\nabla\bu)^{\sf T}$ denotes the transposed matrix of $\nabla\bu$. 
$u_0$ is the given initial velocity. 
In our setting \eqref{eqn;FNS}, we do not take into account the effect 
from  the gravity force or the surface tension.  
\par
Free boundary problems for the incompressible fluids 
are considered by many authors.
The pioneer work was done by Solonnikov \cite{Sol77}, he established 
local well-posedness of \eqref{eqn;FNS} 
whose initial state $\Omega$ is a bounded domain in the frame work 
of H\"older spaces $C^{2+\alpha,1+\alpha/2}$ with $\alpha\in (\frac12, 1)$. 
Solonnikov also proved global well-posedness of \eqref{eqn;FNS} 
in the class of Sobolev space $W^{2,1}_p$ with $n<p<\infty$ when $n=2$, $3$, 
where surface tension is excluded.  
When initial state is {\ppl bounded} and the surface tension is excluded,  
Mucha--Zajaczkowski considered the case where the self-gravitational 
force exists,  they proved in  \cite{MuZa00},  \cite{MuZa} the local 
in time unique solvability in $W^{2,1}_p$ with $n=3$ and 
$3<p<\infty$ for arbitrary initial data. Shibata--Shimizu \cite{SbSz07}, 
\cite{SbSz08} developed the $L^p$-theory for the problem and showed 
global well-posedness of \eqref{eqn;FNS} in the class of Sobolev space 
$W^{2,1}_{q,p}$ with $n<q<\infty$ and $2<p<\infty$ when $n\ge 2$. 
\par
In the case when initial state is {\ppl bounded} and the surface tension is included,  
Solonnikov proved the global in time solvability in $W^{2+\alpha,1+\alpha/ 2}_2$
with $1/2<\alpha<1$ provided that initial data are sufficiently small 
and the initial domain is sufficiently close to a ball. 
There are many other contributions in the case when the effect of  
surface tension is included, for instance \cite{KPW13}, \cite{18}, 
\cite{26}, \cite{27}, \cite{Sol84}-\cite{Sol91} and references therein, 
we do not get involved with the case because in this paper we consider 
without surface tension case.  

Another typical free boundary problem describes the motion of a fluid 
which occupies a semi-infinite domain between the moving upper surface 
and a fixed bottom. 
Beale \cite{B81}, \cite{B84} considered the free surface problem in a three 
dimensional region with a bottom, in the $L^2$-based Bessel potential spaces 
$H^{r/2,r}_2$ where $3<r<7/2$. 
His problem (called as the ocean problem) has a similar setting of 
the following Lagrange coordinate equations and showed that the global 
in time solvability for small initial data in $L_2$ Bessel-potential 
space setting. Since the ocean problem has a finite depth, however, 
the spectral property for the linearized problem is different from the case for 
a domain close to the half-space. 
Pr\"uss--Simonett \cite{PrSi10}, \cite{PrSi16} proved local well-posedness 
of \eqref{eqn;FNS} whose initial state $\Omega_0$ is 
{\ppl close} to the half-space 
$\re^n_+$ in the class of Sobolev space $W^{2,1}_p$ 
with $p>n+2$. 
There are many other contributions on this direction, for instance, 
\cite{Ab}, \cite{BN85}, \cite{BNT20}, \cite{DM09}-\cite{DM15}, \cite{Gui}, 
\cite{Gu-Tc}, \cite{MuZa00}-\cite{NTY04},   
 \cite{PrSi10}-\cite{Sa}, 
\cite{SbSz07},  \cite{Sl-Tn}-\cite{Tn2} 
and reference therein.

Recently, Shibata \cite{Sb16}, \cite{Sb20} considered local and global 
well-posedness on general unbounded domain in the space $W^{2,1}_{q,p}$ 
with $n<q<\infty$ and $2<p<\infty$.

The incompressible Navier--Stokes equations {\ppl are} invariant 
under the following scaling: For all $\lam>0$,
\eqn{ \label{eqn;inv-scaling-NS}
 \left\{
 \spl{
   & \bu(t,y)\to \bu_{\lam}(t,y)\equiv \lam  \bu(\lam^2 t, \lam y), \\
   & \bp(t,y)\to \bp_{\lam}(t,y)\equiv \lam^2\bp(\lam^2 t, \lam y).
  }
 \right.
}
Subsequently, it is well-known that the Cauchy problem of the Navier--Stokes 
equations can be solved globally in time in the invariant Bochner--Sobolev 
space  $L^{\r}\big(\re_+;\dot{H}^{s}_p(\re^n;\re^n)\big)$
\eq{  \label{Serrin}
  \frac{2}{\r}+\frac{n}{p}=1+s,
} 
which is observed in the celebrated result by   
Fujita--Kato \cite{F-K} (see also Prodi \cite{Pr} and Serrin \cite{Sr} 
for the relation between regularity of solutions and the scaling invariance).
When we choose $\r=\infty$, we obtain $s=-1+n/p$ by \eqref{Serrin}, 
and the critical class at $s=0$ is given, 
in particular, by $L^{\infty}(0,T;L^n(\re^n))$, where
Kato \cite{Kt1} considered global well-posedness of the Cauchy problem. 
Such a critical setting for the Cauchy problem is considered by several 
authors in the framework of the scaling critical Besov spaces 
$\dB^{-1+n/p}_{p,\sg}(\re^n)$,  
where $1\le p<\infty$ and $1\le \sg \le \infty$
(\cite{Am00}, \cite{Cn}, \cite{Cn-Pl}, {\ppl \cite{C-L95}},  \cite{K-Y}).
Meanwhile, it is proved ill-posedness 
of the problem in \cite{B-P}, \cite{Wn},  \cite{Yn}, 
namely the continuous dependence on the initial 
data in the classes $u_0\in \dB^{-1}_{\infty,\sg}(\re^n)$,  
$1\le \sg \le \infty$ breaks down. 
In view of those of well-posedness results to the Cauchy problem,
it is natural to ask if the free surface problem can be 
solvable in such a scaling critical function class. 
Our main motivation is to consider the free surface problem \eqref{eqn;FNS} 
near the half-space $\re^n_+$ in the scaling critical function space.
\par
In this paper, we show global in time well-posedness of 
the Lagrangian transformed problem for \eqref{eqn;FNS} 
under small data in the scaling critical Besov space 
$\dB^{-1+n/p}_{p,1}(\re^n_+)$ for all $n\le p< 2n-1$, via maximal 
$L^1$-regularity of the linearized problem associated with \eqref{eqn;FNS}. 
As far as the authors know, there is almost no result of global 
well-posedness to \eqref{eqn;FNS} in the scale critical space whose 
initial state $\Omega_0$ is an unbounded domain except the recent result 
due to Danchin--Hieber--Mucha--Tolksdorf \cite{DHMT20}. 
They consider the analogous problem 
in the scaling critical Besov spaces 
{\gr for $n-1<p<n$}.

Let the half {\ppl Euclidean} space and its boundary be denoted by
\eqn{
 \spl{
\re^n_+\equiv &\{(x', x_n);\ x'\in \re^{n-1},\ x_n>0\},\\ 
\pt\re^n_+\equiv &\re^{n-1}\times\{0\}=\{(x', x_n);\ x'\in \re^{n-1},\ x_n=0\}. 
 }
 }
We also set $\re^n_-$ as the negative part of {\ppl $\re^n$},
i.e., $\re^n=\re^n_+\cup\pt\re^n_+\cup \re^n_-$.
 Aside from the dynamical boundary condition, 
a further kinematic condition 
for the free surface is satisfied which gives 
$\pt\Om_t$ as a set of points $y=y(t,x)$, 
$x\in \pt\Om_0=\pt\re^n_+$, 
where $y(t,x)$ is the solution of the Cauchy problem:    
\begin{equation} \label{def-flow}
  \frac{dy}{dt}=\bu\big(t,y(t)\big),\ \  t>0,
  \qquad y(0)=x. 
\end{equation}
Let the Euler coordinates  $y\in \Om_t$
be transformed into the Lagrangian coordinates 
$x\in\re^n_+$ connected by \eqref{def-flow}. 
If $\bu(t,y)$ is Lipschitz continuous with respect to $y$, 
then \eqref{def-flow} can be 
solved uniquely by  
\begin{equation} \label{def-flow-int}
    y(t,x)=x+\int_0^t \bu\big(s,y(s,x)\big) ds. 
\end{equation}
By the kinematic condition of the original boundary $\Om_t$, 
it is described by the map 
$Y_{\bu}: (t,x)\in \re_+\times \re^n_+\to (t,y)\in \Om_t$, where  
$\Om_t$ is given by 
$$
\Om_t=Y_{\bu}(t,\re_+^n)\equiv 
\big\{(t,y(t,x));\; t>0, \,   \text{ $y(t,x)$ satisfies 
\eqref{def-flow-int} and $x\in \re^n_+$} \big\}.
$$
Setting 
\eqn{
 \left\{
 \spl{
   & u(t,x)\equiv\bu(t,y(t,x)), \\
   & p(t,x)\equiv\bp(t,y(t,x)), 
     }
 \right.
    }
and applying the Lagrangian coordinate to the original problem 
\eqref{eqn;FNS} yields that the system is transformed into 
the following form:
\begin{equation} \label{eqn;NS}
   \left\{
   \begin{split}
    & \pt_t u  - \Del u + \N p =  F_u(u)+F_p(u,p),
        &\quad &t>0,\ \ x\in\re^n_+,\\
    & \qquad \div\, u =  G_{\div}(u),        
        &\quad &t>0,\ \ x\in\re^n_+,\\
    & \big(\N u+(\N u)^{\sf T}-pI)\big)\nu_n
                =H_u(u)+H_p(u,p),        
        &\quad &t>0,\ \ x\in \pt \re^n_+,\\   
    & \qquad  u(0,x) =u_0(x),
        &\quad &\phantom{\ t>0,}\ \ 
               x\in\Bbb R^n_+,
     \end{split}
     \right.   \hfill
\end{equation}
\vskip2mm\noindent
where $\nu_n=(0,\cdots, 0,-1)^{\sf T}$ denotes the outward normal
\footnote{Practically natural setting is 
 $\Om_0=\re^n_-$ under the gravity circumstance.}
 and the nonlinear terms of \eqref{eqn;NS} are given by 
\begin{align}
   F_u (u) 
      \equiv\, &
      \div \Big(J(Du)^{-1}\big(J(Du)^{-1}\big)^{\sf T} \N u -\N u\Big)       
      =\Pi_{{u}}^{2n-2}
           \left(\int_0^t D u\,ds\right)
           D^2 u, 
    \label{eqn;laplace-purterb}
\\
  F_p(u,p)
      \equiv\, &
        -\big(J(D u)^{-1}-I\big)^{\sf T}\,\N p 
       = \Pi_{p}^{n-1}
             \left(\int_0^t D u\,ds\right)
             \N p, 
    \label{eqn;pressure-purterb}
\\
  G_{\div}(u)
     \equiv\, 
     &-\text{tr}\Big(\big(J(D u)^{-1}-I\big)^{\sf T}\; \N u\Big)
         = 
          \text{tr}\left(\Pi_{\div}^{n-1}
          \left(\int_0^t D u\,ds\right)
           D u\right) \notag\\
        =\, &\div\left(
             \Pi_{\div}^{n-1}\left(\int_0^t D u\, ds\right)
                 u\right),
    \label{eqn;divergence-purterb}
\\
  H_{u}(u)\equiv\, & 
    -\Big(\big(J\big(Du\big)^{-1}\big)^{\sf T}\;\N u  
             +(\N u)^{\sf T}\; J\big(Du\big)^{-1}   
     \Big) 
    \big(J(Du)^{-1}-I\big)^{\sf T} \nu_n  
    \notag \\
   &\quad
     -\Big(\big(J\big(Du\big)^{-1}-I\big)^{\sf T}\;\N u  
             +(\N u)^{\sf T}\;\big( J\big(Du\big)^{-1} -I\big)  
      \Big) \nu_n
    \notag \\
 = & \Pi_{bu}^{2n-2}
            \left(\int_0^t D u\,ds\right)
           D u\, \nu_n,
     \label{eqn;boundary-purterb-u} \\
 H_{p}(u,p)\equiv\, &   
     p I \big(J(Du)^{-1}-I\big)^{\sf T} \nu_n
    =   \Pi_{bp}^{n-1}
        \left(\int_0^t D u\,ds\right) p \, \nu_n.
     \label{eqn;boundary-purterb-p}      
\end{align}
Here $J(Du)^{-1}$ denotes the inverse of the Jacobian matrix,  
$I$ denotes the identity matrix,
$(D u)_{ij}=\frac{\pt u_i}{\pt x_j}$ and 
$\Pi^{m}_*(d)$ denote  {\ppl $m$-th order polynomials of}   $d=(d_{jk})_{1\le j,k\le n}$ with 
$$
 d_{jk}
  =\Big(\int_0^t D u(s)ds\Big)_{jk}
 \equiv  \int_0^t\pt_{x_k}  u_j(s)\, ds
$$ 
(here the notation $*$ stands for either $u$, $p$, $\text{div}$ 
or $bu$,  $bp$). Those polynomials  are indeed given by 
the inverse matrix of the Jacobi matrix, $J(Du)^{-1}$ as follows:
\algn{
\Pi_{*}^{m}
           \left(\int_0^t D u\,ds\right)
 =&\sum_{\ell=1}^{m}\prod^{\ell}_{1\le j_{\ell},k_{\ell}\le n}
   \sg_{k_{\ell} j_{\ell}}
   \left(\int_0^t \pt_{x_{k_{\ell}}} u_{j_{\ell}}(s,x)ds\right) 
}
with $\sg_{\ppl k_{\ell} j_{\ell}}$ is either $1$ or $-1$.

By using the Lagrangian transformation, the free surface problem \eqref{eqn;FNS}
can be transformed into the initial-boundary 
value problem in $\re^n_+$ with the fixed boundary $\pt\re^n_+$ and 
the system is transformed into the quasi-linear parabolic equation \eqref{eqn;NS}  
(see e.g., \cite{Sol88}).
\par

Before stating our results, we define the Besov spaces and Lizorkin--Triebel  
spaces in the half-space.
Since the global estimate requires the base space for spatial
variable $x$ in the homogeneous Besov space,  
we introduce the homogeneous Besov space over $\re^n_+$ 
(see for details, Bergh--L\"ofstr\"om \cite{BL}, Lizorkin \cite{Lz}, 
Peetre \cite{Pr75}, \cite{Pr76}, Triebel \cite{Tr73}-\cite{Tr83}).   

\par\vspace{1pc} \noindent
{\it Definition} (The Besov spaces).
Let $s\in \re$,  $1\le p, \sg\le \infty$. 
Let $\{\phi_j\}_{j\in \Z}$ be the Littlewood--Paley dyadic decomposition 
of unity for $x\in \re^n$, 
namely $\widehat \phi$ is the Fourier transform of a smooth radial 
function $\phi$ with 
$\widehat \phi(\xi)\ge 0$ and ${\rm supp}\,
  \widehat\phi\subset \{\xi\in\re^n\mid 2^{-1}<|\xi|<2\}$,  
and
\begin{align} 
    &\widehat\phi_j(\xi)=\widehat\phi(2^{-j}\xi),
     \quad \sum_{j\in\Z}\widehat\phi_j(\xi)=1
     \quad \text{for\ any}\ \xi\in\re^n\setminus\{0\}, 
     \quad j\in \Z  
   \nn\\
    &\text{ and }\quad 
     \widehat{\phi}_{\widehat{0}}(\xi)
    +\sum_{j\ge 1}\widehat \phi_j(\xi)=1\quad 
\text{for\ any}\ \xi\in\re^n,\label{eqn;LP-decomp}
\end{align}
where $\widehat{\phi}_{\widehat{0}}(\xi)
\equiv \widehat{\zeta}(|\xi|)$ with a low frequency cut-off 
\eq{ \label{eqn;low-freq-zeta}
   \widehat{\zeta}(r)
   =
  \begin{cases}
   1,\quad    & 0\le r<1, \\
   \text{decreasing in }& 1\le r< 2,\\
   0,\quad    & 2\le r.
   \end{cases}
 } 
 For $s\in\re$ and $1\le p,\sigma\le \infty$,
 let $\dot{B}^s_{p,\sg}(\re^n)$ be the homogeneous Besov space with norm
\eqn{
 \|\tilde{f}\|_{\dot{B}^s_{p,\sg}}
  \equiv 
  \left\{
  \begin{aligned} 
   & \Bigl(\sum_{j\in \Z}2^{s\sg j}
          \|\phi_j*\tilde{f}\|_p^{\sg}
     \Bigr)^{1/\sg},
   & 1\le \sigma<\infty,\\
   &\, \sup_{j\in \Z} 2^{s j}\|\phi_j*\tilde{f}\|_p, 
   & \phantom{1\le }\sigma=\infty,
  \end{aligned} 
  \right.
}
where $\phi_j * f$ stands for the convolution operation 
with a constant correction $c_n=(2\pi)^{-n/2}$ given by 
\eq{
 \phi_j*f=c_n\int_{\re^n}\phi_j(x-y)f(y)dy
}
for $f\in \S(\re^n)$ and its standard extension to $f\in \S'(\re^n)$.
In what follows, we always regard this correction of the constant 
against the convolution operations for all kinds of the Littlewood--Paley 
decompositions.

Also let ${B}^s_{p,\sg}(\re^n)$ be the inhomogeneous Besov space 
with norm
$$
 \|\tilde{f}\|_{B^s_{p,\sg}}
  \equiv 
  \begin{cases} \dsp
    \Bigl(\|\phi_{\hat{0}}*\tilde{f}\|_p
    +\sum_{j\in \Z}2^{s\sg j}\|\phi_j*\tilde{f}\|_p^{\sg}\Bigr)^{1/\sg},
         & 1\le \sigma<\infty,    \\
    \dsp
    \|\phi_{\hat{0}}*\tilde{f}\|_p
      +\sup_{j\in \Z} 2^{s j}\|\phi_j*\tilde{f}\|_p, 
         & \phantom{1\le }\sigma=\infty.
  \end{cases} 
$$

We define the homogeneous Besov space $\dot{B}^s_{p,\sg}(\re^n_+)$ 
as the set of all {\ppl the restriction $f$ of the distribution $\tilde{f}\in \dB^s_{p,\sg}(\re^n)$, i.e., 
$f=\tilde{f}\big|_{\re^n_+}$  with}
\begin{equation*}
  \|f\|_{\dot{B}^s_{p,\sg}(\re^n_+)}
   \equiv 
   \inf \Bigg\{ 
      \|\tilde f\|_{\dB^s_{p,\sg}(\re^n)}<\infty; 
      \quad\tilde f =\sum_{j\in\Z}\phi_j*\tilde{f} 
              \text{ in } \mathcal{S}',
             \quad   f=\tilde{f}\big|_{\re^n_+}
      \bigg\}.
  \eqntag \label{eqn;Besov-halfspace}
\end{equation*}
Analogously we define the inhomogeneous Besov space ${B}^s_{p,\sg}(\re^n_+)$ 
in a similar manner.

\vskip1mm
\noindent
{\it Definition} (The Bochner--Lizorkin--Triebel  spaces).
Let $s\ge 0$,  $1\le p, \sg\le \infty$ and $X$ be a Banach 
space  with the norm $\|\cdot \|_{X}$.
Let $\{\psi_k\}_{k\in \Z}$ be the Littlewood--Paley dyadic decomposition 
of unity for $t\in \re$. 
For $s\in\re$ and $1\le p< \infty$, 
$\dot{F}^s_{p,\sg}(\re;X)$ be 
the Bochner--Lizorkin--Triebel  space with norm
\eqn{
 \|\tilde{f}\|_{ \dot{F}^s_{p,\sg}(\re;X) }
  \equiv 
  \left\{
  \begin{aligned} 
   & \Big\|\Big(\sum_{k\in \Z}2^{s\sg k}
         \|\psi_k*\tilde{f}(t,\cdot)\|_X^{\sg}
     \Bigr)^{1/\sg}\Big\|_{L^p(\re)},
   & 1\le \sigma<\infty,\\
   &\,\Big\| \sup_{k\in \Z} 
       2^{s k}\|\psi_{\ppl k}*\tilde{f}(t,\cdot)\|_X
      \Big\|_{L^p(\re)}, 
   & \phantom{1\le }\sigma=\infty.
  \end{aligned} 
  \right.
}

Analogously above, we define the Bochner--Lizorkin--Triebel  spaces $\dot{F}^s_{p,\sg}(I;X)$ 
as the set of all {\ppl the restriction $f$ of a distribution $\tilde{f}\in \dot{F}^s_{p,\sg}(\re;X)$ i.e.,
$f=\tilde{f}\big|_I$ } on $X$ with 
$$
  \|f\|_{\dot{F}^s_{p,\sg}(I;X)}
   \equiv \inf\left\{
           \|\tilde f\|_{\dot{F}^s_{p,\sg}(\re;X)}<\infty;\quad
             f=\tilde f\big|_{I}
             \right\},
$$ 
where $I=(0,T)$ denotes the time interval. 
We denote $\re_+=(0,\infty)$ as the half real line and 
$\overline{\re_+}=[0,\infty)$ as its closure. 
\noindent
We note that all those homogeneous spaces are understood as 
the Banach spaces by introducing the quotient spaces identifying 
any difference of polynomials.

\vskip2mm
Let $C_b(I;X)$  be  a set of all bounded continuous 
functions from an interval $I$ to a Banach space $X$.  
We also use the notation $C_v(\re^n_+)$ for a set of all 
continuous functions vanishing at $|x|\to \infty$.
Obviously $C_v(\re^n_+)\subset C_b(\re^n_+)$.
\par
\medskip

\begin{thm}[Global well-posedness under the Lagrangian coordinates] 
\label{thm;main} Let  $n\le p< 2n-1$.  
{\ppl There exists} small $\ep_0>0$ such that if
the initial data $u_0\in \dB^{-1+{\rd n/p}}_{p,1}(\re^n_+)$ 
{\ppl with $\div u_0=0$ in the sense of distribution satisfying}  
\eq{ \label{eqn;ini}
  \|u_0\|_{\dot B^{-1+\frac np}_{p,1}(\re^n_+)}\leq\varepsilon_0, 
}
then \eqref{eqn;NS} admits a unique global solution
\begin{align*}
  &u\in C_b(\overline{\re_+}; \dot B^{-1+\frac np}_{p,1}(\re^n_+))\cap
     \dot W^{1,1}(\re_+; \dot B^{-1+\frac np}_{p,1}(\re^n_+)),\\
  &\Del u,\; \N p \in  L^1(\re_+; \dot B^{-1+\frac np}_{p,1}(\re^n_+)),
     \\
  & p|_{x_n=0}
  \in \dF^{\frac12-\frac{1}{2p}}_{1,1}(\re_+;\dB^{-1+\frac{n}{p}}_{p,1}(\re^{n-1}))
      \cap 
      {\crd L^1}(\re_+; \dot B^{\frac{n-1}p}_{p,1}(\re^{n-1}))
\end{align*}
with the estimate
\eq{\label{eqn;mainest}
 \spl{
  \|\pt_t u &\|_{L^1(\re_+; \dot B^{-1+\frac np}_{p,1}(\re^n_+))} 
   +\|D^2 u\|_{L^1(\re_+; \dot B^{{\crd -1}+\frac np}_{p,1}(\re^n_+))}
   +\|\nabla p\|_{L^1(\re_+; \dot B^{-1+\frac np}_{p,1}(\re^n_+))} \\
   &+\big\| p|_{x_n=0}
     \big\|_{\dF^{\frac12-\frac{1}{2p}}_{1,1}(\re_+; \dot B^{-1+\frac np}_{p,1}(\re^{n-1}))}
    +\big\| p|_{x_n=0}\big\|_{{\crd L^1}(\re_+; \dot B^{\frac{n-1}p}_{p,1}(\re^{n-1}))}
   \le \varepsilon_1, 
}
}
where $D^2 u$ denotes all the second order derivatives of $u$ 
by $x$ and 
 $\varepsilon_1=\ep_1(n,p,\varepsilon_0)$ is a constant.
\end{thm}

\begin{cor}[Global well-posedness]\label  {cor;main} 
Let $n\le p< 2n-1$. For the same $\varepsilon_0$ in Theorem \ref{thm;main} and 
$u_0\in \dot B^{-1+n/p}_{p,1}(\re^n_+)$ with {\ppl $\div u_0=0$ in the sense of distribution} satisfying \eqref{eqn;ini}, 
let $(u,p)$ be the global solution of \eqref{eqn;NS} obtained in Theorem \ref{thm;main}. 
Then the pull-back $(\bu,\bp)$ of $(u,p)$ with the estimate \eqref{eqn;mainest} 
satisfies \eqref{eqn;FNS}. 
\end{cor}

Concerning the half-space problem, Danchin--Mucha \cite{DM09} {\ppl proved} well-posedness of the 
Cauchy--Dirichlet problem of the density-dependent incompressible 
Navier--Stokes equations with the $0$-Dirichlet boundary data. 
The result there is also applicable for the 
incompressible Navier--Stokes equations in the scaling invariant Besov 
spaces $p=n$, namely in $\dB^{0}_{n,1}(\re^n_+)$.  
\par

Let us mention on the two results between Danchin--Hieber--Mucha--Tolksdorf 
\cite{DHMT20} and ours.  Their result based on the abstract 
interpolation spaces based on the original idea that goes back to 
Da Prato--{\ppl Grisvard} \cite{dPG} and based on the result due to 
Danchin--Mucha \cite{DM09} and \cite{DM12}, where the authors considered the 
$0$-Neumann boundary condition for the linearized system of 
the Stokes equation. Our approach is very much different from theirs. 
We  handle 
the boundary potential for non-stress boundary condition directly
in the homogeneous Besov spaces 
and as a result,  our result covers the initial data 
as a class of distributions (negative indices of regularity 
in the scaling critical homogeneous Besov space 
$\dot{B}^{-1+n/p}_{p,1}(\re^n_+)$
with $n\le p< 2n-1$), while they treats the function case {\gr $n-1<p<n$ for 
$\dot{B}^{-1+n/p}_{p,1}(\re^n_+)$ case in \cite{DHMT20}}.  
Theorem \ref{thm;main}
(and hence Corollary \ref{cor;main}) is the 
first result for treating the scaling invariant distribution as 
an initial data for the free boundary value problem, as far as the 
authors can find.

\vskip1mm
Our proof of Theorem \ref{thm;main} is 
heavily depending on the end-point estimate of maximal regularity
for the initial-boundary value problem of the Stokes system in the 
half-space $\re^n_+$.
Many of the existence results are related to the spectral analysis for 
the linearized equation and derive the decay property for the linearized 
Stokes equations.  On the other hand, maximal regularity 
for the parabolic equation gives a suitable estimate for treating
the quasi-linear terms effectively 
(\cite{Am95}, \cite{Am19}). 
In contrast with those results, our method is a direct application 
of maximal $L^1$-regularity for the half time line $\re_+$  
to treat the system under the Lagrange transformations.  
This method enables us to handle 
main terms appearing the quasi-linear perturbations 
 \eqref{eqn;laplace-purterb}--\eqref{eqn;boundary-purterb-p} {\ppl directly and we may
treat them} globally in the transformed problem \eqref{eqn;NS}. 
Namely, to obtain global well-posedness of \eqref{eqn;NS}, 
it is required to treat the terms with
$$
 \tilde{d}_{kj}\equiv 
 \del_{k,j}+ \lim_{t\to \infty}\int_0^t \pt_{x_k} u_j(s,x)ds,
 \qquad  k,j=1,2,\cdots, n. 
$$ 
We then establish maximal $L^1$-regularity for the transformed Stokes system
via maximal regularity estimate for the initial-boundary value problem of 
the heat equations obtained in the previous work of authors \cite{OgSs20-2}
(see for its announce \cite{OgSs20-1}).
Such argument was developed by Danchin \cite{Dc03}, \cite{Dc07} for 
treating global well-posedness for the Cauchy problem of the
compressible or incompressible density dependent Navier--Stokes equations.
The main difference here is to treat the boundary inhomogeneous terms
appearing in $H_b(u)$ by maximal $L^1$-regularity and usage of the sharp
trace estimate of the boundary terms.
Such an estimate is available for analyzing the potential expression of 
the pressure term $p$ for the Stokes system with the free surface boundary 
condition obtained in \cite{SbSz07}.  Maximal regularity and its sharpness 
is obtained by establishing the almost orthogonal estimates for the pressure 
potential and the Littlewood--Paley space-time decompositions of unity that 
defines our sharp function class of the well-posedness.

In order to enlarge the solution class into the critical Besov spaces, 
the divergence free condition is crucial.  In particular to 
enlarge the class for the bilinear estimate remains valid,  
the multiple divergence-rotation-free structure is another crucial point
(cf. \cite{OgSs18}).  This nonlinear structure was partially observed 
by Solonnikov \cite{Sol77}
and Shibata--Shimizu \cite{SbSz03} for treating the terms in the Sobolev spaces.
However in order to apply the bilinear estimate in the critical Besov space,
we need to ensure such a special structure for each decomposition 
steps of sub-matrix expansion of the inverse of Jacobi matrix.
In this stage, we show that  {\it a divergence-curl free structure}
(div-curl structure, in short) holds 
for each step of sub-cofactor of expansion involving the 
{\it null-Lagrangian structure} (cf. Evans \cite{Ev}).
This was shown in \cite{OgSs18} for the initial value problem for 
the Lagrangian coordinate case. 
We develop the analogous estimate  and establish the multiple Besov 
estimate in the half-spaces. 
It is well-known that the convection term $\bar{u}\cdot \N \bar{u}$ 
maintains the div-curl structure and it helps to enlarge the solution class. 
Although the convection term vanishes after the transformation into 
the Lagrangian coordinate,  all the nonlinear terms inherit 
the div-curl structure from the divergence free condition 
and then the solution class can be 
reach the critical homogeneous Besov space.
 
\vskip2mm
We should like to notice that regularity for the solution obtained in 
Theorem \ref{thm;main} is weaker than known results, we do not assume 
the compatibility conditions on the initial and boundary data.
The regularity of solution ensures us that the velocity fields has 
a sufficient regularity $\N u\in L^1((0,\infty);\dB^{n/p}_{p,1}(\re^n_+))$ 
so  that the Lagrange transformation \eqref{def-flow-int} is uniquely 
determined and the inverse of the transformation has meaningful by 
$\dB^{n/p}_{p,1}(\re^n_+)\subset  C_v(\re^n_+)$.  
Thus the original problem \eqref{eqn;FNS} is solvable.
\vskip1mm

The rest of this paper is organized as follows. 
We present a solution formula of the linear problem of 
\eqref{eqn;FNS} in the next section.
Maximal $L^1$-regularity of the Stokes system 
(Theorem \ref{thm;L1MR} stated in Section \ref{Sec;2})  
 is a key estimate for our argument. 
Section \ref{Sec;3} is devoted to prove almost orthogonality 
between the pressure potential and the space-time 
Littlewood--Paley dyadic decomposition, which is {\ppl crucial}  
to prove maximal $L^1$-regularity of the Stokes system.
Using the almost orthogonal estimates, we show maximal regularity for 
the Stokes system in Section \ref{Sec;4}. 
The bilinear estimates as well as the div-curl lemma are discussed 
in Section \ref{Sec;5}, both of them are necessary to treat 
nonlinear equations.  Finally we devote to the proof of  
Theorem \ref{thm;main} in Section \ref{Sec;6}.
Some supplementary estimates are described in the Appendix.

\vskip1mm
Throughout this paper we use the following notations. 
For $x\in \re^n$, $\<x\>\equiv (1+|x|^2)^{1/2}$.  
The transpose of a matrix $A$ is denoted by $A^{\sf T}$.
The Fourier and the inverse Fourier transforms are defined 
with $c_n=(2\pi)^{-n/2}$ by
$$
 \widehat{f}(\xi)=\F[f](\xi)
  \equiv c_n\int_{\re^n} e^{-ix\cdot \xi} f(x)dx,
 \quad
 \F^{-1}[f](x)\equiv c_n\int_{\re^n} e^{ix\cdot \xi} f(\xi)d\xi.
$$
For any functions 
$f=f(t,x',x_n)$ and $g=g(t,x',x_n)$, $f\underset{(t)}{*}g $, 
$f\underset{(t,x')}{*}g $ and $f\underset{(x_n)}{*}g$ stand for 
the convolution between $f$ and $g$ with respect to the variable 
indicated under $*$, respectively. If both $f$ and $g$ are 
vector field functions, 
$f\underset{(t,x')}{\cdot *}g$ denotes the convolution in $x'$ as well as 
the inner-product of $f$ and $g$, i.e.,
\eq{\label{eqn;convolution-innerprod}
 f\underset{(t,x')}{\cdot *}g
 =\sum_{\ell=1}^{n-1}
  \int_{\re}\int_{\re^{n-1}}
     f_\ell(t-s, x'-y')g_\ell(s,y')dy'ds.
}
In the summation 
$\sum_{k\in \Z}$, the parameter $k$ runs for all integers $k\in\Z$ 
and for $\sum_{k\le j}$, $k$ runs for all integers less than or equal to  
$j\in \Z$.   
We denote {\ppl $\mathcal{D}'(\re^n_+)$ the distribution over $\re^n_+$ and} the norm of the Lebesgue space $L^p(\re^{n-1})$ with $x'\in \re^{n-1}$ variable
by $\|\cdot\|_{L^p_{x'}}$.
In the norm for the Bochner 
spaces on $\dF^s_{p,\r}\big(I;X(\re^{n-1})\big)$ we use 
$$
 \|f\|_{\dF^s_{p,\r}(I;X)}
 =\|f\|_{\dF^s_{p,\r}(I;X(\re^{n-1}))}
$$
unless it may cause any confusion.
For the Besov spaces, we abbreviate $\re^n$ for $\dB^s_{p,\sg}
=\dB^s_{p,\sg}(\re^n)$ and its norm $\|\cdot\|_{\dB^s_{p,\sg}}$.
For $a\in \re^n$, we denote $B_R(a)$ as the open ball centered at $a$ 
with its radius $R>0$.  We also denote the compliment of $B_R(0)$ 
by $B_R^c$. $\Gm(\cdot)$ denotes the Gamma function.
Various constants are simply denoted by $C$ unless otherwise stated.

\sect{Maximal $L^1$-regularity for the Stokes equation 
in the half-space}\label{Sec;2} 
\subsection{Maximal $L^1$-regularity for the Stokes flow} 
Maximal $L^1$-regularity in the half-space is considered in 
\cite{OgSs20-1},  \cite{OgSs20-2} (see also Danchin--Mucha 
\cite{DM09} for $0$-boundary data).  Here we develop maximal
$L^1$-regularity for the Stokes system corresponding \eqref{eqn;FNS} and 
\eqref{eqn;NS} with inhomogeneous free stress boundary condition:
\begin{equation} \label{eqn;ST} 
   \left\{
   \begin{aligned}
     &\pt_t u -  \Del u+\N p =f,
        \quad&\ &t>0,\ \ x\in \re^n_+, \\
     &\qquad \div u=g,
        \quad&\  &t>0,\ \ x\in \re^n_+,\\
     &\quad \big( \N u+(\N u)^{\sf T} -pI\big)\, \nu_n =h,
        \quad&\  &t>0,\ \ x\in \pt\re^n_+,\\
     &\qquad \quad
       u(0,x) =u_0(x),
        \quad& &\phantom{t>0,}\ \  x\in \re^n_+,
     \end{aligned}
     \right.  \hfill
\end{equation}
where  $u_0$, $f$, $g$ and $h$ are given initial, external and 
boundary data, respectively and $\nu_n=(0,0\cdots, 0, -1)^{\sf T}$ 
denotes the outer normal on $\pt\re^n_+$. 
The following theorem is the main result of this section. 
\begin{thm}[Maximal $L^1$-regularity]
\label{thm;L1MR}
Let  $1< p< \infty$ and $-1+1/p<s\le 0$.  
The problem \eqref{eqn;ST} admits a unique solution $(u,p)$ with 
\eqn{
 \spl{
  &u\in C_b(\overline{\re_+};\dB^s_{p,1}(\re^n_+))\cap
        \dot{W}^{1,1}(\re_+;\dB^s_{p,1}(\re^n_+)), \\
  & \Del u,\; \N p \in L^{1}(\re_+;\dB^{s}_{p,1}(\re^n_+)) , \\ 
  &p\big|_{x_n=0}
    \in 
    \dF^{\frac12-\frac{1}{2p}}_{1,1}(\re_+;\dB^{s}_{p,1}(\re^{n-1}))
     \cap 
    {\crd L^1}(\re_+;\dB^{s+1-\frac 1p}_{p,1}(\re^{n-1})) 
}
}
if and only if the data in \eqref{eqn;ST} satisfy 
\begin{align*}
  & u_0\in  \dB^{s}_{p,1}(\re^n_+), \, 
  {\ppl 
     \div u_0=g\big|_{t=0} \text{in } \mathcal{D}'(\re^n_+),
  }
\quad
   f\in L^{1}(\re_+;\dB^{s}_{p,1}(\re^n_+)), 
      \\
  &  \N g\in  L^{1}(\re_+;\dB^{s}_{p,1}(\re^n_+)), \quad
     \N(-\Del)^{-1}g\in \dot{W}^{1,1}(\re_+;\dB^{s}_{p,1}(\re^n_+)),
      \\
  & h\in \dF^{\frac12-\frac{1}{2p}}_{1,1}(\re_+;\dB^s_{p,1}(\re^{n-1}))
         \cap 
         {\crd L^1}(\re_+;\dB^{s+1-\frac 1p}_{p,1}(\re^{n-1})),
\end{align*}
{\gr where $(-\Del)^{-1}g$ is given by 
$G*\tilde g|_{\re^n_+}$ with $\tilde{g}$ as the 
{\ppl even} extension of $g$ {\rm (}see \eqref{eqn;Green}  blow{\rm)}.}
Besides the solution $(u,p)$ satisfies the following estimate 
for some constant $C_M>0$ depending only on $p$, $s$ and $n$ 
\begin{align}
    \big\|\pt_t u & \big\|_{L^{1}(\re_+;\dot{B}^s_{p,1}(\re^n_+))}
   +\big\|D^2  u \big\|_{L^{1}(\re_+;\dot{B}^s_{p,1}(\re^n_+))}   
   +\big\|\nabla p \big\|_{L^{1}(\re_+;\dot{B}^s_{p,1}(\re^n_+))} \nn\\
   &+\big\|p |_{x_n=0}
    \big\|_{\dF^{\frac12-\frac{1}{2p}}_{1,1}(\re_+;\dB^{s}_{p,1}(\re^{n-1}))}
   +\big\|p |_{x_n=0}
    \big\|_{{\crd L^1}(\re_+;\dB^{s+1-\frac1{p}}_{p,1}(\re^{n-1}))}
    \nn\\
 \le &
  C_M\Big(\|u_0\|_{\dot B^{s}_{p,1}(\re^n_+)}
         +\|f\|_{L^{1}(\re_+;\dot B^{s}_{p,1}(\re^n_+))} \nn\\
 & \hskip1.5cm
     +\|\N g\|_{L^{1}(\re_+;\dot B^{s}_{p,1}(\re^n_+))}
     +\|\pt_t {\gr \N} (-\Del)^{-1}g\|_{L^{1}(\re_+;\dot B^{s}_{p,1}(\re^n_+))}\nn\\
 & \hskip3cm
     +\| h\|_{\dot{F}^{\frac12-\frac{1}{2p}}_{1,1}(\re_+;\dB^s_{p,1}(\re^{n-1}))}
     +\| h\|_{{\crd L^1}(\re_+;\dB^{s+1-\frac{1}{p}}_{p,1}(\re^{n-1}))}
     \Big).
  \label{eqn;L1MR-estimate}
\end{align}
\end{thm}
\vskip2mm
The general theory of maximal regularity for the parabolic type
partial differential equation is extensively developed in 
the UMD Banach space (see, for instance, 
\cite{DHP03}, \cite{DHP07}, \cite{GS11}, \cite{GS},
\cite{HP97}, \cite{LSU}, \cite{OgSs10}, \cite{PrSi16},
\cite{SbSz05}, \cite{Wd05}, \cite{Ws},  \cite{Z-Z}). 
However the end-point exponent is normally excluded 
in the general theory. If the space is restricted 
the homogeneous Besov space or Fourier transformed measures, 
one can see the end-point estimate holds
as is seen in \cite{BCD}, \cite{Dc03}, \cite{DM12}, \cite{GS11}, 
\cite{OgSs16}, \cite{OgSs20-1}, \cite{OgSs20-2}.

\vskip1mm
To establish maximal regularity of the half-space problem \eqref{eqn;ST}, 
we reduce the problem \eqref{eqn;ST} into the several partial 
components of the data and reduce the problem into the 
inhomogeneous problem with only boundary data.  
At first we remove the divergence data. 
Introducing the {\gr even} extension of divergence data $g$ 
with respect to $x_n$;
$$
 \widetilde{g}(t,x)
 =\begin{cases}
    \, g(t,x', x_n), & x_n>0, \\
    \, {\ppl g}(t,x', -x_n), & x_n<0 
  \end{cases}
$$
for $x'=(x_1,x_2,\cdots, x_{n-1})$,
we consider the problem 
\beq\label{eqn;g}
\left\{
\begin{aligned}
    &-\Del \phi =\widetilde {g}, \qquad &t>0, \ \  &x\in \re^n,  \\
    & \quad  \phi \Big|_{x_n=0}= 0,  \qquad  &t>0, \ \  &x'\in \re^{n-1}.
\end{aligned}     
\right.
\eeq
One of the solution of \eqref{eqn;g} is given by the Newtonian potential 
$\phi=(-\Del)^{-1}\widetilde{g}\equiv G* \widetilde{g}$ with the Newtonian kernel 
$G$ in $\re^n$;
\begin{equation} \label{eqn;Green}
 G(x)=\begin{cases}
         \frac{1}{2\pi}\log |x|^{-1}, & n=2,\\
         ((n-2)\om_n)^{-1} |x|^{-(n-2)}, & n\ge 3, 
        \end{cases}
    \qquad 
    \omega_n=\frac{2\pi^{\frac n2}}{\Gamma\big(\frac{n}{2}\big)}.
\end{equation}
Then the gradient of potential $\N \phi$ satisfies the estimate for 
 $1< p< \infty$ and $-1+1/p <s<1/p$ 
\begin{equation} \label{eqn;potential}
\left\{ 
\begin{split}
   &\|\nabla^3 \phi \|_{L^{1}(\re_+;\dot{B}^s_{p,1}(\re^n_+))}
     \le C\|\N g\|_{L^{1}(\re_+;\dot B^{s}_{p,1}(\re^n_+))}, \\
   &\|\pt_t \N \phi \|_{L^{1}(\re_+;\dot{B}^s_{p,1}(\re^n_+))}
     \le C \|\pt_t{\gr \N}(-\Del)^{-1}g\|_{L^{1}(\re_+;\dot B^{s}_{p,1}(\re^n_+))}. 
\end{split}
\right.
\end{equation}
Indeed the corresponding estimate to \eqref{eqn;potential} in $\re^n$ 
follows directly from the elliptic estimate for the Poisson equation 
(or the Bernstein type estimate) and hence the estimate 
\eqref{eqn;potential} in the half-space naturally follows from 
the definition of the Besov space in $\re^n_+$.
Setting $w=u+\N\phi|_{\re^n_+}$, the pair of functions
 $(w,p)$ satisfy the equations 
\begin{equation} \label{eqn;Stokes-1''}
   \left\{
   \begin{split}
    &\pt_t  w-\Del w +\N p=f+\big(\pt_t \N \phi -\Del \N \phi\big)\big|_{x_n>0}, 
      &\quad &t>0, \ \ x\in \re^n_+,\\
    &\qquad  \qquad \div\, w=0,   
      &\quad  &t>0, \ \ x\in \re^n_+,\\
    &\quad \big(\N w+(\N w)^{\sf T}-pI\big)\, \nu_n 
       =h+\big( \N^2\phi + (\N^2 \phi)^{\sf T} \big)\, \nu_n,
             &\quad &t>0, \ \ x\in \pt\re^n_+,\\
    &\qquad  w(0,x)=u_0(x)+\N \phi(0,x)\big|_{x_n>0}, 
       &\quad &\phantom{t>0,\ }  \ x\in \re^n_+,
 \end{split}
 \right.
\end{equation} 
where $\nu_n$ denotes the outer normal to $\pt\re^n_+$.

In order to exclude the external and initial data, 
we extend them into $\re^n$, more precisely,  
 we extend $f_j$ ($1\le j\le n-1$) by odd functions and  
for the $n$-th component $f_n$,  we employ the even extension 
(we write them $f^o_j$ and $f^e_n$, respectively)
and set  $\overline{f}=(f^o_1,\cdots f^o_{n-1}, f^e_n)^{\sf T}$. 
For the initial data $u_0$,  
we also employ the same extension with respect to $x_n$ and write it 
$\overline{u_0}$ and set
\eq{\label{eqn;data-tilde}
\left\{
 \spl{
    \widetilde{f} =\ &  \overline{f}+(\pt_t \N\phi -\Del \N \phi), 
    & \quad t>0, x\in \re^n, \\
    \widetilde{u_0}(x)  =\ & \overline{u_0}(x)+\N \phi (0,x),
    & \quad x\in \re^n,
 }
\right.
}
and we consider the Cauchy problem:
\eq{ \label{eqn;Stokes-2'} 
 \left\{
  \begin{aligned}
    \pt_t \widetilde{u}-\Del \widetilde{u} +\N \widetilde{p}
      &=\widetilde{f} \quad &t>0,\ \  x\in \re^n,\\
    \qquad \div\, \widetilde{u}&=0,    \quad &t>0,\ \  x\in \re^n,\\
    \qquad \widetilde{u}(0,x)&=\widetilde{u_0}(x), \quad & x\in \re^n.   
  \end{aligned}
 \right.
}
Then it is known that 
the solution  $(\widetilde{u},\widetilde{p})$ of the equation 
\eqref{eqn;Stokes-2'} satisfies 
maximal $L^1$-regularity 
\begin{align*}
  \|\pt_t \widetilde{u} \|_{L^{1}(\re_+;\dot{B}^s_{p,1}(\re^n))}
     +& \|\nabla^2  \widetilde{u} \|_{L^{1}(\re_+;\dot{B}^s_{p,1}(\re^n))} 
     +  \|\N \widetilde{p} \|_{L^{1}(\re_+;\dot{B}^s_{p,1}(\re^n))} \\
 \le& C_M\Big(
          \|\widetilde{u_0}\|_{\dot B^{s}_{p,1}(\re^n)}
         +\|\widetilde{f}\|_{L^{1}(\re_+;\dot B^{s}_{p,1}(\re^n))} \Big)
  \eqntag \label{eqn;UP}         
\end{align*}
for any {\gr $-1/p'<s<1/p$} and {\gr $1< p < \infty$} (see Danchin--Mucha  \cite{DM09} 
and  Ogawa--Shimizu \cite{OgSs16},  see also \cite{Z-Z}).  
By restricting the solution $(\widetilde{u}, \widetilde{p})$
over the half-space $\re^n_+$ (and we denote them in the 
same notation) we directly obtain 
from \eqref{eqn;potential} and \eqref{eqn;UP} that
\begin{align*}
  \|\pt_t \widetilde{u} \|_{L^{1}(\re_+;\dot{B}^s_{p,1}(\re^n_+))}
     +& \|\nabla^2  \widetilde{u} \|_{L^{1}(\re_+;\dot{B}^s_{p,1}(\re^n_+))} 
     +  \|\N \widetilde{p} \|_{L^{1}(\re_+;\dot{B}^s_{p,1}(\re^n_+))} \\
 \le& C_M\Big(
          \|\widetilde{u_0}\|_{\dot B^{s}_{p,1}(\re^n_+)}
         +\|\widetilde{f}\|_{L^{1}(\re_+;\dot B^{s}_{p,1}(\re^n_+))} \Big)\\
  \le& C_M\Big(
          \|u_0\|_{\dot B^{s}_{p,1}(\re^n_+)}
         +\|f\|_{L^{1}(\re_+;\dot B^{s}_{p,1}(\re^n_+))}  \\
     &\qquad
         +\|\N g\|_{L^{1}(\re_+;\dot B^{s}_{p,1}(\re^n_+))}
         +\|\pt_t{\gr \N}(-\Del)^{-1} g
          \|_{L^{1}(\re_+;\dot B^{s}_{p,1}(\re^n_+))}\Big),
  \eqntag \label{eqn;up-tilde}         
\end{align*}
where the inverse operator $(-\Del)^{-1}$ is given by the solution 
operator to the elliptic problem \eqref{eqn;g} and it is realized by the 
Green's function \eqref{eqn;Green}.

Finally we consider the difference between the solutions $(w,p)$ to 
\eqref{eqn;Stokes-1''} and  $(\widetilde{u},\widetilde{p})$  to 
\eqref{eqn;Stokes-2'} restricted in $\re^n_+$. Letting 
$v= w-\widetilde{u}|_{x_n>0}\equiv u+\N \phi|_{x_n>0}-\widetilde{u}|_{x_n>0}$ 
and $q=p-\widetilde{p}|_{x_n>0}$ 
and we reduce the original  problem into 
the following initial boundary value problem for $(v,q)$: 
\begin{equation} \label{eqn;ST2}
   \left\{
   \begin{aligned}
     &\pt_t v  -\Del v +\N q =0,
        \quad&\  &t>0,\ \ x\in  \re^n_+,\\\
     &\qquad \div v=0,
        \quad&\  &t>0, \ \ x\in\re^n_+,\\\
     &\quad 
         \big(\N v+(\N v)^{\sf T}- q I\big)\, \nu_n =H,
        \quad&\  &t>0, \ \ x\in \pt\re^n_+,\\\
     &\qquad v(0,x) =0,
        \quad&\ &\phantom{\ t>0,} \ \ x\in \re^n_+,
     \end{aligned}
     \right.  \hfill
\end{equation}
where we set
\eq{\label{eqn;boundary-data}
 \spl{
  H\equiv &
    \widetilde{h}
     -\Big(\N \widetilde{u}+(\N \widetilde{u})^{\sf T}
           -\widetilde{p}I\Big) \nu_n \\
   = & h-\big( \N^2\phi + (\N^2 \phi )^{\sf T} \big) \nu_n 
      -\Big(\N \widetilde{u}+(\N \widetilde{u})^{\sf T}-\widetilde{p}I
       \Big) \nu_n .
}
}
In order to prove Theorem \ref{thm;L1MR}, it is essential to show 
maximal $L^1$-regularity for \eqref{eqn;ST2}. 

\begin{thm}\label{thm;L1MR2-b}
Let $-1+1/p<s \le  0$ and $1< p< \infty$.  
The problem \eqref{eqn;ST2} admits a unique solution 
\eq{  \label{eqn;Stokes-L1Max-regularity}
\spl{
   &v  \in C_b(\overline{\re_+};\dB^s_{p,1}(\re^n_+))\cap 
           \dot{W}^{1,1}(\re_+;\dB^s_{p,1}(\re^n_+)),\\ 
   &\Del v, \;\N q \in  L^{1}(\re_+;\dB^{s}_{p,1}(\re^n_+)),\\
   &q|_{x_n=0}\in 
   \dF^{\frac12-\frac{1}{2p}}_{1,1}(\re_+;\dB^s_{p,1},(\re^{n-1}))
    \cap 
   {\crd L^1}(\re_+;\dB^{s+1-\frac 1p}_{p,1}(\re^{n-1}))
}
}
if and only if the data in \eqref{eqn;ST2} satisfy
\eq{\label{eqn;boundary-assump}
 H\in \dF^{\frac12-\frac{1}{2p}}_{1,1}(\re_+;\dB^s_{p,1}(\re^{n-1}))
       \cap 
     {\crd L^1}(\re_+;\dB^{s+1-\frac 1p}_{p,1}(\re^{n-1})). 
}
Besides the solution $(v,q)$ satisfies the following estimate 
for some constant $C_M>0$ depending only on $p$, $s$ and $n$ 
\eq{\label{eqn;L1MR2-b-estimate}
 \spl{
 \big\|\pt_t v&\big\|_{L^1(\re_+;\dB^s_{p,1}(\re^n_+))}
    +\big\|D^2 v\big\|_{L^1(\re_+;\dB^s_{p,1}(\re^n_+))}  
    +\big\|\N q\big\|_{L^1(\re_+;\dB^s_{p,1}(\re^n_+))}  \\
  & +\big\|q|_{x_n=0}
     \big\|_{\dF^{\frac12-\frac{1}{2p}}_{1,1}(\re_+;\dB^s_{p,1}(\re^{n-1}))}
    +\big\|q|_{x_n=0}
     \big\|_{{\crd L^1}(\re_+;\dB^{s+1-\frac1{p}}_{p,1}(\re^{n-1}))}
 \\
 \le& C\Big(\|H\|_{\dF^{\frac12-\frac{1}{2p}}_{1,1}(\re_+;\dB^s_{p,1}(\re^{n-1}))} 
           +\|H\|_{{\crd L^1}(\re_+;\dB^{s+1-\frac{1}{p}}_{p,1}(\re^{n-1}))} 
       \Big).
  }
}
\end{thm}

\vskip2mm
Once we obtain corresponding maximal regularity 
Theorem \ref{thm;L1MR2-b} to the solution $(v,q)$ for 
\eqref{eqn;ST2}, we may prove maximal 
$L^1$-regularity for the original Stokes system \eqref{eqn;ST} 
combining with those estimates \eqref{eqn;potential} and \eqref{eqn;up-tilde}
and the relation
\eq{ \label{eqn;up-decomposion}
\left\{
 \spl{
  u(t,x)=&\widetilde{u}(t,x)+v(t,x)+\N \phi(t,x), & t>0,\ \  x\in \re^n_+,\\   
  p(t,x)=&\widetilde{p}(t,x)+q(t,x), & t>0,\ \  x\in \re^n_+,
 }
\right.
}
as well as the following trace estimate (see Appendix below, 
cf. \cite{OgSs20-1}, \cite{OgSs20-2}):
\eq{\label{eqn;up-tilde-trace}
\spl{
   \Big\|\Big(\N  \widetilde{u}+&(\N  \widetilde{u})^{\sf T}\Big)\,\nu_n 
   \Big\|_{\dF^{\frac12-\frac{1}{2p}}_{1,1}(I;\dB^s_{p,1}(\re^{n-1}))} 
  +\Big\|\Big(\N \widetilde{u}+(\N \widetilde{u})^{\sf T}\Big)\, \nu_n 
   \Big\|_{{\crd L^1}(I;\dB^{s+1-\frac{1}{p}}_{p,1}(\re^{n-1}))} \\
  \le &C\Big(\|\pt_t \widetilde{u}\|_{L^1(I;\dB^s_{p,1}(\re^{n}_+))}
            +\|\Del \widetilde{u}\|_{L^1(I;\dB^s_{p,1}(\re^{n}_+))}
       \Big). 
}
}

\vskip1mm
\subsection{Solution formula of the Stokes equation}
We construct the solution formula of \eqref{eqn;ST} according 
to the method by  Shibata--Shimizu \cite{SbSz03} and \cite{SbSz08}.  

Let $H=H(t,x')\equiv (H'(t,x'),H_n(t,x'))$ be the boundary 
data extended into $t<0$ by the  zero extension.
Let $\F$ and $\F^{-1}$ denote the Fourier and the inverse Fourier 
transform with respect to $x'\in \re^{n-1}$,  and $\L$ and $\L^{-1}$ 
denote the Laplace and the inverse Laplace transform
for $t$ and $\t$, respectively.  Namely
\begin{align*}
  &\L \widehat f(\lambda, \xi', x_n)
   =(2\pi)^{-\frac n2}\int_{\re_+}\int_{\re^{n-1}}
        e^{-\lambda t-i x'\cdot\xi'} f(t, x',x_n)\,dx'dt,\\
  &\L^{-1}\F^{-1} f(t, x', x_n)
   =(2\pi)^{-\frac {n-1}2}
    \frac{1}{2\pi i}\int_{\Gm}\int_{\re^{n-1}}
     e^{\lambda t+i x'\cdot\xi'} f({\ppl \lam}, \xi',x_n)\,d\xi'd\lam,
\end{align*}
where $\Gm$ denotes an integral path given for 
{\ppl some $\gm>0$} by 
$\Gm=\{\lam=\gm+i\t;\, \t\in \re\}$.
Applying the Fourier-Laplace transform with respect to 
$(x',t)$ to \eqref{eqn;ST2}, we have the solution formula for the 
$n$-th component of the velocity and the pressure as follows:
\algn{
 \wh{v_n}(\tau,\xi',x_n)
  =& \frac{|\xi'|}{(B(\t,\xi')-|\xi'|)D(\t,\xi')}
     \Big(2i|\xi'|B(\t,\xi') \big(\xi'\cdot \widehat{H'}\big)
           -(|\xi'|^2+B(\t,\xi')^2)\widehat{H_n}
      \Big)e^{-|\xi'|x_n}  \\
   &+\frac{|\xi'|}{(B(\t,\xi')-|\xi'|)D(\t,\xi')}
     \Big(
           -(|\xi'|^2+B(\t,\xi')^2) i\xi'\cdot \widehat{H'}
           +2|\xi'|B(\t,\xi')  \widehat{H_n}
     \Big)e^{-B(\t,\xi')x_n},
   \eqntag \label{eqn;vn} \\
 \widehat{q}(\tau,\xi',x_n) 
  =& \frac{|\xi'|+B(\t,\xi')}{D(\t,\xi')} 
    \Big( 2B(\t,\xi') \big(i\xi'\cdot \wh{H'}\big) 
         -(|\xi'|^2+B(\tau,\xi')^2)\wh{H_n}
    \Big)e^{-|\xi'|x_n}, 
   \eqntag \label{eqn;q}
}
where we have set 
$\widehat{H}=(\widehat{H'}(\t,\xi'),\widehat{H_n}(\t,\xi'))$ 
as the Fourier-Laplace transform of the given boundary data and 
\begin{align}
  &B(\t,\xi')=\sqrt{i\tau+|\xi'|^2},\quad \text{Re}\, B(\t,\xi')\ge 0,
   \label{eqn;B} \\ 
  &D(\tau,\xi')
   =B(\t,\xi')^3+|\xi'|B(\t,\xi')^2+3|\xi'|^2B(\t,\xi')-|\xi'|^3.
    \label{eqn;D}
\end{align}
Hence we see for any smooth rapidly decreasing boundary data 
$(\widehat{H'},\widehat{H_n})$ in both $(\t,\xi')$ variables, 
we see by passing $\gm\to 0$ to obtain 
\alg{
v_n&(t,x',x_n)\nn\\
  =&c_{n+1}\text{p.v.}\iint_{\re^{n}}
               e^{it\t+ix'\cdot \xi'}
               \bigg\{
                \frac{|\xi'|}{(B-|\xi'|)D(\t, \xi')}
                 \Big(
                  2B(i\xi'\cdot \widehat{H}')
                  -(|\xi'|^2+B^2) \widehat{H_n}
                 \Big)e^{-|\xi'|x_n }  \notag\\
    &\hskip3.5cm
      +\frac{|\xi'|}{(B-|\xi'|)D(\t, \xi')}
                 \Big(
                  -(|\xi'|^2+B^2)i\xi' \cdot \widehat{H}')
                  +2|\xi'| B \widehat{H_n}
                 \Big)  e^{-Bx_n }            
               \bigg\}              
               d\t d\xi',
         \label{eqn;Stokes-n-velocity}\\
q&(t,x',x_n) \nn\\
  =&c_{n+1}\text{p.v.}\iint_{\re^{n}}
               e^{it\t+ix'\cdot \xi'}
               \bigg\{
                \frac{|\xi'|+B}{D(\t, \xi')}
                 \Big(
                  2B(i\xi'\cdot \widehat{H}')
                  -(|\xi'|^2+B^2) \widehat{H_n}
                 \Big)
               \bigg\}
               e^{-|\xi'|x_n }
               d\t d\xi',  
        \label{eqn;Stokes-n-pressure}
}
where we take a limit of the integral pass avoiding the 
singularity at $(\t,\xi')=(0,0)$.
All the other components of the velocity fields $v_{\ell}(t,x)$ 
($\ell=1,2,\cdots, n-1$) are given by the above 
two components $(v_n, q)$ and the boundary data $H=(H',H_n)$ 
from the equation \eqref{eqn;ST2} 
(see for the details \cite{SbSz08}).

Our main task is to prove maximal $L^1$-regularity of 
the pressure term  $q$ in \eqref{eqn;ST2} which is 
directly obtained from the inhomogeneous boundary data. 
Then the maximal $L^1$-regularity estimate for the velocity term $v$ 
of \eqref{eqn;ST2} follows from the estimate for  $q$.  
Applying the gradient to the solution formula \eqref{eqn;Stokes-n-pressure}, 
we obtain the explicit expression of  $\N q$ as
\begin{align}
  \N & q(t,x',x_n)\nn\\
  =&c_{n+1}\iint_{\re^{n}}
               e^{it\t+ix'\cdot \xi'}
               (i\xi',-|\xi'|)^{\sf T}
               \bigg\{
                \frac{|\xi'|+B}{D(\t, \xi')}
                 \Big(
                  2B(i\xi'\cdot \widehat{H}')
                  -(|\xi'|^2+B^2) \widehat{H_n}
                 \Big)
               \bigg\}
               e^{-|\xi'|x_n }
               d\t d\xi', \label{eqn;Stokes-n-pressure-grad}
\end{align}
where  $B=B(\t,\xi')$ and $D(\tau, \xi')$ are defined by 
\eqref{eqn;B} and \eqref{eqn;D}, respectively. 
We also set the following Fourier multiplier 
$m(\t,\xi'):\re\times \re^{n-1}\to \re^n$ as
\eq{\label{eqn;singular-int-multi}
\spl{
 m(\t,\xi') 
 = & \big(m'(\t,\xi'), m_n(\t,\xi') \big)\\
  \equiv & \Big(\frac{2B(|\xi'|+B)}{D(\t,\xi')}
                   i\xi',\ 
           -\frac{\big(|\xi'|+B\big)(|\xi'|^2+B^2)}{D(\t,\xi')}
           \Big).
}
} 
%
%
\subsection{The homogeneous Besov spaces on the half-space}\par\label{subsec;2.3}
First we recall the summary for the homogeneous Besov spaces over 
the half Euclidean space $\re^n_+$. We recall the retraction 
and the coretraction defined in the way of Triebel \cite{Tr78}
as follows:

\noindent
{\it Definition} (\cite{Tr78}). Let $A$ and $B$ be Banach spaces and 
let  $R$ and $E$ be linear operators as 
\eq{ 
\spl{
   & R:A\to B \text{ bounded,  } \\
   & E:B\to A \text{ bounded,  }  \\
   & RE=Id:B\to B  \text{ bounded, } 
   }
 }
where $Id$ is the identity operator from $B$ to $B$.
Then {\ppl $R$} is called as {\it retraction} and $E$ is called as {\it coretraction}.
\vskip2mm
 
\vskip2mm
\noindent
{\it Definition.}  Let $1\le p<\infty$ and $1\le \sg <\infty$ with {\ppl $s\in \re$}.
Let 
\begin{align}
   &\overset{\circ\quad}{B^s_{p,\sg}}(\re^n_+)
      \equiv \overline{C_0^{\infty}(\re^n_+)}^{\dB^s_{p,\sg}(\re^n_+)}, \\
   &\overset{\odot\quad}{B^s_{p,\sg}}(\re^n_+)
      \equiv \overline{
        \{f\in \dB^s_{p,\sg}(\re^n); 
            \supp f \subset \re^n_+\}}^{\dB^s_{p,\sg}(\re^n)}.
\end{align}
It is shown that the above defined space coincides {\ppl with}
the space $\dB^s_{p,\sg}(\re^n_+)$ defined by the restriction 
in \eqref{eqn;Besov-halfspace}.
First we observe that the duality is well-defined in certain 
range of the exponent.
\begin{prop}[cf. \cite{Tr78}]
\label{prop;duality} Let $1\le p<\infty$ and  $1\le \sg <\infty$.
For $0\le s<1/p$, it holds 
\eqn{
  \big(\overset{\circ\quad}{B^s_{p,\sg}}(\re^n_+)\big)'
  \simeq \dB^{-s}_{p',\sg'}(\re^n_+).
  }
where $\simeq$ stands for the both spaces being equivalent 
as the normed space.
\end{prop}

The part of the following proposition is shown by Danchin--Mucha \cite{DM09}.
\begin{prop}[\cite{DM09}]
\label{prop;DM-09} Let $1\, {\ppl \le}\, p<\infty$ and  $1\le \sg <\infty$.
For $-1+1/p<s<1/p$,
\algn{
  \overset{\circ\quad}{B^s_{p,\sg}}(\re^n_+)
  \simeq \dB^s_{p,\sg}(\re^n_+), \\
  \overset{\odot\quad}{B^s_{p,\sg}}(\re^n_+)
  \simeq \dB^s_{p,\sg}(\re^n_+),
  }
where $\simeq$ stands for the both spaces being equivalent 
as the normed space.
\end{prop}

\vskip2mm
We consider the restriction operator $R_0$ 
by 
\eq{\label{eqn;restriction}
 R_0f(x)=f(x)\Big|_{x\in \re^n_+}
}
for all $f\in \dB^s_{p,\sg}(\re^n)$ with $s> 0$ and 
it is understood in the sense of distribution for $s\le 0$.

Let $\chi_+$ be a cut-off operation defined 
by multiplying a cut-off function 
\eqn{
 \chi_{\re^n_+}(x)=
  \begin{cases}
   1,    & \text{ in } \re^n_+, \\
   0,    & \text{ in } \overline{\re^n_-}.
  \end{cases}
}
Let the extension operator $E_0$ from 
$\overset{\odot\quad}{B^s_{p,\sg}}(\re^n_+)$ 
given by the zero-extension, i.e.,
for any $f\in \overset{\odot\quad}{B^s_{p,\sg}}(\re^n_+)$, set 
\eq{\label{eqn;extension}
 E_0f=
 \begin{cases}
   f(x), & \text{ in } \re^n_+, \\
   0,    & \text{ in } \overline{\re^n_-}.
 \end{cases}
}

One can find that those operators are basic tool to 
recognize the homogeneous Besov spaces. Using 
Proposition \ref{prop;DM-09},
the following statement is a variant introduced by 
Triebel \cite[p.228]{Tr78} .
 
\begin{prop}\label{prop;retraction-coretraction}
Let $1\le p<\infty$, $1\le \sg<\infty$ and $-1+1/p<s<1/p$,
and let $R_0$ and $E_0$ be operators defined in \eqref{eqn;restriction}
and \eqref{eqn;extension}
It holds that  
\begin{align}
    &R_0:\dB^s_{p,\sg}(\re^n)
        \to\dB^s_{p,\sg}(\re^n_+),  \label{eqn;retraction}\\ 
    &E_0:\dB^s_{p,\sg}(\re^n_+)
        \to \dB^s_{p,\sg}(\re^n), \label{eqn;co-retraction}
\end{align}
are linear bounded operators.  Besides
it holds that
\eq{
    R_0E_0=Id: \dB^s_{p,\sg}(\re^n_+)\to \dB^s_{p,\sg}(\re^n_+),
      }
where $Id$ denotes the identity operator.  
Namely $R_0$ and $E_0$ are a retraction and a coretraction, respectively.
\end{prop}
\vskip2mm
The proof of Proposition \ref{prop;retraction-coretraction} 
is along the same line of the proof in \cite{Tr78}.  Note that the 
spaces are homogeneous Besov spaces and then the arrangement 
appears in Proposition 3 in Danchin--Mucha \cite{DM09} is required.

\vskip2mm
\begin{prf}{Proposition \ref{prop;retraction-coretraction}}
To see the first operator \eqref{eqn;retraction} is bounded,
let $f\in \dB^s_{p,\sg}(\re^n)$ 
and  we show that 
$$
 \|R_0 f\|_{\dB^s_{p,\sg}(\re^n_+)}
 =    \inf \|\widetilde{R_0f}\|_{\dB^s_{p,\sg}(\re^n)}
 \le \| \chi_{\re^n_+}f\|_{\dB^s_{p,\sg}(\re^n)}
 \le \| f\|_{\dB^s_{p,\sg}(\re^n)}
$$
is valid under the restriction $s>0$.
Let $-1/p'<s<0$ and $f\in \dB^{s}_{p,\sg}(\re^n)$.
For any test $\phi \in C_0^{\infty}(\re^n)$,  
\eq{
 \spl{
  \big|\<\chi_{\re^n_+}f,\phi\>\big|
  =  & \big|\<f, \chi_{\re^n_+}\phi\>\big|
  \le \|f\|_{\dB^{s}_{p,\sg}(\re^n)} 
      \|\chi_{\re^n_+}\phi\|_{\dB^{-s}_{p',\sg'}(\re^n)} \\
 \le &\|f\|_{\dB^{s}_{p,\sg}(\re^n)} 
      \|\phi\|_{\dB^{-s}_{p',\sg'}(\re^n)}, 
}}
since $0<-s<1/p'$ and the last inequality follows from 
the pointwise sense.
Thus from the definition of the norm in $\dB^s_{p,\sg}(\re^n_+)$,
it holds similarly to the above that 
\eqn{
    \|R_0 f \|_{\dB^{s}_{p,\sg}(\re^n_+)}
   \le \| \chi_{\re^n_+} f \|_{\dB^{s}_{p,\sg}(\re^n)} 
   = \sup_{\phi\in\dB^{-s}_{p',\sg'}(\re^n)\setminus\{0\} }
    \frac{  \big|\<\chi_{\re^n_+} f,\phi\>\big|  }
         {  \|\phi\|_{\dB^{-s}_{p',\sg'}(\re^n)} }
   \le \|f\|_{\dB^{s}_{p,\sg}(\re^n)}. 
  }
{\ppl
For the second bound \eqref{eqn;co-retraction}, 
see \cite[Proposition 3]{DM09}.
Since the both operators are bounded, 
we see that }
$$
R_0E_0=Id:\dB^s_{p,\sg}(\re^n_+)
        \to \dB^s_{p,\sg}(\re^n_+),
$$
holds by the density argument in Proposition \ref{prop;DM-09}.
\end{prf}
\vskip2mm
\begin{prop}[cf. \cite{DHMT20}, \cite{Tr78}]\label{prop;derivative}
Let $1<p<\infty$ and $-1+1/p<s\le n/p-1$. Then for 
any $f\in \dB^s_{p,1}(\re^n_+)$,
$$
  \|\N f\|_{\dB^s_{p,1}(\re^n_+)}
  \simeq \|f\|_{\dB^{s+1}_{p,1}(\re^n_+)},
$$
where $\simeq$ stands for that the both side of the norm is 
equivalent.
\end{prop}
\vskip2mm\noindent
{\ppl
For the proof, see \cite[Proposition 3.19, Corollary 3.20]{DHMT20}.
}
\vskip2mm

\vskip1mm\noindent
{\bf Remark.}
In what follows, we restrict ourselves to the regularity 
range of the Besov spaces $\dB^s_{p,\sg}(\re^n_+)$ 
in $-1+1/p<s<1/p$ for $1<p<\infty$ 
unless otherwise stated.
According to Proposition \ref{prop;retraction-coretraction},
we can regard that 
any distribution in $\dB^s_{p,\sg}(\re^n_+)$ under such restriction 
on $s$ and $p$ can be extended into a distribution over whole space 
$\re^n$ and 
conversely any distribution in $\dB^{s}_{p,\sg}(\re^n)$ is 
restricted into a distribution over the half-space $\re^n_+$.
We frequently use those facts without noticing for every case
below.

\subsection{The L-P decomposition with separation of variables}\par
In order to split the variables $x'\in \re^{n-1}$ and $x_n\in \re_+$, 
we introduce an $x'$-parallel decomposition and an $x_n$-parallel
decomposition by Littlewood--Paley type. In what follows $\eta\in \re_+$
denotes a parameter for $x_n$-axis in $\re^n_+$. 
We introduce $\{\overline{\Phi_m}\}_{m\in \Z}$ as a Littlewood--Paley dyadic 
frequency decomposition of unity in separated variables $(\x',\x_n)$. 
\vskip2mm
%
%
\begin{center}
  \begin{picture}(300,210)(-50,0)
    \thicklines
    \put(-10,80){\vector(1,0){220}}  
    \put(100,50){\vector(0,1){165}}  
    \thicklines
    {
     \multiput(0,80)(0,2){50}{\line(1,0){1.2} }
     \multiput(200,84)(0,2){50}{\line(1,0){1.2} }
     \multiput(55,80)(0,2){31}{\line(1,0){1.2} }
     \multiput(145,80)(0,2){31}{\line(1,0){1.2} }
     \multiput(0,180)(2,0){100}{\line(0,1){1.2} }
     \multiput(55,140)(2,0){46}{\line(0,1){1.2} }
%
    }
    {
    \multiput(145,160)(2,0){28}{\line(0,1){1.2} }
    \multiput(145,150)(2,0){28}{\line(0,1){1.2} }
    \multiput(145,140)(2,0){28}{\line(0,1){1.2} }
    \multiput(145,130)(2,0){28}{\line(0,1){1.2} }
    \multiput(145,120)(2,0){28}{\line(0,1){1.2} }
    \multiput(145,110)(2,0){28}{\line(0,1){1.2} }
    \multiput(145,100)(2,0){28}{\line(0,1){1.2} }
    \multiput(145,90)(2,0){28}{\line(0,1){1.2} }
    \multiput(145,80)(2,0){28}{\line(0,1){1.2} }
    \multiput(1,160)(2,0){28}{\line(0,1){1.2} }
    \multiput(1,150)(2,0){28}{\line(0,1){1.2} }
    \multiput(1,140)(2,0){28}{\line(0,1){1.2} }
    \multiput(1,130)(2,0){28}{\line(0,1){1.2} }
    \multiput(1,120)(2,0){28}{\line(0,1){1.2} }
    \multiput(1,110)(2,0){28}{\line(0,1){1.2} }
    \multiput(1,100)(2,0){28}{\line(0,1){1.2} }
    \multiput(1,90)(2,0){28}{\line(0,1){1.2} }
    \multiput(1,80)(2,0){28}{\line(0,1){1.2} }
    }
    {
    \multiput(0,140)(0,2){22}{\line(1,0){1.2} }
    \multiput(10,140)(0,2){22}{\line(1,0){1.2} }
    \multiput(20,140)(0,2){22}{\line(1,0){1.2} }
    \multiput(30,140)(0,2){22}{\line(1,0){1.2} }
    \multiput(40,140)(0,2){22}{\line(1,0){1.2} }
    \multiput(50,140)(0,2){22}{\line(1,0){1.2} }
    \multiput(60,140)(0,2){22}{\line(1,0){1.2} }
    \multiput(70,140)(0,2){22}{\line(1,0){1.2} }
    \multiput(80,140)(0,2){22}{\line(1,0){1.2} }
    \multiput(90,140)(0,2){22}{\line(1,0){1.2} }
    \multiput(100,140)(0,2){22}{\line(1,0){1.2} }
    \multiput(110,140)(0,2){22}{\line(1,0){1.2} }
    \multiput(120,140)(0,2){22}{\line(1,0){1.2} }
    \multiput(130,140)(0,2){22}{\line(1,0){1.2} }
    \multiput(140,140)(0,2){22}{\line(1,0){1.2} }
    \multiput(150,140)(0,2){22}{\line(1,0){1.2} }
    \multiput(160,140)(0,2){22}{\line(1,0){1.2} }
    \multiput(170,140)(0,2){22}{\line(1,0){1.2} }
    \multiput(180,140)(0,2){22}{\line(1,0){1.2} }
    \multiput(190,140)(0,2){22}{\line(1,0){1.2} }
    }
{
    \put(75,185){$2^{m+1}$}
    \put(77,160){$2^{m}$}
    \put(77,128){$2^{m-1}$}
    \put(89,65){$0$}
    \put(175,65){$2^m$}
    \put(142,65){$2^{m-1}$}
    \put(200,65){$2^{m+1}$}}
\large 
    \put(80,207){$\xi_n$}   
    \put(235,65){$|\xi'|$}
\end{picture}\\
\vspace{-15mm}
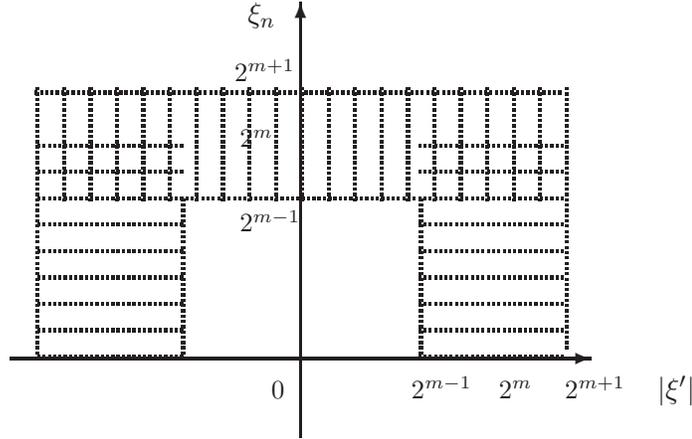
Fig 1:  The support of Littlewood--Paley decomposition 
$\{\overline{\Phi_m}\}_{m\in \Z}$
\end{center}

\noindent{\it Definition} (The Littlewood--Paley decomposition 
of separated variables).
For  $m\in \Z$, let
\eq{\label{eqn;zeta} 
 \begin{aligned}
 &\widehat{\zeta_m}(\xi_n)
   =\left\{
    \begin{aligned}
       1,&\qquad 0\le |\xi_n|\le 2^{m}, \\
       \text{smooth}&, 2^{m}<|\xi_n|<2^{m+1},\\
       0,&\qquad  2^{m+1}\le |\xi_n|, 
    \end{aligned}
    \right. \\
 &\widehat{\zeta_m}(\xi_n) 
    =\widehat{\zeta_{m-1}}(\xi_n)+\widehat{\phi_m}(\xi_n)
 \end{aligned}
}
(one can choose 
$\widehat{\zeta_m}(r)
 =\sum_{\ell\le m-1}\widehat{\phi_{\ell}}(r)+\widehat{\phi_{-\infty}}(r)$
with a correction distribution $\widehat{\phi_{-\infty}}(r)$ 
supported at $r=0$)
and set
\eq{\label{eqn;direct-sum-L-P}
  \widehat{\overline{\Phi_m}}(\xi)
  \equiv\widehat{\phi_m}(|\xi'|)\otimes \widehat{\zeta_{m-1}}(\xi_n)
   +\widehat{\zeta_{m}}(|\xi'|)\otimes \widehat{\phi_{m}}(\xi_n).
}
Then it is obvious from Fig. 2 (restricted on the upper half region in $\re^n$) 
that  
\eq{ \label{eqn;Sum_Phi_is_1}
 \sum_{m\in \Z}\widehat{\overline{\Phi_m}}(\xi)\equiv 1,
 \quad \xi=(\xi',\xi_n) \in \re^n\setminus \{0\}.
}
Indeed, from \eqref{eqn;zeta} and \eqref{eqn;direct-sum-L-P}, 
{\allowdisplaybreaks
\algn{
   &\sum_{m\in \Z}\widehat{\overline{\Phi_m}}(\xi)\\
  =&  \sum_{m\in\Z} \widehat{\phi_m}(|\xi'|)
                    \otimes \sum_{-\infty\le\ell\le m-1}\widehat{\phi_{\ell}}(\xi_n)
    + \sum_{m\in\Z} \sum_{\ell\le m} \widehat{\phi_{\ell}}(|\xi'|)
                    \otimes\widehat{\phi_{m}}(\xi_n)
    + \sum_{m\in\Z} \widehat{\phi_{-\infty}}(|\xi'|)\otimes \widehat{\phi_{m}}(\xi_n) 
  \\
  =&  \sum_{m\in\Z} \widehat{\phi_m}(|\xi'|)\otimes 
                   \Big(\sum_{-\infty\le\ell\le m-1}\widehat{\phi_{\ell}}(\xi_n)
                   + \sum_{\ell\ge m}\widehat{\phi_{\ell}}(\xi_n)\Big) 
      + \widehat{\phi_{-\infty}}(|\xi'|)\otimes \sum_{m\in\Z} \widehat{\phi_{m}}(\xi_n) 
 \\
  =&  \sum_{m\in\Z} \widehat{\phi_m}(|\xi'|)\otimes 
      \sum_{\ell\in\Z\cup\{-\infty\}} \widehat{\phi_{\ell}}(\xi_n)
        + \widehat{\phi_{-\infty}}(|\xi'|)
          \otimes \sum_{\ell\in\Z\cup\{-\infty\}} \widehat{\phi_{m}}(\xi_n) 
        - \widehat{\phi_{-\infty}}(|\xi'|)\otimes \widehat{\phi_{-\infty}}(\xi_n)
 \\
   =& 1- \widehat{\phi_{-\infty}}(|\xi'|)\otimes \widehat{\phi_{-\infty}}(\xi_n).
}
}
\vskip1mm
\noindent
{\it Definition} (Varieties of the Littlewood--Paley dyadic decompositions).
Let $(\t,\xi',\xi_n)\in \re\times \re^{n-1}\times \re$ be 
Fourier adjoint variables corresponding to 
$(t,x',\eta)\in \re\times \re^{n-1}\times \re$.
\begin{enumerate}
\item[$\bullet$] $\{\Phi_m(x)\}_{m\in \Z}$:
      the standard (annulus type) Littlewood--Paley dyadic decomposition by \hfill\break 
      $ x=(x',\eta)\in \re^n$.  
\item[$\bullet$]  $\{ \overline{\Phi_m}(x) \}_{m\in \Z}$:
      the Littlewood--Paley dyadic decomposition  over $x=(x',\eta)\in \re^n$
      given by \eqref{eqn;direct-sum-L-P}.
\item[$\bullet$]  $\{\psi_k(t)\}_{k\in \Z}$:
     the Littlewood--Paley dyadic decompositions  in $t\in \re$.
\item[$\bullet$]  $\{\phi_j(x')\}_{j\in \Z}$ and $\{\phi_j(\eta)\}_{j\in \Z}$:
      the standard (annulus type) Littlewood--Paley dyadic decompositions in $x'\in \re^{n-1}$ 
      and $\eta\in \re$,
      respectively.
\item[$\bullet$]  $\{\zeta_m(x')\}_{m\in \Z}$ and $\{\zeta_m(\eta)\}_{m\in \Z}$:
      the lower frequency smooth cut-off 
      given by \eqref{eqn;zeta}, respectively.
\item[$\bullet$] Let $\widetilde{\phi_j}=\phi_{j-1}+\phi_j+\phi_{j+1}$ be the 
 Littlewood--Paley dyadic decompositions with its $j$-neighborhood to $\phi_j$.
\item[$\bullet$]  All the above defined 
decompositions are even functions. 
\end{enumerate}


Then in view of Proposition \ref{prop;retraction-coretraction}
and the remark at the end of the previous sub-section, 
we see that the norm of the Besov spaces on $\re^n$
defined by $\{\Phi_m\}_m$ is equivalent to the one from the 
Littlewood--Paley decomposition of direct sum type, 
$\{\overline{\Phi_m}\}_m$ over $\re^n$ and hence 
one can identify those norms as it appears the homogeneous 
Besov space over $\re^n_+$ as follows.
Indeed, for any $1<p<\infty$ and $-1+1/p<s<1/p$,
\begin{align*}
\|\N q (t)\|_{\dB^s_{p,1}(\re^n_+)} 
 \le & C\|\N q (t)\|_{\overset{\odot\quad}{B^s_{p,1}}(\re^n_+)} 
   =   C\|E_0\big[\N q (t)\big]\|_{\dB^s_{p,1}(\re^n)}  \\
    =& C\sum_{m\in \Z} 2^{sm}
       \Big\|\Phi_m\underset{(x)}{*}  \sum_{|m-k|\le 1} 
          \overline{\Phi_{k}}\underset{(x)}{*} 
           E_0\big[\N q (t)\big]\Big\|_{L^p(\re^n)} \nonumber\\
 \le & 3C\sum_{m\in \Z} 2^{sm} 
       \big\|\overline{\Phi_m}\underset{(x)}{*} E_0\big[\N q (t)\big]
       \big\|_{L^p(\re^n)} \\
 \le & 3C\sum_{m\in \Z} 2^{sm} 
       \big\|\overline{\Phi_m}\underset{(x)}{*}
         \sum_{|m-k|\le 1} \Phi_{k}\underset{(x)}{*} E_0\big[\N q (t)\big]
       \big\|_{L^p(\re^n)} \\
 \le & 3^2C\sum_{m\in \Z} 2^{sm} 
       \big\|\Phi_{k}\underset{(x)}{*} E_0\big[\N q (t)\big]
       \big\|_{L^p(\re^n)} \\
 \le & 3^2C\big\| \N q (t) \big\|_{\dB^s_{p,1}(\re^n_+)}.   
 \eqntag\label{eqn;potential-besov-1}        
\end{align*}
In what follows, we freely use the retraction and 
coretraction operators as observed above and 
for simplicity we avoid reprised usage of them.

\sect{Almost orthogonality of the pressure potential} \label{Sec;3}
Almost orthogonality is the key lemma to obtain the maximal $L^1$-regularity estimate. 
In this  section,  we derive almost orthogonality concerning the pressure.

\vskip2mm
\noindent
{\it Defintion} (The pressure potentials).
For $j,k\in \Z$, let $\{\psi_k(t)\}_{k\in \Z}$, $\{\phi_j(x')\}_{j\in\Z}$ be 
the Littlewood--Paley decompositions for $t\in \re$ and $x'\in \re^{n-1}$
 valuables, respectively. We set for $\eta=x_n>0$,
\begin{align}
  &\left\{
  \begin{aligned}
   &\pi(t,x',\eta)
   \equiv 
    c_{n+1}\iint_{\re\times \re^{n-1}}
               e^{i t\t+ix'\cdot \xi'}
               (i\xi',-|\xi'|)^{\sf T}
               m(\t,\xi')e^{-|\xi'|\eta }
                d\t d\xi',  \\
   &\pi_{k,j}(t,x',\eta)
    \equiv \psi_k\underset{(t)}{*}\phi_j\underset{(x')}{*}\pi(t,x',\eta) \\
   &\phantom{\pi_{k,j}(t,x',\eta)}
    =\big(\pi_{k,j}'(t,x',\eta),\pi_{n,k,j}(t,x,\eta)\big),
  \end{aligned}
  \right.  \label{eqn;pressre-potential} 
\end{align}
where $m:\re\times \re^{n-1}\to \re^n$ is defined in 
\eqref{eqn;singular-int-multi}.  We extend the potential $\pi(t,x',\eta)$ into
all $\eta\in \re$ by the even extension (i.e. exchange $\eta$ into $|\eta|$).
\vskip2mm

Setting 
$\widetilde{\phi_j}=\phi_{j-1}+\phi_j+\phi_{j+1}$,
$\widetilde{\psi_k}=\psi_{k-1}+\psi_k+\psi_{k+1}$ 
and noting that  
$$
\sum_{k\in \Z}\sum_{j\in\Z}
            {\ppl\widehat\psi_k}(\t)\widehat{\phi_j}(\xi)
 \equiv 1,
 \quad {\ppl \t,\xi\neq 0}, 
$$
we have  for $x_ n>0$ that
{\allowdisplaybreaks
\algn{
  \N q&(t,x',x_n)  \\
  =& c_{n+1}\iint_{\re^{n}}
               e^{i t\t+ix'\cdot \xi'}
               (i\xi',-|\xi'|)^{\sf T}
               \Big(m'(\t,\xi')\cdot \widehat{H'}
                    +m_n(\t,\xi') \widehat{H_n}
               \Big)
                e^{-|\xi'|x_n }
               \sum_{k\in \Z}\sum_{j\in\Z}
          {\ppl\widehat{\psi_k}(\t)}\widehat{\phi_j}(\xi)
                d\t d\xi' 
\\
 \equiv 
   & \sum_{k\in \Z}\sum_{j\in\Z}
     \Big(
       {\pi'}_{k,j}\underset{(t,x')}{\cdot *} 
        \Big(\widetilde{\psi_k} \underset{(t)}{*}   
                  \widetilde{\phi_j}\underset{(x')}{*}H'\Big)
   +   {\pi_n}_{k,j}\underset{(t,x')}{*} 
        \Big(\widetilde{\psi_k} \underset{(t)}{*}
                  \widetilde{\phi_j}\underset{(x')}{*}H_n\Big)
    \Big), 
}
} 
where we use the notion of the inner product-convolution 
\eqref{eqn;convolution-innerprod} and the data is extended 
by the zero extension for $t\le 0$.
We show the almost orthogonality and its variation in the following.

\subsection{The almost orthogonality}
For the symbol of the gradient of the pressure, 
we introduce the useful notation for a part of the symbol 
defined by \eqref{eqn;B}; $B(\t,\xi')=\sqrt{i\tau+|\xi'|^2}$.
\vskip2mm
\noindent
{\it Definition.}  Let $\sg\in \re$
and  $\zeta'\in \re^{n-1}$  
with $1/2<|\sg|,|\zeta'|<2$.
For $a>0$, we set 
\eq{\label{eqn;b-d-def}
 \left\{
 \spl{
 &b(\sg,\zeta',a)=\sqrt{i\sg+a^{-2}|\zeta'|^2},\\
 &d(\sg,\zeta',a)=\sqrt{a^{-2}i\sg+|\zeta'|^2}. 
 }
 \right.
}


\begin{lem}\label{lem;B-bound} Let $\sg\in \re$ and
 $\zeta'\in \re^{n-1}$.\par\noindent
 {\rm (1)} For the time dominated region $k\ge 2j$, 
\eq{ \label{eqn;b-bound-time}
    \frac{1}{\sqrt{2}}\le |b(\sg,\zeta',2^{\frac{k}{2}-j})|\le 20^{1/4},  
   }
in particular, there exist constants $0<c<C$ independent of $j$ and $k$ such that 
\eq{ \label{eqn;B-bound-time}
  c2^{\frac{k}{2}} 
  \le |B(2^k\sg,2^{j}\zeta')|=2^{\frac{k}{2}} |b(\sg,\zeta',2^{\frac{k}{2}-j})|
  \le C 2^{\frac{k}{2}} .
   }
 \par\noindent
 {\rm (2)} For the space dominated region $k< 2j$, 
 $$
  \frac12\le |d(\sg,\eta',2^{j-\frac{k}{2}})|
         \le 20^{1/4}, 
 $$
in particular,  there exist constants $0<c<C$ independent of $j$ and $k$ such that 
 \eq{ \label{eqn;B-bound-space}
  c 2^{j} 
  \le |B(2^k\sg,2^j\zeta')|=2^{j} |d(\sg,\zeta',2^{j-\frac{k}{2}})|
  \le C 2^{j+3}.
 }
\end{lem}

\begin{prf}{Lemma \ref{lem;B-bound}} \ 
(1) In the case when $k\ge 2j$, 
by using $2^{-1}<|\sg|<2$, $2^{-1}<|\zeta'|<2$, it holds that  
\algn{
   B(2^k\sg,2^j\zeta')
   =&2^{\frac{k}{2}}b(\sg,\zeta',a)|_{a=2^{\frac{k}{2}-j}} 
          =2^{\frac{k}{2}} \sqrt{i\sg+(2^{j-\frac{k}{2}})^2|\zeta'|^2} \\
          =&2^{\frac{k}{2}}\cdot  {}^4\sqrt{\sg^2+(2^{j-\frac{k}{2}}|\zeta'|)^4}
             \exp\Big(\frac{i}{2}\tan^{-1}\frac{2^k\sg}{2^{2j}|\zeta'|^2}\Big), 
}
and \eqref{eqn;B-bound-time} follows from 
$$
 2^{-\frac12}2^{\frac{k}{2}}
    \le 2^{\frac{k}{2}}\cdot{}^4\sqrt{\sg^2}
  \le |B(2^k\sg,2^j\zeta')|
      = 2^{\frac{k}{2}}\cdot{}^4\sqrt{\sg^2+2^{4j-2k}|\zeta'|^4}
    \le  20^{1/4}\cdot 2^{\frac{k}{2}}.
$$
(2) In the case when $k< 2j$, it holds that  
\eqn{
 \spl{
    2^{-1} 2^{j}\le 2^{j}\cdot{}^4\sqrt{|\zeta'|^4}
  & \le |B(2^k\sg,2^j\zeta')|= 2^{j}\cdot{}^4\sqrt{2^{2k-4j}\sg^2+|\eta'|^4}
    \le  20^{1/4}\cdot 2^{j}.
}
}
The constants $c$ and $C$ can be taken as $c=1/\sqrt{2}$ and $C=\sqrt{2\sqrt{5}}$.
\end{prf}
\vskip2mm

\begin{lem}[Almost orthogonality I]
\label{lem;pressure-orthogonal}
For $k,j\in \Z$,  let  
$\pi_{k,j}(t, x',\eta)$  be the pressure potentials defined by 
\eqref{eqn;pressre-potential} 
and let 
$\{\psi_k(t)\}_{k\in\Z}$ and $\{\phi_j(x)\}_{j\in \Z}$ 
be the Littlewood--Paley decompositions 
for time and space, respectively. 
\begin{enumerate}
\item[(1)]  
For the time-dominated region $k\ge 2j$,
there exists $C_n>0$ such that for any $\eta\in \re_+$
and $t\in \re$,
\alg{ 
  &\big\|\pi_{k,j}(t,\cdot,\eta)\big\|_{L^1_{x'}}
    \le C_n 2^{j}\big(1+(2^j\eta)^{n+2}\big)
        e^{-2^{(j-1)}\eta}\frac{2^k}{\<2^k t\>^{2}},
    \label{eqn;pressure-orthogonarity-T}
  }
where $\|\cdot\|_{L^1_{x'}}$ denotes 
the $L^1(\re^{n-1})$ norm in $x'$-variable.
\item[(2)] For the space-dominated region $k<2j$, 
there exists $C_n>0$ such that for any $\eta\in \re_+$ 
and $t\in \re$,
\eq{ 
  \Big\|\sum_{k< 2j}\pi_{k,j}(t,\cdot,\eta) \Big\|_{L^1_{x'}} 
   \le C_n 2^{j}\big(1+(2^j\eta)^{n+2}\big)
       e^{-2^{(j-1)}\eta} \frac{2^{2j}}{\<2^{2j} t\>^{2}}.
      \label{eqn;pressure-orthogonarity-S} 
  } 
\end{enumerate}
The estimates are extended to $\eta\in \re$ by the even extensions.
\end{lem}

\begin{prf}{Lemma \ref{lem;pressure-orthogonal}}\ 
(1) In the time-dominated region $k\ge 2j$, by using the expression 
of the fundamental solution and using change of variables $\t=2^k\sg$, 
$\xi'=2^j\zeta'$ and then $x'=2^{-j}y'$, 
we first observe that
{\allowdisplaybreaks  
\algn{ 
   \|\pi_{k,j}&(t,\cdot,\eta)\|_{L^1_{x'}} \\
  =&\left\|
     c_{n+1}
     \int_{\re}\int_{\re^{n-1}}
         e^{it\t+ix'\cdot \xi'}
         (i\xi', -|\xi'|)^{\sf T}m(\t,\xi')
          e^{-|\xi'|\eta}
         \widehat{\psi}(2^{-k}\t)
         \widehat{\phi}(2^{-j}\xi')
     d\xi'd\t   
    \right\|_{L^1_{x'}}
\\
  =&\left\|
     c_{n+1}\int_{\re}\int_{\re^{n-1}}
         e^{i2^k t\sg+i2^jx'\cdot \zeta'}
         (2^ji\zeta',-2^j|\zeta'|)^{\sf T}
         m(2^k\sg, 2^j\zeta')  
         e^{-2^j|\zeta'|\eta}
         \widehat{\psi}(\sg)
         \widehat{\phi}(\zeta')
         2^{(n-1)j}d\zeta'\, 2^kd\sg   
    \right\|_{L^1_{x'}}
 \\
   =& 2^{j+k} e^{-(2^{j-1}\eta)} 
      \left\|
         c_{n+1}\int_{\re}\int_{\re^{n-1}}
            e^{i2^k t\sg+i2^jx'\cdot \zeta'}
            (i\zeta',-|\zeta'|)^{\sf T}
            m(2^k\sg,2^j\zeta') 
      \right.\\
    &\hskip4cm \times \left.
            \exp\big(- 2^j\eta(|\zeta'|-\frac12)\big)
             \widehat{\psi}(\sg)
             \widehat{\phi}(\zeta')
            2^{(n-1)j}d\zeta' d\sg   
      \right\|_{L^1_{x'}}
\\
     =&2^{k+j} e^{-(2^{j-1}\eta)}
         \left\|
            c_{n+1}\int_{\re}\int_{\re^{n-1}}
            e^{i2^k t\sg+iy'\cdot \zeta'}
            (i\zeta',-|\zeta'|)^{\sf T}
          \right. \\
      &\hskip4cm
            \times
            \left.
            m(2^k\sg,2^j\zeta')
            \exp\big(-2^j\eta(|\zeta'|-\frac12)\big)
            \widehat{\psi}(\sg)
            \widehat{\phi}(\zeta')
            d\zeta'\, d\sg   
         \right\|_{L^1_{y'}}.   
     \eqntag \label{eqn;ortho-1}
}
}
Since the Fourier inverse transform of the most right term of 
the above equation contains the Littlewood--Paley cut-off for 
$\sg$ and $\zeta'$, it is integrable absolutely with respect to 
$\sg$ and $\zeta'$. 
If the symbol $m(2^k \t,2^j\zeta')$ is bounded, 
then $\pi_{k,j}(t,2^{-j}y',\eta)$ are integrable with respect to $y'$ when $|y'|<1$. 
Therefore we check the boundedness of the symbol  $m(2^k \sg,2^j\zeta')$. 
Recalling the definition of $m'(\t,\xi')$ in \eqref{eqn;singular-int-multi}  
with using $b(\sg,\zeta,2^{\frac{k}{2}-j})$ in \eqref{eqn;b-d-def} 
and its bound \eqref{eqn;B-bound-time},  it holds that
\algn{
 m'&(2^k \sg, 2^j\zeta') \\
 =&2i\frac{\xi'}{|\xi'|}
    \left.\frac{|\xi'|^2B+|\xi'|B^2}{B^3+|\xi'|B^2+3|\xi'|^2B-|\xi'|^3}
     \right|_{\t=2^k\sg,\xi'=2^j\zeta'} \\
 =&2i\frac{\zeta'}{|\zeta'|}
     \frac{2^{\frac{k}2+2j}|\zeta'|^2b(\sg,\zeta,2^{\frac{k}{2}-j}) 
           +2^{k+j}|\zeta'|b(\sg,\zeta,2^{\frac{k}{2}-j})^2
           }
          {2^{\frac32k}b(\sg,\zeta,2^{\frac{k}{2}-j})^3
           +2^{k+j}|\zeta'|b(\sg,\zeta,2^{\frac{k}{2}-j})^2
           +3\cdot2^{\frac{k}2+2j}|\zeta'|^2b(\sg,\zeta,2^{\frac{k}{2}-j})
           -2^{3j}|\zeta'|^3}\\
 =&2i\frac{\zeta'}{|\zeta'|}
      \frac{ 2^{-2(\frac{k}{2}-j)}|\zeta'|^2 b(\sg,\zeta',2^{2j-k}) 
           + 2^{-(\frac{k}{2}-j)}|\zeta'| b(\sg,\zeta,2^{\frac{k}{2}-j})^2
           }
          {(  b(\sg,\zeta,2^{\frac{k}{2}-j})^3
            +2^{-(\frac{k}{2}-j)}|\zeta'|\,b(\sg,\zeta,2^{\frac{k}{2}-j})^2
            +3\cdot 2^{-2(\frac{k}{2}-j)}|\zeta'|^2\, b(\sg,\zeta,2^{\frac{k}{2}-j})
            -2^{-3(\frac{k}{2}-j)}|\zeta'|^3 }
 }
and thus for the  case  $k -2j>3$, 
it holds $2^{-m(\frac{k}{2}-j)}\le 2^{-\frac{3}{2}m} $ with $m=1,2,3$ and 
from \eqref{eqn;b-bound-time}, we have  
\algn{
| m'&(2^k \sg, 2^j\zeta') |
 \le C\frac{2^{-(\frac{k}{2}-j)} |\zeta'| }{b(\sg,\zeta,2^{\frac{k}{2}-j})}
 \le C,
\eqntag \label{eqn;ortho-2}
}
and otherwise $0\le  k-2j\le 3$ it is obviously bounded from above and below 
since the denominator never vanishes.
Furthermore,
\algn{
 m_n&(2^k \sg,  2^j\zeta') \\
 =&\left.\frac{(|\xi'|+B)(|\xi'|^2+B^2)}{B^3+|\xi'|B^2+3|\xi'|^2B-|\xi'|^3}
   \right|_{\t=2^k\sg,\xi'=2^j\zeta'} \\
 =&\frac{\big( 2^j|\zeta'|+ 2^{\frac{k}{2}}b(\sg,\zeta,2^{\frac{k}{2}-j})\big)
         \big(2^{2j}|\zeta'|^2+2^{k} b(\sg,\zeta,2^{\frac{k}{2}-j})^2\big)}
        {(2^{\frac{3k}{2}} b(\sg,\zeta,2^{\frac{k}{2}-j})^3
           +2^{k+j}|\zeta'| b(\sg,\zeta,2^{\frac{k}{2}-j})^2
           +3\cdot2^{\frac{k}{2}+2j}|\zeta'|^2 b(\sg,\zeta,2^{\frac{k}{2}-j})
           -2^{3j}|\zeta'|^3}\\
  =&\frac{\big( 2^{-(\frac{k}{2}-j)}|\zeta'|+ b(\sg,\zeta,2^{\frac{k}{2}-j})\big)
         \big(2^{-2(\frac{k}{2}-j)}|\zeta'|^2+ b(\sg,\zeta,2^{\frac{k}{2}-j})^2\big)}
        {(  b(\sg,\zeta,2^{\frac{k}{2}-j})^3
            +2^{-(\frac{k}{2}-j)}|\zeta'|\,b(\sg,\zeta,2^{\frac{k}{2}-j})^2
            +3\cdot 2^{-2(\frac{k}{2}-j)}|\zeta'|^2\, b(\sg,\zeta,2^{\frac{k}{2}-j})
            -2^{-3(\frac{k}{2}-j)}|\zeta'|^3}
}
and {\ppl similarly}
\algn{ 
| m_n&(2^k \sg, 2^j\zeta') |
 \lesssim \frac{b(\sg,\zeta,2^{\frac{k}{2}-j})^3 }{b(\sg,\zeta,2^{\frac{k}{2}-j})^3}
 \le C.
 \label{eqn;ortho-3}\eqntag
}
Note that the common {\crd denominator} $B^3+|\xi'|B^2+3|\xi'|^2B-|\xi'|^3$ has  no 
zero point except $(\t,\xi')=(0,0)$ (cf. Lemma 4.4 in \cite{SbSz03}).

For $t<1$, we obtain from \eqref{eqn;ortho-2} and \eqref{eqn;ortho-3} that 
\algn{\label{eqn;ortho-4}\eqntag
  &\|\pi_{k,j}(t,\cdot)
     \|_{ L^1_{x'}(B_{2^{-j}}) } \\
   =&C_n2^{k+j} e^{-(2^{j-1}\eta) } 
     \bigg\|
      \int_{\re}\int_{\re^{n-1}}
         e^{i(2^kt\sg+ y'\cdot \zeta')}\\
    &\phantom{2^{k-2\ell}\|c_n\int_{\re}\int_{\re^n}} 
         \times 
          (i\zeta', -|\zeta'|)^{\sf T} 
          m(2^k\sg, 2^j\zeta')
         \exp\big(- 2^j\eta(|\zeta'|-\frac12)\big) 
         \widehat{\psi}(\sg)\widehat{\phi}(\zeta')
         d\zeta' d\sg  
      \bigg\|_{L^1_{y'}(B_1)}\\
   \le& C_n2^{k+j} e^{-(2^{j-1}\eta)} .  
}
Next we consider the case  when $t>1$. It is important that we gain 
decay of time for $t>1$ by integration by parts.   
Noting that 
\eqn{ 
   e^{i(2^kt\sg+y'\cdot \zeta')}
   =\left(\frac{1}{2^k it}\right)^2
    \pt_\sg^2 e^{i(2^kt\sg+y'\cdot \zeta')}, 
} 
and integrating by parts with respect to $\sg$ twice, we obtain 
\algn{ \eqntag \label{eqn;ortho-6}
  &\|\pi_{k,j}(t,\cdot)
     \|_{ L^1_{x'}(B_{2^{-j}}) } \\
   =& 2^{k+j} e^{-(2^{j-1}\eta) } 
     \bigg\|
      c_{n+1}\int_{\re}\int_{\re^{n-1}}
         \frac{1}{(2^k it)^2} e^{i(2^kt\sg+ y'\cdot \zeta')}\\
    &\hskip3.5cm 
         \times  (i\zeta', -|\zeta'|)^{\sf T} 
         \pt_\sg^2
          \left( 
             m(2^k\sg, 2^j\zeta')
            \widehat{\psi}(\sg)
          \right)         
          \exp\big(- 2^j\eta(|\zeta'|-\frac12)\big) 
          \widehat{\phi}(\zeta')d\zeta' d\sg  
      \bigg\|_{L^1_{y'}(B_1)}.
}
Again we separate the region $k\ge 2j$ and $k<2j$. For $m'$, 
from \eqref{eqn;ortho-2} we use $b=b(\sg,\zeta',2^{\frac{k}{2}-j})$ 
defined in \eqref{eqn;b-d-def} and  
$a=2^{\frac{k}{2}-j}$ to see
\algn{
\label{eqn;ortho-7} \eqntag
\pt_{\sg}&m'(2^k\sg,2^j\zeta') \\
   =&\frac{\pt}{\pt\sg}
     \left(\left.2i\frac{\xi'}{|\xi'|}
     \frac{|\xi'|^2B+|\xi'|B^2}
          {B^3+|\xi'|B^2+3|\xi'|^2B-|\xi'|^3}     
    \right|_{\t=2^k\sg,\xi'=2^j\zeta'} 
    \right)\\
   =& 2i\frac{\zeta'}{|\zeta'|}
        \frac{\pt b}{\pt \sg}
        \frac{\pt}{\pt b}        
        \frac{ |\zeta'|^2 a^{-2} b+ |\zeta'| a^{-1} b^2
             }
             {b^3 +a^{-1}|\zeta'|\, b^2 +3a^{-2}|\zeta'|^2\, b -a^{-3} |\zeta'|^3
             }
\\
   =& \frac{\zeta'}{|\zeta'|b}
       \cdot
        \bigg\{       
        \frac{      a^{-1} |\zeta'|  b^4 
               +  2 a^{-2}|\zeta'|^2\,b^3 
               -  2 a^{-3}|\zeta'|^3\,b^2         
               +  2 a^{-4}|\zeta'|^4b
               +    a^{-5} |\zeta'|^5
             }           
             {\big(
                b^3
              + a^{-1}|\zeta'|\, b^2
              +3a^{-2}|\zeta'|^2\,b
              - a^{-3}|\zeta'|^3
              \big)^2
            }
        \bigg\}. 
}
Analogously
\algn{
\label{eqn;ortho-8} \eqntag
\pt_{\sg}&m_n(2^k\sg,2^j\zeta') \\
   =&\frac{\pt}{\pt \sg}
     \left(\left.
     \frac{(|\xi'|+B)(|\xi'|^2+B^2)}
          {B^3+|\xi'|B^2+3|\xi'|^2B-|\xi'|^3}
   \right|_{\t=2^k\sg,\xi'=2^j\zeta'}
   \right) \\
  =&\frac{\pt b}{\pt \sg}
     \frac{\pt}{\pt b}
     \frac{ (a^{-1}|\zeta'|+ b)
            (a^{-2}|\zeta'|^2+b^2)}
          {\big(b^3
             +a^{-1}|\zeta'|b^2
             +3a^{-2}|\zeta'|^2b
             -a^{-3}|\zeta'|^3\big)}
  \\
  =&\frac{2i}{b}
    \frac{  a^{-2}|\zeta'|^2 b^3
           -a^{-3}|\zeta'|^3 b^2 
           -a^{-4}|\zeta'|^4 b
           -a^{-5}|\zeta'|^5
         }
         {\big(b^3
             + a^{-1}|\zeta'| b^2
             +3a^{-2}|\zeta'|^2b
             - a^{-3}|\zeta'|^3\big)^2
         }.
}
By Lemma \ref{lem;B-bound}, $\pt_{\sg}m'(2^k\sg,2^j\zeta')$ 
and $\pt_{\sg}m_n(2^k\sg,2^j\zeta')$ are bounded 
from above for all $(k,j)$ when $k\ge 2j$. 
Since the denominator does not vanish because it is smooth on the support 
of $\sg$ and $\eta$ such that $|\sg|, |\zeta'|\in  (1/2,2)$, 
 there is no diverging coefficient 
from $\sg$-derivative of $m$, $\pt_{\sg}^2 m$ is bounded on the support 
of $\widehat{\psi}(\sg)$. 
The situation is same for the second derivative with respect to $\sg$. 
Therefore for $t>1$, it holds from \eqref{eqn;ortho-6}-\eqref{eqn;ortho-8} that 
\algn{
  &\|\pi_{k,j}(t,\cdot)
   \|_{ L^1_{x'}(B_{2^{-j}}) } \\
   =& 2^{k+j} e^{-(2^{j-1}\eta)} \frac{1}{(2it)^2} 
     \bigg\|
      c_{n+1}\int_{\re}\int_{\re^{n-1}}
         e^{i(2^kt\sg+ y'\cdot \zeta')}\\
    &\phantom{ 2^{k-2\ell}\|c_n\int_{\re}\int_{\re}\int_{\re^{n-1}}  } 
         \times  (i\zeta', -|\zeta'|)^{\sf T} 
         \pt_{\sg}^2
         \Big(m(2^k\sg, 2^j\zeta')\widehat{\psi}(\sg)
         \Big)
         \exp\big(- 2^j\eta(|\zeta'|-\frac12)\big) 
         \widehat{\phi}(\zeta')
         d\zeta' d\sg  
      \bigg\|_{L^1_{y'}(B_1)}\\
   \le& C2^{j} e^{-(2^{j-1}\eta) } \frac{2^k}{(2^kt)^2}. 
   \eqntag \label{eqn;ortho-10}
}
On the other hand, in the case $y'\in B_1^c(0)$ we differentiate the symbol 
$$
 m(2^k\sg,2^j\zeta')
 (i\zeta',-|\zeta'|)
 e^{-2^j\eta(|\zeta'|-1/2)}\widehat{\phi}(\zeta')
$$
$n$ times with respect to $\zeta'$.  If we obtain boundedness 
of the symbol uniformly $k$ and $j$, 
then we obtain $(t,y')$-decay estimate by using
\eqn{\label{eqn;ortho-12}
   e^{i(2^kt\sg+y'\cdot \zeta')}
   =\left(\frac{1}{|y'|^n}\frac{1}{(2^ki t)^2}\right)
    \pt_\sg^2(-\Del_{\zeta'})^{\frac{n}{2}} 
      e^{i(2^kt\sg+y'\cdot \zeta')}. 
}
In this way it holds that  
\algn{
  \|& \pi_{k,j}(t,\cdot)
     \|_{ L^1_{x'}(B_{2^{-j}}^c) } \\
   =& 2^{k+j} e^{-(2^{j-1}\eta)} 
      \frac{1}{(2^k t)^2} 
     \bigg\|
       \frac{c_{n+1}}{\<y'\>^n} 
       \int_{\re}\int_{\re^{n-1}}
         e^{i(2^kt\sg+ y'\cdot \zeta')}\\
    &\phantom{2} 
         \times 
         (1-\pt_{\sg}^2)\cdot (1-\Del_{\zeta'})^{\frac{n}{2}}
         \Big((i\zeta', -|\zeta'|)^{\sf T} 
               m(2^k\sg, 2^j\zeta') 
            \exp\big(- 2^j\eta(|\zeta'|-\frac12)\big)
            \widehat{\psi}(\sg) 
            \widehat{\phi}(\zeta')
         \Big)
         d\zeta' d\sg  
      \bigg\|_{L^1_{y'}(B_1^c)}\\
   \le&C2^{j}\big(1+(2^{j}\eta)^{n+2}\big) e^{-(2^{j-1}\eta)}  
       \frac{2^k}{(2^kt)^2}.
  \eqntag\label{eqn;ortho-14}
}
Combining with \eqref{eqn;ortho-4}, \eqref{eqn;ortho-10} 
and \eqref{eqn;ortho-14}, we obtain the estimate  
\eqref{eqn;pressure-orthogonarity-T}.
\par

(2)
 In the space-dominated region $k< 2j$, 
the proof is almost the same as in (1).
Using the notation $\widehat{\zeta}(2^{-2j} \t)=\sum_{k<2j}\widehat{\psi}(2^{-k}\t)$
for the Littlewood--Paley decomposition (see \eqref{eqn;low-freq-zeta} and
\eqref{eqn;zeta}), 
and applying the change of variables $\xi'=2^j\zeta'$, $\t=2^{2j}\sg$ and then 
$x'=2^{-j}y'$, we have for 
$\widehat{\zeta}(2^{-2j}\t)=\sum_{k\le 2j}\widehat{\psi}(2^{-k}\t)$ that 
\algn{ 
   \Big\|&\sum_{k<2j}\pi_{k,j}(t,\cdot,\eta)\Big\|_{L^1_{x'}} \\
  =&\left\|
     c_{n+1}\int_{\re}\int_{\re^{n-1}}
         e^{it\t+ix'\cdot \xi'} (i\xi', -|\xi'|) ^{\sf T} 
         m(\t,\xi') e^{-|\xi'|\eta}
         \sum_{k<2j} 
         \widehat{\psi}(2^{-k}\t)
         \widehat{\phi}(2^{-j}\xi')
     d\xi'd\t   
    \right\|_{L^1_{x'}}
\\
  =&\left\|
     c_{n+1}\int_{\re}\int_{\re^{n-1}}
         e^{i2^{2j} t\s+i2^jx'\cdot \zeta'}
         (2^ji\zeta',-2^j|\zeta'|)^{\sf T}
         m(2^{2j}\sg, 2^j\zeta')  
         e^{-2^j|\zeta'|\eta}
         \widehat{\zeta}(2^{-2j}\t)
         \widehat{\phi}(\zeta')
         2^{(n-1)j}d\zeta'\cdot 2^{2j}d\sg   
    \right\|_{L^1_{x'}}
 \\
     =&C_n2^{3j} e^{-(2^{j-1}\eta)}
         \left\|
            \int_{\re}\int_{\re^{n-1}}
            e^{2^{2j}it\s+iy'\cdot \zeta'}
            (i\zeta',-|\zeta'|)^{\sf T}
            m(2^{2j}\s,2^j\zeta')
          \right. \\
      &\hskip4cm
            \times
            \left.       
            \exp\big(-2^j\eta(|\zeta'|-1/2)\big)
            \widehat{\zeta}(\sg)
            \widehat{\phi}(\zeta')
            d\zeta'\, d\sg   
         \right\|_{L^1_{y'}}.   
}
Using $d=d(\sg,\zeta,1)$, we have from \eqref{eqn;B-bound-space} that
\algn{ \label{eqn;ortho-21} \eqntag
  m'(2^{2j}\s,2^j\zeta') 
 =&2i\frac{\zeta'}{|\zeta'|}
     \frac{ 2^{3j}        |\zeta'|^2\sqrt{i\s+|\zeta'|^2} 
           + 2^{3j}       |\zeta'|  (i\s+|\zeta'|^2)}
          {  2^{3j}             \sqrt{i\s+|\zeta'|^2}^3
           + 2^{3j}       |\zeta'| (i\s+|\zeta'|^2)
           +3\cdot  2^{3j}|\zeta'|^2\sqrt{i\s+|\zeta'|^2}
           - 2^{3j}|\zeta'|^3}
     \\
  =&2i\frac{\zeta'}{|\zeta'|}
     \frac{       |\zeta'|^2\sqrt{i\s+|\zeta'|^2} 
           +      |\zeta'|  (i\s+|\zeta'|^2)}
          {                 \sqrt{i\s+|\zeta'|^2}^3
           +      |\zeta'|   (i\s+|\zeta'|^2)
           +3\cdot|\zeta'|^2\sqrt{i\s+|\zeta'|^2}
           -      |\zeta'|^3}
     \\
  =&2i\frac{\zeta'}{|\zeta'|}
     \frac{ |\zeta'|^2 d
           +|\zeta'|   d^2
           }
          {            d^3
           +|\zeta'|   d^2
           +3|\zeta'|^2d
           -|\zeta'|^3}
 \simeq  i\frac{|\zeta'|^2\zeta'}{|\zeta'|^3},
 \\ 
  m_n(2^{2j}\s, 2^j \zeta')
 =&\frac{\big(|\zeta'|+\sqrt{i\sg+|\zeta'|^2}\big)
         \big(|\zeta'|^2+(i\sg+|\zeta'|^2)\big)}
        {\sqrt{i\sg+|\zeta'|^2}^3
           +|\zeta'|    (i\t+2^{2j}|\zeta'|^2)
           +3|\zeta'|^2 \sqrt{i\t+2^{2j}|\zeta'|^2}
           -|\zeta'|^3}\\
 =&\frac{\big( |\zeta'| + d\big)
         \big(|\zeta'|^2+ d^2\big)
        }
        {             d^3
         + |\zeta'|   d^2
         +3|\zeta'|^2 d
         - |\zeta'|^3}
 \quad \simeq C,
}
where $\simeq$ stands for the equivalence with a constant.
Besides by $d\equiv d(\sg,\zeta',1)$ for simplicity, and noting 
$\frac{\pt}{\pt \sg}d(\sg,\zeta',1)=i(2d)^{-1}$, we see 
\algn{
\pt_{\s}m'(2^{2j}\sg,2^j\zeta') 
   =& 2i\frac{\zeta'}{|\zeta'|}
        \frac{\pt d}{\pt \s}
        \frac{\pt}{\pt d}        
        \frac{ |\zeta'|^2 d
              +|\zeta'|  d^2
             }
             { d^3
               +|\zeta'|\, d^2
               +3|\zeta'|^2\, d
               - |\zeta'|^3
             }
 \simeq \frac{C}{d^3(\sg,\zeta',1)} \simeq \frac{C}{\<\sg\>^{\frac32}},
\\
\pt_{\sg}m_n(2^{2j}\s,2^j\zeta') 
  =& \frac{\pt d}{\pt \s}
     \frac{\pt}{\pt d}
     \frac{ |\zeta'|^3+ |\zeta'|^2d
          + |\zeta'| d^2 +d^3
          }
          {\big(d^3
             + |\zeta'|  d^2
             +3|\zeta'|^2d
             - |\zeta'|^3\big)}
\simeq  \frac{C}{d^2} \simeq \frac{C}{\<\sg\>}.
\label{eqn;ortho-22} \eqntag
}
 Similarly one can estimate the second derivative of the symbol 
 $m(2^{2j} \sg, 2^j\zeta')$ 
 and it is now clear that the second derivatives are also bounded 
 over the support of $\widehat{\psi}(\sg)$ and $\widehat{\phi}(\zeta')$.

Hence by the boundedness obtained from \eqref{eqn;ortho-21} and \eqref{eqn;ortho-22}, 
the rest of computation go through along the same line 
to \eqref{eqn;ortho-4}, \eqref{eqn;ortho-10} and \eqref{eqn;ortho-14},
and we conclude  for the ball $B_{2^{-j}}$ with radius $2^{-j}$ 
around the origin that 
\eq{\label{eqn;ortho-26}
 \spl{
  \Big\|\sum_{k<2j}&\pi_{k,j}(t,\cdot)\Big\|_{L^1_{x'}(B_{2^{-j}})} \\
  \le &C_n2^{3j} e^{-(2^{j-1}\eta)} \frac{1}{\<2^{2j}t\>^2} 
     \bigg\|
     \int_{\re}\int_{\re^{n-1}}
         e^{i(2^{2j}t\sg+ y'\cdot \zeta')}
         (i\zeta',-|\zeta'|)^{\sf T}
         (1+\pt_{\sg}^2)
         \Big(m(2^{2j}\sg, 2^j\zeta')\widehat{\zeta}(\sg)
         \Big)\\
    &\hskip5cm 
         \times          
         \exp\big(- 2^j\eta(|\zeta'|-\frac12)\big) 
         \widehat{\phi}(\zeta')
         d\zeta' d\sg  
      \bigg\|_{L^1_{y'}(B_1)}\\
   \le& C2^{j} e^{-(2^{j-1}\eta) } \frac{2^{2j}}{\<2^{2j}t\>^2}. 
   }
}

Very much similar way to the case except the cut off function 
$\widehat{\zeta}(\sg)$ instead of $\widehat{\psi}(\sg)$, we proceed as before that
{\allowdisplaybreaks  
\algn{
  \Big\|\sum_{k< 2j}&\pi_{k,j}(t,\cdot)
  \Big\|_{ L^1_{x'}(B_{2^{-j}}^c) } \\
   =& 2^{3j} e^{-(2^{j-1}\eta)} 
      \frac{1}{(2^{2j} t)^2} 
     \bigg\|
       \frac{c_{n+1}}{\<y'\>^n} 
       \int_{\re}\int_{\re^{n-1}}
         e^{i(2^kt\sg+ y'\cdot \zeta')}\\
    & \times 
         (1-\pt_{\sg}^2)\cdot (1-\Del_{\zeta'})^{\frac{n}{2}}
         \Big((i\zeta',-|\zeta'|)^{\sf T}
               m(2^k\sg, 2^j\zeta')                    
            \exp\big(- 2^j\eta(|\zeta'|-\frac12)\big)
            \widehat{\zeta}(\sg) 
            \widehat{\phi}(\zeta')
         \Big)
         d\zeta' d\sg  
      \bigg\|_{L^1_{y'}(B_1^c)}\\
   \le&C2^{j}\big(1+(2^{j}\eta)^{n+2}\big) e^{-(2^{j-1}\eta)}  
       \frac{2^{2j}}{\<2^{2j}t\>^2}.
  \eqntag\label{eqn;ortho-27}
}
}
From \eqref{eqn;ortho-26} and \eqref{eqn;ortho-27}, 
we conclude the desired estimate. 
\end{prf}
\vskip1mm
\subsection{The second almost orthogonality} 
We consider the  almost orthogonality estimate of second type  
which will be used for the triumphal arch type 
Littlewood--Paley dyadic decomposition.

\begin{lem}[Almost orthogonality II]
\label{lem;pressure-orthogonal2}
Let $k,j,m\in \Z$ and 
$\pi_{k,j}(\t,\xi',\eta)$
be the pressure potential given by 
\eqref{eqn;pressre-potential} and let  
$\{\psi_k(t)\}_{k\in\Z}$ and $\{\phi_j(x)\}_{j\in \Z}$ be 
the Littlewood--Paley decompositions 
for time and space, respectively. Assume that $j\le m$.
Then for any $N\in \Nt$, there exists a constant $C_{n,N}>0$ 
depending on $n$ and $N$ such that 
the following estimates hold:
\begin{enumerate}
\item[(1)]  
For  the time-dominated region $k\ge 2j$,
\alg{ 
  \big\| \phi_m\underset{(\eta)}{*}
         \pi_{k,j}(t,\cdot,\eta)
  \big\|_{L^1_{x'}}
   \le &  C_{n,N}\frac{2^j  2^{-(m-j)} }{\< 2^j \eta\>^N}
         \frac{2^k }{\<2^k t\>^{2}}.
        \label{eqn;pressure-orthogonarity-T-2} 
}
\item[(2)] For the space-dominated region $k< 2j$, it holds that 
\eq{ \label{eqn;pressure-orthogonarity-S-2}
   \Big\|\sum_{k<2j} \phi_m\underset{(\eta)}{*}\pi_{k,j}(t,\cdot,\eta)
   \Big\|_{L^1_{x'}}
  \le C_{n,N} \frac{2^j  2^{-(m-j)} }{\<  2^j \eta\>^N}
          \frac{2^{2j}}{\<2^{2j}t\>^{2}}. 
  }
\end{enumerate}
\end{lem}
\vskip1mm
\begin{prf}{Lemma \ref{lem;pressure-orthogonal2}} 
(1) In the time-dominated region $k\ge 2j$, by using the expression 
of the fundamental solution and using change of variables 
$\t=2^k\sg$, $\xi'=2^j\zeta'$ and then $x'=2^{-j}y'$, we have 
{\allowdisplaybreaks
\algn{ 
   \|\phi_m\underset{(\eta)}{*}&\pi_{k,j}(t,\cdot,\eta)\|_{L^1_{x'}} \\
 =&\left\|
     c_{n+1}\int_{\re}\int_{\re^{n-1}}
         e^{it\t+ix'\cdot \xi'}
         (i\xi',-|\xi'|)^{\sf T}
         m(\t,\xi')
         \big( \phi_m\underset{(\eta)}{*} e^{-|\xi'|\eta}
         \big)
         \widehat{\psi}(2^{-k}\t)
         \widehat{\phi}(2^{-j}\xi')
     d\xi'd\t   
    \right\|_{L^1_{x'}}
 \\
     =&2^{k} \left\|c_{n+1}\int_{\re}\int_{\re^{n-1}}
               e^{i2^k t\sg+iy'\cdot \zeta'}
               2^j(i\zeta',-|\zeta'|)^{\sf T}
               m(2^k\sg,2^j\zeta')
            \big(
                \phi_m\underset{(\eta)}{*} e^{-2^j|\zeta'|\eta}
            \big)
              \widehat{\psi}(\sg)
              \widehat{\phi}(\zeta')
            d\zeta'\, d\sg   
         \right\|_{L^1_{y'}}.
   \eqntag \label{eqn;ortho-30}
}}
Applying 
\eq{ \label{eqn;ortho-5'}
   e^{i(2^kt\sg+y'\cdot \zeta')}
   = \frac{1}{(2^k it)^2}
     \frac{1}{(iy')^n}
     \pt_{\sg}^2\pt_{\zeta'}^n
      e^{i(2^kt\sg+y'\cdot \zeta')}, 
}
and integration by parts in the right hand side 
of \eqref{eqn;ortho-30}, we see that 
{\allowdisplaybreaks   
\algn{
 \|&\phi_m\underset{(\eta)}{*} \pi_{k,j}(t,\cdot,\eta)\|_{L^1_{x'}} \\
 \simeq &\frac{2^k}{\<2^k t\>^2}
    \bigg\|
     \int_{\re}\int_{\re^{n-1}}
         \frac{1}{\<y'\>^n}e^{2^kit\sg+iy'\cdot \zeta'}
         (1-\pt_{\sg}^2)
         \Big(\sum_{|\al_1|+|\al_2|+|\al_3|\le n}C_n
          \pt_{\zeta'}^{\al_1}
          \big(2^j (i\zeta',-|\zeta'|)^{\sf T}
               (\phi_m\underset{(\eta)}{*} e^{-2^j|\zeta'|\eta})
          \big)\nn\\
     &\hskip5cm    
         \times \pt_{\zeta'}^{\al_2}
          \big(  m(2^k\sg,2^j\zeta') 
          \big)
          \pt_{\zeta'}^{\al_3} \widehat{\phi}(\zeta') 
          \Big) \widehat{\psi}(\sg)
     d\xi'd\t   
    \bigg\|_{L^1_{y'}}
 \\
  \simeq &
    \frac{2^k}{\<2^k t\>^2}
    \bigg\| \frac{1}{\<y'\>^n}
     \int_{\re}\int_{\re^{n-1}}
        e^{2^kit\s+iy'\cdot \zeta'}
         2^j(1-\pt_{\sg}^2)
         \bigg(\sum_{|\al_1|+|\al_2|+|\al_3|\le n}C_n
         \Big(
          \phi_m\underset{(\eta)}{*}
          \big(1+|\zeta'|+(2^j\eta)^{\al_1}\big)                 
           e^{-2^j|\zeta'|\eta}  \Big)\\
     &\hskip5cm    
         \times       
         \pt_{\zeta'}^{\al_2}
         \big(m(2^k\s,2^j\zeta') \big)  
         \pt_{\zeta'}^{\al_3}\widehat{\phi}(\zeta')
         \bigg) \widehat{\psi}(\sg)
     d\zeta'd\sg   
    \bigg\|_{L^1_{y'}}.   
  \eqntag \label{eqn;ortho-32}
}
}
Note that the above estimates are also valid even for $\al=\al_1+\al_2+\al_3=0$. 
Since all the functions in the integrand 
involving  $\xi'(=2^j\zeta')$ are spherically  symmetric, we employ 
integration by parts with respect to $|\zeta'|$, then 
$\bar\eta= 2^j|\zeta'|\eta$  appears. 
Namely, if we take integration by parts $|\alpha|$ times, the same estimate holds. 
There is no influence from integration 
by parts with respect to $\sg$.

To consider the effect from the convolution $\phi_m\underset{(\eta)}{*}$, 
we restrict $j\le m$ (cf. \eqref{eqn;M_1-2}). 
Applying the change of variable 
 $\tilde{\th}=2^m\th$, $\teta=2^m\eta$
 and setting 
$$
p_{\al}(\zeta',\eta,\th,\nu)
  =\Big(|\zeta'|\big(1+|\zeta'|+( 2^j(\eta-\nu\th))^{\al}\big)
    +\al\big(2^j(\eta-\nu\th)\big)^{\al-1}\Big)
$$
and noting $\int_{\re}\phi_m(\th)d\th=0$,
it follows that 
\algn{
   \Big|
     \phi_m\underset{(\eta)}{*} &
      \big(1+|\zeta'|+(2^j\eta)^{\al}\big) e^{-2^j|\zeta'|\eta}
   \Big| \\
= & \left|\int_{\re}\phi_m(\th) 
             2^j
             \Big(
             \big(1+|\zeta'|+( 2^j(\eta-\th))^{\al}\big)
             \exp\big(- 2^j|\zeta'||\eta-\th| \big)
         \right.\\
  &\hskip4.5cm \left.
          -  \big(1+|\zeta'|+( 2^j\eta)^{\al}\big)
             \exp\big(- 2^j|\zeta'||\eta| \big)
             \Big)
             d\th
      \right|
       \\
 \le &\int_{\re}
          \left|\phi_m(\th) 
             \Big( 2^j
             \int_0^1\frac{d}{d\nu}
             \big(1+|\zeta'|+( 2^j(\eta-\nu\th))^{\al}\big)
             \exp\big(- 2^j|\zeta'||\eta-\nu\th| \big)
             d\nu\Big)
           \right|
           d\th
     \\
 \le &\int_0^1 \int_{\re}
            \big|\phi_m(\th)\big| 
            2^j
           \Big( 2^j|\th||\zeta'|\big(1+|\zeta'|+( 2^j(\eta-\nu\th))^{\al}\big)
               + 2^j|\th| \al \big( 2^j(\eta-\nu\th)\big)^{\al-1}
             \Big)
      \\       
      &\hskip4.5cm      
       \times \exp\big(-2^j|\zeta'||\eta-\nu\th| \big)        
             d\th
             d\nu
      \\
  = &\int_0^1 \int_{\re}2^{2j}
              | 2^m\phi(2^m\th) |
              |\th| p_{\al}(\zeta',\eta,\th,\nu)
             \exp\big(- 2^j|\zeta'||\eta-\nu\th| \big)             
             d\th d\nu
     \\
   = &  2^j|\zeta'|^{-1}
             \int_0^1 \int_{|\bar\th|> \frac12|\bar\eta|}
             \big|\phi(2^m 2^{-j}|\zeta'|^{-1}\bar\th)\big|
              2^m  2^{-j}|\zeta'|^{-1}|\bar\th|
             \bar{p}_{\al}(\bar\eta,\bar\th)
              \exp\big(-|\bar\eta-\nu\bar\th| \big)            
             d\bar\th
             d\nu \\
     &+ 2^j|\zeta'|^{-1}
             \int_0^1 \int_{|\bar\th|\le \frac12|\bar\eta|}
             \big|\phi(2^m  2^{-j}|\zeta'|^{-1}\bar\th)\big|
                       2^m  2^{-j}|\zeta'|^{-1}|\bar\th|
              \bar{p}_{\al}(\bar\eta,\bar\th)
              \exp\big(-|\bar\eta-\nu\bar\th| \big)             
             d\bar\th
             d\nu\\
   \equiv & I+II,  \label{eqn;ortho-36} \eqntag       
 }
where we set   
\alg{ 
   &
   \bar\eta=2^j|\zeta'|\eta, 
   \quad 
   \bar\th= 2^j|\zeta'|\th, 
   \label{eqn;orth-37} \\
   &\bar{p}_{\al}(\bar\eta,\bar\th)
   =p_{\al}(\zeta',\bar{\eta},\bar{\th},\nu).
   \label{eqn;orth-37-2}
}
By $1/2\le |\zeta'|\le2$ we note that 
\eq{\label{eqn;polynomial-bound}
 \big|\bar{p}_{\al}(\bar\eta,\bar\th)\big| \exp\big(- 2^{-1}|\bar\eta-\nu\bar\th| \big) 
 \le C\Big(1+|\bar\eta-\nu\bar\th|^{\al}\Big)  
      \exp\big(- 2^{-1}|\bar\eta-\nu\bar\th| \big) 
 \le C. 
}
For the first term $I$ of \eqref{eqn;ortho-36}, 
by using the decay property of $\phi\in \mathcal{S}$, 
we know  for sufficiently large $N\in \Nt$ there exists $C_N>0$  such that 
by \eqref{eqn;polynomial-bound} that
\algn{
 I
     \le & 2^j|\zeta'|^{-1}
         \int_0^1 \int_{|\bar\th|> \frac12|\bar\eta|} 
                \frac{C_N 2^m 2^{-j}|\zeta'|^{-1}|\bar\th|}
                     {\< 2^m 2^{-j}|\zeta'|^{-1}\bar\th\>^{2N}}  
             \big|\bar{p}_{\al}(\bar\eta,\bar\th)\big|
             \exp\big(-|\bar\eta-\nu\bar\th|\big)                          
             d\bar\th
             d\nu 
 \\
  = & 2^j|\zeta'|^{-1}
         \frac{C_N (2^{m-j} |\zeta'|^{-1})^{-1}}
              {   \< 2^{m-1-j} |\zeta'|^{-1}\bar\eta\>^N}
        \int_0^1 \int_{|\bar\th|> \frac12|\bar\eta|} 
                      \frac{ (2^{m-j}|\zeta'|^{-1})^2|\bar\th|}
                           {\<2^{m-j}|\zeta'|^{-1}\bar\th\>^N} 
             \big|\bar{p}_{\al}(\bar\eta,\bar\th)\big|
             \exp\big(-|\bar\eta-\nu\bar\th|\big)                          
             d\bar\th
             d\nu
 \\
 \le& 2^j \frac{ C_N2^{-m+j} }{ \< 2^{m-1} \eta\>^N}  
            \int_{\re} 
              \frac{ (2^{m-j}|\zeta'|^{-1})^2|\bar\th|}
                   {\<2^{m-j}|\zeta'|^{-1}\bar\th\>^N} 
              d\bar\th
   \\
   \le&  \frac{C_N  2^{-m+j} 2^j|\zeta'|^{-2}}
              { \<2^m \eta\>^N}
   \le C_N  2^j\frac{ 2^{-m+j} }
           { \< 2^j \eta\>^N}, 
           \eqntag \label{eqn;ortho-38} 
}
where we used $j\le m$.
For an estimate of the second term $II$ of \eqref{eqn;ortho-36}, we use the following lemma. 
\begin{lem}\label{lem;b-bound} For  $N\in \Nt\setminus\{1\}$ and $a>0$, it holds that  
\algn{
{\ppl \int_{-a\le x\le a}}\frac{dx}{(1+|x|^2)^{N/2}}
 \le &\frac{4a}{(1+a^2)^{1/2}}
 \eqntag \label{eqn;almost-orthogonal-A}.}
 \end{lem}
 
\begin{prf}{Lemma \ref{lem;b-bound}}
For any $N\ge 1$ and $a\le 1$
\eqn{
 \int_{-a\le x\le a}\frac{dx}{(1+|x|^2)^{N/2}}\le \int_{|x|\le a}dx=2a
\le \frac{4a}{(1+a^2)^{1/2}},
}
 while for any $N\ge 2$ and $a\ge 1$, 
\eqn{
 \int_{-a\le x\le a}\frac{dx}{(1+|x|^2)^{N/2}}
 \le \int_{|x|\le a}\frac{dx}{1+|x|^2}=2\tan^{-1}a  
 \le \frac{4a}{(1+a^2)^{1/2}}.
}
This shows \eqref{eqn;almost-orthogonal-A}.
\end{prf}

\vskip1mm\noindent
{\bf Proof of Lemma \ref{lem;pressure-orthogonal2}, continued.}\ 
Under the condition $|\bar\th|\le |\bar\eta|/2$, it holds that 
           $|\bar\eta-\nu\bar\th|\ge |\bar\eta|-|\bar\th|
            \ge |\bar\eta|-|\bar\eta|/2=|\bar\eta|/2$. 
By using the above estimate, \eqref{eqn;polynomial-bound} 
and \eqref{eqn;almost-orthogonal-A}, 
by changing the integral variables 
$$
 \bar{\bar\eta}=  2^m 2^{-j}|\zeta'|^{-1}\bar\eta,\qquad 
 \bar{\bar\th} =  2^m 2^{-j}|\zeta'|^{-1} \bar\th,
$$
the second term $II$ of \eqref{eqn;ortho-36} is estimated 
as follows: 
\algn{
 II 
  =& 2^j |\zeta'|^{-1}
             \int_0^1 \int_{|\bar\th|\le \frac12|\bar\eta|}
              \frac{C_N  2^m 2^{-j}|\zeta'|^{-1}|\bar\th|}
                     {\< 2^m 2^{-j}|\zeta'|^{-1}\bar\th\>^{N+1}} 
              \big|\bar{p}_{\al}(\bar\eta,\bar\th)\big|
              \exp\big(-|\bar\eta-\nu\bar\th| \big)             
             d\bar\th
             d\nu
 \\
  =&  C_N 2^j 2^{-m+j}
             \int_0^1 \int_{|\bar\th|\le \frac12|\bar\eta|}
              \frac{( 2^m 2^{-j}|\zeta'|^{-1})^2|\bar\th|}
                   {\< 2^m 2^{-j}|\zeta'|^{-1}\bar\th\>^{N+1}} 
              \big|\bar{p}_{\al}(\bar\eta,\bar\th)\big|
              \exp\big(-|\bar\eta-\nu\bar\th| \big)             
             d\bar\th
             d\nu
 \\
 \le&  C_N2^j 2^{-m+j}
          \int_0^1 \int_{|\bar\th|\le \frac12|\bar\eta|} 
               \frac{( 2^m 2^{-j}|\zeta'|^{-1})^2|\bar\th|}
                    {\<2^{m-j}|\zeta'|^{-1}\bar\th\>^{N+1}} 
               \exp\big(- 2^{-1}|\bar\eta-\nu\bar\th|\big)                          
              d\bar\th
             d\nu
      \\
 \le& C_N2^j 2^{-m+j}\exp\big(-\frac14|\bar\eta|\big) 
        \int_{|\bar{\bar\th}|\le \frac12|\bar{\bar\eta}|} 
                  \frac{\bar{\bar\th}}
                       {\<\bar{\bar\th}\>^{N+1}} 
        d\bar{\bar\th} 
\\
   \le& C_N 2^j 2^{-m+j}\exp\big(-\frac14|\bar\eta|\big)
        \frac{8\cdot 2^{m}2^{-j}|\zeta|^{-1} \bar\eta}
             {\< 2^m2^{-j}|\zeta|^{-1}  \bar\eta\>} 
\\
   \le&   C_N \frac{ 2^j 2^{-m+j} }
              {\< 2^j \eta\>^{N-1}},
  \eqntag \label{eqn;ortho-40}
}
where  we used \eqref{eqn;orth-37} and 
$C_N$ is a constant depending on  sufficiently large $N\in \mathbb{N}$.
Hence the similar estimate as in Lemma \ref{lem;pressure-orthogonal} holds and 
we gain the decay for variables $y'$ and $t$. 
Applying \eqref{eqn;ortho-36}, \eqref{eqn;ortho-38}, \eqref{eqn;ortho-40} into 
\eqref{eqn;ortho-32} and recalling 
\eqref{eqn;ortho-22}, we see that 
\algn{
 \|\phi_m\underset{(\eta)}{*}& \pi_{k,j}(t,\cdot,\eta)\|_{L^1_{x'}} 
\\
  \simeq &
    \frac{2^k}{\<2^k t\>^2}
    \bigg\| \frac{1}{\<y'\>^n}
     \int_{\re}\int_{\re^{n-1}}
        e^{2^kit\s+iy'\cdot \zeta'}
         2^j(1-\pt_{\sg}^2)
         \bigg(\sum_{|\al_1|+|\al_2|+|\al_3|\le n}C_n
         \pt_{\zeta'}^{\al_2}
         \big(m(2^k\s,2^j\zeta')\big) \\
     &\hskip4cm    
         \times \big(
          \phi_m\underset{(\eta)}{*}
          \big(1+|\zeta'|+(2^j\eta)^{\al_1}\big)
           e^{-2^j|\zeta'|\eta}
         \big)         
         \widehat{\psi}(\sg)
         \pt_{\zeta'}^{\al_3}\widehat{\phi}(\zeta')
         \bigg)
     d\zeta'd\sg   
    \bigg\|_{L^1_{y'}}
 \\ 
 \le &
    \frac{2^k}{\<2^k t\>^2}
    \bigg\| \frac{1}{\<y'\>^n}
     \int_{2^{-1}<|\sg|<2}\int_{2^{-1}<|\zeta'|<2}
         2^j
         \big|
          \phi_m\underset{(\eta)}{*}
          \big(1+|\zeta'|+(2^j\eta)^{n}\big)
           e^{-2^j|\zeta'|\eta}
         \big|         
     d\zeta'd\sg   
    \bigg\|_{L^1_{y'}}
 \\ 
 \le &\frac{ C_N  2^j  2^{-(m-j)}  }
           { \< 2^j \eta\>^{N} }
      \frac{2^k}{\<2^k t\>^2}.   
  \eqntag \label{eqn;ortho-42}
}
From the estimate \eqref{eqn;ortho-42}, 
we conclude that \eqref{eqn;pressure-orthogonarity-T-2}  holds.
\par
(2) To see the estimate \eqref{eqn;pressure-orthogonarity-S-2}, 
we recall the low frequency restriction $\zeta$ given by 
\eqref{eqn;low-freq-zeta} and it follows 
\algn{ 
   \Big\|\sum_{k<2j} & \phi_m\underset{(\eta)}{*} \pi_{k,j}(t,\cdot,\eta)
   \Big\|_{L^1_{x'}} \\
 =&\left\|
         c_{n+1}\int_{\re}\int_{\re^{n-1}}
         e^{it\t+ix'\cdot \xi'}
         (i\xi',-|\xi'|)^{\sf T} m(\t,\xi')
         \big( \phi_m\underset{(\eta)}{*} e^{-|\xi'|\eta}
         \big)
         \sum_{k<2j}
         \widehat{\psi}(2^{-k}\t)
         \widehat{\phi}(2^{-j}\xi')
     d\xi'd\t   
    \right\|_{L^1_{x'}}
 \\
 =&\left\|
         c_{n+1}\int_{\re}\int_{\re^{n-1}}
         e^{it\t+ix'\cdot \xi'}
         (i\xi',-|\xi'|)^{\sf T} m(\t,\xi')
         \big( \phi_m\underset{(\eta)}{*} e^{-|\xi'|\eta}
         \big)
         \widehat{\zeta}(2^{-2j}\t)
         \widehat{\phi}(2^{-j}\xi')
     d\xi'd\t   
    \right\|_{L^1_{x'}}
 \\
     =&C_n2^{2j} \left\|\int_{\re}\int_{\re^{n-1}}
               e^{i2^{2j} t\sg+iy'\cdot \zeta'}
               2^j (i\zeta',-|\zeta'|)^{\sf T}
               m(2^{2j}\sg,2^j\zeta')
            \big(
                \phi_m\underset{(\eta)}{*} e^{-2^j|\zeta'|\eta}
            \big)
              \widehat{\zeta}(\sg)
              \widehat{\phi}(\zeta')
            d\zeta'\, d\sg   
         \right\|_{L^1_{y'}}.
   \eqntag \label{eqn;ortho-50}
}
Applying \eqref{eqn;ortho-5'}
and integration by parts in the right hand side 
of \eqref{eqn;ortho-50}, we see by using 
 \eqref{eqn;ortho-36},  \eqref{eqn;ortho-38} and \eqref{eqn;ortho-40} 
 again that 
\algn{
 \Big\|\sum_{k<2j}& \phi_m\underset{(\eta)}{*} \pi_{k,j}(t,\cdot,\eta)
 \Big\|_{L^1_{x'}} \\
  \simeq &
    \frac{2^{2j}}{\<2^{2j} t\>^2}
    \bigg\| \frac{1}{\<y'\>^n}
     \int_{\re}\int_{\re^{n-1}}
        e^{2^{2j} it\s+iy'\cdot \zeta'}
         2^j(1-\pt_{\sg}^2)
         \bigg(\sum_{|\al_1|+|\al_2|+|\al_3|\le n}C_n
         \pt_{\zeta'}^{\al_2}
         \big(m(2^{2j}\sg,2^j\zeta') \big) \\
     &\hskip4cm    
         \times \big(
          \phi_m\underset{(\eta)}{*}
          \big(1+|\zeta'|+(2^j\eta)^{\al_1}\big)
           e^{-2^j|\zeta'|\eta}
         \big)         
         \widehat{\zeta}(\sg)
         \pt_{\zeta'}^{\al_3}\widehat{\phi}(\zeta')
         \bigg)
     d\zeta'd\sg   
    \bigg\|_{L^1_{y'}}
 \\
 \le & \frac{ C_N  2^{j} 2^{-(m-j)} }
           {\< 2^j \eta\>^{N}}
      \frac{2^{2j}}{\<2^{2j} t\>^2}.   
  \eqntag \label{eqn;ortho-52}
}
The estimate 
\eqref{eqn;ortho-52} shows \eqref{eqn;pressure-orthogonarity-S-2}.
This completes the proof of Lemma \ref{lem;pressure-orthogonal2}.
\end{prf}

\sect{Proof of maximal $L^1$-regularity}\label{Sec;4}\par
In this section, we prove maximal $L^1$-regularity 
Theorem \ref{thm;L1MR2-b}.  The key estimate is the 
 bound for the derivative of the pressure term $\N q(t,x)$.
Indeed, once we obtain the required estimate for the 
pressure, then the estimate for the velocity directly 
follows from the  estimate for the heat equation.
Note that the velocity term can be also 
expressed by the potential as is shown in \eqref{eqn;vn}.

\subsection{Maximal regularity for the pressure} 
To show Theorem \ref{thm;L1MR2-b}, we show 
maximal $L^1$-regularity for the pressure term.
We recall the notations for the potential 
\eqref{eqn;pressre-potential}   
for the pressure $\N q$.

\begin{prop}
\label{prop;grad-pressure-bound}
Let $1<p<\infty$ and  $-1+1/p< s < 1/p$. For given data
\eqn{
  H\in \dF^{\frac12-\frac{1}{2p}}_{1,1}(\re_+;\dB^{s}_{p,1}(\re^{n-1}))
        \cap 
       {\crd L^1}(\re_+;\dB^{s+1-\frac{1}{p}}_{p,1}(\re^{n-1})),
}
there exists $C>0$ independent of $H$ 
such that the pressure part $q$ of the problem \eqref{eqn;ST2} satisfies 
the estimate 
\begin{equation}
\spl{
 \big\|\N q \big\|_{ L^1(\re_+; \dB^{s}_{p,1}(\re^{n}_+)) } 
  \le& C\Big(
         \| H \|_{\dF^{\frac12-\frac{1}{2p}}_{1,1}(\re_+;\dB^s_{p,1}(\re^{n-1}))}
        +\| H \|_{{\crd L^1}(\re_+;\dB^{s+1-\frac{1}{p}}_{p,1}(\re^{n-1}))}
       \Big). 
 }
  \label{eqn;grad-pressure-bound}
\end{equation}
\end{prop}
\vskip2mm
To show the pressure estimate \eqref{eqn;grad-pressure-bound},
we use the potential expression  
$\pi(t,x',\eta)$ in \eqref{eqn;Stokes-n-pressure-grad} 
and the Littlewood--Paley decomposition of unity \eqref{eqn;direct-sum-L-P};
 \begin{align} 
  \overline{\Phi_m}&\underset{(x',\eta)}{*}\big(\pi(t,x',\eta)\big) 
  \notag \\
   =& c_{n+1}\overline{\Phi_m}(x',\eta)
                               \underset{(x',\eta)}{*}
            \int_{\re}\int_{\re^{n-1}}  e^{it\t +ix'\cdot \xi'}
               (i\xi', -|\xi'|)^{\sf T} m(\t,\xi')  e^{-|\xi'|\eta}  
           d\xi'd\t\nn\\   
    =&   \zeta_{m-1}(\eta)\underset{(\eta)}{*}
             c_{n+1}\int_{\re^{n-1}}\int_{\re} 
                e^{it\t +ix'\cdot \xi'}  
                \widehat{\phi_m}(|\xi'|)
                (i\xi', -|\xi'|)^{\sf T} m(\t,\xi')  e^{-|\xi'\eta}  
            d\xi'd\t  \nn\\
    &+   \phi_{m}(\eta) \underset{(\eta)}{*}   
         c_{n+1}\int_{\re}\int_{\re^{n-1}}\int_{\re} 
                e^{it\t +ix'\cdot \xi'} \widehat{\zeta_{m}}(|\xi'|)
                (i\xi', -|\xi'|)^{\sf T} m(\t,\xi') e^{-|\xi'|\eta}
                d\xi'd\t \nn\\
  \equiv &\zeta_{m-1}(\eta)\underset{(\eta)}{*}  
           \phi_m(x')\underset{(x')}{*}\pi(t,x',\eta)
         +\phi_{m}(\eta) \underset{(\eta)}{*}  
          \zeta_m(x')\underset{(x')}{*}\pi(t,x',\eta).
           \label{eqn;potential-besov-3} 
\end{align}

\noindent
Concerning the first term of the right-hand side of \eqref{eqn;potential-besov-3},  
we estimate that the convolution with $\zeta_{m-1}(\eta)$ can be treated by 
the Hausdorff--Young inequality in $\eta$-variable.
Note that the potential $\pi(t,x',\eta)$ has the even extension in  
$\eta\in \re$ and hence the $L^p(\re^{n}_+)$ norm of the 
term is estimated as follows: 
\eq{\label{eqn;first-zeta-exclude}
 \spl{
 \Big\|\zeta_{m-1}\underset{(\eta)}{*} &
       \Big(\phi_m(x')\underset{(x')}{*}\pi(t,x',\eta) \Big)
 \Big\|_{L^p(\re_{+,\eta};L^p(\re^{n-1}))} \\
 \le & 
 \Big\|\zeta_{m-1}\underset{(\eta)}{*} 
       \Big(\phi_m(x')\underset{(x')}{*}\pi(t,x',\eta) \Big)
 \Big\|_{L^p(\re_{+,\eta};L^p({\gr \re^{n-1}}))} \\
 \le & \|\zeta_{m-1}\|_{L^1(\re_{+,\eta})}
       \big\|\phi_m(x')\underset{(x')}{*}\pi(t,x',\eta) 
       \big\|_{L^p(\re_{+,\eta};L^p({\gr \re^{n-1}}))} \\
 \le &C\big\|\phi_m(x')\underset{(x')}{*}\pi(t,x',\eta) 
       \big\|_{L^p(\re_{+,\eta};L^p({\gr \re^{n-1}}))}
}}
and we apply Lemma \ref{lem;pressure-orthogonal}.   
Concerning the second term of the right-hand side of \eqref{eqn;potential-besov-3}, 
the number of overlapping  supports of the kernel 
$\zeta_m(x')\underset{(x')}{*}\phi_j(x')$ is infinite, 
i.e., $m$ and $j$ run independently.  
We apply the almost orthogonality of the second type 
stated in Lemma \ref{lem;pressure-orthogonal2}.

\vskip2mm
\begin{prf}{Proposition \ref{prop;grad-pressure-bound}}
Let us recall that the boundary data 
 $H(t,x')=(H'(t,x'),H_n(t.x'))$ is extended into $t<0$ by the zero extension.
By \eqref{eqn;potential-besov-3}, we divide the term into two terms. 
{\allowdisplaybreaks
\algn{
  \N q&(t,x',x_n)  \\
  =& c_{n+1}\iint_{\re^{n}}
               e^{i t\t+ix'\cdot \xi'}(i\xi', -|\xi'|)^{\sf T}
               \Big(m'(\t,\xi')\cdot \widehat{H'}
                    +m_n(\t,\xi') \widehat{H_n}
               \Big) e^{-|\xi'|x_n }
               \sum_{k\in \Z}\sum_{j\in\Z}
               \widehat{\phi_k}(\t)\widehat{\phi_j}(\xi)
                d\t d\xi' 
\\
 \equiv 
   & \sum_{k\in \Z}\sum_{j\in\Z}
     \Big(
       {\pi'}_{k,j}\underset{(t,x')}{\cdot *} 
        \Big(\widetilde{\psi_k} \underset{(t)}{*}   
                  \widetilde{\phi_j}\underset{(x')}{*}H'\Big)
   +   {\pi_n}_{k,j}\underset{(t,x')}{*} 
        \Big(\widetilde{\psi_k} \underset{(t)}{*}
                  \widetilde{\phi_j}\underset{(x')}{*}H_n\Big)
    \Big). 
}
} 
Then observing the estimate \eqref{eqn;first-zeta-exclude},
we see
{\allowdisplaybreaks
\begin{align*}
 \|\N &  q\|_{L^1(\re_+;\dB^{s}_{p,1}(\re^{n}_+))} 
\\
  \le &C \bigg\|\sum_{m\in \Z} 2^{sm}
         \bigg\|\Big(\int_{\re}
                \Big| \zeta_{m-1}(\teta)\underset{(\teta)}{*}  
                   \phi_m(x')\underset{(x')}{*}
                    \pi(t,x',\eta)
                    \underset{(t,x')}{\cdot *} H
                \Big|^p  d\teta 
                \Big)^{1/p}
         \bigg\|_{L^p({\gr \re^{n-1}})}
         \bigg\|_{L^1_t(\re_+)}
\\
       &+ C \bigg\|\sum_{m\in \Z}2^{sm}
          \bigg\|\Big(\int_{\re}
                 \Big| \phi_m(\teta)\underset{(\teta)}{*}  
                       \zeta_{m}(x')\underset{(x')}{*}
                       \pi(t,x',\eta)
                      \underset{(t,x')}{\cdot *} H
                 \Big|^p  d\teta 
                 \Big)^{1/p}
          \bigg\|_{L^p({\gr \re^{n-1}})}
          \bigg\|_{L^1_t(\re_+)}
 \\
  \le &C \bigg\|\sum_{m\in \Z}2^{sm}
         \bigg(\int_{\re} 
         \Big\| \phi_m(x')\underset{(x')}{*}
                \pi(t,x',\eta) 
                \underset{(t,x')}{\cdot *}            
                 H(t,x')   
         \Big\|_{L^p({\gr \re^{n-1}})}^p 
         d\tilde{\eta} 
         \bigg)^{1/p}
         \bigg\|_{L^1_t(\re_+)}
 \\
        &+C \bigg\|\sum_{m\in \Z}2^{sm}
               \bigg(\int_{\re}
               \Big\|
                  \Big(                 
                   \phi_{m}(\teta)\underset{(\teta)}{*}
                   \zeta_{m}(x')\underset{(x')}{*}
                   \pi(t,x',\eta) 
                  \Big) 
                  \underset{(t,x')}{\cdot *} 
                   H(t,x')   
               \Big\|_{L^p(\re^{n-1}_{x'})}^p  d\tilde{\eta} 
               \bigg)^{1/p}
       \bigg\|_{L^1_t(\re_+)}
\\
\equiv &\|P_1(t)\|_{L^1_t(\re_+)}+\|P_2(t)\|_{L^1_t(\re_+)},
\eqntag \label{eqn;potential-besov-9}
\end{align*}
}
where we denote the inner product-convolution  
by \eqref{eqn;convolution-innerprod}.
Noting that the data $H$ is divided into the the time-dominated region 
$k\ge 2j$ and the space-dominated 
region $k<2j$, respectively, as
\eq{\label{eqn;H-decomp0}
 H(t,x')
  =\sum_{k\in\Z} \sum_{2j\le k} H_{k,j}(t,x')  
   +\sum_{k\in\Z} \sum_{2j>k} H_{k,j}(t,x'),
}
where we set 
\eq{ \label{eqn;H-decomp}
 \spl{
 & H_{k,j}(t,x')=\widetilde{\psi_k}(t)\underset{(t)}{*}
                 \widetilde{\phi_j}(x')\underset{(x')}{*}H(t,x'),\\ 
 & H_{j}(t,x')= \widetilde{\phi_j}(x')\underset{(x')}{*}H(t,x'),  
}}
where we use $\widetilde{\phi_j}=\phi_{j-1}+\phi_j+\phi_{j+1}$ and 
$\widetilde{\psi_k}$ similar arrangement.
Then applying $\widetilde{\phi_j}\underset{(x')}{*}\phi_j=\phi_j$
and $\widetilde{\psi_k}\underset{(t)}{*}\psi_k=\psi_k$,
and Proposition \ref{prop;retraction-coretraction},
we  divide $P_1(t)$ into $L_1(t)$ and $L_2(t)$ to have the following: 
\algn{ \label{eqn;potential-besov-10}
 P_1(t) 
 \le& C \sum_{m\in \Z}2^{sm}\bigg\|
        \Big\| \phi_m(x')\underset{(x')}{*}
          \pi(t,x',\eta)
          \underset{(t,x')}{\cdot *}                  
           \sum_{j\in\Z}\sum_{k\ge 2j} H_{k,j}(t,x') 
        \Big\|_{L^p(\re^{n-1}_{x'})}
       \bigg\|_{L^p(\re_{+,\eta})} \\
    &+ C \sum_{m\in \Z}2^{sm}\bigg\|
          \Big\| 
            \phi_m(x')\underset{(x')}{*}
            \pi(t,x',\eta)\underset{(t,x')}{\cdot *}
             \sum_{j\in\Z}\sum_{k<2j}
             H_{k,j}(t,x')   
         \Big\|_{L^p(\re^{n-1}_{x'})}
         \bigg\|_{L^p(\re_{+,\eta})}\\
 =& C \sum_{m\in \Z}2^{sm}
       \bigg( \int_{\re_+}
        \Big\|
          \sum_{k\ge 2m}
          \pi_{k,m}\underset{(t,x')}{\cdot *}        
           H_{k,m}(t,x')  
        \Big\|_{L^p(\re^{n-1}_{x'})}^p
        d\eta
      \bigg)^{1/p} \\
    &+ C \sum_{m\in \Z}2^{sm}
        \bigg(\int_{\re_+}
         \Big\| 
          \sum_{k<2m} 
          \pi_{k,m}\underset{(t,x')}{\cdot *}
           H_{m}(t,x')  
          \big) 
         \Big\|_{L^p(\re^{n-1}_{x'})}^p
         d\eta
        \bigg)^{1/p}\\
  \equiv & L_1(t) + L_2(t),  \eqntag 
}
where $\{\pi_{k,m}\}_{k,m}$ are defined in \eqref{eqn;pressre-potential}.
For the time dominated part $L_1$, since $k\ge 2m$, 
we apply the almost orthogonality estimate 
\eqref{eqn;pressure-orthogonarity-T} in 
Lemma \ref{lem;pressure-orthogonal}, 
by using the change of valuable $2^m\eta=\bar{\eta}$ 
it holds that 
\begin{align*}
 \big\| L_1 \big\|_{L^1(\re_+)}
 \le & C\bigg\|\sum_{m\in\Z} 2^{sm}
       \bigg(\int_{\re_+}
        \Big\{  \sum_{k\ge 2m} 
             \int_{\re}
             \Big\|\pi_{k,m}(t-s,x',\eta)
             \Big\|_{L^1_{x'}}
             \Big\|H_{k,m}(s,x')  \Big\|_{L^p_{x'}}
             ds
        \Big\}^p  
         d\eta
     \bigg)^{1/p}\bigg\|_{L^1_t(\re_+)}
\\
  \le& C\bigg\| 
         \sum_{m\in \Z}2^{sm}
         \bigg(\int_{\re_+}
          \Big\{ \sum_{k\ge 2m}
               \big(2^{m}(1+(2^{m}\eta)^{n+2}\big)
                    e^{-(2^{m-1}\eta)}
        \\
    &\hskip4cm \times
             \int_{\re}
              \frac{2^k}{\<2^k(t-s)\>^2}
              \Big\|
                  \widetilde{\psi_k}\underset{(t)}{*}
                  \widetilde{\phi_m}\underset{(x')}{*}H(s,\cdot)
              \Big\|_{L^p_{x'}} 
              ds 
          \Big\}^pd\eta
          \bigg)^{1/p}
         \bigg\|_{L^1_t(\re_+)} 
 \\
    = & C\bigg\| 
         \sum_{m\in \Z}  2^{sm}
         \bigg(
         \Big\{\sum_{k\ge 2m} 2^{m}
             \int_{\re}
              \frac{2^k}{\<2^k(t-s)\>^2}
              \Big\|
                   \widetilde{\psi_k}\underset{(t)}{*}
                   H_m(s,\cdot)
              \Big\|_{L^p_{x'}} 
              ds 
          \Big\}^p             
         \\
    &\hskip4cm \times
           {\ppl 2^{-m} }
           \Big(\int_{\re_+}         
               \Big((1+\bet^{n+2})e^{-\bet/2}
               \Big)^p                   
               d\bet
          \Big) 
          \bigg)^{1/p}
         \bigg\|_{L^1_t(\re_+)} 
\\
  \le &C\bigg\| 
         \sum_{m\in \Z}  
            2^{sm}
            2^{(1-\frac{1}{p})m}
         \sum_{k\ge 2m} 
          \int_{\re}
            \frac{2^k}{\<2^k(t-s)\>^2}
            \Big\|
              \widetilde{\psi_k}\underset{(t)}{*}H_m(s,\cdot)
            \Big\|_{L^p_{x'}} 
          ds 
         \bigg\|_{L^1_t(\re_+)} 
\\
   \le &C \sum_{k\in\Z} 2^{\frac{k}{2}(1-\frac{1}{p})}
          \bigg\|
          \int_{\re}
            \frac{2^k}{\<2^k(t-s)\>^2}
            \sum_{m\in \Z}  2^{sm}
            \Big\|
             \widetilde{\psi_k}\underset{(t)}{*}H_m(s,\cdot)
            \Big\|_{L^p_{x'}} 
          ds 
         \bigg\|_{L^1_t(\re_+)} 
\\
   \le &C \sum_{k\in\Z} 2^{\frac{k}{2}(1-\frac{1}{p})}
          \bigg\|
           \sum_{m\in \Z}  2^{sm}
            \Big\|
             \widetilde{\psi_k}\underset{(t)}{*}
             {\ppl H_m(t,\cdot)}
            \Big\|_{L^p_{x'}} 
         \bigg\|_{L^1_t(\re_+)} 
\\
   \le &C\bigg\| \sum_{k\in\Z} 2^{\frac{k}{2}(1-\frac{1}{p})} 
            \Big\|
             \psi_k\underset{(t)}{*}{\ppl H_m(t,\cdot)}
            \Big\|_{\dB^{s}_{p,1}(\re^{n-1})}  
         \bigg\|_{L^1_t(\re_+)} 
\\
  \le&C\big\| H 
       \big\|_{\dF^{\frac12-\frac{1}{2p}}_{1,1}(\re_+;\dB^{s}_{p,1}(\re^{n-1}))}.
    \eqntag \label{eqn;L_1-2}  
\end{align*}

On the other hand, when $k<2m$, 
For the space dominated part $L_2$,  
applying the almost orthogonality estimate 
\eqref{eqn;pressure-orthogonarity-S} in 
 Lemma \ref{lem;pressure-orthogonal} with
using the Minkowski inequality, the Hausdorff--Young inequality, 
we obtain 
{\allowdisplaybreaks
\begin{align*}
 \big\|   L_2\|_{L^1_t}  
  \le & \bigg\| 
        \sum_{m\in\Z} 2^{sm}
       \bigg(\int_{\re_+}
        \Big\{ 
             \int_{\re}
             \Big\| 
              \sum_{k<2m}\pi_{k,m}(t-s,x',\eta)\Big\|_{L^1_{x'}}
             \big\|H_{m}(s,\cdot)  \big\|_{L^p_{x'}}
             ds
        \Big\}^p  
         d\eta
     \bigg)^{1/p}
        \bigg\|_{L^1_t(\re_+)}
 \\
  \le& C\bigg\| 
         \sum_{m\in \Z} 2^{sm}
         \bigg(\int_{\re_+}
           \Big\{ 
               \big(2^{m}(1+(2^{m}\eta)^{n+2}) 
               e^{-(2^{m-1}\eta)}\big)
          \nn\\
       &\hskip4cm
              \times 
              \int_{\re}
              \frac{2^{2m}}{\<2^{2m}(t-s)\>^2}
              \big\|H_{m}(s,\cdot)  \big\|_{L^p_{x'}}
              ds 
          \Big\}^p
           d\eta
          \bigg)^{1/p}
         \bigg\|_{L^1_t(\re_+)} 
  \\
   \le& C\bigg\| 
         \sum_{m\in \Z} 2^{sm}
         \bigg( 2^{mp(1-\frac{1}{p})}
             \int_{\re_+}
             \Big((1+\bet^{n+2})
                  e^{-(2^{-1}\bet)}
             \Big)^p
             d\bet \\
      &\hskip7cm       
          \times\Big\{
              \int_{\re}
              \frac{2^{2m}}{\<2^{2m}(t-s)\>^2}
              \big\| H_{m}(s,\cdot)  \big\|_{L^p_{x'}}
              ds 
          \Big\}^p
          \bigg)^{1/p}
         \bigg\|_{L^1_t(\re_+)} 
 \\        
    \le & C 
         \sum_{m\in \Z}
          2^{(s+1-\frac{1}{p})m}    
             \Big\|\frac{2^{2m}}{\<2^{2m}t\>^2}\Big\|_{L_t^1(\re)}                         
             \bigg\|
              \big\|H_{m}(t,\cdot)  \big\|_{L^p_{x'}}
             \bigg\|_{L^1_t(\re_+)} 
 \\
  \le & C\big\|H\big\|_{ L^1 (\re_+;\dB^{s+1-\frac{1}{p}}_{p,1}(\re^{n-1}) )}.     
   \eqntag \label{eqn;L_2-2}
\end{align*}
}

In the same way for $P_1(t)$, we decompose $P_2(t)$ 
as a space-dominated region and a time-dominated region. 
\algn{
 P_2(t)
    \le& C \sum_{m\in \Z}2^{sm}
      \bigg\| 
         \Big\| \phi_m(\eta)\underset{(\eta)}{*}
                \zeta_m(x')\underset{(x')}{*} 
                \pi(t,x',\eta)
                \underset{(t,x')}{\cdot *}        
                \sum_{k\in\Z}\sum_{2j\le k} H_{k,j}(t,x')  
        \Big\|_{L^p(\re^{n-1}_{x'})}
      \bigg\|_{L^p(\re_{+,\eta})} \\
    &+C \sum_{m\in \Z}2^{sm}
       \bigg\|
         \Big\| 
           \phi_m(\eta)\underset{(\eta)}{*} 
           \zeta_m(x')\underset{(x')}{*} 
                \pi(t,x',\eta) \underset{(t,x')}{\cdot *} 
           \sum_{k\in\Z}\sum_{2j>k}
             H_{k,j}(t,x')     
         \Big\|_{L^p(\re^{n-1}_{x'})}
       \bigg\|_{L^p(\re_{+,\eta})}
\\
 \le& C \sum_{m\in \Z}2^{sm}
       \bigg( \int_{\re_+}  
        \Big\|  
          \sum_{k\in\Z}\sum_{2j\le \min(2m,k)}    
          \phi_m(\eta)\underset{(\eta)}{*}
          \pi_{k,j}\underset{(t,x')}{\cdot *}        
          H_{k,j}(t,x')  
        \Big\|_{L^p(\re^{n-1}_{x'})}^p
        d\eta
      \bigg)^{1/p} \\
    &+ C \sum_{m\in \Z}2^{sm}
        \bigg(\int_{\re_+}
         \Big\| 
           \sum_{k\in\Z}\sum_{k<2j\le2m} 
           \phi_m(\eta)\underset{(\eta)}{*}
           \pi_{k,j}\underset{(t,x')}{\cdot *}
           H_{j}(t,x')  
           \big) 
         \Big\|_{L^p(\re^{n-1}_{x'})}^p
         d\eta
        \bigg)^{1/p}
 \\
  \equiv & M_1(t) + M_2(t).
   \eqntag \label{eqn;potential-besov-11}
}
For the time dominated part $ M_1$, 
Setting $h_j\equiv \widetilde{\phi_j}*h$, 
using the Minkowski inequality and the Hausdorff--Young inequality, 
and also using \eqref{eqn;pressure-orthogonarity-T-2} in 
Lemma \ref{lem;pressure-orthogonal2} (1)
(the second almost orthogonality), 
we have
\algn{
 \big\| &M_1\big\|_{L^1_t(\re_+)}  \\
  \le &\bigg\|\sum_{m\in\Z} 2^{sm}
       \Big(\int_{\re_+}
        \Big\{ \sum_{k\in\Z}\sum_{2j\le \min(2m,k)}
             \int_{\re}
             \Big\|\phi_m(\eta)\underset{(\eta)}{*}
                   \pi_{k,j}(t-s,x',\eta)
             \Big\|_{L^1_{x'}}\\
   &\hskip8cm \times              
             \Big\|\psi_k\underset{(s)}{*} H_j(s,x')  
             \Big\|_{L^p_{x'}}
             ds
        \Big\}^p  d\eta
     \Big)^{1/p}
     \bigg\|_{L^1_t(\re_+)}
\\
  \le&C\bigg\|\sum_{m\in\Z} 2^{sm} \bigg(\int_{\re_+}            
             \bigg\{ 
              \sum_{k\in\Z}\sum_{2j\le \min(2m,k)}               
              \frac{2^j  2^{-(m-j)} }{\< 2^j \eta\>^N}
               \\
     &\hskip6cm \times             
              \int_{\re}
              \frac{2^k}{\<2^k(t-s)\>^2}
              \Big\|
               \psi_k \underset{(s)}{*} H_j(s,\cdot)
              \Big\|_{L^p_{x'}}
             ds 
             \bigg\}^pd\eta \bigg)^{1/p} 
     \bigg\|_{L^1_t(\re_+)} 
\\
  =&C\bigg\|\sum_{m\in\Z} 2^{sm}            
              \sum_{k\in\Z}\sum_{2j\le \min(2m,k)}                           
              \\
     &\hskip4cm \times 2^{j}   2^{-(m-j)}          
              \int_{\re}
              \frac{2^k}{\<2^k(t-s)\>^2}
              \Big\|
               \psi_k \underset{(s)}{*} H_j(s,\cdot) 
              \Big\|_{L^p_{x'}}
             ds         
          {\ppl\bigg(
             \int_{\re_+}   
             \frac{1}{\<\bet\>^{pN}}              
              2^{-j}
             d\bet 
             \bigg)^{1/p} 
           }
     \bigg\|_{L^1_t(\re_+)} 
   \\
    =&C\bigg\|\sum_{m\in\Z} 2^{sm} 
              \sum_{k\in\Z}\sum_{2j\le\min(2m,k)}               
               2^{j}  2^{-(m-j)}   2^{-\frac{j}{p}} 
              \int_{\re}
              \frac{2^k}{\<2^k(t-s)\>^2}
              \Big\|
               \psi_k \underset{(s)}{*} H_j(s,\cdot)
              \Big\|_{L^p_{x'}}
              ds 
     \bigg\|_{L^1_t(\re_+)} 
   \\
     {\ppl \le} &C\bigg\| 
              \sum_{k\in\Z}\sum_{2j\le k}   
              \sum_{m\ge j} 2^{sm}  2^{-(m-j)}            
                            2^{(1-\frac{1}{p})j} 
              \int_{\re}
              \frac{2^k}{\<2^k(t-s)\>^2}
              \Big\|
               \psi_k \underset{(s)}{*} H_j(s,\cdot)
              \Big\|_{L^p_{x'}}
              ds 
     \bigg\|_{L^1_t(\re_+)} 
   \\
   &\text{(setting $m=m'+j$ and changing $m\to m'$)}
   \\
   \le &C\bigg\| 
              \sum_{k\in\Z}\sum_{2j\le k}               
              \sum_{m'\ge 0} 
                 2^{s(m'+j)}
                 2^{- m'} 
                 2^{(1-\frac 1p)j}                          
              \int_{\re}
              \frac{2^k}{\<2^k(t-s)\>^2}
              \Big\|
              \psi_k \underset{(s)}{*} H_j(s,\cdot)
              \Big\|_{L^p_{x'}}
              ds 
     \bigg\|_{L^1_t(\re_+)} 
\\
      \le &C \sum_{k\in\Z} \sum_{2j\le k}       
              2^{(1-\frac{1}{p})j} 2^{sj}      
             \bigg\| \int_{\re}
              \frac{2^k}{\<2^k(t-s)\>^2}
              \Big\|
              \psi_k \underset{(s)}{*} H_j(s,\cdot) 
              \Big\|_{L^p_{x'}}
              ds 
     \bigg\|_{L^1_t(\re_+)} 
    \\
      \le &C\bigg\|
            \sum_{k\in\Z}2^{(1-\frac{1}{p})\frac{k}{2}}  
            \sum_{2j\le k}  2^{sj}       
              \Big\|
               \psi_k \underset{(s)}{*} H_j(s,\cdot)
              \Big\|_{L^p_{x'}}
     \bigg\|_{L^1_t(\re_+)} 
   \\
   \le &C\big\|H\big\|_{\dF^{\frac12-\frac{1}{2p}}_{1,1}(\re_+;\dB_{p,1}^{s}(\re^{n-1}))}, 
   \label{eqn;M_1-2} \eqntag
}
where we use  $s< 1$ for the convergence of the 4th line from the bottom.

The space dominated part $ M_2$ is estimated in the similar way as $ M_1$. 
We apply the Minkowski inequality and the Hausdorff--Young inequality
\begin{align*}
    M_2(t)
    \le &\sum_{m\in\Z} 2^{sm}
       \bigg(\int_{\re_+}
       \Big\{\sum_{j\le m}
             \int_{\re}
             \Big\|\phi_m(\eta)\underset{(\eta)}{*}
                   \sum_{k<2j} \pi_{k,j}(t-s,x',\eta)
             \Big\|_{ L^1_{x'} } 
             \big\| H_j(s,x') \big\|_{L^p_{x'}}
             ds
        \Big\}^p  d\eta
     \bigg)^{1/p}.
\end{align*}
Using the almost orthogonality 
\eqref{eqn;pressure-orthogonarity-S-2} 
in Lemma \ref{lem;pressure-orthogonal2} (2) 
for $k< 2j$ we have
\algn{
 \big\| M_2  \big\|_{L^1(\re_+)} 
  \le&C\bigg\|\sum_{m\in\Z} 2^{sm}
              \bigg(\int_{\re_+}            
              \bigg\{ 
               \sum_{j\le m}              
               \frac{C_N 2^{j} 2^{-(m-j)} }{\< 2^j \eta\>^N}
               \int_{\re}
              \frac{2^{2j}}{\<2^{2j}(t-s)\>^2}
              \Big\| H_j(s)\Big\|_{L^p_{x'}}
              ds 
             \bigg\}^pd\eta \bigg)^{1/p} 
     \bigg\|_{L^1_t(\re_+)} 
\\ 
  \le&C\bigg\|\sum_{m\in\Z}  2^{sm}        
              \sum_{j\le m}  2^{j}  2^{-(m-j)}
                \int_{\re}
               \frac{2^{2j}}{\<2^{2j}(t-s)\>^2}
                \Big\| H_j(s)\Big\|_{L^p_{x'}}
                ds 
            {\ppl
              \bigg(
              \int_{\re_+} 
              \frac{1}{\< \widetilde\eta\>^{pN}}
               2^{-j} d\widetilde\eta 
             \bigg)^{1/p}
            } 
     \bigg\|_{L^1_t(\re_+)} 
\\
  \le&C\bigg\|\sum_{m\in\Z} 2^{sm}          
              \sum_{j\le m} 
               2^{j} 
               2^{-\frac{j}{p}}
               2^{-(m-j)} 
              \Big(
               \int_{\re}
              \frac{2^{2j}}{\<2^{2j}(t-s)\>^2}
               \Big\| H_j(s)\Big\|_{L^p_{x'}}
               ds 
              \Big)
     \bigg\|_{L^1_t(\re_+)} 
\\
  \le&C\sum_{m\in\Z}         
       \sum_{j\le m}     
         2^{-(m-j)} 2^{sm}
         2^{(1-\frac{1}{p})j}
     \bigg\| \int_{\re}
             \frac{2^{2j}}{\<2^{2j}(t-s)\>^2}
             \Big\| H_j(s)\Big\|_{L^p_{x'}}
             ds 
     \bigg\|_{L^1_t(\re_+)} 
\\
  \le&C\bigg\|\sum_{j\in\Z}  2^{(s+1-\frac{1}{p})j}       
              \sum_{m'\ge 0}    
               2^{-m'+sm'}                
               \Big\|\phi_j\underset{(x')}{*} H(t)\Big\|_{L^p_{x'}}
       \bigg\|_{L^1_t(\re_+)} 
\\
  =&C \big\| H 
      \big\|_{L^1(\re_+;\dB^{s+1-\frac{1}{p}}_{p,1}(\re^{n-1}))}, 
  \eqntag\label{eqn;M_2-3}   
}
where we use $m-j=m'\ge 0$,  $m\to m'$ in the third line from the bottom, 
and  $s< 1$ to obtain the convergence of  $\sum_{m'\ge 0}2^{-m'+sm'}$. 
In the last line, we enter the $t$-integral in the sum of of $\ell'$, 
and delete the convolution with respect to $t$ by using the the sum of of $\ell'$. 
Combining \eqref{eqn;potential-besov-9}, 
\eqref{eqn;potential-besov-10}--\eqref{eqn;M_2-3},  
we obtain  the estimate \eqref{eqn;grad-pressure-bound}.
The restriction on the regularity exponent $s$ stems from 
the {\crd structure} of the homogeneous Besov space stated in 
Propositions \ref{prop;duality}--\ref{prop;derivative}.

This complete the proof.
\end{prf}
\vskip2mm

The following estimate is required for showing maximal regularity 
for the velocity part of the Stokes equation.

\vskip2mm
\begin{prop}\label{prop;pressure-bound} 
Let $1\le p< \infty$ and  $s\in \re$. Given boundary data
$$
H\in \dF^{\frac12-\frac{1}{2p}}_{1,1}(\re_+;\dB^s_{p,1}(\re^{n-1})),
$$ 
let $q$ be the pressure term defined by 
\eqref{eqn;Stokes-n-pressure}. Then there exists a constant 
$C>0$ such that the following estimates hold:
\begin{align}
 &\big\|q|_{x_n=0}\big\|_{\dF^{\frac12-\frac{1}{2p}}_{1,1}(\re_+;\dB^s_{p,1}(\re^{n-1}))}
  \le C\big\|H
        \big\|_{\dF^{\frac12-\frac{1}{2p}}_{1,1}(\re_+;\dB^s_{p,1}(\re^{n-1}))}. 
  \label{eqn;pressure-bound-1}
\end{align}
\end{prop}

\begin{prf}{Proposition \ref{prop;pressure-bound}}
Let $\{\psi_k\}_{k\in \Z}$ and $\{\phi_j\}_{j\in \Z}$
be the Littlewood--Paley dyadic decomposition of the unity in 
$t\in \re$ and  $x'\in \re^{n-1}$ variables, respectively. 
For simplicity, we assume that $q\in \S_0(\re^{n-1})$ 
and show the estimates \eqref{eqn;pressure-bound-1}. 
The results follows by 
the density $\S_{0}(\re^{n-1})\subset \dB^s_{p,1}(\re^{n-1})$, 
where $\S_{0}(\re^{n-1})$
denotes the rapidly decreasing functions 
with vanishing at the origin of their Fourier images.
Then the resulting estimates follows from the following bounds.
\alg{ \label{eqn;pressure-key}
 &\Big\| \| \psi_k\underset{(t)}{*}\phi_j\underset{(x')}{*}q\big|_{x_n=0}
        \|_{L^p(\re^{n-1})}
  \Big\|_{L^1_t(\re_+)}
 \le C  \Big\| \| \psi_k\underset{(t)}{*}\phi_j\underset{(x')}{*}H
               \|_{L^p(\re^{n-1})}
        \Big\|_{L^1_t(\re_+)} .
}
Indeed, admitting the above estimate \eqref{eqn;pressure-key},
the Minkowski inequality yields

\algn{
 \big\|q\big|_{x_n=0}\big\|_{\dF^{\frac12-\frac{1}{2p}}_{1,1}(\re_+;\dB^s_{p,1}(\re^{n-1}))}
 \le &C\sum_{k\in \Z} 2^{(\frac12-\frac{1}{2p})k } 
           \sum_{j\in \Z} 2^{sj}
           \left\|            
              \| \psi_k\underset{(t)}{*}\phi_j\underset{(x')}{*}q\big|_{x_n=0}
              \|_{L^p(\re^{n-1})}    
           \right\|_{L^1_t(\re_+)} \\
 \le &C\sum_{k\in \Z} 2^{(\frac12-\frac{1}{2p})k } 
           \sum_{j\in \Z} 2^{sj}
           \left\|            
              \| \psi_k\underset{(t)}{*}\phi_j\underset{(x')}{*}H
              \|_{L^p(\re^{n-1})}    
           \right\|_{L^1_t(\re_+)} \\
 \le&C \left\| \sum_{k\in \Z} 2^{(\frac12-\frac{1}{2p})k } 
           \sum_{j\in \Z} 2^{sj}
          \| \psi_k\underset{(t)}{*}\phi_j\underset{(x')}{*} H
           \|_{L^p(\re^{n-1})}    
      \right\|_{L^1_t(\re_+)} \\
 \le& C \big\|H
        \big\|_{\dF^{\frac12-\frac{1}{2p}}_{1,1}(\re_+;\dB^s_{p,1}(\re^{n-1}))}, 
}
which implies \eqref{eqn;pressure-bound-1}.

To see \eqref{eqn;pressure-key}, from \eqref{eqn;Stokes-n-pressure}, it follows 
\eq{\label{eqn;Stokes-n-pressure-0}
\spl{
\psi_k\underset{(t)}{*}&\phi_j\underset{(x')}{*}q(t,x',x_n)\big|_{x_n=0}
\\
  =&c_{n+1}\iint_{\re^{n}}
               e^{it\t+ix'\cdot \xi'}
               \bigg\{
                \frac{ |\xi'|+B}{D(\t, \xi')}
                 \Big(
                  2B(i\xi'\cdot \widehat{H}')
                  -(|\xi'|^2+B^2) \widehat{H_n}
                 \Big)
               \bigg\}
               \widehat{\psi_k}(\t)
               \widehat{\phi_j}(\xi')
               d\t d\xi'
  } 
}
and the support of the symbol on the right hand side is
in an annulus domain and hence there is no singular point in both 
$\t$, $|\xi'|$- variables and it gives a smooth symbol.
Besides, for the time-like region $k\ge 2j$, by Lemma \ref{lem;B-bound},
\eqref{eqn;B-bound-time}, \eqref{eqn;ortho-2} and 
\eqref{eqn;ortho-3}  implies that
\eqn{
\spl{
 &\left|\frac{2^j|\zeta'|+\sqrt{2^ki\sg+2^{2j}|\zeta'|^2}}
             {D(2^k\sg,2^j\zeta')} 
                 \Big(
                  -2\cdot2^j|\zeta'| \sqrt{2^ki\sg+2^{2j}|\zeta'|^2}
                  i\frac{\zeta'}{|\zeta'|}
                  \Big)
  \right| 
  = \; O(1).
}}
Analogously for the space-like region, we see from 
\eqref{eqn;B-bound-space} that
\eqn{
 \spl{
 &\left|\frac{2^j|\zeta'|+\sqrt{2^ki\sg+2^{2j}|\zeta'|^2}}
             {D(2^k\sg,2^j\zeta')} 
        \Big((i2^k\sg+2\cdot 2^{2j}|\zeta'|^2) \Big)
  \right| 
   = \; O(1). 
}
}
Those bounds enable us to treat the operator 
given by \eqref{eqn;Stokes-n-pressure-0} is $L^p(\re^{n-1})$ 
bounded in $x'$ and $L^1$ bound in $t$-variable. 
Thus the estimate \eqref{eqn;pressure-key} holds  
for all $1\le p\le \infty$. 
\end{prf}

\vskip3mm
\begin{prop}
\label{prop;grad-pressure-trace} Let  $1< p \, <\, \infty$
and $-1+1/p\,  <  s\,\,  {\mt <1/p} $.
There exists $C>0$  such that for any 
$\N(-\Del)^{-1}f\in C(\overline{\re_+};\dB^s_{p,1}(\re^n_+) \cap 
\dot{W}^{1,1}(\re_+; \dB^{s}_{p,1}(\re^{n}_+))$, 
$\N f\in  L^1(\re_+; \dB^{s}_{p,1}(\re^{n}_+))$, it holds that 
\begin{align} 
  \sup_{x_n\in \re_+}
  \Big( \big\| f(\cdot,\cdot,x_n) &
        \big\|_{ \dF^{\frac12-\frac{1}{2p}}_{1,1}(\re_+;\dB^{s}_{p,1}(\re^{n-1})) }
      + \big\| f(\cdot,\cdot,x_n) 
        \big\|_{{\crd L^1}(\re_+;\dB^{s+1-\frac{1}{p}}_{p,1}(\re^{n-1})) }
  \Big)
  \notag\\
 \le & C\Big(
        \big\|\N f\big\|_{ L^1(\re_+; \dB^{s}_{p,1}(\re^{n}_+)) },
       +\big\|\pt_t \N(-\Del)^{-1} f
        \big\|_{ L^1(\re_+; \dB^{s}_{p,1}(\re^{n}_+)) }
        \Big),
   \label{eqn;pressure-trace}
\end{align}
In particular,  for $1\le p<\infty$,
\begin{align} 
  \sup_{x_n\in \re_+}
   \big\| f(\cdot,\cdot,x_n) 
        \big\|_{ L^1(\re_+;\dB^{s+1-\frac{1}{p}}_{p,1}(\re^{n-1})) }
 \le & C\big\|\N f\big\|_{ L^1(\re_+; \dB^{s}_{p,1}(\re^{n}_+)) }.
   \label{eqn;pressure-trace-2}
\end{align}
\end{prop}

\vskip2mm
The proof of the trace estimate \eqref{eqn;pressure-trace} 
is along the line of proof for the trace estimate \eqref{eqn;sharp-trace-N} 
in Theorem \ref{thm;sharp-boundary-trace-N} shown in 
Appendix  and we show it in subsection \ref{Sec;7-2} of the Appendix below.

\subsection{Estimate for the velocity}
Once we obtain the estimates for the pressure $\N q$ to \eqref{eqn;ST2}, 
the required estimates for the velocities $v$ of the solution to \eqref{eqn;ST2} 
can be obtained by applying maximal regularity for the initial boundary 
value of the heat equations:
\eq{ \label{eqn;heat-N}
   \left\{
  \begin{aligned} 
    &\pt_t u -\Del u=f, 
     &\qquad t>0,\ \  &x\in \re^n_+,\\
   &\quad \left.\pt_n u(t,x',x_n)\right|_{x_n=0}=h(t,x'), 
     &\qquad t>0,\ \  &x'\in \re^{n-1}, \\
   &\quad u(t,x)\big|_{t=0}=u_0(x),
      &\qquad  &x\in \re^n_+,
   \end{aligned}     
    \right.
}
where $x=(x',x_n)\in \re^n_+$ and $\pt_n$ denotes the 
normal derivative $\pt/\pt x_n$ at any boundary point of $\re^n_+$.

\vskip2mm
\begin{prop}[Maximal $L^1$-regularity \cite{OgSs20-1}, \cite{OgSs20-2}]
\label{prop;MaxReg-Neumann} 
Let $1< p< \infty$ and  $-1+1/p<s \le 0$.
The problem \eqref{eqn;heat-N} admits a unique solution 
\eqn{
 \spl{
   &u\in C_b(\overline{\re_+};\dot{B}^s_{p,1}(\re^n_+))
         \cap \dot{W}^{1,1}(\re_+;\dot{B}^s_{p,1}(\re^n_+)),\\
   &\Del u\in L^1(\re_+;\dot{B}^{s}_{p,1}(\re^n_+)) 
}}
if and only if the external, the initial and the boundary data in 
\eqref{eqn;heat-N} satisfy
\begin{align*}
    &u_0\in \dB^{s}_{p,1}(\re^n_+),  \,
     \quad
    f\in L^1(\re_+;\dot{B}^s_{p,1}(\re^n_+)),  
     \\
    &h\in \dot{F}^{\frac12-\frac{1}{2p}}_{1,1}(\re_+;\dB^s_{p,1}(\re^{n-1}))
          \cap 
          L^1(\re_+;\dB^{s+1-\frac1p}_{p,1}(\re^{n-1})),
\end{align*}  
respectively.  Moreover following the maximal $L^1$-regularity
estimate holds:  
\eq{ \label{eqn;heat-L1MR}
  \begin{aligned}
   & \|\pt_t u\|_{L^1(\re_+;\dot{B}^s_{p,1}(\re^{n}_+))}
    +\|D^2 u \|_{L^1(\re_+;\dot{B}^s_{p,1}(\re^{n}_+))} \\
    \le& C\Big(
        \|u_0\|_{\dot{B}^s_{p,1}(\re^{n}_+)}
        +\|f\|_{L^1(\re_+;\dot{B}^s_{p,1}(\re^{n}_+))}
        +\|h\|_{\dot{F}^{\frac12-\frac{1}{2p}}_{1,1}(\re_+;\dB^s_{p,1}(\re^{n-1}))}
        +\|h\|_{L^1(\re_+;\dB^{s+1-\frac1p}_{p,1}(\re^{n-1}))}
         \Big), 
    \end{aligned}
  }
where $C$ is depending only on $p$, $s$ and $n$. 
\end{prop}
\vskip2mm

For the proof of Proposition \ref{prop;MaxReg-Neumann}, 
see \cite{OgSs20-1} and \cite{OgSs20-2}. 

\par\vspace{0.5pc}
\noindent
\begin{prf}{Theorem \ref{thm;L1MR2-b}} 
Let the boundary data satisfy the regularity assumption 
\eqref{eqn;boundary-assump}.
First we consider the $n$-th component of the unknown velocity that satisfies 
the the corresponding system 
{\bl
(cf. Shibata-Shimizu \cite[(4.24)]{SbSz03}, 
and  \cite[(5.19)]{SbSz08}).
}
Namely $(v_n,q)$ is given by the expressions 
\eqref{eqn;Stokes-n-velocity} and \eqref{eqn;Stokes-n-pressure}.
{\bl
In particular from \eqref{eqn;Stokes-n-velocity}, we see that 
\begin{align}
 \pt_t & v_n(t,x',x_n)\nn\\
    =& c_{n+1}\text{p.v.}\iint_{\re^{n}}
               e^{it\t+ix'\cdot \xi'}
               \bigg\{
                \frac{(B+|\xi'|)|\xi'|}{D(\t, \xi')}
                 \Big(
                  2B(i\xi'\cdot \widehat{H}')
                  -(|\xi'|^2+B^2) \widehat{H_n}
                 \Big)e^{-|\xi'|x_n }
                   \notag\\
    &\hskip3.5cm
      +
       \frac{(B+|\xi'|)|\xi'|}{D(\t, \xi')}
                 \Big(
                  -(|\xi'|^2+B^2)(i\xi'/|\xi'|\cdot \widehat{H}')
                  +2|\xi'|^2  \widehat{H_n}
                 \Big)  e^{-Bx_n }            
               \bigg\}              
               d\t d\xi'.
         \label{eqn;Stokes-n-velocity-del-t}
\end{align}
Via a very much similar argument for the pressure estimate in 
Proposition \ref{prop;grad-pressure-bound}, we may derive 
the estimates for $\pt_t v_n$, $\Del v_n$. 
Namely the $n$-component of the velocity} fulfills the estimate:
\eq{\label{eqn;MR-velocity-n-comp}
\spl{
    \|\pt_t v_n&\|_{L^1(\re_+;\dot{B}^s_{p,1}(\re^{n}_+))}
     +\|D^2  v_n \|_{L^1(\re_+;\dot{B}^s_{p,1}(\re^{n}_+))} 
 \\
   \le& C\big(
         \|H\|_{\dF^{\frac12-\frac{1}{2p}}_{1,1}(\re_+;\dB^s_{p,1}(\re^{n-1}))}
        +\|H\|_{L^1 (\re_+;\dB^{s+1-\frac{1}{p}}_{p,1}(\re^{n-1}))}
         \big). 
}}
{\bl
Note that the first term of the Fourier image of $\pt_t v_n$ in \eqref{eqn;Stokes-n-velocity-del-t} 
is indeed expressed by the pressure 
and the rest of the symbol which is the parabolic part involving the symbol   
$$
 \widetilde{m}(\t,\xi')=\frac{(B(\t,\xi')+|\xi'|)}{D(\t,\xi')}
  \Big(-(|\xi'|^2+B(\t,\xi')^2)\frac{i\xi'}{|\xi'|}, 2|\xi'|^2\Big)
$$ 
and the above symbol denotes the singular integral part and 
it is analogous to $m(\t,\xi)$ in \eqref{eqn;singular-int-multi} so that the estimate 
\eqref{eqn;MR-velocity-n-comp} 
follows from the estimate of the pressure term and maximal 
regularity for the parabolic part with quite similar argument 
found in the previous work, in particular
using the \cite[Lemma 6.5]{OgSs20-2} with a modification involving 
$\widetilde{m}$ as is shown in \eqref{eqn;ortho-2}, \eqref{eqn;ortho-3}, 
\eqref{eqn;ortho-7} and \eqref{eqn;ortho-8}.
}
{\ppl
Hence the maximal regularity estimate for the $n$-th component of the 
velocity as well as the pressure follows from 
the estimate \eqref{eqn;MR-velocity-n-comp}, 
Proposition \ref{prop;grad-pressure-bound} and
Proposition \ref{prop;pressure-bound} and we obtain that
\eq{\label{eqn;MR-velocity-n-comp+pressure}
\spl{
    \|\pt_t v_n&\|_{L^1(\re_+;\dot{B}^s_{p,1}(\re^{n}_+))}
     +\|D^2  v_n \|_{L^1(\re_+;\dot{B}^s_{p,1}(\re^{n}_+))}  \\
     & + \|\N q \|_{ L^1(\re_+; \dB^{s}_{p,1}(\re^{n}_+)) }
       + \big\|q|_{x_n=0}\big\|_{\dF^{\frac12-\frac{1}{2p}}_{1,1}(\re_+;\dB^s_{p,1}(\re^{n-1}))}
 \\
   \le& C\big(
         \|H\|_{\dF^{\frac12-\frac{1}{2p}}_{1,1}(\re_+;\dB^s_{p,1}(\re^{n-1}))}
        +\|H\|_{L^1 (\re_+;\dB^{s+1-\frac{1}{p}}_{p,1}(\re^{n-1}))}
         \big). 
}}
The trace estimate \eqref{eqn;pressure-trace-2}
in Proposition \ref{prop;grad-pressure-trace} enable us to control 
the term 
$\big\|q|_{x_n=0}\big\|_{L^1(\re_+;\dB^{s+1-\frac1{p}}_{p,1}(\re^{n-1}))}$.
}

The other components of the velocity fields 
$v'=(v_1(t,x), v_2(t,x), \cdots, v_{n-1}(t,x))$  satisfy the 
initial boundary value problem of the heat equations 
as the pressure and the $n$-th component velocity as the 
external force and boundary condition as follows:
For $\ell=1,2\cdots,n-1$,
\begin{equation}  \label{eqn;Stokes-ell-component}
   \left\{
   \begin{aligned}
     \pt_t v_\ell -\Delta v_\ell &=-\pt_{\ell} q,
        \quad&\ t>0,\  x\in \re^n_+,\\
     \pt_n v_\ell &=-H_{\ell}-\pt_{\ell} v_n,
        \quad&\ t>0,\  x\in \pt\re^n_+,\\
     v_\ell(0,x) &=0,
        \quad& \phantom{\ t>0.}\ x\in\re^n_+.
     \end{aligned}
     \right.  \hfill
\end{equation}
Similarly to the above estimate, we have from 
Proposition \ref{prop;grad-pressure-bound}, 
Proposition \ref{prop;grad-pressure-trace},
Proposition \ref{prop;MaxReg-Neumann}
and the estimate \eqref{eqn;MR-velocity-n-comp}
that the solution $v_\ell(t,x)$ to the problem 
\eqref{eqn;Stokes-ell-component} has the estimate 
\eq{\label{eqn;MR-velocity-ell-comp}
\spl{
     \|\pt_t v_{\ell} &\|_{L^1(\re_+;\dot{B}^s_{p,1}(\re^{n}_+))}
    +\|D^2  v_{\ell}  \|_{L^1(\re_+;\dot{B}^s_{p,1}(\re^{n}_+))} \\
  \le& C\big(
         \|\pt_\ell q\|_{L^1(\re_+;\dot{B}^s_{p,1}(\re^{n}_+))}
        +\|H_{\ell}\|_{\dot{F}^{\frac12-\frac{1}{2p}}_{1,1}(\re_+;\dB^s_{p,1}(\re^{n-1}))}
        +\|H_{\ell}\|_{{\crd L^1}(\re_+;\dB^{s+1-\frac1p}_{p,1}(\re^{n-1}))} \\
     &\qquad \quad
        +\big\|\pt_{\ell} v_n\big|_{x_n=0}
         \big\|_{\dot{F}^{\frac12-\frac{1}{2p}}_{1,1}(\re_+;\dB^s_{p,1}(\re^{n-1}))}
        +\big\|\pt_{\ell} v_n\big|_{x_n=0}
         \big\|_{{\crd L^1}(\re_+;\dB^{s+1-\frac1p}_{p,1}(\re^{n-1}))}
         \big)
 \\
   \le& C\big(
         \|\N q\|_{L^1(\re_+;\dot{B}^s_{p,1}(\re^{n}_+))}
        +\|H_{\ell}\|_{\dot{F}^{\frac12-\frac{1}{2p}}_{1,1}(\re_+;\dB^s_{p,1}(\re^{n-1}))}
        +\|H_{\ell}\|_{{\crd L^1}(\re_+;\dB^{s+1-\frac1p}_{p,1}(\re^{n-1}))} \\
      &\qquad \quad 
        +\|\pt_t v_n\|_{L^1(\re_+;\dot{B}^s_{p,1}(\re^{n}_+))}
        +\|\N^2  v_n\|_{L^1(\re_+;\dot{B}^s_{p,1}(\re^{n}_+))} 
         \big)
 \\
   \le& C\big(
         \|H\|_{\dot{F}^{\frac12-\frac{1}{2p}}_{1,1}(\re_+;\dB^s_{p,1}(\re^{n-1}))}
        +\|H\|_{{\crd L^1}(\re_+;\dB^{s+1-\frac{1}{p}}_{p,1}(\re^{n-1}))}
         \big). 
}}
Combining the estimates \eqref{eqn;MR-velocity-n-comp} and 
\eqref{eqn;MR-velocity-ell-comp} for all $\ell=1,2,\cdots, n-1$ as well as 
the pressure estimate \eqref{eqn;grad-pressure-bound} in Proposition 
\ref{prop;grad-pressure-bound}, we conclude that the desired estimate 
\eqref{eqn;L1MR2-b-estimate} holds.

Conversely, if the solution $(v,q)$ to the problem \eqref{eqn;ST2} 
exists, then it holds by letting $f$ by $v$ in the trace estimate 
\eqref{eqn;sharp-trace-N}  of Theorem \ref{thm;sharp-boundary-trace-N}  in Appendix 
and Proposition \ref{prop;grad-pressure-trace} that
\eq{\label{eqn;L1MR2-trace-estimate}
 \spl{
      \|H&\|_{\dF^{\frac12-\frac{1}{2p}}_{1,1}(\re_+;\dB^s_{p,1}(\re^{n-1})))} 
     +\|H \|_{L^1(\re_+;\dB^{s+1-\frac{1}{p}}_{p,1}(\re^{n-1}))} 
\\
\le &  2\|\N v\|_{\dF^{\frac12-\frac{1}{2p}}_{1,1}(\re_+;\dB^s_{p,1}(\re^{n-1})))} 
      +2\|\N v\|_{L^1(\re_+;\dB^{s+1-\frac{1}{p}}_{p,1}(\re^{n-1}))}   \\
    &\quad
      +\big\|q\,|_{x_n=0}
        \big\|_{\dF^{\frac12-\frac{1}{2p}}_{1,1}(\re_+;\dB^s_{p,1}(\re^{n-1})))} 
      +\big\|q\,|_{x_n=0}
       \big\|_{{\crd L^1}(\re_+;\dB^{s+1-\frac{1}{p}}_{p,1}(\re^{n-1}))}   
\\
  \le& C \Big(\|\pt_t v\|_{L^1(\re_+;\dB^s_{p,1}(\re^n_+))}
             +\|\N^2 v\|_{L^1(\re_+;\dB^s_{p,1}(\re^n_+))} \\
     &\qquad
             +\|\N q\|_{L^1(\re_+;\dB^s_{p,1}(\re^n_+))}
             +\big\|q \,|_{x_n=0}
              \big\|_{\dF^{\frac12-\frac{1}{2p}}_{1,1}(\re_+;\dB^s_{p,1}(\re^{n-1})))} 
             +\big\|q\,|_{x_n=0}
              \big\|_{{\crd L^1}(\re_+;\dB^{s+1-\frac{1}{p}}_{p,1}(\re^{n-1}))} 
         \Big).
}
}
This shows regularity for the boundary data is necessary.
\vskip1mm
{\bl 
Concerning the uniqueness, we invoke the argument employed  
 in \cite[Theorem 4.3 and 5.7]{SbSz08} for the half space.
Let $(v,q)$ be a solution of the Stokes system \eqref{eqn;ST} with 
vanishing data {\mt and satisfying the} regularity Theorem \ref{thm;L1MR}. 
For any $\phi\in C^\infty_0(\re\times \re^n_+)$ and for any 
{\mt $T>0$ with}
$\phi=0$ for  {\mt  $t\in (-\infty, -1)\cup (T/2,\infty)$,  set} $\phi_T(t,x)=\phi(T-t, x)$. 
Let $(\Phi,\theta)$ be the solution of \eqref{eqn;ST} 
with vanishing data except the external force $f=\phi_T$ and set 
$v_*(t,x)=\Phi(T-t, x)$ and $q_*(t,x)=\theta(T-t,x)$
 and arrange its support into the time interval 
$(-1,T/2)$.
Then $(v_*,q_*)$ solves the adjoint Stokes system 
except the pressure sign
in the subset of the dual space $L^{\infty}(I;\dot{H}^{-s,p'}(\re^n_+))
\subset L^{\infty}(I;\dB^{-s}_{p',\infty}(\re^n_+))$.
If we choose {\mt $\r>2p/(p+1)$}, then 
$v_*\in W^{1,\r}(I;L^{p'}(\re^n_+))
\cap L^{\r}(I;H^{2,p'}(\re^n_+)) \subset L^{\infty}(I;\dH^{-s,p'}(\re^n_+))$.
Here we note that 
$\dot{H}^{-s,p'}(\re^n_+)\subset 
\dB^{-s}_{p',\infty}(\re^n_+)=(\dB^s_{p,1}(\re^n_+))^*$,
where $0\le -s < 1/p'$.
Let $\tilde\chi(r)$ be a smooth cut-off function of $r>0$
over the annulus $B_2(0)\setminus \bar{B}_1(0)$ with 
$\chi(x)\equiv \tilde\chi(|x|)$  and
set $\chi_R(x)\equiv R^{-1}\chi(R^{-1}x)$ for any $R>0$ 
and $D_R=\supp \chi_R(x)\equiv \{x\in\re^n_+;\ R\le |x|\le 2R\}$. 
By the Poincar\'e--Wirtinger inequality for $I=(-2,T)$, there exists $\theta\in (1,\infty)$ such that
\begin{align*}
        &\Big|\int_I\int_{\re^n_+}q(t,x)\chi_R(x) v_*(t,x)\,dx\,dt\Big|\\
  \le &C \Big(\|q|_{x_n=0}\|_{L^{1}(I;\dot{B}^{s+1-\frac{1}{p}}_{p,1}(\re^{n-1}))}
               +\| \N q\|_{L^{1}(I;\dot{B}^{s}_{p,1}(\re^n_+))}  
         \Big)
            \|\chi_R\,  v_*\|_{C_b(I; \dH^{ -s, p'}(D_R))}, 
     \eqntag\label{eqn;boundary-coupling-1}
  \\
  &\Big|\int_I\int_{\re^n_+}q_*(t,x)\chi_R(x) v(t,x)\,dx\,dt\Big|\\
  \le &C \Big( \| q_*|_{x_n=0}
              \|_{L^{\th}(I; \dot{W}_{p'}^{-s +1-\frac 1{p'}}(\re^{n-1}))}
             +\|\N q_*\|_{L^{\th}(I; \dot{H}^{-s, p'}(\re^n_+))} \Big)
              \|\chi_R\,v\|_{C_b(I;\dB^s_{p,1}(\re^n_+))}.
     \eqntag\label{eqn;boundary-coupling-2}
\end{align*}
By passing $R\to \infty$ 
the both terms in the right hand side of \eqref{eqn;boundary-coupling-1} and \eqref{eqn;boundary-coupling-2} 
are vanishing (cf. by the bilinear estimate \eqref{eqn;bilinear1.0-ap} and $\|\chi_R\|_{\dB^{n/p}_{p,1}(\re^n_+)}= O(R^{-1})$), which justify the integration by parts 
(cf. \cite{SbSz03} and \cite{SbSz08}
for the details). Since the external force is smooth, the exceptional regularity can be avoided. 
The analogous estimate above also justify the dual couplings at the boundary.   
Those observations ensure that the following argument remains valid; 
\begin{align*}
\langle v, \phi\rangle_{\re\times\re^n_+}
&=\langle v, -\pt_t v_*-\Delta v_*+\N q_*\rangle_{\re\times\re^n_+}\\
&=\langle \pt_t v, v_*\rangle_{\re\times\re^n_+} +\langle \N v, \N v_*+(\N v_*)^{\sf T}-q_*I\rangle_{\re\times\re^n_+}-(v|_{x_n=0},T(v_*, q_*)\nu_n|_{x_n=0}\rangle_{\re\times\re^{n-1}}\\
&=\langle \pt_t v, v_*\rangle_{\re\times\re^n_+} +\langle \N v+(\N v)^{\sf T}-qI, \N v_*\rangle_{\re\times\re^n_+}\\
&=\langle \pt_t v, v_*\rangle_{\re\times\re^n_+} +\langle -\Delta v+\N q, v_*\rangle_{\re\times\re^n_+}
+\langle T(v, q)\nu_n|_{x_n=0}, v_*\rangle_{\re\times\re^{n-1}}\\
&=\langle \pt_t v-\Delta v +\N q, v_*\rangle_{\re\times\re^n_+}=0,
\end{align*}
from which we conclude $v=0$ by the arbitrariness of $\phi$, and hence $q=0$ by $\N q=0$ in $\re^n_+$ 
and $q(\cdot, 0)=0$ by \eqref{eqn;ST2}.  
}
This proves Theorem \ref{thm;L1MR2-b}. 
\end{prf}
\vskip 3mm

\begin{prf}{Theorem \ref{thm;L1MR}}
Applying the maximal $L^1$-regularity result 
to the initial-boundary value problem of the Stokes equations with the
boundary condition, 
we obtain end-point maximal $L^1$-maximal regularity from 
\eqref{eqn;MR-velocity-n-comp}, \eqref{eqn;MR-velocity-ell-comp}
and 
Since by \eqref{eqn;up-decomposion} 
$u=\widetilde{u}+v-\N\phi |_{x_n>0}$ and $p=\widetilde{p}+q$ is the solution to \eqref{eqn;ST}.  Hence by combining 
\eqref{eqn;g}-\eqref{eqn;potential},
\eqref{eqn;data-tilde},
\eqref{eqn;up-tilde},
\eqref{eqn;boundary-data},
\eqref{eqn;L1MR2-b-estimate} in Theorem \ref{thm;L1MR2-b},
 we obtain \eqref{eqn;L1MR-estimate}.
 
Conversely, by using
\eqref{eqn;up-tilde-trace},
\eqref{eqn;pressure-trace} 
in Proposition \ref{prop;grad-pressure-trace}
as well as 
\eqref{eqn;sharp-trace-N} in 
Theorem \ref{thm;sharp-boundary-trace-N}, 
we conclude that regularity for data is necessary for 
the existence of the solution $(u,p)$ for the Stokes 
system \eqref{eqn;ST}.
This completes the proof of Theorem \ref{thm;L1MR}.
\end{prf}

\sect{Multiple div-curl structure and critical multi-linear 
estimates}\label{Sec;5}

In this section, we show {\it the multiple divergence-free-curl-free 
structure} related to Jacobi matrix of transformation 
from the Euler coordinates to the Lagrange coordinate, 
which is essential to obtain global well-posedness in the 
critical Besov space $\dB^{-1+n/p}_{p,1}(\re^n_+)$. 
The single divergence structure was firstly 
pointed out by Solonnikov \cite{Sol88}, 
and it is applied by Shibata--Shimizu \cite{SbSz07} for 
the free boundary value problem.
Our case is a multiple extension from those 
divergence structure.  Namely in order to apply the 
bilinear estimate in the critical Besov spaces, 
we need to ensure that the divergence-free, rotation-free 
structure for every step when we apply the bilinear estimate.
Namely for multi-linear case, we need to make it clear that the 
nonlinear terms in the equation maintains the multiple 
div-curl free structures.  This was shown in \cite{OgSs18} for 
the initial value problem for the Lagrangian coordinate case.  
We develop the analogous estimate  and establish 
 the multiple Besov estimate in the half-spaces.

\subsection{Multiple div-curl structure}
We show that the inverse matrix of Jacobian for the Lagrangian transform 
and consequently the perturbation terms $F_u$, $F_p$ and $G_{\div}$ have 
a special divergence structure.  We call a inner product of two vector fields 
$f\cdot g$ for $f, g\in \D^*$
maintains the divergence free 
rotation-free structure (in short div-curl structure) if 
$\rot f=0$ and $\div g=0$ and the multiple-div-curl structure for 
$f\cdot \Pi(g_1\otimes g_2\otimes \cdots \otimes g_{\ell})$
if $\rot f=0$ and $\div \Pi=0$ and $\Pi$ can be decomposed into 
lower order component consisting of div-curl structure.
Such kind of structure easily yields us the original term can 
be expressed in the divergence form and the bilinear estimate can be 
enlarged in the critical Besov framework.

We show such a structure holds for each of the semi-linear terms
of the system \eqref{eqn;NS}.

\vskip3mm
\begin{prop}[Multiple divergence structure] 
\label{prop;multiple div-structure} Let $I=(0,T)$ with 
$T\le \infty$ and suppose that $(u,p)$ has the following regularities;
\eqn{
 \spl{
  \N u\in  L^1(I;\dB^{\frac{n}{p}}_{p,1}(\re^n_+)\big), \quad
  \N p\in  L^1(I;\dB^{-1+\frac{n}{p}}_{p,1}(\re^n_+)\big).
 }
}
Let $F_p(u,p)$, $G_{\div}(u)$ 
be the polynomials of $d_{jk}$ of order $n-1$ defined in
\eqref{eqn;pressure-purterb},
\eqref{eqn;divergence-purterb}. 
Then the terms
are subject to the multiple div-curl structure.
Namely every component of those polynomial consist of the 
inner products of divergence free vector and rotation free vectors.
\end{prop}
\vskip3mm

Before going into details of the proof, we introduce several 
notations. Let
\eq{ \label{eqn;element}
 {\td}_{ij}=\del_{ij}+d_{ij}
      =\del_{ij}+\int_0^t \pt_j u_i(s)ds,
}
then the Jacobi matrix is written as 
\eqn{
\spl{
  J(Du)=&\begin{pmatrix}
      1+d_{11} & d_{12}   &\cdots  &d_{1n}\\
      d_{21}   & 1+d_{22} &\cdots  &d_{2n}\\
      \vdots   & \vdots   &\ddots  & \vdots\\
      d_{n1}   & d_{n2}   &\cdots  &1+d_{nn}
     \end{pmatrix} 
  =  \begin{pmatrix}
      \td_{11} & \td_{12} &\cdots &\td_{1n}\\
      \td_{21} & \td_{22} &\cdots &\td_{2n}\\
      \vdots   & \vdots   &\ddots &\vdots\\
      \td_{n1} & \td_{n2} &\cdots & \td_{nn}
     \end{pmatrix}
}}
and set 
\eqn{
  J(Du)^{-1}
  =  \begin{pmatrix}
      b_{11} & b_{12} &\cdots  & b_{1n}\\
      b_{21} & b_{22} & \cdots & b_{2n}\\
      \vdots & \vdots & \ddots & \vdots\\
      b_{n1} & b_{n2} & \cdots & b_{nn}
     \end{pmatrix}.
}
For any $1\le \ell\le n-1$, let us denote 
{\ppl an $\ell\times \ell$-submatrix} of $J(Du)$ as 
\eq{\label{eqn;inverse-Jacobi-matrix}
  J(Du)^{[\ell]}
   \equiv
     \begin{pmatrix}
      \td_{\sg_1\t_1}   &\cdots &\td_{\sg_1\t_\ell} \\
       \vdots            &\ddots &\vdots      \\
      \td_{\sg_\ell\t_1}   &\cdots &\td_{\sg_\ell\t_\ell}
     \end{pmatrix},                                   
  }
where $(\sg_1,\sg_2,\cdots,\sg_\ell)$ and $(\t_1,\t_2,\cdots,\t_\ell)$
are any combination of ordered sub-factor from $(1,2,\cdots, n)$, 
namely $1\le \sg_1<\sg_2<\cdots <\sg_\ell\le n$ and $1\le \t_1<\t_2<\cdots {\ppl <}\t_\ell\le n$.
We notice that 
\eqn{ 
 b_{kj}=(-1)^{k+j}\det J(Du)^{[n-1]}_{jk}. 
}

We prove {\it the multiple div-curl structure} 
for those \eqref{eqn;pressure-purterb}--\eqref{eqn;divergence-purterb} by induction. 
The proof for the case of $F_p(u,p)$ in \eqref{eqn;pressure-purterb},
$G_{\div}(u)$ in \eqref{eqn;divergence-purterb} 
can be shown along a similar way. Hence we mainly show \eqref{eqn;pressure-purterb} 
for the case $F_p(u,p)$.
In the case of $n=2$ and $n=3$, such a  structure is shown in  an explicit way 
for $\re^n$ case (see \cite{OgSs18}). 
It is easy to show \eqref{eqn;pressure-purterb} in the case of $n=2$.

According to Evans \cite{Ev} (section 8.1),
we recall the null Lagrangian structure for the Jacobian of a Lipschitz 
continuous function  $u$. 
Let $A$ be a $n\times n$ matrix and consider its $\ell\times \ell$ sub-matrix 
$A^{[\ell]}$.

\vskip3mm
\begin{lem}[Divergence free for sub-cofator]
\label{lem;subcofactor-divfree}
Let  $u:\re^n\to \re^n$ be a Lipschitz continuous function 
and for any $1\le \ell\le n$, let 
$J(Du)^{[\ell]}$ be {\ppl an} $\ell\times \ell$-submatrix of the Jacobi 
matrix $J(Du)$ and 
$\cof(J(Du)^{[\ell]})_{kj}=(-1)^{\sg_k}\det J(Du)^{[\ell-1]}_{kj}$ 
be the cofactor matrix of the {\crd sub-matrix} $J(Du)^{[\ell]}$. Then 
for $(k,j)$ component of $\cof(J(Du)^{[\ell]})$, it holds
$$
 \div_{\hspace{-1mm}j}\, \cof(J(Du)^{[\ell]}\big)_{k j}=0
$$
for any point $x\in \re^{\ell}$ with  $\det(J(Du)^{[\ell]})(x)\neq0$.
\end{lem}
\vskip2mm

We show the proof of Lemma \ref{lem;subcofactor-divfree} in 
the Appendix (Lemma \ref{lem;cofactor-divfree}) below.

\vskip3mm
\begin{prf}{Proposition \ref{prop;multiple div-structure}}
Assume that $\N u\in L^1\big(I;\dB^{n/p}_{p,1}(\re^n_+)\big)$.
Since $\dB^{n/p}_{p,1}(\re^n_+)\subset C_v(\re^n_+)$
for almost every $t\in I=(0,T)$, we regard that $u$ is the Lipschitz continuous
 function in $x\in \re^n_+$ for almost all $t\in I$.

Each component of $J(Du)^{-1}$ can be realized by 
the cofactor expansion by $J(Du)$, namely
\eqn{
 b_{kj}=(\det J(Du))^{-1}\sum_{j=1}^n (-1)^{k+j}\td_{kj}
         \det(J(Du)^{[n-1]})_{kj}
       =\sum_{j=1}^n \td_{kj} \cof\big(J(Du)\big)_{kj},
}
where we recall that $\div \bar{u}=0$ implies $\det (J(D\bu))=\det(J(Du))=1$.
Then it is easy to see that from Lemma \ref{lem;subcofactor-divfree} and 
for each $k$,
$$
 \td_{kj}=\del_{kj}+\int_0^t \pt_j u_k(s)ds
$$
is a rotation-free vector and hence the each component of 
cofactor of $J(Du)$, namely $J(Du)^{-1}$ has the div-curl structure 
and this structure can be decomposed into any order of its sub-factor 
by expanding the determinant of cofactor matrices
$J(Du)_{kj}^{[n-1]}$. It can be realized by the form
\algn{
 \det \big(J(Du)^{[\ell]}\big)
  =&\sum_{j=1}^{\ell}(-1)^{k+j}\td_{kj}\, \det\big(J(Du)^{[\ell-1]}\big)_{kj}
  =\sum_{j=1}^{\ell}\td_{kj}\, \cof\big(J(Du)^{[\ell]}\big)
}
for all $\ell=1,2,\cdots{\ppl ,} n-1$. It is clear from Lemma \ref{lem;subcofactor-divfree}
that the above expression also maintains div-curl structure
unless $\det \big(J(Du)^{[\ell]}\big)=0$, since 
$\td_{k}\equiv (\td_{k1},\td_{k2},\cdots{\ppl ,} \td_{kn})$ is a 
rotation free vector for each $k=1,2,\cdots{\ppl ,} n$.

Now we finalized the proof to see that 
\eqn{
 \spl{
 \big(F_p(u,p)\big)_k= \Big(\big(J(Du)^{-1}-I\big)^{\sf T} \N p\Big)_k
          =\sum_{j=1}^n
           \cof\big( J(Du)^{[n-1]} \big)_{kj}\pt_j p -\pt_k p
}
}
with observing that the first term is an inner product of the rotation free vector $\N p$ 
and divergence free element $\cof\big(J(Du)_{kj}^{[n-1]}\big)$
and the second is also with a trivial curl-free element $\del_{kj}$.

The proof for $G_{\ppl \rm div}(u)$ goes almost the same way since each component of 
trace has the div-curl structure as the above.
The boundary term $H_u(u)$ is also decomposed into the div-curl free structure 
before taking the inner product with $\nu_n$.
This completes the proof.
\end{prf}

%
\subsection{Bilinear estimate for div-curl structure}\label{Sec6}
In general, the following bilinear estimate does not hold in the Besov space 
over $\re^n$: 
$$
  \|fg\|_{\dot{B}^{-1+\frac np}_{p,1}}
  \le C\|f\|_{\infty}\|g\|_{\dot{B}^{-1+\frac np}_{p,1}}. 
$$
However it is possible to change the norm of $L^\infty$ into slightly stronger norm of 
$\dot{B}^{n/q}_{q,1}$, which have the same scaling invariance with $L^{\infty}$.  
The following bilinear estimate is essentially obtained 
by Abidi-Paicu \cite{AP07} (cf. \cite{OgSs16},  \cite{OgSs18}) in $\re^n$.

\begin{prop}
\label{prop;abidi-paicu-halfspace}
Let $1\le p,p_1,p_2\le \infty$ and $1/p=1/p_1+1/p_2$.  
Under the assumption $s+s'>1$, in particular, for $1\le p<2n$, 
for all 
$f\in \dB^{s-1}_{p_1,1}(\re^n_+)$ and 
$g\in  L^{p_2}(\re^n_+)\cap \dot B^{s'}_{p_2,\infty}(\re^n_+)$, 
there exists $C>0$ independent of  $f$, $g$ such that the following estimate holds:
\eq{\label{eqn;bilinear1.1}
  \|f g\|_{\dB^{-1+s}_{p,1}(\re^n_+)}
   \le    C \|f\|_{\dot B^{-1+s}_{p_1,1}(\re^n_+)}
            \big(  \|g\|_{L^{p_2}(\re^n_+)} 
                 + \|g\|_{\dot B^{s'}_{p_2,\infty}(\re^n_+)}\big).
}
In particular for $g\in  \dot B^{n/p}_{p,1}(\re^n_+)$
\eq{  \label{eqn;bilinear1.2}
  \|f g\|_{\dot{B}^{-1+\frac{n}{p}}_{p,1}(\re^n_+)}
  \le
    C\|f\|_{\dot{B}^{-1+\frac{n}{p}}_{p,1}(\re^n_+)} 
     \|g\|_{\dot{B}^{\frac{n}{p}}_{p,1}(\re^n_+)}. 
    }
\end{prop}
\vskip2mm
\begin{prf}{Proposition \ref{prop;abidi-paicu-halfspace}}
Since the bilinear estimates in \eqref{eqn;bilinear1.1}
and \eqref{eqn;bilinear1.2} are established in the whole space, 
we merely show the case for the half-space for \eqref{eqn;bilinear1.2}.
Let 
$f\in \dB^{-1+n/p}_{p,1}(\re^n_+)$ and 
$g\in \dB^{n/p}_{p,1}(\re^n_+)$ then 
from the definition of the Besov space in the half-space,
for any $\ep>0$ there exists 
$\tilde{f}\in \dB^{-1+n/p}_{p,\sg}(\re^n)$ and 
$\tilde{g}\in\dB^{n/p}_{p,1}(\re^n)$
such that $f=\tilde{f}$ in {\ppl $\D'(\re^n_+)$}  and $g=\tilde{g}$ 
over $\re^n_+$ and
\eqn{
 \spl{
 &\|\widetilde{f}\|_{\dB^{-1+\frac{n}{p}}_{p,1}}
  \le \|f\|_{\dB^{-1+\frac{n}{p}}_{p,1}(\re^n_+)}+\ep,\\
 &\|\widetilde{g}\|_{\dB^{\frac{n}{p}}_{p,1}}
  \le \|g\|_{\dB^{\frac{n}{p}}_{p,1}(\re^n_+)}+\ep.
}
}
Then the corresponding
estimate \eqref{eqn;bilinear1.1} in $\re^n$ now implies
\eqn{
 \spl{
  \|fg\|_{\dB^{-1+\frac{n}{p}}_{p,1}(\re^n_+)}
  \le & \Big\|\widetilde{f g}\big|_{\re^n_+}
        \Big\|_{\dB^{s-1}_{p,1}}
  \le   \Big\|\widetilde{f} \widetilde{g}
        \Big\|_{\dB^{s-1}_{p,1}}
  \le C \|\widetilde{f}\|_{\dB^{-1+\frac{n}{p}}_{p,1}}
        \|\widetilde{g}\|_{\dB^{\frac{n}{p}}_{p,1}} 
\\
 \le &C\|f\|_{\dB^{-1+\frac{n}{p}}_{p,1}(\re^n_+)}
       \|g\|_{\dB^{\frac{n}{p}}_{p,1}(\re^n_+)}
      +\ep.
}
}
Since $\ep>0$ is arbitraly, this proves the estimate 
 \eqref{eqn;bilinear1.2} in $\re^n_+$.
The estimate  \eqref{eqn;bilinear1.1} follows in a similar 
way.
\end{prf}
\vskip2mm
\vskip3mm
\begin{prop}[Bilinear estimate under 
div-curl structure]\label{prop;Besov-bilinear-div-curl-half}
Let $1\le p<\infty$.  
For any vector valued functions $f \in  \dB^{-1+n/p}_{p,1}(\re^n_+)$ 
 and $g\in  \dB^{n/p}_{p,1}(\re^n_+)$ with   $\div f=0$ and $\rot g=0$ 
in the distribution sense,  
there exists a constant $C>0$ independent of $f$, $g$ such that 
\eq{\label{eqn;bilinear3}
  \|f\cdot g\|_{\dot B^{-1+\frac{n}{p}}_{p, 1}(\re^n_+)} 
      \le C \|f\|_{\dot B^{-1+\frac np}_{p,1}(\re^n_+)}
            \|g\|_{\dB^{\frac np}_{p,1}(\re^n_+)}.
}
If $f$ satisfies $\rot f=0$ and $g$ satisfies $\div g=0$ in ${\mathcal D}'$, 
then \eqref{eqn;bilinear3} also holds. 
\end{prop}

\begin{prf}{Proposition \ref{prop;Besov-bilinear-div-curl-half}} 
The corresponding estimate to \eqref{eqn;bilinear3} in the whole 
space is shown in Proposition \ref{prop;Besov-bilinear-div-curl} 
in Appendix below (cf. \cite{OgSs18}).
To show the half-space case, the argument of proof of 
Proposition \ref{prop;abidi-paicu-halfspace} works as well and 
this shows the proof.
\end{prf}

\vskip3mm
\subsection{Multi-linear estimate under the div-curl structure}
The perturbation terms for the Navier--Stokes equations in the 
Lagrangian coordinate exhibit the multiple-div-structure.  
In this case, we use the improved bilinear 
estimate in Proposition \ref{prop;Besov-bilinear-div-curl-half} 
in the critical Besov space 
for the nonlinear term(cf. for the whole space case \cite{OgSs18},
Proposition 4.5--4.6).

\vskip3mm
\begin{prop}[Multiple div-curl estimates 1]
\label{prop;multiple-bilinear-1}
Let $n\ge 2$, $1\le p<\infty$. 
For $\pt_t u$, $D^2 u$ and $\N p\in L^1(\re_+;\dB^{-1+n/p}_{p,1}(\re^n_+))$,
let $F_{p}(u,p)$ and $G_{\div}(u)$ be the terms defined 
in \eqref{eqn;pressure-purterb} and \eqref{eqn;divergence-purterb}, respectively. 
Then the following estimates hold:
\begin{align}
  \label{eqn;m-bilinear-apriori-p}
  \|F_{p}(u,p)\|_{L^1(\re_+;\dB^{-1+\frac{n}{p}}_{p,1}(\re^n_+))}
    \le& C\sum_{k=1}^{n-1}
            \|D^2 u\|_{L^1(\re_+;\dB^{-1+\frac{n}{p}}_{p,1}(\re^n_+))}^k
            \|\N p\|_{L^1(\re_+;\dB^{-1+\frac{n}{p}}_{p,1}(\re^n_+))}, 
\\
  \label{eqn;m-bilinear-apriori-div}
  \|\N G_{\div}(u)\|_{L^1(\re_+;\dB^{-1+\frac{n}{p}}_{p,1}(\re^n_+))}
    \le& C\sum_{k=1}^{n-1}
           \|D^2 u\|_{L^1(\re_+;\dB^{-1+\frac{n}{p}}_{p,1}(\re^n_+))}^{k+1}. 
\end{align}
In particular for $1\le p<2n$, it holds that
\begin{align}
 \label{eqn;m-bilinear-apriori-div-dt}
 \|\pt_t (-\Del)^{-1}\N G_{\div}(u)
 \|_{L^1(\re_+;\dB^{-1+\frac{n}{p}}_{p,1}(\re^n_+))}
    \le& C\sum_{k=1}^{n-1}
          \|  D^2 u\|_{L^1(\re_+;\dB^{-1+\frac{n}{p}}_{p,1}(\re^n_+))}^k
          \|\pt_t u\|_{L^1(\re_+;\dB^{-1+\frac{n}{p}}_{p,1}(\re^n_+))}.
\end{align}
\end{prop}
\vskip2mm

\begin{prf}{Proposition \ref{prop;multiple-bilinear-1}}
First we show  the estimate for $F_p$:
\eqn{
 \spl{
    F_p(u,p)
        = -\big(J(D u)^{-1}-I\big)^{\sf T}\;\N p
       \equiv& \Pi_{p}^{n-1}
             \left(\int_0^t D u\,ds\right)
             \N p .
}}
Here $\Pi_{ p}^{n-1}(\cdot)$ is the $n-1$-th order polynomial 
of the component of inverse of Jacobi matrix $J(Du)^{-1}$ without a constant 
term: 
$$
  \Pi_{ p}^{n-1}
           \Big(\int_0^t D u\, ds\Big)
   =\sum_{k=1}^{n-1}c_k\prod_{m,\ell\le n}^k
    \Big(\int_0^t \pt_{\ell} u_m\,ds\Big). 
$$

By using \eqref{eqn;bilinear3} with $1\le  p<\infty$ repeatedly, we have 
by inductively that 
%
%
{\allowdisplaybreaks
\algn{
  \|F_p&(u,p)\|_{L^1(\re_+;\dB^{-1+\frac{n}{p}}_{p,1}(\re^n_+))}\\
  \le & C\sup_{t>0}
         \Big\|\int_0^t D u\,ds\Big\|_{\dB^{\frac{n}{p}}_{p,1}(\re^n_+)}
          \bigg(\|\N p\|_{L^1(\re_+;\dB^{-1+\frac{n}{p}}_{p,1}(\re^n_+))}\\
     &\hskip4.5cm   +
           \Big\|\sum_{k=2}^{n-2}c_k\prod_{m,\ell\le n}^k
                \Big(\int_0^t \pt_{\ell} u_m\,ds\Big) 
                \N p        
          \Big\|_{L^1(\re_+;\dB^{-1+\frac{n}{p}}_{p,1}(\re^n_+))}
         \bigg)
   \\
  \le & C\left\|\int_0^t D u\, ds\right\|_{L^{\infty}(\re_+;\dB^{\frac{n}{p}}_{p,1}(\re^n_+))}
         \|\N p\|_{L^1(\re_+;\dB^{-1+\frac{n}{p}}_{p,1}(\re^n_+))} \\
      &\hskip1cm
         +\Big(\sup_{t>0}
                \Big\|\int_0^t D u\, ds\Big\|_{\dB^{\frac{n}{p}}_{p,1}(\re^n_+)}
          \Big)^2
           \Big\|\sum_{k=0}^{n-3}c_k\prod_{m,\ell\le n}^k
                \Big(\int_0^t \pt_{\ell} u_m\,ds\Big) 
                \N p        
          \Big\|_{L^1(\re_+;\dB^{-1+\frac{n}{p}}_{p,1}(\re^n_+))}
   \\
   \le & C\sum_{k=1}^{n-1}
           \Big\| \int_0^t D u\, ds
           \Big\|_{L^{\infty}(\re_+;\dB^{\frac{n}{p}}_{p,1}(\re^n_+))}^k
           \| \N  p\|_{L^1(\re_+;\dB^{-1+\frac{n}{p}}_{p,1}(\re^n_+))}
\\
   \le & C\sum_{k=1}^{n-1}
           \| D^2 u\|_{L^1(\re_+;\dB^{-1+\frac{n}{p}}_{p,1}(\re^n_+))}^k
           \| \N p\|_{L^1(\re_+;\dB^{-1+\frac{n}{p}}_{p,1}(\re^n_+))}.
}
}
This proves the estimate \eqref{eqn;m-bilinear-apriori-p}.

\vskip2mm\noindent
Next we show \eqref{eqn;m-bilinear-apriori-div}. The estimate of $\N G_{\div}$:
$$
  \N G_{\div}(u)
    \equiv \N \text{tr}\left(\Pi_{\div}^{n-1}
              \Big(\int_0^t D u\, ds\Big)D u\right). 
$$
Here $\Pi_{\div}^{n-1}(\cdot)$ is the $(n-1)$-th order polynomial 
without a constant term:
$$
  \Pi_{\div}^{n-1}
           \Big(\int_0^t D u\, ds\Big)
   =\sum_{k=1}^{n-1}c_k\prod_{m,\ell\le n}^k
    \Big(\int_0^t \pt_{\ell} u_mds\Big). 
$$
By using \eqref{eqn;bilinear3}, we have 
{\allowdisplaybreaks
\algn{
 & \|\N G_{\div}(u)\|_{L^1(\re_+;\dB^{-1+\frac{n}{p}}_{p,1}(\re^n_+))}\\
 \le &   \Big\|\Pi_{\div}^{n-1}
                \Big(\int_0^t D u\,ds\Big)   
                 D^2 u      
         \Big\|_{L^1(\re_+;\dB^{-1+\frac{n}{p}}_{p,1}(\re^n_+))} 
     +   \Big\|\Big( D\, \Pi_{\div}^{n-1}
                  \Big(\int_0^t D u\,ds\Big)
                \Big)
                \div\, u         
         \Big\|_{L^1(\re_+;\dB^{\frac{n}{p}}_{p,1}(\re^n_+))}
  \\
  \le & C\sup_{t>0}
         \Big\|\int_0^t D u\,ds\Big\|_{\dB^{\frac{n}{p}}_{p,1}(\re^n_+)}
         \Big\|\Pi_{\div}^{n-2}
                \Big(\int_0^t D u\,ds\Big)
                 D^2\, u         
         \Big\|_{L^1(\re_+;\dB^{-1+\frac{n}{p}}_{p,1}(\re^n_+))}\\
     & +C\sup_{t>0}
         \Big\|\int_0^t D u\,ds\Big\|_{\dB^{\frac{n}{p}}_{p,1}(\re^n_+)}
         \Big\| D\,\Pi_{\div}^{n-2}
                \Big(\int_0^t D u\,ds\Big)
                \div\, u         
         \Big\|_{L^1(\re_+;\dB^{-1+\frac{n}{p}}_{p,1}(\re^n_+))}
   \\
  \le & C\left\| D^2 u\right\|_{L^1(\re_+;\dB^{-1+\frac{n}{p}}_{p,1}(\re^n_+))}\\
      &\times
         \bigg(
         \Big\|\Pi_{\div}^{n-2}
                \Big(\int_0^t D u\,ds\Big)
                 D^2\, u         
         \Big\|_{L^1(\re_+;\dB^{-1+\frac{n}{p}}_{p,1}(\re^n_+))}\\
     &\phantom{\times \bigg(\|}
         +\Big\|\Pi_{\div}^{n-3}
                \Big(\int_0^t D u\,ds\Big)
                 D\,\Big(\int_0^t D u\,ds\Big)
                \div\, u         
         \Big\|_{L^1(\re_+;\dB^{-1+\frac{n}{p}}_{p,1}(\re^n_+))}
        \bigg)
   \\
   \le &C\sum_{k=1}^{n-2}
         \left\|D^2 u\right\|_{L^1(\re_+;\dB^{-1+\frac{n}{p}}_{p,1}(\re^n_+))}^k
         \bigg(
         \Big\|\Big(\int_0^t D u\,ds\Big) D \div\, u       
         \Big\|_{L^1(\re_+;\dB^{-1+\frac{n}{p}}_{p,1}(\re^n_+))} \\
       &\qquad\qquad\qquad
         + \big\|\div\, u
           \big\|_{L^1(\re_+;\dB^{\frac{n}{p}}_{p,1}(\re^n_+))}
           \Big\| \int_0^t D^2 u\,ds                        
           \Big\|_{L^{\infty}(\re_+;\dB^{-1+\frac{n}{p}}_{p,1}(\re^n_+))}
         \bigg)
   \\
   \le & C\sum_{k=1}^{n-1}
           \| D^2 u\|_{L^1(\re_+;\dB^{-1+\frac{n}{p}}_{p,1}(\re^n_+))}^{k+1}.
    \eqntag \label{eqn;G_div-nonlinear}
}
}
This shows \eqref{eqn;m-bilinear-apriori-div}.

\vskip2mm\noindent
Finally we show \eqref{eqn;m-bilinear-apriori-div-dt}.
$\Pi_{\div}^{n-1}(\cdot)$ is the $n-1$-th order polynomial of 
the component of inverse of Jacobi matrix $J(Du)^{-1}$ without a constant term.   
By the div-curl structure it holds that
\algn{ \eqntag\label{eqn;div-curl-G_div}
  G_{\div}(u)
   = &
   \Pi_{\div}^{n-1}\Big(\int_0^t D u\,ds\Big)\div u
   =\div_k\bigg(
     \sum_{k=1}^{n-1}c_k\prod_{m,\ell\le n}^k
     \Big(\int_0^t \pt_{\ell} u_mds\Big)
     u_k\bigg) 
  =\div \overline{G_{\div}}(u).
}

%
%
%
%
%
%
%
%

By using \eqref{eqn;bilinear3} several times, we have
{ \allowdisplaybreaks
 \algn{
  &\|\pt_t (-\Del)^{-1}\N G_{\div}(u)\|_{L^1(\re_+;\dB^{-1+\frac{n}{p}}_{p,1}(\re^n_+))} \\
  \le & \|\pt_t  \overline{G_{\div}}(u)\|_{L^1(\re_+;\dB^{-1+\frac{n}{p}}_{p,1}(\re^n_+))} \\
  \le & C\left\|\Pi_{\div}^{n-1}
         \Big(\int_0^t D u\, ds\Big)                
                \pt_t  u 
         \right\|_{L^1(\re_+;\dB^{-1+\frac{n}{p}}_{p,1}(\re^n_+))} 
        +C\left\|\Pi_{\div}^{n-2}
                \big(\int_0^t D u\,ds\Big)         
                \N u \otimes  u 
         \right\|_{L^1(\re_+;\dB^{-1+\frac{n}{p}}_{p,1}(\re^n_+))}   
\\
  \le & C\sup_{t>0}
         \Big\|\int_0^t D u\,ds\Big\|_{\dB^{\frac{n}{p}}_{p,1}(\re^n_+)}
         \Big\|\Pi_{\div}^{n-2}
                \Big(\int_0^t D u\,ds\Big)
                 \pt_t\, u         
         \Big\|_{L^1(\re_+;\dB^{-1+\frac{n}{p}}_{p,1}(\re^n_+))}\\
     & +C\sup_{t>0}
         \Big\|\int_0^t D u\,ds\Big\|_{\dB^{\frac{n}{p}}_{p,1}(\re^n_+)}
         \Big\|\Pi_{\div}^{n-2}
                \Big(\int_0^t D u\,ds\Big)
                \N u \otimes u        
         \Big\|_{L^1(\re_+;\dB^{-1+\frac{n}{p}}_{p,1}(\re^n_+))}
   \\
   \le & C\sum_{k=1}^{n-1}
          \left\| D^2  u \right\|_{L^1(\re_+;\dB^{-1+\frac{n}{p}}_{p,1}(\re^n_+))}^k
          \| \pt_t  u \|_{L^1(\re_+;\dB^{-1+\frac{n}{p}}_{p,1}(\re^n_+))}\\
      &+C\sum_{k=0}^{n-2}
          \left\| D^2  u \right\|_{L^1(\re_+;\dB^{-1+\frac{n}{p}}_{p,1}(\re^n_+))}^k
          \| u\,  D u     \|_{L^1(\re_+;\dB^{-1+\frac{n}{p}}_{p,1}(\re^n_+))}
 \\
   \le & C\sum_{k=1}^{n-1}
          \left\| D^2  u \right\|_{L^1(\re_+;\dB^{-1+\frac{n}{p}}_{p,1}(\re^n_+))}^k
          \| \pt_t  u \|_{L^1(\re_+;\dB^{-1+\frac{n}{p}}_{p,1}(\re^n_+))}\\
      &+C\sum_{k=0}^{n-2}
          \left\| D^2  u \right\|_{L^1(\re_+;\dB^{-1+\frac{n}{p}}_{p,1}(\re^n_+))}^k
          \left\|\N    u \right\|_{L^1(\re_+;\dB^{\frac{n}{p}}_{p,1}(\re^n_+))}
          \sup_{t>0} \|  u \|_{\dB^{-1+\frac{n}{p}}_{p,1}(\re^n_+)}
\\
   \le & C\sum_{k=1}^{n-1}
          \left\| D^2  u \right\|_{L^1(\re_+;\dB^{-1+\frac{n}{p}}_{p,1}(\re^n_+))}^k
        \Big(  \| \pt_t  u \|_{L^1(\re_+;\dB^{-1+\frac{n}{p}}_{p,1}(\re^n_+))}
                +\sup_{t>0} 
                 \left\| \int_t^{\infty}\pt_t  u ds
                 \right\|_{\dB^{-1+\frac{n}{p}}_{p,1}(\re^n_+)}\Big)
\\
   \le & C\sum_{k=1}^{n-1}
          \left\| D^2  u \right\|_{L^1(\re_+;\dB^{-1+\frac{n}{p}}_{p,1}(\re^n_+))}^k
          \| \pt_t  u \|_{L^1(\re_+;\dB^{-1+\frac{n}{p}}_{p,1}(\re^n_+))}.
           \eqntag \label{eqn;G_div-nonlinear2}
 }}
\end{prf}
\vskip3mm
For the quasi-linear term $F_u(u)$ associated with the Laplace operator
defined in \eqref{eqn;laplace-purterb}, the structure is slightly 
different from the others since total {\it multiple-div-structure} 
fails in the following whole form:
\eqn{
 F_u(u)= \div \Big(\Big(J(Du)^{-1}\big(J(Du)^{-1}\big)^{\sf T}-I\Big)\N u\Big).
}
Indeed, the {\it multiple-div-curl structure} remains valid for  
$\big(J(Du)^{-1}-I\big)^{\sf T}\N u$ partially,  as well as  \break
for $\div \Big(J(Du)^{-1}\; F\Big)$ 
with any vector field $F$.  However since the coefficient function is 
$J(Du)^{-1}$ which is adjoint of $\big(J(Du)^{-1}\big)^{\sf T}$, 
the derivatives for the divergence `$\div$'  outside does not 
commute with $J(Du)^{-1}$ and the whole {\it multiple div-curl structure}
does not hold.  To recover this difficulty, we use the  bilinear estimate 
in Proposition \ref{prop;abidi-paicu-halfspace}.

\vskip3mm

\begin{prop}[Multiple div-curl estimate 2]
\label{prop;multiple-bilinear-2}
Let $n\ge 2$, $1\le p<\infty$.  
For $D^2 u\in L^1(\re_+;\dB^{-1+n/p}_{p,1}(\re^n_+))$, 
let $F_{u}(u)$ be  the terms defined 
by \eqref{eqn;laplace-purterb}. Then the following 
estimate holds:
\begin{equation} \label{eqn;m-bilinear-apriori-2nd}
  \|F_{u}(u)\|_{L^1(\re_+;\dB^{-1+\frac{n}{p}}_{p,1})}
    \le C\sum_{k=1}^{2n-2}
         \|D^2 u\|_{L^1(\re_+;\dB^{-1+\frac{n}{p}}_{p,1}(\re^n_+))}^{k+1}.
\end{equation}
\end{prop}
\vskip2mm

\begin{prf}{Proposition \ref{prop;multiple-bilinear-2}}
To show the estimate \eqref{eqn;m-bilinear-apriori-2nd},
we first decompose the terms as
\algn{
 F_u(u)=\, & \div \Big(J(Du)^{-1}\big(J(Du)^{-1}\big)^{\sf T} \N u -\N u\Big) 
\\
    =\, & \div \Big(J(Du)^{-1}\big(J(Du)^{-1}-I\big)^{\sf T}
        \; \N u\Big)
      +\div \Big(\big(J(Du)^{-1}-I\big)\N u\Big)
\\
    =\, 
     & \div \Big(\big(J(Du)^{-1}-I\big)
        \big(J(Du)^{-1}-I\big)^{\sf T}\; \N u\, \Big)
      +\div \Big(\big(J(Du)^{-1}-I\big)^{\sf T}\N u\Big) \\
     &+\div \Big(\big(J(Du)^{-1}-I\big)\N u\Big) \\
 \equiv\, & F_u^1(u)+F_u^2(u)+F_u^3(u).    
}
Here $\div E$ stands for $[\N^{\sf T}E]^{\sf T}$, where $E$ denotes 
the $n\times n$-matrix valued function.
To show the estimate  $F_u^1$, one can use the {\it div-structure}
up to estimate for the terms for $\div\Big(\big(J(Du)^{-1}-I\big)F\Big)$
for the vector field $F\equiv  (J(Du)^{-1}-I\big)^{\sf T}\N u$ since 
it maintains the {\it multiple-div-structure}.
Namely since $\div\Big(\big(J(Du)^{-1}-I\big) F\Big)$ is 
the adjoint operator of $\big(J(Du)^{-1}-I\big)^{\sf T}\N F$, 
it maintains the structure and it follows 
that for any $1\le p<\infty$, 
\eqn{
 \spl{
     \Big\|\div\Big(\big(J(Du)^{-1}-I\big)\;F\Big)
     \Big\|_{\dot{B}^{-1+\frac{n}{p}}_{p,1}(\re^n_+)}
     \le C \sum_{k=1}^{n-1}
           \|\int_0^t \N u ds\|_{\dot{B}^{\frac{n}{p}}_{p,1}(\re^n_+)}^k
           \|\div F\|_{\dot{B}^{-1+\frac{n}{p}}_{p,1}(\re^n_+)}
 }
}
and 
\eq{ \label{eqn;quasi-linear-estimate-1}
 \spl{
    \| 
      \div\Big(\big( &J(Du)^{-1}-I\big)\;F\Big)
    \|_{L^1(\re_+;\dot{B}^{-1+\frac{n}{p}}_{p,1}(\re^n_+))}
\\
   \le& C\sum_{k=1}^{n-1}c_k
          \|\N u\|_{L^1(\re_+;\dB^{\frac{n}{p}}_{p,1}(\re^n_+))}^k
          \|\div F\|_{L^1(\re_+;\dB^{-1+\frac{n}{p}}_{p,1}(\re^n_+))} 
\\
   \le& C\sum_{k=1}^{n-1}
         \|D^2  u\|_{L^1(\re_+;\dB^{-1+\frac{n}{p}}_{p,1}(\re^n_+))}^{k}
         \|\div F\|_{L^1(\re_+;\dB^{-1+\frac{n}{p}}_{p,1}(\re^n_+))}.
  }
}
Letting $F\equiv (J(Du)-1)^{\sf T}\N u$ and using Proposition 
\ref{prop;abidi-paicu-halfspace} and the div-curl bilinear estimate \eqref{eqn;bilinear3} in 
Proposition \ref{prop;Besov-bilinear-div-curl-half} for $n-2$-times,  
we see for the last term of the right hand 
side of \eqref{eqn;quasi-linear-estimate-1} that 
\algn{ 
     \|\div F\|_{L^1(\re_+;\dB^{-1+\frac{n}{p}}_{p,1}(\re^n_+))}
 =&
     \Big\|
        \div \Big(\big( J(Du)^{-1}-I\big)^{\sf T}\;\N u\Big)
     \Big\|_{L^1(\re_+;\dB^{-1+\frac{n}{p}}_{p,1}(\re^n_+))}
\\
    \le C & \big\| \big(J(Du)^{-1}-I\big)^{\sf T}\N u 
            \big\|_{ L^1(\re_+;\dB^{\frac{n}{p}}_{p,1}(\re^n_+)) } 
\\
    \le C & \sum_{k=1}^{n-1}c_k\Big\|\int_0^t \N u ds 
            \Big\|_{L^{\infty}(\re_+;\dB^{\frac{n}{p}}_{p,1}(\re^n_+) ) }^{k} 
            \|\N u\|_{L^1(\re_+;\dB^{\frac{n}{p}}_{p,1}(\re^n_+))} 
\\
    \le C & \sum_{k=1}^{n-1}
            \big\|D^2 u \big\|_{ L^1(\re_+;\dB^{-1+\frac{n}{p}}_{p,1}(\re^n_+) )}^{k} 
            \|D^2 u\|_{L^1(\re_+;\dB^{-1+\frac{n}{p}}_{p,1}(\re^n_+))}.   
        \eqntag\label{eqn;quasi-linear-estimate-2}
}
Combining  
\eqref{eqn;quasi-linear-estimate-1},  
\eqref{eqn;quasi-linear-estimate-2}, the estimates for 
$F_u^1(u)$ and $F_u^2(u)$ are proven.
It is then easy to see that a similar argument of 
\eqref{eqn;quasi-linear-estimate-2} can be 
applicable for estimating the last term $F_u^3(u)$ and we conclude that
\eqn{
 \spl{
  \|F_u(u)\|_{L^1(\re_+;\dB^{-1+\frac{n}{p}}_{p,1}(\re^n_+))} 
  \le & C\sum_{k=1}^{2n-2}
         \|D^2 u\|_{L^1(\re_+;\dB^{-1+\frac{n}{p}}_{p,1}(\re^n_+))}^{k+1}.
 }}
\end{prf}
\vskip 3mm

We finally treat the boundary nonlinearity as follows.
\vskip2mm
\begin{prop}[Multiple estimates for boundary nonlinearity]
\label{prop;multiple-bilinear-3}
Let $n\ge 2$,  $1< p<2n-1$ 
{\rd and} assume that functions $u$ and $p$ {\crd satisfy} $\pt_t u$, $D^2 u$, 
$\N p \in L^1(\re_+;{\rd \dB^{-1+n/p}_{p,1}}(\re^n_+))$ 
with  
$p|_{x_n=0}
 \in 
 \dF^{1/2-1/2p}_{1,1}(\re_+;\dB^{-1+n/p}_{p,1}(\re^{n-1}))
 \cap
 {\crd L^1}(\re_+;{\rd \dB^{(n-1)/p}_{p,1}}(\re^{n-1}))$.
For the boundary terms  $H_{b}(u)$ and $H_p(u,p)$ defined 
by \eqref{eqn;boundary-purterb-u}  and \eqref{eqn;boundary-purterb-p}, 
respectively. 
Then the following estimates hold:
\algn{
   \|H_{p}(u,p)&\|_{\dF^{\frac12-\frac{1}{2p}}_{1,1}(\re_+;\dB^{-1+\frac{n}{p}}_{p,1}(\re^{n-1}))} \\
    \le&
       C\Big(
        \big\| p|_{x_n=0}
        \big\|_{\dF^{\frac12-\frac{1}{2p}}_{1,1}(\re_+;\dB^{-1+\frac{n}{p}}_{p,1}(\re^{n-1}))}
 {\crd+ \big\| p|_{x_n=0}
        \big\|_{L^1(\re_+;\dB^{-1+\frac{n}{p}}_{p,1}(\re^{n-1}))}
  }
        \Big) \\
     &\hskip2cm\times
        \sum_{k=1}^{n-1}    
         \Big(\|\pt_t u\|_{L^1(\re_+;\dB^{-1+\frac{n}{p}}_{p,1}(\re^n_+))}
                +
             \|D^2 u\|_{L^1(\re_+;\dB^{-1+\frac{n}{p}}_{p,1}(\re^n_+))}
         \Big)^k, 
       \eqntag \label{eqn;m-bilinear-apriori-b-T-p}
\\
  \|H_{p}(u,p)&\|_{{\crd L^1}(\re_+; \dB^{\frac{n-1}{p}}_{p,1}(\re^{n-1}))}
    \le 
       C\big\| p|_{x_n=0}
         \big\|_{{\crd L^1}(\re_+;\dB^{\frac{n-1}{p}}_{p,1}(\re^{n-1}))}
        \sum_{k=1}^{n-1}    
         {\crd
              \|D^2   u\|_{L^1(\re_+;\dB^{-1+\frac{n}{p}}_{p,1}(\re^n_+))}^k,
         }
     \eqntag  \label{eqn;m-bilinear-apriori-b-S-p} 
}
\algn{ 
  \|H_{u}(u)\|_{\dF^{\frac12-\frac{1}{2p}}_{1,1}\re_+;\dB^{-1+\frac{n}{p}}_{p,1}(\re^{n-1}))} 
    \le&C
       \sum_{k=2}^{2n-1}
         \Big(\|\pt_t u\|_{L^1(\re_+;\dB^{-1+\frac{n}{p}}_{p,1}(\re^n_+))}
             +\|D^2 u  \|_{L^1(\re_+;\dB^{-1+\frac{n}{p}}_{p,1}(\re^n_+))}
         \Big)^k,
       \eqntag \label{eqn;m-bilinear-apriori-b-T-u}
\\
  \| H_{u}(u)\|_{{\crd L^1}(\re_+; \dB^{\frac{n-1}{p}}_{p,1}(\re^{n-1}))}  
  \le& C 
       \sum_{k=2}^{2n-1}
         {\crd 
          \| D^2 u \|_{L^1(\re_+;\dB^{-1+\frac{n}{p}}_{p,1}(\re^n_+))}^k.
         }
       \eqntag \label{eqn;m-bilinear-apriori-b-S-u}
 }
\end{prop}
\vskip2mm

%
%
\begin{prf}{Proposition \ref{prop;multiple-bilinear-3}} 
From \eqref{eqn;boundary-purterb-u} and  \eqref{eqn;boundary-purterb-p} 
and from the regularity assumptions;
we notice that the sharp trace estimate implies 
\algn{
 &Du\,\big|_{x_n=0},\; p\,\big|_{x_n=0} \in 
  \dF^{\frac12-\frac{1}{2p}}_{1,1}(\re_+;\dB^{-1+\frac{n}{p}}_{p,1}(\re^{n-1})) 
  \cap
  {\crd 
  L^1(\re_+;\dB^{\frac{n-1}{p}}_{p,1}(\re^{n-1}).
  }
  \eqntag\label{eqn;p-reg-1}
 }

We first see the estimate \eqref{eqn;m-bilinear-apriori-b-T-p}.
Since $J(Du)^{-1}-I$ consists of a polynomial of $\int_0^t Du\,ds$ with 
its order up to $n-1$,
we show that 
\eq{\label{eqn;crutial-nonlinear-est2}
  H_p(u,p)
   = 
  \Pi_{bp}^{n-1}\left(\int_0^t D u\,ds\right) p\, \nu_n
   \in \dF^{\frac12-\frac{1}{2p}}_{1,1}
       \big(\re_+;\dB^{-1+\frac{n}{p}}_{p,1}(\re^{n-1})\big).
}
\vskip1mm
{\crd 
We introduce auxiliary norms of Chemin--Lerner type {\bl (cf. \cite{C-L95})} for the proof of
Proposition \ref{prop;multiple-bilinear-3}.
\vskip2mm
\noindent
{\it Definition.} For $1\le p, \r\le \infty$ and $r,s\in \re$, 
the Bochner--Besov spaces of
Chemin--Lerner type  
$\widetilde{\dB^r_{\r,1}\big(\re_+};\dB^{s}_{p,1}(\re^{n-1})\big)$
and
$\widetilde{L^{\rho}\big(\re_+};\dB^{s}_{p,1}(\re^{n-1})\big)$ 
are defined by the following norms:
\eq{ \label{eqn;chemin-lerner-norms}
 \spl{
  \big\|f\big\|_{\widetilde{\dB^r_{\r,1}\big(\re_+};\dB^{s}_{p,1}(\re^{n-1})\big)}
  \equiv &
   \sum_{k\in \Z} 2^{rk} \sum_{j\in\Z} 2^{sj}  
        \big\| 
           \psi_k\underset{(t)}{*}
           \phi_j\underset{(x')}{*}
           f(t,x')
       \big\|_{L^{\rho}_t(\re_+;L^p({\gr \re^{n-1}}))}, 
\\
 \big\|f\big\|_{\widetilde{L^{\r}\big(\re_+};\dB^{s}_{p,1}(\re^{n-1})\big)}
  \equiv &
  \sum_{j\in \Z}2^{sj}
       \big\| 
          {\bl \phi_j\underset{(x')}{*} }\,
           f(t,x')
       \big\|_{L^{\rho}_t(\re_+;L^p({\gr \re^{n-1}}))}.
  } 
}
\vskip2mm
\begin{lem}[Multiple estimates for boundary nonlinearity]
\label{lem;space-time-bilinear-0}
Let $n\ge 2$,  $1< p<2n-1$ and assume that functions 
$F$ and $G$ over  $\re_+\times \re^{n-1}$ satisfy 
$
F\in 
 \dF^{1/2-1/2p}_{1,1}(\re_+;\dB^{-1+n/p}_{p,1}(\re^{n-1}))
 \cap
 L^1(\re_+;\dB^{(n-1)/p}_{p,1}(\re^{n-1}))$
and 
$
 G\in 
 \widetilde{\dB^{1/2-1/2p}_{\infty,1}(\re_+};\dB^{-1+n/p}_{p,1}(\re^{n-1}))
 \cap
 \widetilde{L^{\infty}(\re_+};\dB^{(n-1)/p}_{p,1}(\re^{n-1}))$.
Then the following estimate holds:
\algn{
   \|F\, G\|_{\dF^{\frac12-\frac{1}{2p}}_{1,1}(\re_+;\dB^{-1+\frac{n}{p}}_{p,1}(\re^{n-1}))}
    \le& 
      C\bigg(
        \big\| F
        \big\|_{\dF^{\frac12-\frac{1}{2p}}_{1,1}(\re_+;\dB^{-1+\frac{n}{p}}_{p,1}(\re^{n-1}))}
     +  \big\| F
        \big\|_{L^1(\re_+;\dB^{\frac{n-1}{p}}_{p,1}(\re^{n-1}))}
       \bigg)\\
    &\hskip5mm\times
    \bigg(
        \big\|G    
        \big\|_{\widetilde{\dB^{\frac12-\frac{1}{2p}}_{\infty,1}
           (\re_+;}\dB^{-1+\frac{n}{p}}_{p,1}({\gr \re^{n-1}}))}       
      + \big\|
           G
        \big\|_{\widetilde{L^{\infty}(\re_+;}\dB^{\frac{n-1}{p}}_{p,1}(\re^{n-1}))}
   \bigg).
       \eqntag \label{eqn;m-bilinear-b-T-p}
 }
\end{lem}
\vskip2mm
The proof of Lemma \ref{lem;space-time-bilinear-0} directly follows from Proposition 
\ref{prop;double-bony-bilinear} with $\r=1$ and \eqref{eqn;F-B-space-equiv-1}--\eqref{eqn;F-B-space-equiv-2} in Appendix below.

\vskip3mm
Now we set 
\algn{
 &F(t,x') \equiv p(t,x',x_n)|_{x_n=0}, \\
 &G(t,x') \equiv \Pi^{n-1}_{bp}\Big(\int_0^tDu(s, x', x_n)ds\Big)\Big|_{x_n=0}
}
in Lemma \ref{lem;space-time-bilinear-0} with regarding 
\eqref{eqn;p-reg-1} to find that 
\algn{
   \|H_p(u,p)&\|_{\dF^{\frac12-\frac{1}{2p}}_{1,1}(\re_+;\dB^{-1+\frac{n}{p}}_{p,1}(\re^{n-1}))}   \\
    \le& 
      C\bigg(
        \big\| p|_{x_n=0}
        \big\|_{\dF^{\frac12-\frac{1}{2p}}_{1,1}(\re_+;\dB^{-1+\frac{n}{p}}_{p,1}(\re^{n-1}))}
     +  \big\| p|_{x_n=0}
        \big\|_{L^1(\re_+;\dB^{\frac{n-1}{p}}_{p,1}(\re^{n-1}))}
       \bigg)\\
    &\hskip5mm\times
    \bigg(
        \Big\|\Pi^{n-1}_{bp}\Big(\int_0^tDu(s, x', x_n)ds\Big)\Big|_{x_n=0}   
        \Big\|_{\widetilde{\dB^{\frac12-\frac{1}{2p}}_{\infty,1}
           (\re_+;}\dB^{-1+\frac{n}{p}}_{p,1}({\gr \re^{n-1}}))}   \\
    &\hskip12mm     
      + \Big\|
           \Pi^{n-1}_{bp}\Big(\int_0^tDu(s, x', x_n)ds\Big)\Big|_{x_n=0}
        \Big\|_{\widetilde{L^{\infty}(\re_+;}\dB^{\frac{n-1}{p}}_{p,1}(\re^{n-1}))}
   \bigg).
       \eqntag \label{eqn;Hp(u,p)-1}
 }

\vskip1mm
To complete the estimate we use the following lemma.
\vskip1mm

\vskip2mm
\noindent
\begin{lem}\label{lem;product-estimates}
For any $1\le p<\infty$,
\algn{
  \Big\| \int_0^t Du(s)\,ds\,\Big|_{x_n=0}    
  \Big\|_{\widetilde{\dB^{\frac12-\frac{1}{2p}}_{\infty,1}
           (\re_+;}\dB^{-1+\frac{n}{p}}_{p,1}({\gr \re^{n-1}}))}
   \le&\, C\big\|Du|_{x_n=0}\big\|_{\dF^{\frac12-\frac{1}{2p}}_{1,1}
           (\re_+;\dB^{-1+\frac{n}{p}}_{p,1}({\gr \re^{n-1}}))}, 
    \eqntag \label{eqn;Du-1}
 \\
  \Big\| \int_0^t Du(s)\,ds\, \Big|_{x_n=0}    
  \Big\|_{\widetilde{L^{\infty}
           (\re_+;}{\bl \dB^{\frac{n-1}{p}}_{p,1}}({\gr \re^{n-1}}))}
   \le&\, C\big\|Du|_{x_n=0}\big\|_{L^1
           (\re_+;{\bl \dB^{\frac{n-1}{p}}_{p,1}}({\gr \re^{n-1}}))}.
    \eqntag \label{eqn;Du-2}
  }
\end{lem}
\vskip2mm

\begin{prf}{Lemma \ref{lem;product-estimates}}
The first estimate \eqref{eqn;Du-1} follows by 
using $\widetilde{\psi}_k(t)=\psi_{k-1}(t)+\psi_k(t)+\psi_{k+1}(t)$ 
and noticing $\|\pt_t^{-1}\psi_k\|_{L^{\infty}(\re_+)}
\le \|\psi_k\|_{L^1(\re_+)}$, where $\pt_t^{-1}\psi_k$ is defined as in 
\eqref{eqn;int-psi-k} below that
\algn{
  \Big\| \int_0^t & Du(s)\,ds\,\Big|_{x_n=0}   
  \Big\|_{\widetilde{\dB^{\frac12-\frac{1}{2p}}_{\infty,1}
           (\re_+;}\dB^{-1+\frac{n}{p}}_{p,1}({\gr \re^{n-1}}))}  \\
   =& \sum_{j\in \Z} 2^{(-1+\frac{n}{p})j} 
       \sum_{k\in \Z} 2^{(\frac12-\frac{1}{2p})k}   
     \bigg\|  
         \Big\|  
          \psi_{k}\underset{(t)}{*} 
          \phi_{j}\underset{(x')}{*}
          \Big(\int_0^t Du(s)\,ds\,\Big|_{x_n=0}\Big)
          \Big\|_{L^{p}({\gr \re^{n-1}})}         
     \bigg\|_{L^{\infty}_t(\re_+)} 
\\
  \le& \sum_{j\in \Z} 2^{(-1+\frac{n}{p})j} 
       \sum_{k\in \Z} 2^{(\frac12-\frac{1}{2p})k}   
     \bigg\|   
       \Big\|
          (\pt_t^{-1}\psi_{k})\underset{(t)}{*}
          \widetilde{\psi_{k}}\underset{(t)}{*} 
          \phi_{j}\underset{(x')}{*}
            \pt_t\Big(\int_0^t Du(s)\,|_{x_n=0}\,ds\Big)
        \Big\|_{L^{p}({\gr \re^{n-1}})}
      \bigg\|_{L^{\infty}_t(\re_+)} 
 \\
  \le& \sum_{j\in \Z} 2^{(-1+\frac{n}{p})j} 
       \sum_{k\in \Z} 2^{(\frac12-\frac{1}{2p})k}   
     \big\|   
          \pt_t^{-1}\psi_{k}
     \big\|_{L^{\infty}_t(\re_+)}
     \bigg\|
     \Big\|
          \widetilde{\psi_{k}}\underset{(t)}{*} 
          \phi_{j}\underset{(x')}{*}
           Du\,|_{x_n=0}
      \Big\|_{L^{p}({\gr \re^{n-1}})}
     \bigg\|_{L^{1}_t(\re_+)} 
 \\
   \le& \sum_{j\in \Z} 2^{(-1+\frac{n}{p})j} 
       \sum_{k\in \Z} 2^{(\frac12-\frac{1}{2p})k}   
     \big\|   
          \psi_{k}
     \big\|_{L^{1}_t(\re_+)}
     \bigg\|
     \Big\|
          \widetilde{\psi_{k}}\underset{(t)}{*} 
          \phi_{j}\underset{(x')}{*}
           Du\,|_{x_n=0}
      \Big\|_{L^{p}({\gr \re^{n-1}})}
     \bigg\|_{L^{1}_t(\re_+)} 
  \\
   \le&C     
    \bigg\|
        \sum_{j\in \Z} 2^{(-1+\frac{n}{p})j} 
        \sum_{k\in \Z} 2^{(\frac12-\frac{1}{2p})k}   
     \Big\|
          \widetilde{\psi_{k}}\underset{(t)}{*} 
          \phi_{j}\underset{(x')}{*}
           Du\,|_{x_n=0}
      \Big\|_{L^{p}({\gr \re^{n-1}})}
     \bigg\|_{L^{1}_t(\re_+)}
 \\
 \le&C\big\|
           Du\,|_{x_n=0}
    \big\|_{\dF^{\frac12-\frac{1}{2p}}_{1,1}(\re_+;\dB^{-1+\frac{n}{p}}_{p,1}({\gr \re^{n-1}}))}.
}
The second inequality \eqref{eqn;Du-2} follows from the following estimate:
\algn{
 \Big\|\int_0^t & Du(s)\,ds\,\Big|_{x_n=0}
 \Big\|_{\widetilde{L^{\infty}(\re_+;}\dB^{\frac{n-1}{p}}_{p,1}({\gr \re^{n-1}}))}  
  \le 
   \sum_{j\in\Z} 2^{\frac{n-1}{p}j}  
       \Big\|           
        \big\|
            \int_0^t \phi_j\underset{(x')}{*} Du(s)\,|_{x_n=0} ds
        \big\|_{L^p({\gr \re^{n-1}})}
       \Big\|_{L^{\infty}_t(\re_+)} 
 \\
 \le & \sum_{j\in\Z} 2^{\frac{n-1}{p}j}  
       \Big\|           
           \|\phi_j\underset{(x')}{*}Du\,|_{x_n=0}\|_{L^p({\gr \re^{n-1}})}
       \Big\|_{L^1_t(\re_+)} 
 \\
  = & \Big\| \sum_{j\in\Z} 2^{\frac{n-1}{p}j}          
           \|\phi_j\underset{(x')}{*}Du\,|_{x_n=0}\|_{L^p({\gr \re^{n-1}})}
       \Big\|_{L^1_t(\re_+)} 
 =\big\| Du\,|_{x_n=0} \big\|_{L^1(\re_+;\dB^{\frac{n-1}{p}}_{p,1}(\re^{n-1})}.
 }
\end{prf}

\vskip2mm

\noindent
\begin{lem}\label{lem;Banach-algebra-2}
$\widetilde{L^{\infty}(\re_+;}\dB^{(n-1)/p}_{p,1}(\re^{n-1}))$
is the Banach algebra, namely for any 
$f, g\in\widetilde{L^{\infty}(\re_+;}\dB^{(n-1)/p}_{p,1}(\re^{n-1}))$ it holds
\algn{
 \big\|f\, g
 \big\|_{\widetilde{L^{\infty}(\re_+;}\dB^{\frac{n-1}{p}}_{p,1}({\gr \re^{n-1}}))}  
 \le&
     C  \big\|f
        \big\|_{\widetilde{L^{\infty}(\re_+;}\dB^{\frac{n-1}{p}}_{p,1}(\re_+)} 
        \big\|
             g
        \big\|_{\widetilde{L^{\infty}(\re_+;}\dB^{\frac{n-1}{p}}_{p,1}(\re_+)}. 
  \eqntag\label{eqn;banach-algebra}
}
In particular for 
$Du|_{x_n=0}\in L^1(\re_+;\dB^{\frac{n-1}{p}}_{p,1}(\re^{n-1}))$,
\alg{
\Big\|\Pi^{n-1}_{bp}\Big(\int_0^t Du(s)\, ds\Big)\,\Big|_{x_n=0}
\Big\|_{\widetilde{L^{\infty}(\re_+;}\dB^{\frac{n-1}{p}}_{p,1}({\gr \re^{n-1}}))}  
 \le& \sum_{k=1}^{n-1}
      \big\| Du\,|_{x_n=0} 
      \big\|_{L^1(\re_+;\dB^{\frac{n-1}{p}}_{p,1}(\re^{n-1})}^k. 
 \label{eqn;Hp(u,p)-2}
 }
\end{lem}

\begin{prf}{Lemma \ref{lem;Banach-algebra-2}}
To see that \eqref{eqn;banach-algebra} holds, we start from 
Bony's paraproduct decomposition:
Setting 
$P_m\underset{(x')}{*}F\equiv \sum_{m'\le m}\phi_{m'}\underset{(x')}{*}F$
and  
$Q_\ell\underset{(t)}{*}F\equiv \sum_{\ell'\le \ell}\psi_{\ell'}\underset{(t)}{*}F$,
it follows that 
\algn{
  \big\|f&\, g 
  \big\|_{\widetilde{L^{\infty}(\re_+;}\dB^{\frac{n-1}{p}}_{p,1}({\gr \re^{n-1}}))}   
  \\
  \le&\sum_{j\in \Z} 2^{\frac{n-1}{p}j}
        \Big\| \|\phi_j\|_1
        \big\| 
              \widetilde{\phi}_{j}\underset{(x')}{*}f
        \big\|_{L^p(\re^{n-1}_{x'})}
        \big\|
             P_{j-2}\underset{(x')}{*}g
        \big\|_{L^{\infty}({\gr \re^{n-1}})}
       \Big\|_{L^{\infty}_t(\re_+)}    \\
    & +\sum_{j\in \Z} 2^{\frac{n-1}{p}j}
       \Big\|\|\phi_j\|_1
        \big\| 
             P_{j-2}\underset{(x')}{*}f
        \big\|_{L^{\infty}({\gr \re^{n-1}})}
        \big\|
             \widetilde{\phi}_{j}\underset{(x')}{*}g
        \big\|_{L^{p}({\gr \re^{n-1}})}
       \Big\|_{L^{\infty}_t(\re_+)}   \\ 
     &\quad +\sum_{j\in \Z} 2^{\frac{n-1}{p}j}
       \Big\|\|\phi_j\|_r
       \sum_{m\ge j-2}
        \big\| 
             \phi_{m}\underset{(x')}{*} f
        \big\|_{L^p({\gr \re^{n-1}})}
        \big\|
             \widetilde{\phi}_{m}\underset{(x')}{*} g
        \big\|_{L^{r'}({\gr \re^{n-1}})}
       \Big\|_{L^{\infty}_t(\re_+)}    
\\ 
  \le&\sum_{j\in \Z} 2^{\frac{n-1}{p}j}
        \Big\|  
          \big\|
          \widetilde{\phi}_{j}\underset{(x')}{*} f
          \big\|_{L^{p}({\gr \re^{n-1}})}
        \Big\|_{L^{\infty}_t(\re_+)} 
        \Big\|
        \big\| P_{j-2}\big\|_{L^{1}(\re^{n-1})}
        \big\|
            g
        \big\|_{L^{\infty}({\gr \re^{n-1}})}
       \Big\|_{L^{\infty}_t(\re_+)}    \\
    & +\sum_{j\in \Z} 2^{\frac{n-1}{p}j}
       \Big\|\big\| P_{j-2}\big\|_{L^{1}(\re^{n-1})}
        \big\| 
           f
        \big\|_{L^{\infty}({\gr \re^{n-1}})}
        \Big\|_{L^{\infty}_t(\re_+)}
        \Big\|
        \big\|
             \widetilde{\phi}_{j}\underset{(x')}{*} g
        \big\|_{L^{p}({\gr \re^{n-1}})}
       \Big\|_{L^{\infty}_t(\re_+)}   \\ 
   &\quad +\sum_{j\in \Z} 2^{\frac{n-1}{p}j} \sum_{m\ge j-2} 2^{\frac{n-1}{r'}j} \\
   &\hskip1.3cm \times
       \Big\|       
        \big\| 
             \phi_{m}\underset{(x')}{*}f
        \big\|_{L^p({\gr \re^{n-1}})}
        \Big\|_{L^{\infty}_t(\re_+)}
        \Big\|
        \big\|
             \widetilde{\phi}_{m}\underset{(x')}{*}g
        \big\|_{L^{r'}({\gr \re^{n-1}})}
       \Big\|_{L^{\infty}_t(\re_+)}      
\\ 
  \le&
       \sum_{j\in \Z} 2^{\frac{n-1}{p}j} \Big\|
           \big\|\widetilde{\phi}_{j}\underset{(x')}{*}
                 f
           \big\|_{L^p({\gr \re^{n-1}})}
        \Big\|_{L^{\infty}_t(\re_+)} 
        \Big\|
        \big\|
            g
        \big\|_{L^{\infty}({\gr \re^{n-1}})}
       \Big\|_{L^{\infty}_t(\re_+)}    \\
    & +
       \Big\|
        \big\| 
            f
        \big\|_{L^{\infty}({\gr \re^{n-1}})}
        \Big\|_{L^{\infty}_t(\re_+)}
        \sum_{j\in \Z} 2^{\frac{n-1}{p}j}
        \Big\|
        \big\|
             \widetilde{\phi}_{j}\underset{(x')}{*}g
        \big\|_{L^{p}({\gr \re^{n-1}})}
       \Big\|_{L^{\infty}_t(\re_+)}   \\ 
   &\quad +\sum_{m\in \Z}  \sum_{j\le m+2} 2^{\frac{2(n-1)}{p}j}
        \Big\|       
        \big\| 
            \phi_{m}\underset{(x')}{*}f
        \big\|_{L^p({\gr \re^{n-1}})}
        \Big\|_{L^{\infty}_t(\re_+)}
        \Big\|
        \big\|
             \widetilde{\phi}_{m}\underset{(x')}{*}g
        \big\|_{L^{p}({\gr \re^{n-1}})}
       \Big\|_{L^{\infty}_t(\re_+)}      
\\ 
  \le&
        \big\|f
        \big\|_{\widetilde{L^{\infty}(\re_+;}\dB^{\frac{n-1}{p}}_{p,1}(\re_+)} 
        \big\|
            g
        \big\|_{L^{\infty}(\re_+;L^{\infty}({\gr \re^{n-1}}))}    
     +
       \big\|
            f
       \big\|_{L^{\infty}(\re_+;L^{\infty}({\gr \re^{n-1}}))} 
       \big\|
            g
       \big\|_{\widetilde{L^{\infty}(\re_+;}\dB^{\frac{n-1}{p}}_{p,1}(\re_+)}  \\ 
   &\quad +\sum_{m\in \Z}  2^{\frac{(n-1)}{p}m}
        \Big\|       
        \big\| 
             \phi_{m}\underset{(x')}{*}f
        \big\|_{L^p({\gr \re^{n-1}})}
        \Big\|_{L^{\infty}_t(\re_+)}
        \sup_{m\in \Z} 2^{\frac{(n-1)}{p}m}
        \Big\|
        \big\|
             \widetilde{\phi}_{m}\underset{(x')}{*}g
        \big\|_{L^{p}({\gr \re^{n-1}})}
        \Big\|_{L^{\infty}_t(\re_+)}      
\\ 
  \le&
     C  \big\|f
        \big\|_{\widetilde{L^{\infty}(\re_+;}\dB^{\frac{n-1}{p}}_{p,1}(\re_+)} 
        \big\|
            g
        \big\|_{L^{\infty}(\re_+;L^{\infty}({\gr \re^{n-1}}))}    
    + C \big\|
            f
       \big\|_{L^{\infty}(\re_+;L^{\infty}({\gr \re^{n-1}}))} 
       \big\|
            g
       \big\|_{\widetilde{L^{\infty}(\re_+;}\dB^{\frac{n-1}{p}}_{p,1}(\re_+)}  \\ 
   &\quad +       
      C \big\| 
             f
        \big\|_{\widetilde{L^{\infty}(\re_+;}\dB^{\frac{n-1}{p}}_{p,1}(\re_+)}
        \big\|
             g
        \big\|_{\widetilde{L^{\infty}(\re_+;}\dB^{\frac{n-1}{p}}_{p,\infty}(\re_+)}     
\\ 
  \le&
     C  \big\|f
        \big\|_{\widetilde{L^{\infty}(\re_+;}\dB^{\frac{n-1}{p}}_{p,1}(\re_+)} 
        \big\|
             g
        \big\|_{\widetilde{L^{\infty}(\re_+;}\dB^{\frac{n-1}{p}}_{p,1}(\re_+)},  
  }
which shows \eqref{eqn;banach-algebra}.
The estimate \eqref{eqn;Hp(u,p)-2} follows from the estimate \eqref{eqn;Du-2}.
\end{prf}
\vskip2mm
The polynomial term can be estimated as the following way:
From Proposition \ref{prop;double-bony-bilinear}  in Appendix 
{\bl with  $\r=\infty$}, 
Lemma \ref{lem;Banach-algebra-2} 
and {\bl \eqref{eqn;Du-1} in} Lemma \ref{lem;product-estimates}
we see that 
\algn{
  \Big\|\Pi^{n-1}_{bp}\Big( \int_0^t & Du(s)ds \Big)\, \Big|_{x_n=0}   
  \Big\|_{\widetilde{\dB^{\frac12-\frac{1}{2p}}_{\infty,1}
           (\re_+;}\dB^{-1+\frac{n}{p}}_{p,1}({\gr \re^{n-1}}))}  \\
 \le &
  C\sum_{k=1}^{n-1}
    \bigg(\Big\|\int_0^t  Du(s)ds\, \Big|_{x_n=0}    
          \Big\|_{\widetilde{\bl\dB^{\frac12-\frac{1}{2p}}_{\infty,1}
               (\re_+;}\dB^{-1+\frac{n}{p}}_{p,1}({\gr \re^{n-1}}))}
        +\Big\|\int_0^t  Du(s)ds\, \Big|_{x_n=0}    
         \Big\|_{\widetilde{L^{\infty}
               (\re_+;}{\bl \dB^{\frac{n-1}{p}}_{p,1}}({\gr \re^{n-1}}))}
    \bigg)^k
 \\
 \le &C\sum_{k=1}^{n-1} 
 \bigg(
 \big\| Du \,|_{x_n=0}   
 \big\|_{\dF^{\frac12-\frac{1}{2p}}_{1,1}
               (\re_+;\dB^{-1+\frac{n}{p}}_{p,1}({\gr \re^{n-1}}))}
+\big\| Du \,|_{x_n=0}    
 \big\|_{L^{1}(\re_+;{\bl \dB^{\frac{n-1}{p}}_{p,1}}({\gr \re^{n-1}}))}
 \bigg)^k
\\
\le & C\sum_{k=1}^{n-1} 
 \bigg(
 \big\| \pt_t u   
 \big\|_{L^1(\re_+;\dB^{-1+\frac{n}{p}}_{p,1}({\gr \re^{n-1}}))}
+\big\| \Del u    
 \big\|_{L^{1}(\re_+;\dB^{-1+\frac{n}{p}}_{p,1}({\gr \re^{n-1}}))}
 \bigg)^k
\eqntag \label{eqn;Hp(u,p)-3}
}
by the sharp trace estimate Proposition \ref{prop;grad-pressure-trace}.
Combining the estimates 
\eqref{eqn;Hp(u,p)-1}%
-\eqref{eqn;Hp(u,p)-3},
we obtain \eqref{eqn;m-bilinear-apriori-b-T-p}.
}

\vskip2mm
To show the estimate \eqref{eqn;m-bilinear-apriori-b-S-p},
{\crd
we notice that $L^{\infty}(\re_+;\dB^{(n-1)/p}_{p,1}(\re^{n-1}))$ 
is the Banach algebra and from} the sharp trace estimate 
{\crd \eqref{eqn;pressure-trace-2} }, 
it follows from \eqref{eqn;pressure-trace-2} that 
{\crd
\algn{
 \|H_{p}&(u,p)\|_{L^1(\re_+;\dB^{\frac{n-1}{p}}_{p,1}(\re^{n-1}))} \\
    \le &C 
       \|p\,|_{x_n=0} \|_{L^1(\re_+;\dB^{\frac{n-1}{p}}_{p,1}(\re^{n-1}))}
       \Big\|\Pi^{n-1}_{bp}\Big(\int_0^t Du\,ds \Big)\,\Big|_{x_n=0} 
       \Big\|_{L^{\infty}(\re_+;\dB^{\frac{n-1}{p}}_{p,1}(\re^{n-1}))}
 \\
    \le &C 
       \|p \,|_{x_n=0} \|_{L^1(\re_+;\dB^{\frac{n-1}{p}}_{p,1}(\re^{n-1}))}
       \Big\|\int_0^t Du\,ds \,\Big|_{x_n=0} 
       \Big\|_{L^{\infty}(\re_+;\dB^{\frac{n-1}{p}}_{p,1}(\re^{n-1}))}
       \Big\|\Pi^{n-2}_{bp}\Big(\int_0^t Du\,ds \Big)\,\Big|_{x_n=0} 
       \Big\|_{L^{\infty}(\re_+;\dB^{\frac{n-1}{p}}_{p,1}(\re^{n-1}))}
 \\
     \le &C 
       \|p\,|_{x_n=0} \|_{L^1(\re_+;\dB^{\frac{n-1}{p}}_{p,1}(\re^{n-1}))}
       \sum_{k=1}^{n-1}
       \big\| Du\,|_{x_n=0} 
       \big\|_{L^{1}(\re_+;\dB^{\frac{n-1}{p}}_{p,1}(\re^{n-1}))}^k
 \\
 \le &C\|p\,|_{x_n=0} 
       \|_{L^1(\re_+;\dB^{\frac{n-1}{p}}_{p,1}(\re^{n-1}))}
       \sum_{k=1}^{n-1}
          \| D^2 u
          \|_{L^1(\re_+;\dB^{-1+\frac{n}{p}}_{p,1}(\re^n_+))}^k.
 \eqntag \label{eqn;boundary-nonlinear-5'}
}
}
\vskip2mm
On the other hand, for the {\crd estimate} 
\eqref{eqn;m-bilinear-apriori-b-T-u} of the velocity boundary 
term,   we split $H_u(u)$ into two parts 
as 
\algn{
  H_u(u) 
   =& \Big((J(Du)^{-1}\big)^{\sf T}\N u
               +(\N u)^{\sf T}J(Du)^{-1}
      \Big)\big((J(Du)^{-1})^{\sf T}-I\big) \nu_n \\
    &+\Big(\big(J(Du)^{-1}-I\big)^{\sf T}\N u
              +\big(\big(J(Du)^{-1}-I\big)^{\sf T}\N u\big)^{\sf T}
      \Big)\nu_n \\
  \equiv& H_u^1(u)+H^2_u(u).
  \eqntag \label{eqn;boundary-nonlinear-7}
}
{\crd
By setting 
\algn{
 &F(t,x') 
  \equiv \Big((J(Du)^{-1}\big)^{\sf T}\N u+(\N u)^{\sf T}J(Du)^{-1}
         \Big)\,\Big|_{x_n=0}, \\
 &G(t,x') \equiv \Pi^{n-1}_{bu}\Big(\int_0^tDu(s, x', x_n)ds\Big)\Big|_{x_n=0}
}
in Lemma \ref{lem;space-time-bilinear-0} with 
$$
 \big((J(Du)^{-1})^{\sf T}-I\big)\Big|_{x_n=0}
 =  \Pi^{n-1}_{bu}\Big(\int_0^tDu(s, x', x_n)ds\Big)\Big|_{x_n=0},
$$
we find that 
\algn{
   \|H_u^1&\|_{\dF^{\frac12-\frac{1}{2p}}_{1,1}(\re_+;\dB^{-1+\frac{n}{p}}_{p,1}(\re^{n-1}))}  
 \\
  \le& 
      C\bigg(
        \Big\|\Big((J(Du)^{-1}\big)^{\sf T}\N u+(\N u)^{\sf T}J(Du)^{-1}
        \Big)\,\Big|_{x_n=0}
        \Big\|_{\dF^{\frac12-\frac{1}{2p}}_{1,1}(\re_+;\dB^{-1+\frac{n}{p}}_{p,1}(\re^{n-1}))}  
        \\
     &\qquad
      +  \Big\| \Big((J(Du)^{-1}\big)^{\sf T}\N u+(\N u)^{\sf T}J(Du)^{-1}
         \Big)\,\Big|_{x_n=0}
        \Big\|_{L^1(\re_+;\dB^{\frac{n-1}{p}}_{p,1}(\re^{n-1}))}
       \bigg)\\
    &\hskip3mm\times
    \bigg(
        \Big\|\Pi^{n-1}_{bu}\Big(\int_0^tDu(s, x', x_n)ds\Big)\Big|_{x_n=0}   
        \Big\|_{\widetilde{\dB^{\frac12-\frac{1}{2p}}_{\infty,1}
           (\re_+;}\dB^{-1+\frac{n}{p}}_{p,1}({\gr \re^{n-1}}))}   \\
    &\hskip12mm     
      + \Big\|
           \Pi^{n-1}_{bu}\Big(\int_0^tDu(s, x', x_n)ds\Big)\Big|_{x_n=0}
        \Big\|_{\widetilde{L^{\infty}(\re_+;}\dB^{\frac{n-1}{p}}_{p,1}(\re^{n-1}))}
   \bigg).
      \eqntag \label{eqn;boundary-nonlinear-8}
 }
}  
The first term of the right hand side of \eqref{eqn;boundary-nonlinear-8}
is estimated by applying the sharp trace estimate \eqref{eqn;pressure-trace}
as well as a similar way in \eqref{eqn;G_div-nonlinear}, \eqref{eqn;G_div-nonlinear2} 
to obtain
\algn{
    \Big\| \big(J(Du)^{-1}\big)^{\sf T}\N u
            +&(\N u)^{\sf T}J(Du)^{-1}\Big|_{x_n=0}
    \Big\|_{\dF^{\frac12-\frac{1}{2p}}_{1,1}
            (\re_+;\dB^{-1+\frac{n}{p}}_{p,1}(\re^{n-1})) }
  \\
 \le & 
     C\bigg(
      \Big\|\pt_t(-\Del)^{-1}\N 
                  \Big(\big(J(Du)^{-1}\big)^{\sf T}\N u
                  +(\N u)^{\sf T}J(Du)^{-1}\Big)
      \Big\|_{L^1(\re_+;\dB^{-1+\frac{n}{p}}_{p,1}(\re^n_+)) } \\
   &\qquad
     + \Big\| \N \Big(\big(J(Du)^{-1}\big)^{\sf T}\N u
                  +(\N u)^{\sf T}J(Du)^{-1}\Big)
       \Big\|_{L^1(\re_+;\dB^{-1+\frac{n}{p}}_{p,1}(\re^n_+)) }
    \bigg)
 \\
  \le & 
     C\bigg(
      \Big\|\Big(\big(J(Du)^{-1}\big)^{\sf T} \pt_t  u
                  +(\pt_t u)^{\sf T}J(Du)^{-1}\Big)
      \Big\|_{L^1(\re_+;\dB^{-1+\frac{n}{p}}_{p,1}(\re^n_+)) } \\
   &\qquad
     +\Big\|\sum_{k=1}^{n-1}\sg_{k}\Big(\int_0^t Du(s)\,ds \Big)^{k-1} Du \times u      
      \Big\|_{L^1(\re_+;\dB^{-1+\frac{n}{p}}_{p,1}(\re^n_+)) } \\
   &\quad\qquad
     + \Big\|\big(J(Du)^{-1}\big)^{\sf T}\N u
                  +(\N u)^{\sf T}J(Du)^{-1}
       \Big\|_{L^1(\re_+;\dB^{\frac{n}{p}}_{p,1}(\re^n_+)) }
    \bigg)
\\
  \le & 
     C\bigg(
     \sum_{k=0}^{n-1} 
      \Big\|\int_0^t Du(s)\,ds 
      \Big\|_{L^{\infty}(\re_+;\dB^{\frac{n}{p}}_{p,1}(\re^n_+)) }^k
      \big\| \pt_t  u
      \big\|_{L^1(\re_+;\dB^{-1+\frac{n}{p}}_{p,1}(\re^n_+)) } \\
   &\qquad 
    +\sum_{k=1}^{n-1} 
      \Big\| \int_0^t Du(s)\,ds      
      \Big\|_{L^{\infty}(\re_+;\dB^{\frac{n}{p}}_{p,1}(\re^n_+)) }^{k-1} 
      \big\|  Du  
      \big\|_{L^1(\re_+;\dB^{\frac{n}{p}}_{p,1}(\re^n_+)) }
      \big\|  u  
      \big\|_{L^{\infty}(\re_+;\dB^{-1+\frac{n}{p}}_{p,1}(\re^n_+)) }
      \\
   &\quad\qquad
     +\sum_{k=0}^{n-1}
       \Big\|\int_0^t Du(s) ds
       \Big\|_{L^{\infty}(\re_+;\dB^{\frac{n}{p}}_{p,1}(\re^n_+)) }^k
       \big\|\N u
       \big\|_{L^{1}(\re_+;\dB^{\frac{n}{p}}_{p,1}(\re^n_+)) }
   \bigg)
 \\
 \le &
   C\sum_{k=0}^{n-1}
           \|  D^2 u\|_{L^1(\re_+;\dB^{-1+\frac{n}{p}}_{p,1}(\re^n_+))}^k
          \Big(
               \|\pt_t u\|_{L^1(\re_+;\dB^{-1+\frac{n}{p}}_{p,1}(\re^n_+))}
             + \|  D^2 u\|_{L^1(\re_+;\dB^{-1+\frac{n}{p}}_{p,1}(\re^n_+))}
          \Big),
  \eqntag \label{eqn;boundary-nonlinear-9}
}
{\crd 
while 
\eqn{
    \Big\| \big(J(Du)^{-1}\big)^{\sf T}\N u
          +(\N u)^{\sf T}J(Du)^{-1}\Big|_{x_n=0}
    \Big\|_{L^1(\re_+;\dB^{\frac{n-1}{p}}_{p,1}(\re^{n-1}))}
    }
is estimated by the right hand side of \eqref{eqn;boundary-nonlinear-9} 
in much simpler way as in \eqref{eqn;Hp(u,p)-2}, 
since $\dB^{(n-1)/p}_{p,1}(\re^{n-1})$ is the Banach algebra.
}
By {\crd \eqref{eqn;boundary-nonlinear-8}, \eqref{eqn;boundary-nonlinear-9}
and the estimates \eqref{eqn;Hp(u,p)-2} and \eqref{eqn;Hp(u,p)-3} 
}
with the sharp trace estimate imply 
\algn{
  \Big\| 
     H^1_u(u) &
  \Big\|_{ \dF^{\frac12-\frac{1}{2p}}_{1,1}(\re_+;\dB^{-1+\frac{n}{p}}_{p,1}(\re^{n-1}))  } 
    \\
 \le &
   C\sum_{k=2}^{2n-1}
          \Big(
               \|\pt_t u\|_{L^1(\re_+;\dB^{-1+\frac{n}{p}}_{p,1}(\re^n_+))}
             + \|  D^2 u\|_{L^1(\re_+;\dB^{-1+\frac{n}{p}}_{p,1}(\re^n_+))}
          \Big)^k.
  \eqntag \label{eqn;boundary-nonlinear-10}
}
The estimate for the term $H_u^2(u)$ can be shown in the same way as is shown in 
\eqref{eqn;boundary-nonlinear-9},
which shows that \eqref{eqn;m-bilinear-apriori-b-T-u} holds.

For the proof of \eqref{eqn;m-bilinear-apriori-b-S-u},
we notice that  $\dB^{(n-1)/p}_{p,1}(\re^{n-1})$ is the Banach algebra, 
the same argument to \eqref{eqn;boundary-nonlinear-5'} shows for $1\le p<\infty$ that
\algn{
 \|H_{u}(u)&\|_{{\crd L^1}(\re_+;\dB^{\frac{n-1}{p}}_{\crd p,1}(\re^{n-1}))} 
\\
    \le &C 
    {\crd
       \Big\| \Big(\big(J(Du)^{-1}\big)^{\sf T} \N u-(\N u)^{\sf T}J(Du)^{-1}
              \Big) \,\Big|_{x_n=0}
       \Big\|_{ L^1 (\re_+;\dB^{\frac{n-1}{p}}_{p,1}(\re^{n-1}))}
    } \\
    &\hskip2cm {\crd \times 
       \Big\|\Pi^{n-1}_{bu}\Big(\int_0^t Du\,ds \Big)\,\Big|_{x_n=0}
       \Big\|_{ L^{\infty} (\re_+;\dB^{\frac{n-1}{p}}_{p,1}(\re^{n-1}))}
    } \\
    &+C 
    {\crd
       \Big\| \Big(\big(J(Du)^{-1}-I\big)^{\sf T} \N u 
                  -(\N u)^{\sf T}\big(J(Du)^{-1}-I\big)\Big)
            \,\Big|_{x_n=0}
       \Big\|_{ L^1 (\re_+;\dB^{\frac{n-1}{p}}_{p,1}(\re^{n-1}))}
    }
\\
 \le &C\| Du  {\crd \,|_{x_n=0} }
       \|_{ L^1 (\re_+;\dB^{\frac{n-1}{p}}_{p,1}(\re^{n-1}))}
       \sum_{k=1}^{2n-2}          
       \Big\|\int_0^t  D u(s)\, ds \,{\crd \Big|_{x_n=0}}
       \Big\|_{ L^{\infty} (\re_+;\dB^{\frac{n-1}{p}}_{p,1}(\re^{n-1}))} ^k
\\
 \le &C\sum_{k=2}^{2n-1}          
       {\crd
       \big\| D u  \,|_{x_n=0}
       \big\|_{ L^1(\re_+;\dB^{\frac{n-1}{p}}_{p,1}(\re^{n-1}))}^k
       }
 \le C\sum_{k=2}^{2n-1}          
      {\crd
       \|D^2 u
       \|_{L^1(\re_+;\dB^{-1+\frac{n}{p}}_{p,1}(\re^n_+))}^k.
      }
}
This show the estimate \eqref{eqn;m-bilinear-apriori-b-S-u}.
\end{prf}

\sect{The proof of main theorem}\label{Sec;6}
\par\vspace{0.5pc}
\begin{prf}{Theorem \ref{thm;main}}
We define the complete metric space 
\begin{align*}
 X=\left\{ {\ppl (u,p)} :
  \begin{aligned}
      &u\in C(\overline{\re_+};\dot B^{-1+\frac np}_{p,1}(\re^n_+))
            \cap \dot{W}^{1,1}(\re_+;\dot B^{-1+\frac np}_{p,1}(\re^n_+)), \\
      &(\Del u, \N p)\in  L^1(\re_+;\dot B^{-1+\frac np}_{p,1}(\re^n_+)) ,\\
      &p|_{x_n=0}\in 
              \dF^{\frac12-\frac{1}{2p}}_{1,1}(\re_+;\dB^{-1+\frac{n}{p}}_{p,1}(\re^{n-1}))
              \cap
              {\crd L^1}(\re_+;\dB^{\frac{n-1}{p}}_{p,1}(\re^{n-1})), 
  \end{aligned}
  \hskip2mm  \|(u,p)\|_X\le M 
  \right\}, 
\end{align*}
where 
\eqn{
 \spl{
 \|(u,p)\|_X
   \equiv& \|\pt_t u\|_{L^1(\re_+;\dB^{-1+\frac np}_{p,1}(\re^n_+))}
          +\|D^2 u\|_{L^1(\re_+;\dB^{-1+\frac np}_{p,1}(\re^n_+))}   
          +\|\N p\|_{L^1(\re_+;\dB^{-1+\frac np}_{p,1}(\re^n_+))}  \\
         &+\big\| p|_{x_n=0}\big\|_{\dF^{1/2-1/2p}_{1,1}%
                                    (\re_+;\dB^{-1+\frac np}_{p,1}(\re^{n-1}))}
          +\big\| p|_{x_n=0}\big\|_{{\crd L^1}%
                                    (\re_+;\dB^{\frac{n-1}p}_{p,1}(\re^{n-1}))}.             
}
}
The constant $M>0$ is chosen to be small enough depending on the norm of 
the initial data.
Given $(\tu,\tp)\in X$, we consider the liner inhomogeneous 
initial boundary value problem:
\begin{equation} \label{eqn;NS-L_p-2}
   \left\{
   \begin{split}
     &\pt_t u  - \Del u + \N p
        =  F_u(\tu)+F_p(\tu,\tp),
        &\qquad t>0,\ x\in\Bbb R^n_+,\\
     &\qquad \qquad \div\, u =G_{\div}(\tu).        
        &\qquad t>0,\ x\in\Bbb R^n_+,\\
     &\quad \big(\nabla u+(\nabla u)^{\sf T}-p I\big)\cdot \nu_n 
          = H_u(\tu)  +H_p(\tu,\tp).        
        &\qquad t>0,\ x\in\pt\Bbb{R}_+^n,\\   
     &\qquad \qquad u(0,x) =u_0(x),
        &\qquad \phantom{\ t>0,}
               x\in\Bbb R^n_+,
     \end{split}
     \right.   \hfill
\end{equation}
where  $\nu_n=(0,0,\cdots, 0,-1)^{\sf T}$ denotes 
the outer normal and we set
\begin{align}
   F_u(\tu)=&\div\Big(
         J(D \tu)^{-1}\big(J(D \tu)^{-1}\big)^{\sf T}\N \tu-\N \tu\Big)
           = \Pi_{{\tu}}^{2n-2}
           \left(\int_0^t D\tu\; ds\right)
           D^2 \tu, 
    \label{eqn;laplace-purterb-tilde}\\
   F_p(\tu,\tp)&=-\big(J(D \tu)^{-1}-I\big)\, \N \tp 
           =\Pi_{p}^{n-1}
             \left(\int_0^t D \tu\; ds\right)
             \N\tp, 
    \label{eqn;pressure-purterb-tilde}
\\
  G_{\div}(\tu)=&-\text{tr}\Big(\big(J(D \tu)^{-1}-I\big)\; \N \tu\Big)
           = 
          \text{tr}\left(\Pi_{\div}^{n-1}
          \left(\int_0^t D \tu\; ds\right)
           D \tu\right) \notag\\
       =&\, \div\left(
             \Pi_{\div}^{n-1}\left(\int_0^t D \tu\, ds\right)
                 \tu\right),
    \label{eqn;divergence-purterb-tilde}
 \\
   H_{b}(\tu) =& 
        -\Big( \big( J\big(D(\tu)\big)^{-1}\big)^{\sf T}\;\N \tu
                +(\N \tu)^{\sf T} J\big(D(\tu)\big)^{-1}
         \Big)\;\big( J\big(D(\tu)\big)^{-1}-I\big)^{\sf T}\nu_n  \notag\\
      &\quad 
         -\Big( \big( J\big(D(\tu)\big)^{-1}-I\big)^{\sf T}\;\N \tu
                +(\N \tu)^{\sf T}\big( J\big(D(\tu)\big)^{-1}-I\big)
         \Big)\nu_n  
      \notag \\
     =& \Pi_{bu}^{2n-2}
            \left(\int_0^t D \tu\; ds\right)
           D \tu \;\nu_n,  
   \label{eqn;boundary-u-purterb-tilde}     
 \\
 H_{p}(\tu,\tp)=&        
          \tp  \big(J(D\tu)^{-1}-I\big)^{\sf T}\nu_n  
   =   \Pi_{bp}^{n-1}
            \left(\int_0^t D \tu\; ds\right)
             p \;\nu_n.
     \label{eqn;boundary-p-purterb-tilde}      
\end{align}

We define the map 
$$
\Phi:  X\to X
$$
by 
$$
 (\tu,\tp) \to  (u,p)\equiv \Phi\big[\tu, \tp\big] 
$$
and prove that $\Phi$ is contraction on $X$. 

\vskip1mm
First we show that a priori estimate of $\Phi[{u}, {p}]$ 
in $L^1(\re_+;\dB^{-1+n/p}_{p,1}(\re^n_+))$. 
Let $({u}, {p})$ solve \eqref{eqn;NS-L_p-2}. 
Applying Theorem \ref{thm;L1MR} to the equation \eqref{eqn;NS-L_p-2}, 
we have by \eqref{eqn;L1MR-estimate}, 
Propositions \ref{prop;multiple-bilinear-1}--\ref{prop;multiple-bilinear-3}
to the nonlinear terms 
\eqref{eqn;laplace-purterb-tilde}--\eqref{eqn;boundary-p-purterb-tilde} that 
{\allowdisplaybreaks
\algn{
   \|\pt_t u &\|_{L^1(\re_+;\dB^{-1+\frac np}_{p,1}(\re^n_+))}
   +\|D^2 u  \|_{L^1(\re_+;\dB^{-1+\frac np}_{p,1}(\re^n_+))} \\
 & +\|\N   p  \|_{L^1(\re_+;\dB^{-1+\frac np}_{p,1}(\re^n_+))}
   +\big\| p|_{x_n=0}\big\|_{\dF^{1/2-1/2p}_{1,1}(\re_+;\dB^{-1+\frac np}_{p,1}(\re^{n-1}))}
   +\big\| p|_{x_n=0}\big\|_{{\crd L^1}(\re_+;\dB^{\frac{n-1}p}_{p,1}(\re^{n-1}))}
\\
\le & 
  C_M\Big(\|u_0 \|_{\dB^0_{p,1}(\re^n_+)}
         +\|F_u(\tu) \|_{L^1(\re_+;\dB^{-1+\frac np}_{p,1}(\re^n_+))} 
         +\|F_p(\tu,\tp)\|_{L^1(\re_+;\dB^{-1+\frac np}_{p,1}(\re^n_+))}  \\
   &\quad
         +\|\N G_{\div}(\tu)\|_{L^1(\re_+;\dB^{-1+\frac np}_{p,1}(\re^n_+))}
         +\|\pt_t(-\Del)^{-1}\N {G_{\div}}(\tu)
          \|_{L^1(\re_+;\dB^{-1+\frac np}_{p,1}(\re^n_+))}\\
    &\qquad
      +\| H_u(\tu) \|_{\dF^{\frac12-\frac{1}{2p}}_{1,1}(\re_+;\dot B^{-1+\frac np}_{p,1}(\re^{n-1}))}
      +\| H_u(\tu) \|_{{\crd L^1}(\re_+;\dot B^{\frac np-\frac 1p}_{p,1}(\re^{n-1}))}    \\
   &\qquad\quad 
      +\|H_p(\tu,\tp)\|_{\dF^{\frac12-\frac{1}{2p}}_{1,1}(\re_+;\dot B^{-1+\frac np}_{p,1}(\re^{n-1}))}
      +\|H_p(\tu,\tp)\|_{{\crd L^1}(\re_+;\dot B^{\frac np-\frac 1p}_{p,1}(\re^{n-1}))}     
     \Big)
\\
\le & 
  C\bigg(\|u_0 \|_{\dB^{-1+\frac np}_{p,1}(\re^n_+)}
    +\sum_{k=1}^{2n-2}
     \|D^2 \tu\|_{L^1(\re_+;\dB^{-1+\frac np}_{p,1}(\re^n_+))}^{k+1}\\
 &\hskip1cm        
    +\sum_{k=1}^{2n-1}
     \|D^2 \tu\|_{L^1(\re_+;\dB^{-1+\frac np}_{p,1}(\re^n_+))}^k
     \Big( \|\pt_t \tu\|_{L^1(\re_+;\dB^{-1+\frac np}_{p,1}(\re^n_+))} 
           +\|D^2  \tu\|_{L^1(\re_+;\dB^{-1+\frac np}_{p,1}(\re^n_+))} 
     \Big)
   \\
 &\hskip1.5cm
  + \sum_{k=1}^{n-1}
         \Big(
            \|\pt_t \tu\|_{L^1(\re_+;\dB^{-1+\frac np}_{p,1}(\re^n_+))}
          + \|D^2   \tu\|_{L^1(\re_+;\dB^{-1+\frac np}_{p,1}(\re^n_+))}
          \Big)^k  \\
       &\hskip2cm\times
         \Big(
          \|\N  \tp\|_{L^1(\re_+;\dB^{-1+\frac np}_{p,1}(\re^n_+))}    
         +\|\tp|_{x_n=0}\|_{\dF^{\frac12-\frac{1}{2p}}_{1,1}(\re_+;\dB^{-1+\frac np}_{p,1}(\re^{n-1}))}  
         +\|\tp|_{x_n=0}\|_{{\crd L^1}(\re_+;\dB^{\frac {n-1}p}_{p,1}(\re^{n-1}))}  
          \Big)
  \\
  &\hskip2cm
  + \sum_{k=2}^{2n-1}
        \Big(
           \|\pt_t \tu\|_{L^1(\re_+;\dB^{-1+\frac np}_{p,1}(\re^n_+))}
          +\|D^2   \tp\|_{L^1(\re_+;\dB^{-1+\frac np}_{p,1}(\re^n_+))} 
         \Big)^k 
     \bigg).
  \eqntag \label{eqn;L1MR-Stokes-10}
}}
By \eqref{eqn;L1MR-Stokes-10}, it holds that 
\eq{ \label{eqn;apriori-bound}
 \|\Phi[\tu,\tp]\|_X
  \le C_1\bigg(\|u_0\|_{\dB^{-1+\frac np}_{p,1}(\re^n_+)}
               + \sum_{k=1}^{2n-1}M^{k+1}\bigg). 
}
Therefore if we choose the initial data small enough
$$
\|u_0\|_{\dB^{-1+\frac np}_{p,1}(\re^n_+)}\le \frac{1}{2C_1} M<\frac 14
$$
 then we obtain from \eqref{eqn;apriori-bound} that
$$
  \|\Phi[\tu,\tp]\|_X\le M. 
$$
Moreover, for all $(u_1,p_1)$, $(u_2,p_2)\in X$, we know that the difference 
 $$
  w=u_1-u_2, 
  \qquad
  q=p_1-p_2
 $$
satisfy the same estimate \eqref{eqn;L1MR-Stokes-10}: 
\eqn{
 \spl{
    \|\pt_t w &\|_{L^1(\re_+;\dB^{-1+\frac{n}{p}}_{p,1}(\re^n_+))}
   +\|D^2   w  \|_{L^1(\re_+;\dB^{-1+\frac{n}{p}}_{p,1}(\re^n_+))}
   +\|\N  q    \|_{L^1(\re_+;\dB^{-1+\frac{n}{p}}_{p,1}(\re^n_+))} \\
  &\qquad
   +\big\| q|_{x_n=0}
    \big\|_{\dF^{\frac12-\frac{1}{2p}}_{1,1}(\re_+;\dB^{-1+\frac np}_{p,1}(\re^{n-1}))}
   +\big\| q|_{x_n=0}\big\|_{{\crd L^1}(\re_+;\dB^{\frac{n-1}p}_{p,1}(\re^{n-1}))}
\\
\le & 
  C_2\max_{i=1,2}\sum_{k=1}^{2n-2}\Big(
          \|\pt_t u_i\|_{L^1(\re_+;\dB^{-1+\frac{n}{p}}_{p,1}(\re^n_+))}
      +   \|D^2   u_i\|_{L^1(\re_+;\dB^{-1+\frac{n}{p}}_{p,1}(\re^n_+))}
     \Big)^k
\\
   &\quad\times
     \Big(\|\pt_t w\|_{L^1(\re_+;\dB^{-1+\frac{n}{p}}_{p,1}(\re^n_+))}
         +\|D^2 w\|_{L^1(\re_+;\dB^{-1+\frac{n}{p}}_{p,1}(\re^n_+))}
         +\|\N q\|_{L^1(\re_+;\dB^{-1+\frac{n}{p}}_{p,1}(\re^n_+))} \\
       &\qquad 
        + \big\| q|_{x_n=0}
          \big\|_{\dF^{\frac12-\frac{1}{2p}}_{1,1}(\re_+;\dB^{-1+\frac np}_{p,1}(\re^{n-1}))}
        + \big\| q|_{x_n=0}
          \big\|_{{\crd L^1}(\re_+;\dB^{\frac{n-1}p}_{p,1}(\re^{n-1}))}
     \Big).
   }
}
Therefore if we choose
$$
  C_2\max_{i=1,2}\sum_{k=1}^{2n-2}\Big(
          \|\pt_t u_i\|_{L^1(\re_+;\dB^{-1+\frac{n}{p}}_{p,1}(\re^n_+))}
      +   \|D^2 u_i\|_{L^1(\re_+;\dB^{-1+\frac{n}{p}}_{p,1}(\re^n_+))}
     \Big)^k
   \le C_2 \sum_{k=1}^{2n-2} M^k\le  \frac12, 
$$
then it holds that 
$$
  \| \Phi[w,q] \|_X\le  \frac12 \|(w,q)\|_X, 
$$
which shows the map 
$$
\Phi :  X\to X
$$
is contraction.  By the fixed point theorem of Banach-Cacciopolli, 
there exists a unique fixed point $(u,p)$ of the map $\Phi$ in $X$.

We finally confirm that the boundary equation in \eqref{eqn;NS-L_p-2} is fulfilled.
Let the difference between the 
solution and the date as 
 $$
  \widetilde{w}=\widetilde{u}-u.
 $$
The sharp trace estimate Proposition \ref{prop;grad-pressure-trace} ensure 
that 
{\allowdisplaybreaks
\algn{ 
  \big\|\N \widetilde{w}&|_{x_n=0} 
  \big\|_{\dF^{\frac12-\frac{1}{2p}}_{1,1}(\re_+;\dB^{-1+\frac{n}{p}}_{p,1}(\re^{n-1}))} 
 +\big\|\N \widetilde{w} |_{x_n=0}
  \big\|_{{\crd L^1}(\re_+;\dB^{\frac{n}{p}-\frac1p}_{p,1}(\re^{n-1}))} \\
\le &C\Big( \|\pt_t \widetilde{w}\|_{L^1(\re_+;\dB^{-1+\frac{n}{p}}_{p,1}(\re^n_+))}
           +\|\Del  \widetilde{w}\|_{L^1(\re_+;\dB^{-1+\frac{n}{p}}_{p,1}(\re^n_+))}
      \Big)  
   \eqntag \label{eqn;boundary-convergence}
}}
and the right hand side of \eqref{eqn;boundary-convergence} converges to 
$0$ as the iterative process.
Then the unique fixed point $(u,p)$ satisfies \eqref{eqn;NS-L_p-2} 
with the all right members changed into 
$(u,p)$ and it is a time global strong solution of \eqref{eqn;NS}. 
This completes the proof of Theorem \ref{thm;main}. 
\end{prf}

\sect{Appendix}\label{Sec;7}
\subsection{The optimal boundary trace}\label{Sec;7-2}
The proof of Proposition \ref{prop;grad-pressure-trace} is based 
on the following trace estimate (cf. \cite{MV14}).
\vskip2mm

\begin{thm}[Sharp boundary derivative trace \cite{OgSs20-2}]
\label{thm;sharp-boundary-trace-N} 
For  $1< p< \infty$ and {\bl $-1+1/p<s$}, 
there exists a constant $C>0$ 
such that for any function $f=f(t,x',\eta)$ with  
$f\in C_b(\overline{\re_+};\dB^s_{p,1}(\re^{n}_+))
\cap \dot{W}^{1,1}(\re_+; \dB^s_{p,1}(\re^{n}_+))$, 
${\crd D^2} f \in  L^1(\re_+; \dB^{s}_{p,1}(\re^{n}_+))$, 
it holds for all $\pt_{\ell}=(\pt_{x_{\ell}}, \pt_{\eta})$ with 
$1\le \ell \le n-1$ that 
\begin{align}
   \sup_{\eta\in\re_+}
   \Bigl(\big\| \pt_{\ell} f(\cdot,\cdot, \eta)
        &\big\|_{\dF^{\frac12-\frac{1}{2p}}_{1,1}(\re_+; \dB^s_{p,1}(\re^{n-1}))}
        +\big\| \pt_{\ell} f(\cdot,\cdot, \eta)
         \big\|_{{\crd L^1} (\re_+; \dB^{s+1-\frac1p}_{p,1}(\re^{n-1}))}
  \Bigr)  \nonumber\\
   & \qquad\qquad  \le C\Bigl(
         \|\pt_t f\|_{L^1(\re_+;\dB^s_{p,1}(\re^{n}_+))}
       + \|{\crd D^2}  f\|_{L^1(\re_+;\dB^s_{p,1}(\re^{n}_+))}
       \Bigr) .\label{eqn;sharp-trace-N}
\end{align}
\end{thm}

\begin{prf}{Theorem \ref{thm;sharp-boundary-trace-N}}
For $1< p< \infty$ and {\bl $-1+1/p<s$}, assume 
$f\in C_b(\overline{\re_+};\dB^s_{p,1}(\re^{n}_+))
 \cap \dot{W}^{1,1}(\re_+; \dB^s_{p,1}(\re^{n}_+))$, 
$D^2 f \in  L^1(\re_+; \dB^{s}_{p,1}(\re^{n}_+))$.
{\bl Then by the definition \eqref{eqn;Besov-halfspace}, 
for any $\ep>0$, there exists $\widetilde{f}
\in  C_b(\overline{\re_+};\dB^s_{p,1}(\re^{n}))
\cap \dot{W}^{1,1}(\re_+; \dB^s_{p,1}(\re^{n}))
\cap L^1(\re_+; \dB^{s+2}_{p,1}(\re^{n}))$
such that 
\algn{
  \|\widetilde{f}\|_{ L^{\infty}(\re_+;\dB^s_{p,1}(\re^n))
   \cap 
   \dot{W}^{1,1}(\re;\dB^s_{p,1}(\re^n)) }
  \le  &
  \|f\|_{ L^{\infty}(\re_+;\dB^s_{p,1}(\re^n_+))
    \cap  \dot{W}^{1,1}(\re_+;\dB^s_{p,1}(\re^n_+)) }+\ep,
  \\
   \|D^2\widetilde{f}\|_{ L^{1}(\re_+;\dB^s_{p,1}(\re^n)) }
  \le  &
  \|D^2 f\|_{ L^{1}(\re_+;\dB^s_{p,1}(\re^n_+)) }+\ep.
}
} 
{\ppl 
We then extend  $\widetilde{f}$ into $t<0$ as an 
even extension.
} 
{\bl For simplicity, we denote $\widetilde{f}$ as $f$ in the following.}
It directly follows that 
\eq{\label{eqn;infinity-vanishing}
\pt_{\ell}f \in L^{\ppl 2}({\ppl \re}; \dB^{s+1}_{p,1}({\bl \re^{n}})).
}
Then
\begin{align*}
 \big\|\pt_{\ell}f(\cdot,\cdot, \eta) 
  \big\|_{\dF^{\frac12-\frac{1}{2p}}_{1,1}
           \big(\re_+;\dB^{s}_{p,1}({\gr \re^{n-1}})\big)} 
  \le  & {\ppl
        \big\|\pt_{\ell}f(\cdot,\cdot, \eta) 
        \big\|_{\dF^{\frac12-\frac{1}{2p}}_{1,1}
           \big(\re;\dB^{s}_{p,1}(\re^{n-1})\big)}
        } 
\\
  \le & 
  \Big\| \sum_{j\in\Z} 2^{sj}
         \sum_{k\ge 2j}  
              2^{(\frac12-\frac{1}{2p})k}
              \big\|
                 \psi_k\underset{(t)}{*}
                 \phi_j\underset{(x')}{*}
                 \pt_{\ell}f(t,\cdot,\eta)
              \big\|_{L^p({\gr \re^{n-1}})}
       \Big\|_{L^1_t({\ppl \re})} 
 \\
  &+
    \Big\| \sum_{j\in\Z} 2^{sj}
           \sum_{k\le 2j}  
            2^{(\frac12-\frac{1}{2p})k}
              \big\|
                 \psi_k\underset{(t)}{*}
                 \phi_j\underset{(x')}{*}
                 \pt_{\ell}f(t,\cdot,\eta)
              \big\|_{L^p({\mt \re^{n-1}})}
       \Big\|_{L^1_t({\ppl \re})} 
 \\
 \equiv & I+II.
  \addtocounter{equation}{1}
 \tag{\theequation} \label{eqn;time-like-N0}
 \end{align*}
From \eqref{eqn;infinity-vanishing} {\bl and $1/p<s+1$}, 
{\ppl one can approximate $\pt_{\ell} f$ by a function satisfying }  
\eq{ \label{eqn;compatibility-condition-N}
  \lim_{t\to {\ppl \pm}\infty}\pt_{\ell}f(t,x',\eta)=0,
  \quad \text{ a.a. } (x',\eta)\in \re^{n-1}\times{\ppl \re},
}
and using the assumption \eqref{eqn;compatibility-condition-N}
{\ppl and noting $\widehat{\psi_k}(0)=\int_{\re}\psi_k(t)dt=0$,} 
\begin{align*}
 \psi_k(t)\underset{(t)}{*}\pt_{\ell} f(t,x',\eta)
  =&\int_{\ppl\re}\psi_k(s) 
                \, \pt_{\ell} f(t-s,x',\eta)ds \\
 =&-\int_{\ppl\re}\pt_s\big(\int_s^{\infty}\psi_k(r)dr\big) 
                \pt_{\ell} f(t-s,x',\eta)ds \\
 =&-\Bigl[\Big(\int_s^{\infty} \psi_k(r)dr\Big) 
          \pt_{\ppl \ell} f(t-s,x',\eta)
    \Bigr]_{s=-\infty}^{\infty}
   +\int_{\ppl\re}\Big(\int_s^{\infty}\psi_k(r)dr\Big)
       \pt_s \pt_{\ell} f(t-s,x',\eta)ds  \\
 =&- \int_{\ppl \re}\pt_s^{-1}\psi_k(s) \pt_s \pt_{\ell} f(t-s,x',\eta)ds \\
 =&\, \pt_t^{-1}\psi_k(t)\underset{(t)}{*}  \pt_t \pt_{\ell} f(t,x',\eta),
 \label{eqn;time-like-int-N}\eqntag
\end{align*}
where we set
\eq{\label{eqn;int-psi-k}
 {\ppl
  \pt_t^{-1}\psi_k(s)
  \equiv -\int_s^{\infty}\psi_k(r)dr.
 }
}
Then $\pt_t^{-1}\psi_k(t)=2^{-k}\big(\pt_t^{-1}\psi_0\big)_k(t)$.
Hence from \eqref{eqn;compatibility-condition-N} and 
\eqref{eqn;time-like-int-N} and using the Hausdorff--Young 
inequality, it follows for $\pt_\ell=\pt_\eta$ (i.e., $\ell=n$) that
{\allowdisplaybreaks
\begin{align*}
  I 
 \le &C\Big\| {\bk\sum_{j\in\Z} 2^{sj}
                  \sum_{k\ge 2j}  2^{-(\frac12+\frac1{2p})k}
              }
               \int_{\re}              
                 |(\pt_t^{-1}\psi_0)_k(t-s)|  
                \big\|
                    \phi_j\underset{(x')}{*}
                    \pt_t\pt_{\eta} f(s,x',\eta)
                \big\|_{L^p({\gr \re^{n-1}})}                 
                ds
          \Big\|_{L^1_t({\ppl \re})}
 \\
   \le &  C \sum_{j\in \Z} 2^{sj}
              \sum_{k\ge 2j} 2^{-(\frac12+\frac1{2p})k} 
          \Big\|\int_{\re}
                 \frac{2^k}{\<2^k(t-s)\>^2}
                 \big\|
                 \phi_j\underset{(x')}{*}
                 \pt_t\pt_{\eta}  f({\ppl s},x',\eta) 
                \big\|_{L^p({\gr \re^{n-1}})}               
                ds
        \Big\|_{L^1_t({\ppl \re})}
  \\
   \le & C{\bk \sum_{j\in \Z}  2^{sj}
              \sum_{k\ge 2j}  2^{-(\frac12+\frac1{2p})k} }
          \Big\|\frac{2^k}{\<2^k t\>^2}\Big\|_{L^1_t({\ppl \re})}
          \Big\|
                 \big\|
                 \sum_{m\in \Z} \overline{\Phi_m}\underset{(x',\eta)}{*} 
                 \phi_j\underset{(x')}{*}
                 \pt_t\pt_{\eta} f(t,x',\eta)  
                \big\|_{L^p({\gr \re^{n-1}})}
          \Big\|_{L^1_t(\ppl \re)}
  \\
   \le &C\sum_{j\in \Z} 2^{sj}
          \Big\| {\bk \sum_{k\ge 2j}  2^{-(\frac12+\frac1{2p})k} } \\
            &\hskip24mm\times 
                 \big\|
                 \sum_{|m-j|\le 1} \overline{\Phi_m}\underset{(x',\eta)}{*} 
                  \phi_j\underset{(x')}{*}
                 \pt_{\eta}{\bl \zeta_{j+1} }(\eta)\underset{(\eta)}{*}
                 \zeta_{\bl j+2}(\eta)\underset{(\eta)}{*} 
                 \pt_t f(t,x',\eta)  
                \big\|_{L^p({\gr \re^{n-1}})}                
          \big\|_{L^1_t(\re_+)}
\\
    \le &C\sum_{j\in \Z} 2^{sj}
          \Big\| {\bk \sum_{k\ge 2j}  2^{-(\frac12+\frac1{2p})k} } \\
            &\hskip24mm\times 
                 \int_{\re} 
                    \big| \pt_{\eta}{\bl \zeta_{j+1}}\underset{(\eta)}{*}
                     \zeta_{\bl j+2}(\eta-\xi)\big|
                   \Big\| \sum_{|m-j|\le 1} 
                    \overline{\Phi_m}\underset{(x',\xi)}{*} 
                    \phi_j\underset{(x')}{*}                          
                     \pt_t f(t,x',\xi) \Big\|_{L^p({\gr \re^{n-1}})} 
                  d\xi                               
          \Big\|_{L^1_t(\re_+)}
    \\
   \le &C\sum_{j\in \Z} 2^{sj}
          \Big\|{\bk \sum_{k\ge 2j} 2^{-(\frac12+\frac1{2p})k}  }
                 \|\pt_{\eta}{\bl \zeta_{j+1}}(\eta)
                 \|_{L^{p'}(\re_{\eta})}\\
          &\qquad\qquad\qquad\qquad \times       
                 \Big\|\big\|
                     \overline{\Phi_j}\underset{(x',\eta)}{*} 
                 \phi_j\underset{(x')}{*}               
                \pt_t f(t,x',\eta)  
                \big\|_{L^p({\gr \re^{n-1}})}
                \Big\|_{L^p(\re_{\eta})}                
          \Big\|_{L^1_t(\re_+)}
   \\
      \le &C\Big\|\sum_{j\in \Z}  2^{sj}
                 \big\| 
                     \overline{\Phi_j}\underset{(x',\eta)}{*} 
                 \pt_t f(t,x',\eta)  
                \big\|_{L^p(\re^{n}_{(x',\eta)})}              
          \Big\|_{L^1_t(\re_+)} 
    \le C \|\pt_t f\|_{ L^1(\re_+;\dB^s_{p,1}(\re^n_{+}) )}
         +\tilde\ep.
     \addtocounter{equation}{1}
 \tag{\theequation} \label{eqn;time-like-N1}
\end{align*}
}

For the second term of  \eqref{eqn;time-like-N0}, 
we use the Minkowski inequality to see
{\allowdisplaybreaks
\begin{align*}
  II 
   \le & C {\bk  \sum_{j\in\Z}  2^{sj} \sum_{k\le 2j} }
                2^{(\frac12 -\frac{1}{2p})k}
        \Big\| {\bk |\psi_k|\underset{(t)}{*}}
             \big\|
                   \phi_j\underset{(x')}{*}
                    \pt_{\eta}f(t,x',\eta)
             \big\|_{L^p({\gr \re^{n-1}})}
        \Big\|_{L^1_t({\ppl \re})}
\\
   \le & C {\bk  \sum_{j\in\Z} 2^{sj} \sum_{k\le 2j} 2^{(\frac12 -\frac{1}{2p})k}} \\
       &\hskip20mm\times            
        \Big\|
             \big\|
              \sum_{m\in \Z} \overline{\Phi_m}\underset{(x',\eta)}{*} 
              {\bl \zeta_{j+1}}(\eta)\underset{(\eta)}{*} 
              {\bl \zeta_{j+2}}(\eta)\underset{(\eta)}{*} 
              \phi_j\underset{(x')}{*}
              \pt_{\eta}f(t,x',\eta)
             \big\|_{L^p({\gr \re^{n-1}})}
        \Big\|_{L^1_t(\re_+)}
\\
     \le & C\Big\| {\bk  \sum_{j\in\Z}  2^{(s+1-\frac{1}{p})j} }              
             \big\|
              \sum_{|m-j|\le1} \overline{\Phi_m}\underset{(x',\eta)}{*} 
              \pt_{\eta}\zeta_{\bl j+1}(\eta)\underset{(\eta)}{*} 
              \zeta_{\bl j+2}(\eta)\underset{(\eta)}{*} 
              \phi_j\underset{(x')}{*}
              f(t,x',\eta)
             \big\|_{L^p({\gr \re^{n-1}})}
        \Big\|_{L^1_t(\re_+)}
\\
       \le & C\Big\|   \sum_{j\in\Z}    2^{(s+1-\frac{1}{p})j}  
              \|\pt_{\eta}\widetilde{\zeta_{j}}\|_{L^{p'}(\re_{\eta})}             
             \big\|
               \overline{\Phi_j}\underset{(x',\eta)}{*} 
               f(t,x',\eta)
             \big\|_{L^p( \re^{n})}
        \Big\|_{L^1_t(\re_+)}
\\
       \le & C\Big\|   \sum_{j\in\Z}    2^{(s+2)j} {\bl 2^{-2j}} 
             \big\|
              \overline{\Phi_j}\underset{(x',\eta)}{*} 
              {\crd D^2} f(t,x',\eta)
             \big\|_{L^p(\re^{n})}
        \Big\|_{L^1_t(\re_+)}
\\
    \le &  C\big\|{\crd D^2} f \big\|_{ L^1(\re_+;\dB^s_{p,1}(\re^n)) }
    \le    C\big\|{\crd D^2} f \big\|_{ L^1(\re_+;\dB^s_{p,1}(\re^n_+)) }
            {\bl +\tilde{\ep} },
     \addtocounter{equation}{1}
\tag{\theequation} \label{eqn;time-like-N2}
\end{align*}
}
where we used $1<p$ and $\tilde\ep= C\ep$.

{\crd 
The estimate for the spatial trace term is shown 
along the following way:
\begin{align*}
   \big\|\pt_{\eta}f(\cdot,\cdot, \eta) 
  &\big\|_{{\crd L^1}(\re_+;\dB^{s+1-\frac{1}{p}}_{p,1}({\gr \re^{n-1}}))} 
 \\
  = & 
  \Big\| \sum_{j\in\Z} 2^{sj}  2^{(1-\frac{1}{p})j}               
              \big\|
                 \phi_j\underset{(x')}{*}
                 \pt_{\eta}f(t,\cdot,\eta)
              \big\|_{L^p({\gr \re^{n-1}})}
       \Big\|_{L^1_t(\re_+)} 
\\
     \le & C\Big\|  \sum_{j\in\Z}  
              2^{sj}  
              2^{(1 -\frac{1}{p})j}             
             \big\|
              \sum_{m\in \Z} 
              \overline{\Phi_m}\underset{(x',\eta)}{*} 
              \phi_j\underset{(x')}{*}
              \pt_{\eta}f(t,x',\eta)
             \big\|_{L^p({\gr \re^{n-1}})}
        \Big\|_{L^1_t(\re_+)}
\\
     \le & C\Big\|  \sum_{j\in\Z}  
              2^{sj}  
              2^{(1 -\frac{1}{p})j}             
             \big\|
              \sum_{m\in \Z} 
              \overline{\Phi_m}\underset{(x',\eta)}{*} 
              \pt_{\eta}{\bl \zeta_{j+1}}(\eta)\underset{(\eta)}{*} 
              {\bl \zeta_{j+2}}(\eta)\underset{(\eta)}{*} 
              \phi_j\underset{(x')}{*}
              f(t,x',\eta)
             \big\|_{L^p({\gr \re^{n-1}})}
        \Big\|_{L^1_t(\re_+)}
\\
       \le & C\Big\|   \sum_{j\in\Z} 2^{(s+2)j}\, {\bl 2^{-2j}} 
                \big\|\overline{\Phi_j}\underset{(x',\eta)}{*} 
                  D^2 f(t,x',\eta)
                \big\|_{L^p(\re^{n})}
              \Big\|_{L^1_t(\re_+)}
\\
    \le &  C\big\|D^2 f \big\|_{ L^1(\re_+;\dB^s_{p,1}(\re^n)) }
    \le    C\big\|D^2 f \big\|_{ L^1(\re_+;\dB^s_{p,1}(\re^n_+)) }
            {\bl +\tilde{\ep}}.
     \addtocounter{equation}{1}
     \tag{\theequation} \label{eqn;time-like-N2}
\end{align*}
}
{\bl 
Since $\tilde{\ep}=C\ep>0$ is arbitrary small, we conclude our 
result.
}

The other cases $1\le \ell \le n-1$ can be shown similar way by changing 
$\pt_{\eta}\zeta_j$ into $\zeta_j$ and $\phi_j$ 
into $\pt_{\ell}\widetilde{\phi_j}*\phi_j$.
This completes the proof of Theorem \ref{thm;sharp-boundary-trace-N}.

\end{prf}
\vskip3mm

\begin{prf}{Proposition \ref{prop;grad-pressure-trace}}
The proof of the trace estimate \eqref{eqn;pressure-trace} 
is almost the same line of \eqref{eqn;sharp-trace-N} 
in Theorem \ref{thm;sharp-boundary-trace-N} except the 
regularity.  Hence we show  {\bl an outlined } proof.
For {\gr $1< p< \infty$} and {\bl $-1+1/p<s<1/p$}, assume 
$\N(-\Del)^{-1}f\in {\ppl C_b}(\overline{\re^n_+};\dB^s_{p,1}(\re^n_+))
\cap \dot{W}^{1,1}(\re_+;\dB^{s}_{p,1}(\re^n_+))$
and $\N f\in  L^1(\re_+;\dB^{s}_{p,1}(\re^n_+))$. 
{\bl We employ the similar argument and use the {\crd extended} element 
$\widetilde{f}$  associated with $f$ as above.}
{\bl We regard $f$ as the extended element $\widetilde{f}$. From}
$
 {\bl 
 f\in L^{\ppl 2}({\ppl \re};\dB^{s+1}_{p,1}(\re^n)) 
 }
$
like in  \eqref{eqn;infinity-vanishing}, we may {\ppl assume} that 
\eq{ \label{eqn;infinity-vanishing-2}
  \lim_{t\to {\ppl \pm}\infty}f(t,x',\eta)=0,
  \quad \text{ a.a. } (x',\eta)\in \re^{n-1}\times\re,
}
 we have like \eqref{eqn;time-like-int-N} that
\begin{align*}
 \psi_k(t)\underset{(t)}{*} f(t,x',\eta)
  =&\pt_t^{-1}\psi_k(t)\underset{(t)}{*}  \pt_t f(t,x',\eta),
 \label{eqn;time-like-int-p}\eqntag
\end{align*}
with \eqref{eqn;int-psi-k}.
Then it follows 
\begin{align*}
   \big\|f(\cdot,\cdot, \eta) 
   \big\|_{\dF^{\frac12-\frac{1}{2p}}_{1,1}\big(\re_+;\dB^{s}_{p,1}({\gr \re^{n-1}})\big)} 
  \le & 
  \Big\| \sum_{j\in\Z} 2^{sj}
         \sum_{k\ge 2j}  
              2^{(\frac12-\frac{1}{2p})k}
              \big\|
                 \psi_k\underset{(t)}{*}
                 \phi_j\underset{(x')}{*}
                 f(t,\cdot,\eta)
              \big\|_{L^p({\gr \re^{n-1}})}
       \Big\|_{L^1_t({\ppl \re})} \\
  &+\Big\| \sum_{j\in\Z} 2^{sj}
           \sum_{k\le 2j}  
            2^{(\frac12-\frac{1}{2p})k}
              \big\|
                 \psi_k\underset{(t)}{*}
                 \phi_j\underset{(x')}{*}
                 f(t,\cdot,\eta)
              \big\|_{L^p({\gr \re^{n-1}})}
       \Big\|_{L^1_t({\ppl \re})} 
 \\
 \equiv & I+II.
  \addtocounter{equation}{1}
 \tag{\theequation} \label{eqn;time-like-p0}
 \end{align*}
Hence from \eqref{eqn;infinity-vanishing-2} and 
\eqref{eqn;time-like-int-p} and using the Hausdorff--Young 
inequality

\begin{align*}
  I = & \Big\|{\bk \sum_{j\in\Z} 2^{sj}
                   \sum_{k\ge 2j}  }
                    2^{(\frac12-\frac{1}{2p})k}
              \big\|
                    \pt_t^{-1}\psi_k\underset{(t)}{*}
                    \phi_j\underset{(x')}{*}
                    \pt_t f(t,x',\eta)
             \big\|_{L^p({\gr \re^{n-1}})}
        \Big\|_{L^1_t({\ppl \re})}
  \\
   \le &C\sum_{j\in \Z} 2^{sj}
          \Big\| {\bk \sum_{k\ge 2j}  2^{-(\frac12+\frac1{2p})k} } \\
            &\hskip24mm\times 
                 \big\|
                 \sum_{|m-j|\le 1} \overline{\Phi_m}\underset{(x',\eta)}{*} 
                 \phi_j\underset{(x')}{*}
                 {\bl \zeta_{j+1}}(\eta)\underset{(\eta)}{*}
                 {\bl \zeta_{j+2}}(\eta)\underset{(\eta)}{*} 
                 \pt_t f(t,x',\eta)  
                \big\|_{L^p({\gr \re^{n-1}})}                
          \big\|_{L^1_t(\re_+)}
   \\
      \le &C\Big\|\sum_{j\in \Z}  2^{(s-1)j}
                 \big\| 
                     \overline{\Phi_j}\underset{(x',\eta)}{*} 
                 \pt_t f(t,x',\eta)  
                \big\|_{L^p(\re^{n}_{x',\eta})}              
          \Big\|_{L^1_t(\re_+)} 
 \\
   \le &C\Big\|\sum_{j\in \Z}  2^{sj}
                 \big\| \N(-\Del)^{-1}
                        \overline{\Phi_j}\underset{(x',\eta)}{*} 
                 \pt_t f(t,x',\eta)  
                \big\|_{L^p(\re^{n}_{x',\eta})}              
          \Big\|_{L^1_t(\re_+)} 
 \\
    \le &C \|\pt_t\N (-\Del)^{-1} f                                                                                                         
           \|_{ L^1(\re_+;\dB^s_{p,1}(\re^n)) }
    \le  C \|\pt_t\N (-\Del)^{-1} f                                                                                                         
           \|_{ L^1(\re_+;\dB^s_{p,1}(\re^n_{+})) }
           {\bl + C\ep}.
     \addtocounter{equation}{1}
 \tag{\theequation} \label{eqn;time-like-p1}
\end{align*}

For the second term of  \eqref{eqn;time-like-p0}, 
we use the Minkowski inequality to see
\begin{align*}
  II 
   \le & C\sum_{j\in\Z} 2^{sj} \sum_{k\le 2j} 2^{(\frac12-\frac{1}{2p})k} \\
       &\hskip20mm\times            
        \Big\|
             \big\|
              \sum_{m\in \Z} \overline{\Phi_m}\underset{(x',\eta)}{*} 
         {\bl \zeta_{j+1}}(\eta)\underset{(\eta)}{*} 
         {\bl \zeta_{j+2}}(\eta)\underset{(\eta)}{*} 
              \phi_j\underset{(x')}{*}
              f(t,x',\eta)
             \big\|_{L^p({\gr \re^{n-1}})}
        \Big\|_{L^1_t(\re_+)}
\\
     \le & C\Big\| {\bk  \sum_{j\in\Z}  2^{(s+1-\frac1p)j} }              
             \big\|
              \sum_{|m-j|\le1} \overline{\Phi_m}\underset{(x',\eta)}{*} 
         {\bl \zeta_{j+1}}(\eta)\underset{(\eta)}{*} 
         {\bl \zeta_{j+2}}(\eta)\underset{(\eta)}{*} 
              \phi_j\underset{(x')}{*}
              f(t,x',\eta)
             \big\|_{L^p({\gr \re^{n-1}})}
        \Big\|_{L^1_t(\re_+)}
\\
       \le & C\Big\|   \sum_{j\in\Z}    2^{(s+1)j} 
             \big\|
               |\N|^{-1}\overline{\Phi_j}\underset{(x',\eta)}{*} 
              \N f(t,x',\eta)
             \big\|_{L^p(\re^{n}_{x',\eta})}
        \Big\|_{L^1_t(\re_+)}
    \le  C\big\|\N f \big\|_{ L^1(\re_+;\dB^s_{p,1}(\re^n_+)) }
         {\bl + C\ep}.
     \addtocounter{equation}{1}
\tag{\theequation} \label{eqn;time-like-p2}
\end{align*}

Combining the estimates 
\eqref{eqn;time-like-p0}, 
\eqref{eqn;time-like-p1}, 
\eqref{eqn;time-like-p2}, we obtain 
the first part of the left hand side of \eqref{eqn;pressure-trace}. 
The estimate for the spatial direction 
\eqref{eqn;pressure-trace-2} is slightly simpler.
For $1\le p<\infty$, 
{\crd
$\nabla f \in L^1(\re_+;\dB^{s}_{p,1}(\re^n_+))$ and 
noting \eqref{eqn;Sum_Phi_is_1},
we obtain that
\begin{align*}
  \big\|f(\cdot, &\cdot, \eta) 
  \big\|_{ L^1(\re_+;\dB^{s+1-\frac{1}{p}}_{p,1}({\gr \re^{n-1}})) } \\
   \le &C\Big\|\sum_{j\in \Z}
                 2^{(s+1-1/p)j}          
                \big\|\phi_j\underset{(x')}{*}
                \sum_{m\in\Z}
                \overline{\Phi_m}(x',\eta)\underset{(x',\eta)}{*}
                f(t,x,\eta)
                \big\|_{ L^p({\gr \re^{n-1}}) }  
        \Big\|_{L^1_t(\re_+)}
  \\
   \le &C\Big\|\sum_{j\in \Z} 2^{(s+1-\frac{1}{p})j}     
              \|\zeta_{\bl j+1}\|_{ L^{p'}(\re_{\eta}) }   
              \big\|
                   \phi_j\underset{(x')}{*}                   
              {\bl \zeta_{j+2}}(\eta)\underset{(\eta)}{*} 
               \overline{\Phi_j}(x',\eta)\underset{(x',\eta)}{*}
               f(t,x,\eta)
              \big\|_{ L^p(\re^{n}_{x',\eta}) }
        \Big\|_{L^1_t(\re_+)}
  \\
     \le &C\Big\|\sum_{j\in \Z} 2^{(s+1)j} 
            \big\| \overline{\Phi_j}(x',\eta)\underset{(x',\eta)}{*}
                   f(t,\cdot,\eta)
            \big\|_{L^p(\re^{n}_{x',\eta})}
           \Big\|_{L^1_t(\re_+)} 
     \le C\|\N f\|_{L^1(\re_+;\dB^{s}_{p,1}(\re^n_+))}{\bl + C\ep}.
\end{align*}
}

{\bl
For almost every $t\in \re_+$ and $\eta\in \re$,
\begin{align*}
  \big\|f(t, &\cdot, \eta) 
  \big\|_{ \dB^{s+1-\frac{1}{p}}_{p,1}({\gr \re^{n-1}}) } \\
   \le &C\sum_{j\in \Z}
                 2^{(s+1-1/p)j}          
                \big\|\phi_j\underset{(x')}{*}
                \sum_{m\in\Z}
                \overline{\Phi_m}(\cdot,\eta)\underset{(x',\eta)}{*}
                f(t,\cdot,\eta)
                \big\|_{ L^p({\gr \re^{n-1}}) }  
  \\
   \le &C\sum_{j\in \Z} 2^{(s+1-\frac{1}{p})j}     
              \|\zeta_{\bl j+1}\|_{ L^{p'}(\re_{\eta}) }   
              \big\|
                   \phi_j\underset{(x')}{*}                   
              {\bl \zeta_{j+2}}(\eta)\underset{(\eta)}{*} 
               \overline{\Phi_j}(\cdot,\eta)\underset{(x',\eta)}{*}
               f(t,\cdot,\eta)
              \big\|_{ L^p(\re^{n}_{x',\eta}) }
  \\
     \le &C\sum_{j\in \Z} 2^{(s+1)j} 
            \big\| \overline{\Phi_j}(x',\eta)\underset{(x',\eta)}{*}
                   f(t,\cdot,\cdot)
            \big\|_{L^p(\re^{n}_{x',\eta})}. 
\end{align*}
After taking supremum in $\eta>0$ and integrate it over $\re_+$, 
we obtain that
\begin{align*}
 \Big\|\sup_{\eta>0} \big\|f(\cdot, &\cdot, \eta) 
  \big\|_{\dB^{s+1-\frac{1}{p}}_{p,1}({\gr \re^{n-1}}) } 
 \Big\|_{L^1(\re_+)}
\\
      \le &C\Big\|\sum_{j\in \Z} 2^{(s+1)j} 
            \big\| \overline{\Phi_j}(x',\eta)\underset{(x',\eta)}{*}
                   f(t,\cdot,\eta)
            \big\|_{L^p(\re^{n}_{x',\eta})}
           \Big\|_{L^1_t(\re_+)} 
\\
     \le & C\|\N f\|_{L^1(\re_+;\dB^{s}_{p,1}(\re^n_+))}{\bl + C\ep}.
 \eqntag\label{eqn;pressure-trace-3}
\end{align*}
The estimate \eqref{eqn;pressure-trace-2} follows immediately from  \eqref{eqn;pressure-trace-3}.
}
This completes the proof of Proposition \ref{prop;grad-pressure-trace}.
\end{prf}

%
%
\subsection{Null-Lagrangian structure}
According to Evans \cite[section 8.1]{Ev},
we recall the null Lagrangian structure for the Jacobian of a Lipschitz 
continuous function  $u$. 

For $n\in \Nt$, let $A$ be a $n\times n$ matrix whose components are denoted by
$\{a_{kj}\}$ and consider its $\ell\times \ell$ sub-matrix 
$A^{[\ell]}$ given by 
\eq{\label{eqn;submatrix-ap}
A^{[\ell]}
=\begin{pmatrix}
      a_{\sg_1\t_1}      &\cdots &a_{\sg_1\t_\ell} \\
       \vdots            &\ddots &\vdots      \\
      a_{\sg_\ell\t_1}   &\cdots &a_{\sg_\ell\t_\ell}
 \end{pmatrix},
}
where $\sg_k,\t_j\in \{1,2,\cdots, n\}$ with 
$1\le \sg_1<\sg_2<\cdots,<\sg_{\ell}\le n$ and 
$1\le \t_1<\t_2<\cdots,<\t_{\ell}\le n$

\vskip3mm
\begin{lem}[Evans \cite{Ev}]
\label{lem;cofactor-divfree}
Let $1\le \ell\le n$ and let 
$u:\re^n\to \re^n$ be a Lipschitz continuous function  and  
$J(Du)^{[\ell]}$ denotes the $\ell\times\ell$
 sub-matrix of the Jacobi matrix
defined by \eqref{eqn;submatrix-ap} 
{\rm(cf.} \eqref{eqn;inverse-Jacobi-matrix}{\rm)}, 
$\cof(J(Du)^{[\ell]})_{kj}$ denotes 
the $(k,j)$- cofactor and $\cof(J(Du)^{[\ell]})$ be the cofactor matrix.
Then for any  $x\in \re^n$ with $\det(J(Du)(x)\neq0$ it holds that
$$
 \div_j\big(\cof(J(Du)^{[\ell]})\big)_{kj}=0.   
$$
\end{lem}
\vskip3mm

\begin{prf}{Lemma \ref{lem;cofactor-divfree}}  
Let $P=[p_{ij}]_{1\le i,j\le \ell}$  be a matrix whose $(i,j)$-components 
are $p_{ij}$ and its cofactor matrix be 
$\cof(P)$.  Then $(k,j)$-component of cofactor  is given by 
$$
 \cof(P)_{kj}=(-1)^{k+j}\det P^{[\ell-1]}_{kj}.
$$
Let $I$ be the $\ell\times \ell$ unit matrix and by 
$
 I=P^{\sf T} \, (P^{-1})^{\sf T},
$ 
it follows
\eq{ \label{eqn;matrix-expression}
  \det P\cdot \del_{ij}= \big(P^{\sf T}\, (\cof(P))\big)_{ij}
  =\sum_{k=1}^n(P^{\sf T})_{ik}\cof(P)_{kj}
  =\sum_{k=1}^np_{ki}\cof(P)_{kj}.    
}
Taking the partial derivative of the both side of 
\eqref{eqn;matrix-expression} by  $p_{km}$, 
the component  $p_{kj}$ is missing in $\cof(P)_{kj}$
\eq{\label{eqn;matrix-expression3}
 \frac{\pt}{\pt p_{km}}\det P =\cof(P)_{km}.   
}
Choose $P$ as the sub-matrix of the Jacobian $J(Du)^{[\ell]}$, 
the relation \eqref{eqn;matrix-expression} is now reduced into
\eq{ \label{eqn;matrix-expression2}
  \det J(Du)^{[\ell]}\cdot \del_{ij} 
   =\sum_{k=1}^{\ell} \bd_{ki}
    \cof(J(Du)^{[\ell]})_{kj},
}
where $\bd$ denotes the component of the Jacobi matrix $J(Du)$
defined by \eqref{eqn;element}.
Taking divergence for $j$-raw in 
\eqref{eqn;matrix-expression2} and noting 
\eqref{eqn;matrix-expression3},
\eqn{
 \spl{
  \sum_{j=1}^{\ell}\pt_j(\det J(Du)^{[\ell]})\cdot \del_{ij}
     =&\sum_{j=1}^{\ell} \pt_j \Big(\sum_{k=1}^{\ell} \bd_{ki}
                         \cof(J(Du)^{[\ell]})_{kj}\Big), 
\\
   \sum_{k=1}^{\ell}\sum_{m=1}^{\ell}
        \pt_i \bd_{km} \cdot\cof(J(Du)^{[\ell]})_{km}
    =&\sum_{j=1}^{\ell} \sum_{k=1}^{\ell}\pt_j  \bd_{ki}
       \cdot \cof(J(Du)^{[\ell]})_{kj}
     +\sum_{j=1}^{\ell} \sum_{k=1}^{\ell} \bd_{ki}
       \cdot \pt_j \cof(J(Du)^{[\ell]})_{kj}
 }
}
and thus we obtain
\eqn{
\sum_{k=1}^{\ell}\bd_{ki}\cdot 
\sum_{j=1}^{\ell}  \pt_j \cof(J(Du))_{kj}=0.    
}
Rewrite the above as 
$$
0=J(Du)^{[\ell]} \div_j( \cof(J(Du)^{[\ell]})_{kj})
$$ 
and multiplying the both side by $(J(Du)^{[\ell]})^{-1}$ 
at the point $x_0$ with 
$\det J(Du)^{[\ell]}(x_0)\neq 0$, it follows that 
\eqn{
  \div_j \cof(J(Du)^{[\ell]})_{kj}=0.
} 
\end{prf}
\vskip3mm

%
%

\subsection{Bilinear estimates}

The following bilinear estimate is well-known:

\begin{lem} \label{lem:para-product}
Let $1 \le p \le \infty$, $1 \le \sg \le \infty$. \par
\noindent
 If $s >0$ then for all $f\in L^{q_2}(\re^n)\cap \dot B^s_{r_1, \sg}(\re^n)$ and
$g\in L^{r_2}(\re^n)\cap \dot B^s_{q_1, \sg}(\re^n)$, 
\begin{equation}
   \|f g\|_{\dot B^s_{p, \sg}}
    \le   C\big(\|f\|_{\dB^s_{r_1, \sg}}\|g\|_{r_2}
             +  \|f\|_{q_2} \|g\|_{\dB^s_{q_1, \sg}}\big),
    \label{eqn;bilinear-1}
\end{equation}
where 
$$
  \frac{1}{p}=\frac{1}{r_1}+\frac{1}{r_2}=\frac{1}{q_1}+\frac{1}{q_2}   
$$
{\ppl and} $C>0$ is independent of $f$ and $g$.
\end{lem}

\begin{prf}{Lemma \ref{lem:para-product}}
Let $P_k  g = \dsp{\sum_{\ell=-\infty}^{k-3}\phi _{\ell}\ast g} =
\psi_{2^{-(k-3)}}\ast g$
denotes the low frequency part of the Littlewood--Paley decomposition of $g$.
Then by the para-product decomposition of the product of $f$ and $g$,
\begin{eqnarray}
f\cdot g
  &=&  \sum_{k\in\Z}(\phi_k\ast f)(P_kg)
     +\sum_{k\in\Z}(P_kf)(\phi_k\ast g) 
     + \sum_{k\in\Z}\sum_{|l-k| \le 2}
     (\phi _k\ast f)(\phi _l\ast g) \nonumber  \\
&\equiv& h_1 + h_2 + h_3. \label{eqn;para-1}
\end{eqnarray}
 Since
\begin{eqnarray*}
&& \mbox{supp }\F\left((\phi _k\ast f)(P_kg)\right) \subset \{\xi \in
\re^n; 2^{k-2} \le |\xi| \le 2^{k+2}\}, 
\end{eqnarray*}
we have by the Young inequality that
for
$\dsp
  \frac{1}{p}=\frac{1}{r_1}+\frac{1}{r_2},
$
\begin{equation}  \label{eqn;para-2}
\begin{split}
  \|h_1\|_{\dot B^s_{p, \sg}}
    \le & \left\{ \sum_{j\in\Z}\left( 2^{sj}
          \|\phi_j\|_{1}\|\widetilde{\phi_{j}}\ast f\|_{r_1}
         \|P_{j+2}g\|_{r_2}\right)^{\sg}\right\}^{1/\sg} \\
    \le & C\| g\|_{r_2}
         \left\{ \sum_{j\in\Z}\left( 2^{sj}
          \|\widetilde{\phi_j}\ast f\|_{r_1}\right)^{\sg}\right\}^{1/\sg} \\
    \le & C\|f\|_{\dB^{s}_{r_1,\sg}}\|g\|_{r_2},
\end{split}
\end{equation}
where $C=\|\F^{-1}\phi\|_{1}\|\psi\|_{1}$. 
By replacing the role of $f$ and $g$ with that of $g$ and $f$,
respectively, we see that the second term can be handled in the
similar way as above.  Hence there holds
\begin{equation}\label{eqn;para-3}
   \|h_2\|_{\dot B^s_{p, \sg}}
     \le C\|f\|_{q_2}\|g\|_{\dot B^s_{q_1, \sg}} ,
\quad \frac{1}{p}=\frac{1}{q_1}+\frac{1}{q_2}.
\end{equation}

To deal with the third term, we should notice that
$$
\mbox{supp }\F(\phi_k\ast f\cdot \phi_l\ast g) \subset \{\xi \in
\re^n; |\xi| \le 2^{\max\{k, l\}+2}\},
$$
so there holds
$$
\phi_j \ast((\phi_k\ast f)(\phi_l\ast g) )= 0
\quad\mbox{for $\max\{k, l\} \le j -3$}.
$$
{\gr 
Let $r_1$ and $r_2$ satisfy 
\begin{equation}
  \frac 1p
   =\frac 1{r_1}+\frac 1{r_2},
   \label{eqn;exponent1}
\end{equation}
}
then 
{\gr
\begin{align*}
\|h_3\|_{\dot B^{s}_{p, \sg}}
   =&  \Big(\sum_{j\in \Z} 
        \Big(2^{sj}
        \Big\|\sum_{k\ge j-2}
           \phi_j\ast( (\phi_k* f)(\widetilde{\phi_{k}}* g) )
        \Big\|_{p}
        \Big)^{\sg}\Big)^{1/\sg}\\
\le & \Big(\sum_{j\in \Z} 
           \Big(
            2^{sj}
            \sum_{k\ge j-2}
            \|\phi_j\|_{1} \|\phi_{k}\ast f\|_{r_1}
            \|\widetilde{\phi}_{k}\ast g\|_{r_2}
           \Big)^{\sg}
      \Big)^{1/\sg}
\\
 \le & C
      \Big(
      \sum_{j\in \Z}
      \Big(
      \sum_{k\ge j-2} 2^{ sj } 
         \|\phi_{k}\ast f\|_{r_1}
      \Big)^{\sg}
      \Big)^{1/\sg}
      \sup_{k\in \Z}  \|\widetilde{\phi}_{k}\ast g\|_{r_2}
      \\
  &\text{(changing $k'=k-j$ to see)}
 \\
 \le & C
      \Big(
      \sum_{j\in \Z}
      \Big(
      \sum_{k'\ge -2} 2^{ s(k'+j) } 2^{-sk'}  
         \|\phi_{k'+j}\ast f\|_{r_1}
      \Big)^{\sg}
      \Big)^{1/\sg}
      \sup_{k\in \Z}  \|\widetilde{\phi}_{k}\ast g\|_{r_2}
      \\
 \le & C\sum_{k'\ge -2} 2^{-sk'} 
      \Big(
      \sum_{j\in \Z} 
      \Big(2^{ s(k'+j)}            
           \|\phi_{k'+j}\ast f\|_{r_1}^{\sg}
      \Big)^{\sg}
      \Big)^{1/\sg}
      \sup_{k\in \Z}  \|\widetilde{\phi}_{k}\ast g\|_{r_2}
\\
 \le & C\Big(
      \sum_{j\in \Z}  
          2^{sj\sg}  
         \|\phi_{j}\ast f\|_{r_1}^{\sg}
       \Big)^{1/\sg}
       \| g\|_{r_2}
\\
  \le  & C \|f\|_{\dot B^{s}_{r_1,\sg}}
         \|g\|_{L^{r_2}},
 \eqntag \label{eqn;para-5}
\end{align*} 
}
 where we use $s>0$. 
The estimate \eqref{eqn;bilinear-1} 
follow from  \eqref{eqn;para-1},  \eqref{eqn;para-2},
\eqref{eqn;para-3} 
and \eqref{eqn;para-5}. 
 
\end{prf}
\vskip3mm

The following bilinear estimates over the whole space $\re^n$ are
 obtained by Abidi--Paicu 
\cite{AP07} (cf. \cite{OgSs18}).

\vskip3mm
\noindent
\begin{prop}[\cite{AP07}] 
\label{prop;abidi-paicu2}
Let $1\le p,\, p_1,\, p_2,\, \sg,\, \lam_1,\, \lam_2 \le \infty$,
$1/p\le 1/p_1+1/p_2$, $p_1\le  \lam_2$,
$p_2\le \lam_1$ and
$$
 s_1+s_2+n\inf\bigg(0,1-\frac{1}{p_1}-\frac{1}{p_2}\bigg)>0,\quad 
 \frac{1}{p}\le \frac{1}{p_1}+\frac{1}{\lam_1}\le 1,\quad
 \frac{1}{p}\le \frac{1}{p_2}+\frac{1}{\lam_2}\le 1.
$$
\begin{itemize}
\item[{\rm (1)}]
If  
$s_1+\frac{n}{\lam_2}<\frac{n}{p_1}$ and 
$s_2+\frac{n}{\lam_1}<\frac{n}{p_2}$, 
then there exists $C>0$ such that 
for all $f\in \dB^{s_1}_{p_1,\sg}$ and 
$g\in \dot B^{s_2}_{p_2,\infty}$, 
the following estimate holds 
\eqn{
  \|fg\|_{\dB^{s_1+s_2-n\big(\frac{1}{p_1}+\frac{1}{p_2}-\frac{1}{p}\big)}_{p,\sg}}
  \le C \|f\|_{\dB^{s_1}_{p_1,\sg}} 
        \|g\|_{\dB^{s_2}_{p_2,\infty}} .  
 }
\item[{\rm (2)}]If $s_1+\frac{n}{\lam_2}=\frac{n}{p_1}$ and 
$s_2+\frac{n}{\lam_1}=\frac{n}{p_2}$, then 
there exists $C>0$ such that 
for any 
$f\in \dB^{s_1}_{p_1,1}$ and 
$g\in \dB^{s_2}_{p_2,1}$
the following estimate holds 
\eq{\label{eqn;bilinear1.0-ap}
  \|fg\|_{\dB^{s_1+s_2-n\big(\frac{1}{p_1}+\frac{1}{p_2}-\frac{1}{p}\big)}_{p,1}}
  \le C \|f\|_{\dB^{s_1}_{p_1,1}} 
        \|g\|_{\dB^{s_2}_{p_2,1}} .  
 }
 \item[{\rm (3)}] In particular, if $s_1=-1+n/p$, $s_2=n/p$ and 
$-1+n/p+n/p>\inf(0,n-2n/p)$  in (2), i.e., $1\le p<2n$  then
there exists $C>0$ such that the following estimate holds 
\eq{\label{eqn;bilinear1.1-ap}
  \|fg\|_{\dB^{-1+\frac{n}{p}}_{p,1}}
  \le C \|f\|_{\dB^{-1+\frac{n}{p}}_{p, 1}} 
        \|g\|_{\dB^{\frac{n}{p}}_{p,1}}.
 }
\end{itemize}
\end{prop}
Since Danchin--Mucha \cite{DM12} treats the equations 
depending on the density, the restriction on the exponent 
$p$ in the solution space $\dB^{-1+n/p}_{p,1}(\re^n)$ stems 
from the restriction on $1\le p<2n$ for the above bilinear 
estimate \eqref{eqn;bilinear1.1-ap}. 
One may improve the restriction by using the divergence free - 
curl free structure of nonlinear terms.

The bilinear estimates as above hold for the case when the  
two functions $f:\re^n\to \re^n$ and $g: \re^n\to \re^n$ {\ppl have} the 
 divergence structure condition:
\eqn{  
 f\cdot D_x g =D_x (f\cdot g),
 } 
where $D_x$ denotes any combination of partial derivatives by 
$x=(x_1,x_2,\cdots, x_n)$ of the first order.
A typical case is given by the form when $f$ and $g$ satisfies 
{\it divergence free-rotation free structure} as 
$\div f=0$ and $\rot g=0$.

\vskip3mm
\begin{prop}[Bilinear estimate under divergence structure]
\label{prop;Besov-bilinear-div-curl}
Let $1\le p < \infty$ and 
$f\in  \dB^{-1+n/p}_{p,1}$ 
and $g\in   \dB^{n/p}_{p,1}$. \par

\noindent{\rm (1)} If 
 there exists $F=F(x)$  such that  $f\cdot g =D_{x}(F\cdot g)$ 
with $f=D_x F(x)$ in the sense of distribution, 
where $D_x$ is any combination of the first differentiation in $x$. Then
\eq{ \label{eqn;bilinear5}
  \|f\cdot g\|_{\dot{B}^{-1+\frac{n}{p}}_{p,1}}
  \le
    C\|f\|_{\dot{B}^{-1+\frac{n}{p}}_{p,1}} 
     \|g\|_{\dot{B}^{\frac{n}{p}}_{p,1}}.
}
\noindent{\rm (2)} In particular
with additional conditions  $\div f=0$, $\rot g=0$ in the distribution sense,
then  
\eq{\label{eqn;bilinear6}
  \|f\cdot g\|_{\dot{B}^{-1+\frac{n}{p}}_{p,1}}
  \le
    C\|f\|_{\dot{B}^{-1+\frac{n}{p}}_{p,1}} 
     \|g\|_{\dot{B}^{\frac{n}{p}}_{p,1}}.
}
\end{prop}

\vskip 2mm 
\begin{prf}{Proposition \ref{prop;Besov-bilinear-div-curl}}
Here we show the case when $n<p<\infty$ since the other cases 
$1\le p\le n$  are already proved in Proposition \ref{prop;abidi-paicu2}. 
The second estimate is directly obtained from the first part 
by observing that $f$, $g$ are both vector-valued functions and satisfy 
$\div f=0$, $\rot g=0$ in the sense of distribution.  
Then
$$
 \div (F\wedge g)=\rot F\cdot g-F\cdot (\rot g), 
$$
it holds that 
$f\cdot g=(\rot F) \cdot g =\div (F\wedge g)$
which represents the divergence form structure.

Hence we assume that there exists a function $G$ such that 
$g=D_x G$ and $f\cdot g=D_x(f G)$.
Using \eqref{eqn;bilinear1.1-ap}, it follows that
\eqn{
 \spl{
       \|f g\|_{\dot{B}^{-1+\frac{n}{p}}_{p,1}}
   = & \|D_x(F g)\|_{\dot{B}^{-1+\frac{n}{p}}_{p,1}}
   \le  C\|Fg\|_{\dot{B}^{\frac{n}{p}}_{p,1}} \\
   \le &
        C\|F\|_{\dot{B}^{\frac{n}{p}}_{p,1}} 
         \|g\|_{\dot{B}^{\frac{n}{p}}_{p,1}} 
   \le 
        C\|f\|_{\dot{B}^{-1+\frac{n}{p}}_{p,1}} 
         \|g\|_{\dot{B}^{\frac{n}{p}}_{p,1}}.
 }}
 In particular \eqref{eqn;bilinear5} and hence \eqref{eqn;bilinear6} holds.
\end{prf}

\vskip3mm
{\crd 
\noindent
\begin{prop}[The space-time bilinear estimate]\label{prop;double-bony-bilinear} 
Let $1\le \r\le \infty $  and $1<p<2n-1$. Then for 
$F\in \widetilde{\dB^{1/2-1/(2p)}_{\r,1}(\re_+;}\dB^{-1+n/p}_{p,1}(\re^{n-1}))
\cap \widetilde{L^{\r}(\re_+;}\dB^{(n-1)/p}_{p,1}(\re^{n-1}))$
and 
$G\in \widetilde{\dB^{1/2-1/(2p)}_{\infty,1}
           (\re_+;}\dB^{-1+n/p}_{p,1}(\re^{n-1}))
      \cap \widetilde{L^{\infty}(\re_+;}\dB^{(n-1)/p}_{p,1}(\re^{n-1}))$, 
it holds that
\algn{
 \big\| F\, G &
 \big\|_{\widetilde{\dB^{\frac12-\frac{1}{2p}}_{\r,1}
          (\re_+};\dB^{-1+\frac{n}{p}}_{p,1}(\re^{n-1}))
        }
 \\
\le & 
    C\big\| F(t)
     \big\|_{\widetilde{L^{\r}(\re_+;}\dB^{\frac{n-1}{p}}_{p,1}(\re^{n-1}))}
     \big\| G(t)     
     \big\|_{\widetilde{\dB^{\frac12-\frac{1}{2p}}_{\infty,1}
           (\re_+;}\dB^{-1+\frac{n}{p}}_{p,1}(\re^{n-1}))}
\\
  & 
  + C\Big(
      \big\| F(t)
      \big\|_{\widetilde{\dB^{\frac12-\frac{1}{2p}}_{\r,1}(\re_+;}\dB^{-1+\frac{n}{p}}_{p,1}(\re^{n-1}))} 
      +
      \big\| F(t)
      \big\|_{\widetilde{L^{\r}(\re_+;}\dB^{\frac{n-1}{p}}_{p,1}(\re^{n-1}))}
     \Big) 
       \big\| G(t)          
       \big\|_{\widetilde{L^{\infty}(\re_+;}\dB^{\frac{n-1}{p}}_{p,1}(\re^{n-1}))} 
\\
 \le & {\bl 
  C\Big(
      \big\| F(t)
      \big\|_{\widetilde{\dB^{\frac12-\frac{1}{2p}}_{\r,1}(\re_+;}\dB^{-1+\frac{n}{p}}_{p,1}(\re^{n-1}))} 
      +
      \big\| F(t)
      \big\|_{\widetilde{L^{\r}(\re_+;}\dB^{\frac{n-1}{p}}_{p,1}(\re^{n-1}))}
     \Big) 
     }
  \\
   &\hskip1cm {\bl
   \times
     \Big(
       \big\| G(t)     
       \big\|_{\widetilde{\dB^{\frac12-\frac{1}{2p}}_{\infty,1}
           (\re_+;}\dB^{-1+\frac{n}{p}}_{p,1}(\re^{n-1}))}
    +
       \big\| G(t)          
       \big\|_{\widetilde{L^{\infty}(\re_+;}\dB^{\frac{n-1}{p}}_{p,1}(\re^{n-1}))} 
     \Big),
   }
}
where the norms are defined in \eqref{eqn;chemin-lerner-norms}.
\end{prop}

We should like to note that when $\r=1$, the following spaces are norm-equivalent;
\begin{gather}  
 \widetilde{\dB^{\frac12-\frac{1}{2p}}_{1,1}
          (\re_+};\dB^{-1+\frac{n}{p}}_{p,1}(\re^{n-1}))
  \simeq \dB^{\frac12-\frac{1}{2p}}_{1,1}
          (\re_+;\dB^{-1+\frac{n}{p}}_{p,1}(\re^{n-1}))
  \simeq \dF^{\frac12-\frac{1}{2p}}_{1,1}
          (\re_+;\dB^{-1+\frac{n}{p}}_{p,1}({\gr \re^{n-1}})),
 \label{eqn;F-B-space-equiv-1}
\\
 \widetilde{L^{1}(\re_+;}\dB^{\frac{n-1}{p}}_{p,1}(\re^{n-1}))
  \simeq  L^1(\re_+;\dB^{\frac{n-1}{p}}_{p,1}(\re^{n-1})).
  \label{eqn;F-B-space-equiv-2}
\end{gather}

\vskip3mm
\begin{prf}{Proposition \ref{prop;double-bony-bilinear} }
We employ the doubled Bony paraproduct decomposition in both 
space and time direction:
\begin{align*}
   \big\| F\, G
   \big\|_{\widetilde{\dB^{\frac12-\frac{1}{2p}}_{\r,1}
           \big(\re_+;}\dB^{-1+\frac{n}{p}}_{p,1}(\re^{n-1})\big)} 
  \le & 
     \sum_{j\in\Z} 2^{(-1+\frac{n}{p})j}
     \sum_{k\ge 2j}  
              2^{(\frac12-\frac{1}{2p})k}
        \Big\| 
               \big\|
                 \psi_k\underset{(t)}{*}
                 \phi_j\underset{(x')}{*}
                 \big(F\, G\big)
              \big\|_{L^p(\re^{n-1})}
       \Big\|_{L^{\r}_t(\re_+)} \\
  &+ \sum_{j\in\Z} 2^{(-1+\frac{n}{p})j}
     \sum_{k\le 2j}  
            2^{(\frac12-\frac{1}{2p})k}
        \Big\| 
              \big\|
                 \psi_k\underset{(t)}{*}
                 \phi_j\underset{(x')}{*}
                 \big(F\, G\big)
              \big\|_{L^p(\re^{n-1})}
       \Big\|_{L^{\r}_t(\re_+)} 
 \\
 \equiv & I+II.
  \addtocounter{equation}{1}
 \tag{\theequation} \label{eqn;time-like-N00}
 \end{align*}

\vskip2mm
The estimate for the second term of the right hand side of 
\eqref{eqn;time-like-N00} is straightforward:
\algn{
  II
   \le & \sum_{j\in\Z} 2^{(-1+\frac{n}{p})j}
           \sum_{k\le 2j}  
            2^{(\frac12-\frac{1}{2p})k}
       \Big\|| \psi_k|\underset{(t)}{*}
              \big\|                 
                 \phi_j\underset{(x')}{*}
                 \big(F\, G\big)
              \big\|_{L^p(\re^{n-1})}
       \Big\|_{L^{\r}_t(\re_+)} 
 \\
    \le &C\sum_{j\in\Z} 2^{(-1+\frac{n}{p})j}
           \sum_{k\le 2j}  
            2^{(\frac12-\frac{1}{2p})k}
       \Big\|\big\|                 
                 \phi_j\underset{(x')}{*}
                 \big(F\, G\big)
              \big\|_{L^p(\re^{n-1})}
       \Big\|_{L^{\r}_t(\re_+)} 
 \\
    \le &C\sum_{j\in\Z} 2^{(-1+\frac{n}{p})j}
            2^{(1-\frac{1}{p})j}
       \Big\|
          \big\|                 
                \phi_j\underset{(x')}{*}
                \big(F\, G\big)
          \big\|_{L^p(\re^{n-1})}
       \Big\|_{L^{\r}_t(\re_+)} 
 \\
    \le &{\bl
       C\sum_{j\in\Z} 2^{\frac{n-1}{p}j}
       \Big\| \|\phi_j\|_{L^1(\re^{n-1})}     
        \Big(
          \big\| (\widetilde{\phi}_{j}\underset{(x')}{*}F)
                 (P_{j}\underset{(x')}{*} G)
          \big\|_{L^p(\re^{n-1})}
        + \big\| (P_{j}\underset{(x')}{*} F)
                 (\widetilde{\phi}_{j}\underset{(x')}{*}  G)
          \big\|_{L^p(\re^{n-1})} 
       }
          \\
      &\hskip15mm
      {\bl
        + \big\| \sum_{\ell \ge j-2}
                 (\phi_{\ell}\underset{(x')}{*}F) 
                 (\widetilde{\phi}_{\ell}\underset{(x')}{*}G)
          \big\|_{L^p(\re^{n-1})}
        \Big)
       \Big\|_{L^{\r}(\re_+)} 
      }
 \\
     \le &{\bl
       C\sum_{j\in\Z} 2^{\frac{n-1}{p}j}
       \Big\|    
          \big\| \widetilde{\phi}_{j}\underset{(x')}{*}F
          \big\|_{L^p(\re^{n-1})}
          \big\| G
          \big\|_{L^{\infty}(\re^{n-1})}
        + \big\|  F
          \big\|_{L^{\infty}(\re^{n-1})}
          \big\| \widetilde{\phi}_{j}\underset{(x')}{*}  G
          \big\|_{L^p(\re^{n-1})} 
       }
          \\
      &\hskip15mm
      {\bl
        + \sum_{\ell \ge j-2}
          \big\| \phi_{\ell}\underset{(x')}{*}F
          \big\|_{L^\infty(\re^{n-1})}
          \big\| \widetilde{\phi}_{\ell}\underset{(x')}{*}G
          \big\|_{L^p(\re^{n-1})}
        \Big)
       \Big\|_{L^{\r}(\re_+)} 
      }
 \\
       \le &C\Big(
          \big\| F
          \big\|_{\bl L^{\r}(\re_+;L^{\infty}(\re^{n-1}))} 
          \big\| G
          \big\|_{\widetilde{L^{\infty}(\re_+;}\dB^{\frac{n-1}{p}}_{p,1}(\re^{n-1}))}
       + 
          \big\| G
          \big\|_{{\bl L^{\infty}(\re_+;}L^{\infty}(\re^{n-1}))} 
          \big\| F
          \big\|_{\widetilde{L^{\r}(\re_+;}\dB^{\frac{n-1}{p}}_{p,1}(\re^{n-1}))}
          \Big) \\
      & {\bl
        + C\sum_{j\in\Z} 2^{\frac{n-1}{p}j}
         \Big\|\sum_{\ell \ge j-2} 2^{\frac{n-1}{p}\ell}
           \big\| \phi_{\ell}\underset{(x')}{*}F
           \big\|_{L^p(\re^{n-1})}
           \big\| \widetilde{\phi}_{\ell}\underset{(x')}{*}G
           \big\|_{L^p(\re^{n-1})}
         \Big\|_{L^{\r}(\re_+)} 
      }
 \\
    \le &C\Big(
          \big\| F
          \big\|_{{\bl L^{\r}(\re_+;}L^{\infty}(\re^{n-1}))} 
          \big\| G
          \big\|_{\widetilde{L^{\infty}(\re_+;}\dB^{\frac{n-1}{p}}_{p,1}(\re^{n-1}))}
       + 
          \big\| G
          \big\|_{{\bl L^{\infty}(\re_+;}L^{\infty}(\re^{n-1}))} 
          \big\| F
          \big\|_{\widetilde{L^{\r}(\re_+;}\dB^{\frac{n-1}{p}}_{p,1}(\re^{n-1}))}
          \Big) \\
      & {\bl
        + C\sum_{\ell\in\Z} 2^{\frac{n-1}{p}\ell}
          \sum_{j\le \ell+2} 2^{\frac{n-1}{p}j}
          \Big\|
           \big\| \phi_{\ell}\underset{(x')}{*}F
           \big\|_{L^p(\re^{n-1})}
          \Big\|_{L^{\r}(\re_+)} 
          \Big\|
           \big\| \widetilde{\phi}_{\ell}\underset{(x')}{*}G
           \big\|_{L^p(\re^{n-1})}
          \Big\|_{L^{\infty}(\re_+)} 
      }
\\
    \le &C\big\| F
          \big\|_{\widetilde{L^{\r}(\re_+;}\dB^{\frac{n-1}{p}}_{p,1}(\re^{n-1}))} 
          \big\| G
          \big\|_{\widetilde{L^{\infty}(\re_+;}\dB^{\frac{n-1}{p}}_{p,1}(\re^{n-1}))}. 
   }
For the estimate of $I$ in \eqref{eqn;time-like-N00}, 
we employ the double Bony decomposition in both time and space regions: 
Set 
$P_m\underset{(x')}{*}F\equiv \sum_{m'\le m}\phi_{m'}\underset{(x')}{*}F$
and  
$Q_\ell\underset{(t)}{*}F\equiv \sum_{\ell'\le \ell}\psi_{\ell'}\underset{(t)}{*}F$.
Then the Bony decomposition in space direction gives 
 \algn{
 I \le &
    \sum_{j\in \Z} 2^{(-1+\frac{n}{p})j}
    \sum_{k\ge 2j}  2^{(\frac12-\frac{1}{2p})k}   
    \Big\|     \big\|
            \psi_{k}\underset{(t)}{*}
            \phi_{j}\underset{(x')}{*}  \\
       &\hskip2cm \times           
              \Big(\sum_{m\in \Z}
                \big(
                  \phi_{m}\underset{(x')}{*}F(t)
                \big)\cdot
                \big(  
                  P_{m}\underset{(x')}{*}
                  G(t)\big)
            \Big)
        \big\|_{L^p(\re^{n-1})}
     \Big\|_{L^{\r}_t(\re_+)} \\
    &+  
     \sum_{j\in \Z} 2^{(-1+\frac{n}{p})j}
     \sum_{k\ge 2j}  2^{(\frac12-\frac{1}{2p})k}   
     \Big\|    
       \big\|
             \psi_{k}\underset{(t)}{*}
             \phi_{j}\underset{(x')}{*}  \\
       &\hskip2cm \times                  
              \Big(\sum_{m \in \Z}
                \big(
                  P_{m}\underset{(x')}{*}F(t)
                \big)\cdot
                \big(
                 \phi_{m}\underset{(x')}{*}
                 G(t)\big)
             \Big)
        \big\|_{L^p(\re^{n-1})}
  \Big\|_{L^{\r}_t(\re_+)} \\
  &\quad+ 
    \sum_{j\in \Z} 2^{(-1+\frac{n}{p})j}
        \sum_{k\ge 2j}  2^{(\frac12-\frac{1}{2p})k}   
     \Big\|    \big\|
             \psi_{k}\underset{(t)}{*}
             \phi_{j}\underset{(x')}{*}   \\
       &\hskip2cm \times           
              \Big(\sum_{m\ge j-2}
                \big(
                 \phi_{m}\underset{(x')}{*}F(t)
                \big)\cdot
                \big( 
                 \widetilde{\phi}_{m}\underset{(x')}{*}
                 G(t)\big)
             \Big)
        \big\|_{L^p(\re^{n-1})}
  \Big\|_{L^{\r}_t(\re_+)} 
  \\
\equiv & h_1+h_2+h_3
\le  \sum_{\t=1}^3h_1^{\t}
    +\sum_{\t=1}^3h_2^{\t}
     +\sum_{\t=1}^3h_3^{\t},
\eqntag \label{eqn;space-Bony-decomp0}
}
where $h_1$ can be decomposed by 
the Bony decomposition in time direction such as the following:
\algn{
  h_1 
  \le&  \sum_{j\in \Z} 2^{(-1+\frac{n}{p})j}
        \sum_{k\ge 2j}  2^{(\frac12-\frac{1}{2p})k}   \\
     &\hskip1cm \times 
       \Big\| 
          \big\|
            \psi_{k}\underset{(t)}{*}
            \phi_{j}\underset{(x')}{*}   
            \Big(\sum_{\ell\in \Z}
                \big(
                  \widetilde{\psi}_{\ell}\underset{(t)}{*}
                  \widetilde{\phi}_{j}\underset{(x')}{*}F(t)
                \big)\cdot
                \big(  
                  Q_{\ell-2}\underset{(t)}{*}
                  P_{j-2}\underset{(x')}{*}
                  G(t)\big)
            \Big)
        \big\|_{L^p(\re^{n-1})}
     \Big\|_{L^{\r}_t(\re_+)}   
 \\
   &+ \sum_{j\in \Z} 2^{(-1+\frac{n}{p})j}
      \sum_{k\ge 2j}  2^{(\frac12-\frac{1}{2p})k}    \\
     &\hskip1cm \times 
      \Big\| 
        \big\|
            \psi_{k}\underset{(t)}{*}
            \phi_{j}\underset{(x')}{*}   
            \Big(\sum_{\ell\in \Z}
                \big(
                  Q_{\ell-2}\underset{(t)}{*}
                  \widetilde{\phi}_j\underset{(x')}{*}
                  F(t)
                \big)\cdot
                \big(  
                  \widetilde{\psi}_{\ell}\underset{(t)}{*}
                  P_{j-2}\underset{(x')}{*}
                  G(t)\big)
            \Big)
        \big\|_{L^p(\re^{n-1})}
     \Big\|_{L^{\r}_t(\re_+)}    \\
    &+\sum_{j\in \Z} 2^{(-1+\frac{n}{p})j}
      \sum_{k\ge 2j}  2^{(\frac12-\frac{1}{2p})k}    \\
     &\hskip1cm \times 
        \Big\| \big\|
            \psi_{k}\underset{(t)}{*}
            \phi_{j}\underset{(x')}{*}   
            \Big(\sum_{\ell\ge k-2}
                \big(
                  \psi_{\ell}\underset{(t)}{*}
                  \widetilde{\phi}_{j}\underset{(x')}{*}F(t)
                \big)\cdot
                \big(  
                 \widetilde{\psi}_{\ell}\underset{(t)}{*}
                  P_{j-2}\underset{(x')}{*}
                  G(t)\big)
            \Big)
        \big\|_{L^p(\re^{n-1})}
     \Big\|_{L^{\r}_t(\re_+)}   
  \\
  \equiv & h_1^1+h_1^2+h_1^3.
}
The estimates for the terms $h_1^1$, $h_1^2$ are  straightforward.  
For instance, 
\algn{
h_1^2
 \le &\sum_{j\in \Z} 2^{(-1+\frac{n}{p})j}
        \sum_{k\ge 2j}  2^{(\frac12-\frac{1}{2p})k} 
 \\
    &\hskip5mm \times  
      \Big\|  |\psi_{k}|\underset{(t)}{*}  \big\|           
            \phi_{j}\underset{(x')}{*}   
            \Big(
                \big(
                  Q_{k-2}\underset{(t)}{*}
                  \widetilde{\phi}_{j}\underset{(x')}{*}
                  F(t)
                \big)\cdot
                \big(  
                  \widetilde{\psi}_{k}\underset{(t)}{*}
                  P_{j-2}\underset{(x')}{*}
                  G(t)\big)
            \Big)
        \big\|_{L^p(\re^{n-1})}
     \Big\|_{L^{\r}_t(\re_+)}    
\\
 \le &C \sum_{j\in \Z} 2^{(-1+\frac{n}{p})j}
        \sum_{k\ge 2j}  2^{(\frac12-\frac{1}{2p})k} \\
     &\hskip1cm \times  
       \Big\|  \big\|
            \Big(
                \big(
                  Q_{k-2}\underset{(t)}{*}
                  \widetilde{\phi}_{j}\underset{(x')}{*}
                  F(t)
                \big)\cdot
                \big(  
                  \widetilde{\psi}_{k}\underset{(t)}{*}
                   P_{j-2}\underset{(x')}{*}
                  G(t)\big)
            \Big)
        \big\|_{L^p(\re^{n-1})}
     \Big\|_{L^{\r}_t(\re_+)}    
 \\
    \le &C\sum_{k\in \Z} 2^{(\frac12-\frac{1}{2p})k}
          \sum_{2j\le k} 2^{(-1+\frac{n}{p})j}  \\
     &\hskip1cm \times  
       \Big\| 
         \big\|
               Q_{k-2}\underset{(t)}{*}
               \widetilde{\phi}_{j}\underset{(x')}{*}
               F(t)
        \big\|_{L^{p}(\re^{n-1})}
        \big\| 
             \widetilde{\psi}_{k}\underset{(t)}{*}
             \sum_{m\le j-2} \phi_m\underset{(x')}{*}
             G(t)\big)
        \big\|_{L^{\infty}(\re^{n-1})}
     \Big\|_{L^{\r}_t(\re_+)}
 \\
  \le &C\sum_{k\in \Z} 2^{(\frac12-\frac{1}{2p})k}
        \sum_{2m\le k}\sum_{j\ge m+2} 2^{{\rd(-1+\frac{1}{p})}j} 
                    \\
     &\hskip1cm \times  
       \Big\|  \sup_{j\in \Z}\big\| 2^{\frac{n-1}{p}j} 
               Q_{k-2}\underset{(t)}{*}
               \widetilde{\phi}_{j}\underset{(x')}{*}
               F(t)
        \big\|_{L^{p}(\re^{n-1})}
        \big\| 
             \widetilde{\psi}_{k}\underset{(t)}{*}
             \phi_m\underset{(x')}{*}
             G(t)\big)
        \big\|_{L^{\infty}(\re^{n-1})}
     \Big\|_{L^{\r}_t(\re_+)}
\\
      \le &C\sum_{k\in \Z} 2^{(\frac12-\frac{1}{2p})k}
            \sum_{2m\le k}2^{(-1+\frac{1}{p})m}         
                   \\
     &\hskip1cm \times  
        \Big\| \sup_{j\in \Z}\big\| 2^{\frac{n-1}{p}j} 
               Q_{k-2}\underset{(t)}{*}
               \widetilde{\phi}_{j}\underset{(x')}{*}
               F(t)
        \big\|_{L^{p}(\re^{n-1})}
        \big\| 
             \widetilde{\psi}_{k}\underset{(t)}{*}
             \phi_m\underset{(x')}{*}
             G(t)\big)
        \big\|_{L^{\infty}(\re^{n-1})}
     \Big\|_{L^{\r}_t(\re_+)}
 \\
 \le &C\sup_{k\in \Z}
        \sup_{j\in \Z}2^{\frac{n-1}{p}j} 
        \Big\| 
        \big\| Q_{k-2}\underset{(t)}{*}
               \widetilde{\phi}_{j}\underset{(x')}{*}
               F(t)
        \big\|_{L^{p}(\re^{n-1})}
        \Big\|_{L^{\r}_t(\re_+)}\\
     &\hskip1cm \times  
      \sum_{k\in \Z} 2^{(\frac12-\frac{1}{2p})k}
      \sum_{2m\le k}2^{(-1+\frac{1}{p})m}                            
      \Big\|
        \|\widetilde{\phi}_m\|_{p'}
        \big\| 
             \widetilde{\psi}_{k}\underset{(t)}{*}
             \phi_m\underset{(x')}{*}
             G(t)\big)
        \big\|_{L^{p}(\re^{n-1})}
     \Big\|_{L^{\infty}_t(\re_+)}
\\
  \le& C\big\|F(t)
       \big\|_{\widetilde{L^{\r}_t(\re_+};\dB^{\frac{n-1}{p}}_{p,\infty}(\re^{n-1}))}
       \big\| G(t)
       \big\|_{\widetilde{\dB^{\frac12-\frac{1}{2p}}_{\infty,1}
              (\re_+};\dB^{-1+\frac{n}{p}}_{p,1}(\re^{n-1}))}.
}
Here we need $1<p$.
{\bl While the diagonal term $h_1^3$ can be estimated as
follows.
\algn{
  h_1^3
 \le &
      \sum_{j\in \Z} 2^{(-1+\frac{n}{p})j}
      \sum_{k\ge 2j}  2^{(\frac12-\frac{1}{2p})k}  \\
      &\hskip5mm \times   
       \Big\|  \big\|
            \phi_{j}\underset{(x')}{*}   
            \Big(\sum_{\ell\ge k-2}
                \big(
                  \psi_{\ell}\underset{(t)}{*}
                  \widetilde{\phi}_{j}\underset{(x')}{*}F(t)
                \big)\cdot
                \big(  
                 \widetilde{\psi}_{\ell}\underset{(t)}{*}
                  P_{j-2}\underset{(x')}{*}
                  G(t)\big)
            \Big)
        \big\|_{L^p({\gr \re^{n-1}})}
     \Big\|_{L^{\r}_t(\re_+)}  
 \\
 \le &
      \sum_{j\in \Z} 2^{(-1+\frac{n}{p})j}
      \sum_{k\ge 2j}  2^{(\frac12-\frac{1}{2p})k}  \\
      &\hskip1mm \times   
       \Big\|  \big\|
            \phi_{j}
        \big\|_{L^1({\gr \re^{n-1}})}   
        \sum_{\ell\ge k-2}
        \big\|
              \psi_{\ell}\underset{(t)}{*}
              \widetilde{\phi}_{j}\underset{(x')}{*}F(t)
        \big\|_{L^p({\gr \re^{n-1}})}
        \sup_{\ell,j}
        \Big\| 
             \widetilde{\psi}_{\ell}\underset{(t)}{*}
             P_{j-2}\underset{(x')}{*}
             G(t)
        \Big\|_{L^{\infty}({\gr \re^{n-1}})}
     \Big\|_{L^{\r}_t(\re_+)}  
\\
 \le &C
     \sum_{j\in \Z} 2^{(-1+\frac{n}{p})j}
     \sum_{k\ge 2j}  2^{(\frac12-\frac{1}{2p})k}  \\
      &\hskip1mm \times   
        \sum_{\ell\ge k-2}\Big\|
        \big\|
            \psi_{\ell}\underset{(t)}{*}
            \widetilde{\phi}_{j}\underset{(x')}{*}F(t)
        \big\|_{L^p({\gr \re^{n-1}})}
    \Big\|_{L^{\r}_t(\re_+)}  
    \sup_{\ell,j}
         \Big\|
         \Big\| 
             \widetilde{\psi}_{\ell}\underset{(t)}{*}
             P_{j-2}\underset{(x')}{*}
             G(t)
        \Big\|_{L^{\infty}({\gr \re^{n-1}})}
     \Big\|_{L^{\infty}_t(\re_+)}  
\\
 \le &C
   \sum_{j\in \Z} 2^{(-1+\frac{n}{p})j}
   \sum_{\ell\ge 2j-2} 
        \Big(
         \sum_{k\le \ell+2}
          2^{(\frac12-\frac{1}{2p})k}  
        \Big)
   \Big\|\big\|
          \psi_{\ell}\underset{(t)}{*}
          \widetilde{\phi}_{j}\underset{(x')}{*}F(t)
        \big\|_{L^p({\gr \re^{n-1}})}
    \Big\|_{L^{\r}_t(\re_+)}  
    \Big\|
       G(t)
    \Big\|_{L^{\infty}_t(\re_+;L^{\infty}({\gr \re^{n-1}}))} 
 \\
 \le &C
   \sum_{j\in \Z} 2^{(-1+\frac{n}{p})j}
   \sum_{\ell\ge 2j-2}  2^{(\frac12-\frac{1}{2p})\ell}  
    \Big\| 
        \big\|
              \psi_{\ell}\underset{(t)}{*}
                  \widetilde{\phi}_{j}\underset{(x')}{*}F(t)
        \big\|_{L^p({\gr \re^{n-1}})}
    \Big\|_{L^{\r}_t(\re_+)}  
    \Big\|
      G(t)
    \Big\|_{L^{\infty}_t(\re_+;L^{\infty}({\gr \re^{n-1}}))} 
\\
    \le & 
     C\|F(t)
      \|_{\widetilde{{\gr \dB^{\frac12-\frac{1}{2p}}_{\r,1}}
         (\re_+;}\dB^{-1+\frac{n}{p}}_{p,1}(\re^{n-1}))}
      \big\| G(t)
      \big\|_{\widetilde{L^{\infty}(\re_+;}\dB^{\frac{n-1}{p}}_{p,1}({\gr \re^{n-1}}))}.
 }
There is no restriction except $1<p$. 

}
For the second term $h_2$, we decompose by the time direction and 
typical term can be estimated as follows:
\algn{
  h_2^1
\le & \sum_{j\in \Z} 2^{(-1+\frac{n}{p})j}
      \sum_{k\ge 2j}  2^{(\frac12-\frac{1}{2p})k}   
      \Big\|     
          |\psi_{k}|\underset{(t)}{*}
          \big\| 
            \phi_{j}
          \big\|_{L^1(\re^{n-1})}  
          \big\|
              \widetilde{\psi}_{k}\underset{(t)}{*}
              P_{j-2}\underset{(x')}{*}F(t)
          \big\|_{L^{\infty}(\re^{n-1})}
       \\
    &\hskip4cm \times 
          \big\|
               Q_{k-2}\underset{(t)}{*}
               \phi_j\underset{(x')}{*}
               G(t)
           \big\|_{L^{p}(\re^{n-1})}
     \Big\|_{L^{\r}_t(\re_+)}  
\\  
\le & C\sum_{j\in \Z} 2^{(-1+\frac{n}{p})j}
        \sum_{k\ge 2j}  2^{(\frac12-\frac{1}{2p})k}   
      \Big\| \big\|\sum_{m\le j-2} 
                \phi_m\underset{(x')}{*}
                \widetilde{\psi}_{k}\underset{(t)}{*}
                F(t)
          \big\|_{L^{\infty}(\re^{n-1})}
       \\
    &\hskip4cm \times 
          \big\|
              Q_{k-2}\underset{(t)}{*}
              \phi_j\underset{(x')}{*}
              G(t)
          \big\|_{L^{p}(\re^{n-1})}
     \Big\|_{L^{\r}_t(\re_+)}  
\\
\le & C \sum_{j\in \Z} 2^{(-1+\frac{n}{p})j}  
        \sum_{m\le j-2}
        \sum_{k\ge 2j}  2^{(\frac12-\frac{1}{2p})k}         
       \Big\|  \big\| \phi_m\underset{(x')}{*}
                 \widetilde{\psi}_{k}\underset{(t)}{*}F(t)
         \big\|_{L^{\infty}(\re^{n-1})}
       \Big\|_{L^{\r}_t(\re_+)} \\
   &\hskip4cm\times
        \Big\| \big\|\phi_{j}\underset{(x')}{*}
               G(t)
           \big\|_{L^{p}(\re^{n-1})}
     \Big\|_{L^{\infty}_t(\re_+)} 
 \\
 \le & C  \sum_{m\in \Z}   
          \sum_{j\ge m+2}2^{{\rd(-1+\frac{1}{p})}j}
          \sum_{k\ge 2j}  2^{(\frac12-\frac{1}{2p})k}                      
     \Big\|      
        \big\| \phi_m\underset{(x')}{*}
               \widetilde{\psi}_{k}\underset{(t)}{*}F(t)
        \big\|_{L^{\infty}(\re^{n-1})}
       \Big\|_{L^{\r}_t(\re_+)} \\
   &\hskip4cm\times
     {\bl  \sup_{j\in \Z}}
     \Big\| 2^{\frac{n-1}{p}j} \big\|\phi_{j}\underset{(x')}{*}
            G(t)
           \big\|_{L^{p}(\re^{n-1})}
     \Big\|_{L^{\infty}_t(\re_+)} 
 \\
  \le & C\sum_{k\in\Z}  2^{(\frac12-\frac{1}{2p})k}   
         \sum_{m\in \Z} 2^{(-1+\frac{1}{p})m}2^{\frac{n-1}{p}m}      
         \Big\| 
          \big\| \phi_m\underset{(x')}{*}
                 \widetilde{\psi}_{k}\underset{(t)}{*}F(t)
          \big\|_{L^{p}(\re^{n-1})}
         \Big\|_{L^{\r}_t(\re_+)} \\
   &\hskip4cm\times
     {\bl \sup_{j\in \Z} }
     \Big\|  2^{\frac{n-1}{p}j}
           \big\|\phi_{j}\underset{(x')}{*}
                 G(t)
           \big\|_{L^{p}(\re^{n-1})}
     \Big\|_{L^{\infty}_t(\re_+)} 
 \\   
   \le& C\big\|F(t)
         \big\|_{\widetilde{\dB^{\frac12-\frac{1}{2p}}_{\r,1}
                   (\re_+;}\dB^{-1+\frac{n}{p}}_{p,1}(\re^{n-1}))}
       \big\| G(t)           
       \big\|_{\widetilde{L^{\infty}
               (\re_+;}\dB^{\frac{n-1}{p}}_{p,\infty}(\re^{n-1}))},
  }
where we require the {\crd restriction}  $1<p$ again.
The other terms $h_2^2$ can be treated in similar manner.
{\bl 
Space off diagonal and time diagonal term term $h_2^3$ can be dominated by 
\algn{
  h_2^3
 \le &
     \sum_{j\in \Z} 2^{(-1+\frac{n}{p})j}
     \sum_{k\ge 2j}  2^{(\frac12-\frac{1}{2p})k}  \\
      &\hskip5mm \times   
      \Big\|   \sum_{\ell\ge k-2}
        \big\|
                \big(
                  \psi_{\ell}\underset{(t)}{*}
                  P_{j-2}\underset{(x')}{*}F(t)
                \big)
        \big\|_{L^{\infty}({\gr \re^{n-1}})}
        \Big\| 
             \widetilde{\psi}_{\ell}\underset{(t)}{*}
             \widetilde{\phi}_j\underset{(x')}{*}
             {\gr G(t) }
        \Big\|_{L^{p}({\gr \re^{n-1}})}  
      \Big\|_{L^{\r}_t(\re_+)} 
\\
 \le &  \sum_{j\in \Z} 2^{(-1+\frac{n}{p})j}
        \sum_{k\ge 2j} \sum_{\ell\ge k-2}
         2^{(\frac12-\frac{1}{2p})k} 
        \Big\|  
        \big\|
              \psi_{\ell}\underset{(t)}{*}
              P_{j-2}\underset{(x')}{*}F(t)
        \big\|_{L^{\infty}({\gr \re^{n-1}})} 
        \Big\| 
             \widetilde{\psi}_{\ell}\underset{(t)}{*}
             \widetilde{\phi}_j\underset{(x')}{*}
             {\gr G(t)}
        \Big\|_{L^{p}({\gr \re^{n-1}})} 
      \Big\|_{L^{\r}_t(\re_+)} 
\\
 \le &C\sup_{\ell\in \Z}\sup_{j\in \Z}
   \Big\|
   \big\|
         \psi_{\ell}\underset{(t)}{*}
         P_{j-2}\underset{(x')}{*}F(t)
   \big\|_{L^{\infty}({\gr \re^{n-1}})}
   \Big\|_{L^{\r}_t(\re_+)}  \\
  &\hskip1cm \times
        \sum_{j\in \Z} 2^{(-1+\frac{n}{p})j}
        \sum_{k\ge 2j} \sum_{\ell\ge k-2}
         2^{(\frac12-\frac{1}{2p})k} 
    \Big\|
    \Big\| 
           \widetilde{\psi}_{\ell}\underset{(t)}{*}
           \widetilde{\phi}_j\underset{(x')}{*}
           G(t)
     \Big\|_{L^{p}({\gr \re^{n-1}})}
    \Big\|_{L^{\infty}_t(\re_+)} 
 \\
   \le &C\sup_{\ell\in \Z}\sup_{j\in \Z}
   \Big\|
   \big\| 
         \psi_{\ell}\underset{(t)}{*}
         P_{j-2}\underset{(x')}{*}F(t)
   \big\|_{L^{\infty}({\gr \re^{n-1}})}
   \Big\|_{L^{\r}_t(\re_+)}  \\
  &\hskip1cm \times
        \sum_{j\in \Z} 2^{(-1+\frac{n}{p})j}
        \sum_{\ell\ge 2j-2} \sum_{k\le \ell+2} 
         2^{(\frac12-\frac{1}{2p})k} 
    \Big\|
     \Big\|
           \widetilde{\psi}_{\ell}\underset{(t)}{*}
           \widetilde{\phi}_j\underset{(x')}{*}
           G(t)
     \Big\|_{L^{p}({\gr \re^{n-1}})}
    \Big\|_{L^{\infty}_t(\re_+)} 
\\
 \le &C\big\|
         F(t)
      \big\|_{\gr L^{\r}(\re;L^{\infty}({\gr \re^{n-1}}))} 
        \sum_{j\in \Z} 2^{(-1+\frac{n}{p})j}
        \sum_{\ell\ge 2j-2}  
         2^{(\frac12-\frac{1}{2p})\ell}  
      \Big\|                      
          \Big\|
              \widetilde{\psi}_{\ell}\underset{(t)}{*}
              \widetilde{\phi}_j\underset{(x')}{*}
              G(t)      
          \Big\|_{L^{p}({\gr \re^{n-1}})}
    \Big\|_{L^{\infty}_t(\re_+)} 
\\
 \le &C\big\|
         F(t)
       \big\|_{\widetilde{L^{\r}(\re_+;}\dB^{\frac{n-1}{p}}_{p,1}({\gr \re^{n-1}}))}  
       \Big\|                      
          G(t)      
       \Big\|_{\gr\widetilde{\dB^{\frac12-\frac{1}{2p}}_{\infty,1}(\re_+:}
               \dB^{-1+\frac{n}{p}}_{p,1}({\gr \re^{n-1}})}.
}

We estimate for the {\gr third term $h_3$ of right hand side 
in \eqref{eqn;space-Bony-decomp0}.  It}  can be 
dominated by setting $1/p=1/r+2/p -1$ with $r=p'$, 
\algn{
 h_3^1 
\le &\sum_{j\in \Z} 2^{(-1+\frac{n}{p})j}
     \sum_{k\ge 2j} 2^{(\frac12-\frac{1}{2p})k}    \\
     &\hskip1cm \times
      \Big\| 
         \|\phi_{j}\|_{L^r({\gr \re^{n-1}})}
         \Big\|
            \sum_{m\ge j-2}
              \big(
                  \widetilde{\psi}_{k}\underset{(t)}{*}
                \phi_{m}\underset{(x')}{*}F(t)
                \big)\cdot
                \big(  
                  Q_{k}\underset{(t)}{*}
                  \widetilde{\phi}_{m}\underset{(x')}{*}
                 G(t)\big)            
        \Big\|_{L^{\frac{p}{2}}({\gr \re^{n-1}})}
     \Big\|_{L^{\r}_t(\re_+)}
\\
\le &C
      \sum_{j\in \Z}  2^{(-1+\frac{n}{p}+\frac{n-1}{r'})j}
      \sum_{k\ge 2j}  2^{(\frac12-\frac{1}{2p})k}    \\
   &\hskip2cm\times
      \bigg\|    \sum_{m\ge j-2}
         \Big\|
            \widetilde{\psi}_{k}\underset{(t)}{*}
            \phi_{m}\underset{(x')}{*}F(t)
          \Big\|_{L^{r'}({\gr \re^{n-1}})}
          \Big\|
             Q_{k}\underset{(t)}{*}
             \widetilde{\phi}_{m}\underset{(x')}{*}
             G(t)           
        \Big\|_{L^{p}({\gr \re^{n-1}})}
     \bigg\|_{L^{\r}_t(\re_+)}
\\%
\le &C\sum_{m\in \Z} 
      \sum_{j-2\le m} 2^{(-1+\frac{n}{p}+\frac{n-1}{p})j}
      \sum_{k\ge 2j}  2^{(\frac12-\frac{1}{2p})k}         
      \bigg\| 
         \Big\|            
            \widetilde{\psi}_{k}\underset{(t)}{*}
            \phi_{m}\underset{(x')}{*}F(t)
          \Big\|_{L^{p}({\gr \re^{n-1}})}\\
    &\hskip6cm\times         
        \Big\|
          Q_{k}\underset{(t)}{*}
          \widetilde{\phi}_{m}\underset{(x')}{*}
          G(t)        
        \Big\|_{L^{p}({\gr \re^{n-1}})}
     \bigg\|_{L^{\r}_t(\re_+)}
\\
\le &C
      \sum_{m\in \Z}  2^{(-1+\frac{n}{p}+\frac{n-1}{p})m} 2^{-\frac{n-1}{p}m}
      \sum_{k\in \Z}  2^{(\frac12-\frac{1}{2p})k}         
      \bigg\|    
        \Big\|            
            \widetilde{\psi}_{k}\underset{(t)}{*}
            \phi_{m}\underset{(x')}{*}F(t)
        \Big\|_{L^{p}({\gr \re^{n-1}})}
     \bigg\|_{L^{\r}_t(\re_+)} \\
     &\hskip4cm \times
     \sup_{\gr k\in \Z}
     \sup_{m\in \Z}   2^{\frac{n-1}{p}m}
     \bigg\|
     \Big\|
          Q_{k}\underset{(t)}{*}
          \widetilde{\phi}_{m}\underset{(x')}{*}
          G(t)          
     \Big\|_{L^{p}({\gr \re^{n-1}})}
     \bigg\|_{L^{\infty}_t(\re_+)} 
 \\
    \le &C 
    \|F\|_{\gr \widetilde{\dB^{\frac12-\frac{1}{2p}}_{\r,1}
           (\re_+;}\dB^{-1+\frac{n}{p}}_{p,1}(\re^{n-1}))}     
     \sup_{m\in \Z}
     \bigg\|          
         2^{\frac{n-1}{p}m}
       \big\|\widetilde{\phi}_{m}\underset{(x')}{*} 
             G(t)      
       \big\|_{L^{p}({\gr \re^{n-1}})}          
     \bigg\|_{L^{\infty}_t(\re_+)}
\\
   \le &C 
    \|F\|_{\gr\widetilde{\dB^{\frac12-\frac{1}{2p}}_{\r,1}
           (\re_+;}\dB^{-1+\frac{n}{p}}_{p,1}(\re^{n-1}))}         
    \big\|
        G(t)     
    \big\|_{\widetilde{L^{\infty}(\re_+;}\dB^{\frac{n-1}{p}}_{p,1}({\gr \re^{n-1}}))},
}
where we used  $p<2n-1$.
The second term can be estimated by a very similar way:
\noindent
\algn{
h_3^2
\le &C 
      \sum_{j\in \Z} 2^{(-1+\frac{n}{p}+\frac{n-1}{p})j}
      \sum_{k\ge 2j}  2^{(\frac12-\frac{1}{2p})k}    \\
   &\hskip2cm\times
      \bigg\|
         \sum_{m\ge j-2}
         \Big\|
             Q_{k}\underset{(t)}{*}
            \phi_{m}\underset{(x')}{*}F(t)
          \Big\|_{L^{p}({\gr \re^{n-1}})}
          \Big\|            
             \widetilde{\psi}_{k}\underset{(t)}{*}
             \widetilde{\phi}_{m}\underset{(x')}{*}
             G(t)           
        \Big\|_{L^{p}({\gr \re^{n-1}})}
     \bigg\|_{L^{\r}_t(\re_+)}
\\
\le &C 
      \sum_{j\in \Z} 
      \sum_{m\ge j-2} 2^{(-1+\frac{n}{p}+\frac{n-1}{p})j}
      \sum_{k\ge 2j}  2^{(\frac12-\frac{1}{2p})k}         
       \bigg\|  \Big\|            
             Q_{k}\underset{(t)}{*}
            \phi_{m}\underset{(x')}{*}F(t)
          \Big\|_{L^{p}({\gr \re^{n-1}})}\\
    &\hskip6cm\times         
        \Big\|
             \widetilde{\psi}_{k}\underset{(t)}{*}
             \widetilde{\phi}_{m}\underset{(x')}{*}
             G(t)           
        \Big\|_{L^{p}({\gr \re^{n-1}})}
     \bigg\|_{L^{\r}_t(\re_+)}
\\
  \le &C 
       \sum_{m\in \Z}  2^{-\frac{n-1}{p}m}         
       \sum_{j\le m+2} 2^{(-1+\frac{n}{p}+\frac{n-1}{p})j}               
       \sum_{k\in \Z}  2^{(\frac12-\frac{1}{2p})k}         
        \\
    &\hskip2cm \times   
     \bigg\|   2^{\frac{n-1}{p}m}
        \Big\|  
            Q_{k}\underset{(t)}{*}      
            \phi_{m}\underset{(x')}{*}
            F(t)
        \Big\|_{L^{p}({\gr \re^{n-1}})} 
    \bigg\|_{L^{\r}_t(\re_+)}
    \bigg\| 
        \Big\|  
          \widetilde{\psi}_{k}\underset{(t)}{*}
          \widetilde{\phi}_{m}\underset{(x')}{*}
          G(t)           
        \Big\|_{L^{p}({\gr \re^{n-1}})}
     \bigg\|_{L^{\infty}_t(\re_+)} 
\\
   \le &C 
       \sum_{m\in \Z }
          2^{-\frac{n-1}{p}m}         
          2^{(-1+\frac{n}{p}+\frac{n-1}{p})m}               
     \sum_{k\in \Z}  
          2^{(\frac12-\frac{1}{2p})k}          
        \\
    &\hskip2cm \times   
     \bigg\|   2^{\frac{n-1}{p}m}
        \Big\|  
           \phi_{m}\underset{(x')}{*}F(t)
        \Big\|_{L^{p}({\gr \re^{n-1}})} 
    \bigg\|_{L^{\r}_t(\re_+)}
    \bigg\| 
        \Big\|  
          \widetilde{\psi}_{k}\underset{(t)}{*}
          \widetilde{\phi}_{m}\underset{(x')}{*}
          G(t)           
        \Big\|_{L^{p}({\gr \re^{n-1}})}
     \bigg\|_{L^{\infty}_t(\re_+)} 
 \\
  \le &C 
      \sup_{m\in \Z}
            2^{\frac{n-1}{p}m}
     \bigg\|                 
        \Big\|  
            \phi_{m}\underset{(x')}{*}F(t)
        \Big\|_{L^{p}({\gr \re^{n-1}})}  
    \bigg\|_{L^{\r}_t(\re_+)} \\
    &\hskip2cm \times
          \sum_{k\in \Z}  2^{(\frac12-\frac{1}{2p})k}   
          \sum_{m\in \Z}    2^{(-1+\frac{n}{p})m}
    \bigg\|  
        \Big\|    
          \widetilde{\psi}_{k}\underset{(t)}{*}       
          \widetilde{\phi}_{m}\underset{(x')}{*}
          G(t)           
        \Big\|_{L^{p}({\gr \re^{n-1}})}
     \bigg\|_{L^{\infty}_t(\re_+)} 
 \\
    \le &C 
    \|F(t)\|_{\widetilde{L^{\r}(\re_+;}\dB^{\frac{n-1}{p}}_{p,1}(\re^{n-1}))}      
    \big\| G(t)    
    \big\|_{\gr\widetilde{\dB^{\frac12-\frac{1}{2p}}_{\infty,1}
           (\re_+;}\dB^{-1+\frac{n}{p}}_{p,\infty}({\gr \re^{n-1}}))}.
}
This case we again  need the restriction $1<p<2n-1$. 
}
\algn{
h_3^3
\le &C
     \sum_{j\in \Z} 2^{(-1+\frac{n}{p})j}
     \sum_{k\ge 2j}  2^{(\frac12-\frac{1}{2p})k}   \\
     &\hskip3mm\times
       \|\phi_j\|_{L^{p'}(\re^{n-1})}  
     \Big\| 
       |\psi_{k}|\underset{(t)}{*}
       \Big\|  
            \Big(\sum_{\ell\ge k-2}
                 \sum_{m\ge j-2}
                \big(
                 \psi_{\ell}\underset{(t)}{*}
                 \phi_{m}\underset{(x')}{*}F(t)
                \big)\cdot
                \big(  
                 \widetilde{\psi}_{\ell}\underset{(t)}{*}
                 \widetilde{\phi}_{m}\underset{(x')}{*}
                 G(t)\big)
            \Big)
        \Big\|_{\bl L^{\frac{p}{2}}(\re^{n-1})}
     \Big\|_{L^{\r}_t(\re_+)} 
\\
\le &C \sum_{j\in \Z} 2^{(-1+\frac{n}{p}+\frac{n-1}{p})j}
       \sum_{k\ge 2j}  2^{(\frac12-\frac{1}{2p})k}   \\
     &\hskip3mm\times
     \bigg\|       
       \sum_{\ell\ge k-2}
       \sum_{m\ge j-2}
       \Big\|
          \psi_{\ell}\underset{(t)}{*}
          \phi_{m}\underset{(x')}{*}F(t)
        \Big\|_{L^{p}(\re^{n-1})}
        \Big\| 
          \widetilde{\psi}_{\ell}\underset{(t)}{*}
          \widetilde{\phi}_{m}\underset{(x')}{*}
          G(t)
        \Big\|_{L^p(\re^{n-1})}
     \bigg\|_{L^{\r}_t(\re_+)}
\\
\le &C\sum_{m\in \Z} 
      \sum_{j\le m+2} 
        2^{(-1+\frac{n}{p}+\frac{n-1}{p})j}
      \sum_{\ell\in \Z}\sum_{2j\le k\le \ell+2} 
        2^{(\frac12-\frac{1}{2p})k}        
         \\
     &\hskip1cm\times
     \bigg\|             
       \Big\|
          \psi_{\ell}\underset{(t)}{*}
          \phi_{m}\underset{(x')}{*}F(t)
        \Big\|_{L^{p}(\re^{n-1})}
      \bigg\|_{L^{\r}(\re_+)}
      \bigg\|
        \Big\| 
          \widetilde{\psi}_{\ell}\underset{(t)}{*}
          \widetilde{\phi}_{m}\underset{(x')}{*}
          G(t)
        \Big\|_{L^p(\re^{n-1})}
     \bigg\|_{L^{\infty}(\re_+)}
 \\
 \le &C\sum_{m\in \Z} 
       \sum_{j\le m+2} 
        2^{(-1+\frac{n}{p}+\frac{n-1}{p})j}
        2^{-\frac{n-1}{p}m}
       \sum_{\ell\in \Z}   
        2^{(\frac12-\frac{1}{2p})\ell}  
         \\
     &\hskip1cm\times
     \bigg\|             
       \Big\|
          \psi_{\ell}\underset{(t)}{*}
          \phi_{m}\underset{(x')}{*}F(t)
        \Big\|_{L^{p}(\re^{n-1})}
      \bigg\|_{L^{\r}(\re_+)}
      \bigg\|  2^{\frac{n-1}{p}m}
        \Big\| 
          \widetilde{\psi}_{\ell}\underset{(t)}{*}
          \widetilde{\phi}_{m}\underset{(x')}{*}
          G(t)
        \Big\|_{L^p(\re^{n-1})}
     \bigg\|_{L^{\infty}(\re_+)}
 \\
  \le &C\sum_{m\in \Z}    2^{(-1+\frac{n}{p})m}
        \sum_{\ell\in \Z} 2^{(\frac12-\frac{1}{2p})\ell}     
      \bigg\|      \Big\|
          \psi_{\ell}\underset{(t)}{*}
          \phi_{m}\underset{(x')}{*}F(t)
        \Big\|_{L^{p}(\re^{n-1})}
      \bigg\|_{L^{\r}(\re_+)}  \\
   &\hskip4cm \times
    \sup_{\ell\in \Z}
      \bigg\|
      \sup_{m\in \Z} 2^{\frac{n-1}{p}m}
        \Big\| 
          \widetilde{\psi}_{\ell}\underset{(t)}{*}
          \widetilde{\phi}_{m}\underset{(x')}{*}
          G(t)
        \Big\|_{L^p(\re^{n-1})}
     \bigg\|_{L^{\infty}_t(\re_+)}
 \\
 \le &C 
    \big\| F(t) 
    \big\|_{\widetilde{\dB^{\frac12-\frac{1}{2p}}_{\r,1}
            (\re_+;}\dB^{-1+\frac{n}{p}}_{p,1}(\re^{n-1}))}
    \big\| G(t)
    \big\|_{\widetilde{L^{\infty}(\re_+;}\dB^{\frac{n-1}{p}}_{p,1}(\re^{n-1}))},
}
where we used the condition $p<2n-1$.
The other terms $h_3^1$ and $h_3^2$ can be estimated in a similar way.
\end{prf}

}


\vskip2mm
{\bf Acknowledgments}
The authors would like to thank the referee for careful reading and valuable 
comments.
The first author is partially supported by JSPS grant-in-aid for 
Scientific Research (S) \#19H05597 and Challenging Research 
(Pioneering) \#20K20284. 
The second author is partially supported by JSPS grant-in-aid 
for Scientific Research (B) \#21H00992 and by 
JSPS Fostering Joint Research Program (B) \#18KK0072. 
\par
Conflict of interest: The authors declare that there is no conflict of interest.


\end{document}